\documentclass{amsart} \usepackage{amsmath} \usepackage{amscd}
\usepackage{amsxtra}
\usepackage{verbatim}
\usepackage{multicol}
\usepackage[inline]{enumitem}
\usepackage{mathtools}
\newlength\NotationWidth
\makeatletter
\newcommand\tabfill[1]{%
  \parbox[t]\NotationWidth{%
    \raggedright\leftskip=2em\hspace*{-2em}#1\ifhmode\unskip\nobreak\strut\fi}%
}
\makeatother

\usepackage{amssymb}
\usepackage[mathscr]{eucal}	
\DeclareSymbolFont{pssymbols}     {OMS}{ztmcm}{m}{n}
\DeclareSymbolFontAlphabet{\mathpsscr}   {pssymbols}

\usepackage{amsthm}
\swapnumbers
\theoremstyle{plain}
\newtheorem{thm}[subsection]{Theorem}
\newtheorem{cor}[subsection]{Corollary}
\newtheorem{prop}[subsection]{Proposition}
\newtheorem{lem}[subsection]{Lemma}
\newtheorem{sublem}[subsection]{Sublemma}
\newtheorem{BasicLemma}[subsection]{The Basic Lemma}
\newtheorem*{thm*}{Theorem}
\newtheorem*{cor*}{Corollary}
\newtheorem*{prop*}{Proposition}
\newtheorem*{lem*}{Lemma}
\newtheorem*{RapoportConjecturethm}{Theorem \ref{ssectRapoportConjecture}}
\newtheorem*{vanthm}{Theorem \ref{ssectGlobalVanishing}}
\newtheorem*{SSWHthm}{Theorem \ref{ssectWeightCohomologyMicroSupport}}
\newtheorem*{SSIHcor}{Corollary \ref{ssectIHMicroPurityCorollary}}
\newtheorem*{SSFibercor}{Corollary
\ref{ssectRestrictMicroSupportToFiberCorollary}}
\newtheorem*{GoreskyHarderMacPhersonthm}{Theorem \ref{ssectGoreskyHarderMacPhersonTheorem}}
\newtheorem*{corAlternateDegreeRanges}{Corollary to Proof of Proposition~\ref{ssectAlternateMicroSupport}}

\theoremstyle{definition}
\newtheorem{defn}[subsection]{Definition}
\newtheorem*{defn*}{Definition}

\theoremstyle{remark}
\newtheorem{rem}[subsection]{Remark}

\newtheorem*{rem*}{Remark}
\newtheorem*{rems*}{Remarks}

\newtheorem*{note*}{Note}

\setlist[enumerate]{font=\normalfont}
\setlist[enumerate,1]{label=\textup{(\roman*)}}
\newcommand{\itemref}[1]{\ref{#1}}
\renewcommand{\theenumi}{\roman{enumi}}
\renewcommand{\labelenumi}{(\theenumi)}

\numberwithin{equation}{subsection}
\newcounter{saveenum}           

\usepackage[arrow,matrix,tips,frame,curve,ps]{xy}
\SelectTips{cm}{}
\UseTips
\CompileMatrices

\newdir^{ (}{{}*!/-5pt/@^{(}}
\newdir_{ (}{{}*!/-5pt/@_{(}}

\renewcommand{\AA}{{\mathbb A}}
\newcommand{\CC}{{\mathbb C}}
\newcommand{\RR}{{\mathbb R}}
\newcommand{\ZZ}{{\mathbb Z}}
\newcommand{\NN}{{\mathbb N}}
\newcommand{\QQ}{{\mathbb Q}}
\newcommand{\EE}{{\mathbb E}}
\newcommand{\EEtilde}{{\widetilde{\EE}}}

\newcommand{\HH}{{\mathbb H}}
\newcommand{\VV}{{\mathbb V}}

\DeclareMathOperator{\Ad}{Ad}
\DeclareMathOperator{\codim}{codim}

\newcommand{\CoDifferential}{\delta}
\newcommand{\ModTwo}{\delta}
\DeclareMathOperator{\dom}{Dom}
\newcommand{\Dom}{\dom}
\DeclareMathOperator{\rank}{rank}
\DeclareMathOperator{\Int}{int}
\DeclareMathOperator{\Ker}{Ker}
\renewcommand{\ker}{\Ker}
\renewcommand{\Im}{\operatorname{Im}}
\newcommand{\PP}{\mathbb P}
\DeclareMathOperator{\Range}{Range}
\DeclareMathOperator{\range}{\Range}
\DeclareMathOperator{\QQrank}{\QQ-rank}
\DeclareMathOperator{\CCrank}{\CC-rank}
\DeclareMathOperator{\Span}{Span}
\newcommand{\dbar}{\overline{d}}
\DeclareMathOperator{\supp}{supp}
\newcommand{\geo}{\mathbin{\mathbf o}}
\DeclareMathOperator{\Vol}{Vol}
\DeclareMathOperator{\mS}{SS}
\DeclareMathOperator{\emS}{\mS_{ess}}
\newcommand{\essc}{c_{\textup{ess}}}
\newcommand{\essd}{d_{\textup{ess}}}
\DeclareMathOperator{\Res}{Res}

\DeclareMathOperator{\pr}{pr}
\DeclareMathOperator{\Cl}{cl}
\newcommand{\cl}[1]{\Cl(#1)}
\newcommand{\dd}{d}		
\renewcommand{\l}{\ell}         

\DeclareMathOperator{\restr}{\rho}	
\newcommand{\lvvv}{\lVert}				
\newcommand{\rvvv}[1][]{\rVert_{\if!#1!b\else#1,b\fi}}	

\newcommand{\lsb}[1]{{}_{#1}}
\newcommand{\lsp}[1]{{}^{#1}\!}

\newcommand{\tildearrow}{\xrightarrow{\,\sim\,}}
\newcommand{\qiarrow}{\xrightarrow{\text{q.i.}}}
\newcommand{\longtildearrow}{\tildearrow}
\newbox\arrowbox
\setbox\arrowbox=\hbox{\lower .5ex\rlap{$\longrightarrow$}%
 \raise .5ex\hbox{$\longleftarrow$}}

\newcommand{\G}{\Gamma}

\DeclareMathOperator{\GL}{GL}
\DeclareMathOperator{\SO}{SO}
\DeclareMathOperator{\Sympl}{Sp}
\DeclareMathOperator{\Gal}{Gal}

\DeclareMathOperator{\Hom}{Hom}
\newcommand{\Pl}{\mathscr P}	

\newcommand{\Dbar}{\overline{D}}
\newcommand{\Xbar}{\overline{X}}
\newcommand{\Xhat}{\widehat{X}}
\newcommand{\back}{\backslash}
\newcommand{\Dstar}{D^*}
\newcommand{\Xstar}{X^*}
\newcommand{\Abar}{{\bar{A}}} 
\newcommand{\Ybar}{\overline{Y}}
\newcommand{\Ftilde}{\widetilde{F}}
\newcommand{\X}{\mathcal X}

\renewcommand{\i}{i}		
\newcommand{\ihat}{{\hat{\imath}}}
\renewcommand{\j}{j}		

\newcommand{\jhat}{{\hat{\jmath}}}
\DeclareMathOperator{\id}{id}
\newcommand{\nil}{p}	
\newcommand{\coh}{b}  

\newcommand{\etahat}{{\hat{\eta}}}

\newcommand{\Wbar}{{\overline{W}}}
\newcommand{\What}{{\widehat{W}}}

\newcommand{\IpC}{{\mathcal I_p\mathcal C}}
\renewcommand{\H}{\mathcal H}

\newcommand{\Sheaf}{\mathcal S}
\newcommand{\A}{\mathcal A}
\newcommand{\Asp}{\A_{\,\textup{sp}}}
\newcommand{\Acomb}{\A_{\,\textup{comb}}}
\renewcommand{\sp}{_{\,\textup{sp}}}
\newcommand{\inv}{_{\textup{inv}}}
\newcommand{\WnC}{\mathcal W^\eta \mathcal C}

\renewcommand{\L}{\mathscr L}
\newcommand{\M}{\mathcal M}

\newcommand{\IH}{I_{p_w}H}

\newcommand{\sa}{{\mathfrak a}}
\newcommand{\h}{{\mathfrak h}}
\newcommand{\hb}{{\mathfrak b}}
\renewcommand{\k}{{\mathfrak k}}
\newcommand{\n}{{\mathfrak n}}
\newcommand{\m}{{\mathfrak m}}
\newcommand{\levi}{{\mathfrak l}}
\newcommand{\p}{{\mathfrak p}}
\renewcommand{\r}{{\mathfrak r}}

\newcommand{\al}{\alpha}
\newcommand{\atilde}{{\widetilde{\alpha}}}
\renewcommand{\b}{\beta}
\newcommand{\hsr}{\rho}
\renewcommand{\d}{\delta}
\newcommand{\dtilde}{{\widetilde{\d}}}
\newcommand{\D}{\Delta}
\newcommand{\Dhat}{{\widehat{\D}}}

\newcommand{\DRR}{\lsb\RR\D}
\newcommand{\e}{\varepsilon}
\newcommand{\g}{\gamma}	
\renewcommand{\v}{\nu}

\renewcommand{\t}{\tau}
\newcommand{\U}{\varUpsilon}
\newcommand{\Utilde}{\widetilde{\U}}
\newcommand{\Th}{\varTheta}
\newcommand{\Phitilde}{\widetilde{\Phi}}
\newcommand{\kap}{\kappa}
\renewcommand{\o}{\omega}

\makeatletter
\def\sphat{^{\mathchoice{}{}%
 {\,\,\smash[b]{\hbox{\lower4\ex@\hbox{$\m@th\widehat{\null}$}}}}%
 {\,\smash[b]{\hbox{\lower3\ex@\hbox{$\m@th\hat{\null}$}}}}}\,}
\makeatother

\newcommand{\Cat}{\mathscr C}	

\newcommand{\Rep}{\operatorname{\mathfrak M\mathfrak o\mathfrak d}}
\newcommand{\R}{\Rep}
\newcommand{\IrrRep}{\operatorname{\mathfrak I\mathfrak r\mathfrak r}}
\newcommand{\Sh}{\operatorname{\mathfrak S\mathfrak h}}
\newcommand{\Derived}{\operatorname{\mathbf D}}

\newcommand{\Complex}{\operatorname{\mathbf C}}
\newcommand{\K}{\operatorname{\mathbf K}}
\newcommand{\Graded}{\operatorname{\mathbf G\mathbf r}}

\makeatletter
\def\prime{{\null\prime@\null}}
\mathchardef\prime@="0230
\makeatother

\begin{document}


\author{Leslie Saper}
\address{Department of Mathematics\\ Duke University\\ Box 90320\\ Durham,
NC 27708\\U.S.A.}
\email{saper@math.duke.edu}
\urladdr{http://www.math.duke.edu/faculty/saper}
\title{$\mathscr L$-modules and micro-support}
\thanks{
This research was supported in part by the National Science Foundation
under grants DMS-8957216, DMS-9100383, and DMS-9870162, a grant from The
Duke Endowment, an Alfred P. Sloan Research Fellowship, a grant from the
Maximilian-Bickhoff-Stiftung to visit the Katholischen Universit\"at
Eichst\"att as Hermann-Minkowski Gastprofessur, and a grant from the
Institut des Hautes \'Etudes Scientifiques.  The author
wishes to thank these organizations for their hospitality and support.}
\thanks{The original manuscript was prepared with the \AmS-\LaTeX\ macro
system and the \Xy-pic\ package.}
\begin{abstract}
$\mathscr L$-modules are a combinatorial analogue of constructible complexes of sheaves
on the reductive Borel-Serre compactification of a locally symmetric space.
We define the \emph{micro-support} of an $\mathscr L$-module, a certain set
of irreducible modules for the Levi quotients of the parabolic
$\QQ$-subgroups associated to the strata.  We prove a vanishing theorem for
the global cohomology of a bounded $\mathscr L$-module in terms of its
micro-support.  We calculate the micro-support of the middle weight profile
weighted cohomology and the middle perversity intersection cohomology
$\mathscr L$-modules. (For intersection cohomology we must assume the
$\QQ$-root system has no component of type $D_n$, $E_n$, or $F_4$.)
Finally we prove a functoriality theorem concerning the behavior of
micro-support upon restriction of an $\mathscr L$-module to the pre-image
of a Satake stratum.  As an application we settle a conjecture made
independently by Rapoport and by Goresky and MacPherson, namely, that the
intersection cohomology (for either middle perversity) of the reductive
Borel-Serre compactification of a Hermitian locally symmetric space is
isomorphic to the intersection cohomology of the Baily-Borel-Satake
compactification.  We also obtain a new proof of a result of
Goresky, Harder, and MacPherson on weighted cohomology as well as
generalizations of both of these results to general Satake
compactifications with equal-rank real boundary components.
\end{abstract}
\maketitle


\setcounter{tocdepth}{1}
\tableofcontents


\section{Introduction}
\label{sectIntro}
Let $X=\Gamma\back D$ be an arithmetic quotient of the symmetric space $D$
associated to a connected reductive algebraic group $G$ defined over $\QQ$.
(Here $D= G(\RR)/K A_G$, where $K\subset G(\RR)$ is a maximal compact
subgroup, $A_G=S_G(\RR)^0$, and $S_G$ is the maximal $\QQ$-split torus in
the center of $G$.)  If $X$ is not compact, Borel and Serre
\cite{refnBorelSerre} construct a compact real analytic manifold with
corners $\Xbar$ whose interior is $X$.  The faces of $\Xbar$ are indexed by
$\Pl$, the $\G$-conjugacy classes of parabolic $\QQ$-subgroups $P$ of $G$.
Each such face $Y_P$ is fibered over an arithmetic quotient
$X_P=\G_{L_P}\back D_P$ of the symmetric space $D_P$ associated to the Levi
quotient $L_P$ of $P$.  The fibers are compact nilmanifolds (an arithmetic
quotient of the real points of the unipotent radical of $P$) and Zucker
\cite{refnZuckerWarped} collapses these fibers (for all $Y_P$) to obtain a
new compactification $\Xhat$ of $X$, the \emph{reductive Borel-Serre
compactification}.

The reductive Borel-Serre compactification is important because it is
canonical, its singularities are reasonably easy to understand, and the
more singular Satake compactifications $\Xstar$ (such as the
Baily-Borel-Satake compactification in the case that $D$ is Hermitian
symmetric) may be realized as quotients of it
\cite{refnZuckerSatakeCompactifications}.  It would thus be desirable to
transfer cohomological calculations from $\Xstar$ to $\Xhat$.  In this
paper we develop a combinatorial tool, the theory of $\L$-modules and their
micro-support, to study the cohomology of constructible complexes of
sheaves on $\Xhat$.  As one application, we settle a 40 year old conjecture
on the intersection cohomology made independently by Rapoport
\cite{refnRapoportLetterBorel} and by Goresky and MacPherson
\cite{refnGoreskyMacPhersonWeighted}.  In fact we prove a generalization:

\begin{RapoportConjecturethm}[Rapoport/Goresky-MacPherson Conjecture]
Let $\Xstar$ be a Satake compactification of $X$ for which all real
boundary components $D_{R,h}$ are equal-rank, and let $\pi\colon \Xhat \to
\Xstar$ be the projection from the reductive Borel-Serre compactification.
Let $E$ be an algebraic representation of $G$ \textup(that is, a finite-dimensional representation given by a regular morphism $G\to\GL(E)$\textup) and let $\EE$ be the
corresponding locally constant sheaf on $X$.  Let $p$ be a middle
perversity.  Then there is a natural quasi-isomorphism of intersection
cohomology sheaves $\pi_*\IpC(\Xhat;\EE) \cong \IpC(\Xstar;\EE)$.
\end{RapoportConjecturethm}

Recall that $D$ is said to be \emph{equal-rank} if $\CCrank \lsp0 G =
\rank K$; here $\lsp 0G$ is the intersection of the kernels of the squares
of all characters of $G$ defined over $\QQ$.  The original conjecture was
for $D$ Hermitian symmetric and $\Xstar$ the Baily-Borel-Satake
compactification.  Previously Zucker had observed that the conjecture held
for $G=\Sympl(4)$ with $\EE=\CC$.  Furthermore, Goresky and MacPherson
\cite{refnGoreskyMacPhersonWeighted} announced a proof for $G=\Sympl(4)$,
$G=\Sympl(6)$, and (for $\EE=\CC$) $G=\Sympl(8)$.  A proof in the case
where $\QQrank G=1$ was given by Saper and Stern
\cite[Appendix]{refnRapoport}.

The theory of $\L$-modules applies to all $D$ and to sheaves other than
intersection cohomology, so for future applications we work in a more
general context and do not assume $D$ to be Hermitian or equal-rank except
as needed.

The paper is divided into four parts: I. $\L$-modules and Micro-support;
II. A Global Vanishing Theorem for $\L$-modules; III. Micro-support
Calculations; and IV. Satake Compactifications and Functoriality of
Micro-support.  After Part I, Parts II-IV may be read in any order provided
one is willing to follow up on a few references; the last section of the
paper, however, uses everything that precedes it.

We now briefly outline
the main results of the paper, after which we summarize the notation we
will use; see also the semi-expository article \cite{refnSaperIHP}.

\subsection{$\L$-modules on the reductive Borel-Serre
  compactification}
(\S\S\ref{sectGeodesicActionBundles}--\ref{sectWeightCohomologyLmodule})\ \
The construction of the reductive Borel-Serre compactification $\Xhat=
\coprod_P X_P$ is briefly indicated in
\S\S\ref{sectGeodesicActionBundles}--\ref{sectReductiveBorelSerre}.  Let
$\i_P$ and $\ihat_P$ denote respectively the inclusion of a stratum $X_P$
and its closure $\Xhat_P$ into $\Xhat$.  Let $\Derived_\X(\Xhat)$ be the
derived category whose objects are complexes of sheaves $\Sheaf$ on
$\Xhat$ whose cohomology sheaves $\H(\i_P^!\Sheaf)$ supported on each
stratum $X_P$ are locally constant, say induced from graded representations
$\mathscr E_P$ of $\G_{L_P}$.  (The focus on $\H(\i_P^!\Sheaf)$ instead of $\H(\i_P^*\Sheaf)$ is crucial in the theory of $\L$-modules.)  We would like to consider such objects
equipped with the extra structure that $\mathscr E_P$ is induced by
restriction from a graded representation $E_P$ of $L_P$.

In fact it is simpler to start with the family of graded representations
$E_P$ of $L_P$ and construct from it an object of $\Derived_\X(\Xhat)$.
However we need something to ``glue'' these representations together.
Consider two strata $X_P$ and $X_Q$ such that $X_P\subset \Xhat_Q$, the
closure of $X_Q$; this means $P\subset Q$ (after perhaps replacing $Q$ by
a $\G$-conjugate).  Let $N_P\supset N_Q$ be the corresponding unipotent
radicals.  The intersection of the topological link of the stratum $X_P$ with the stratum $X_Q$ is homotopically equivalent to an arithmetic quotient of $N_P^Q(\RR)
= N_P(\RR)/N_Q(\RR)$.  A theorem of Nomizu and van Est states that the
cohomology of this nilmanifold (with coefficients in $\EE_Q$, the locally
constant sheaf associated to $E_Q$) is isomorphic to the Lie algebra
cohomology $H(\n_P^Q; E_Q)$; furthermore this is an isomorphism of
$\G_{L_P}$-modules.

Thus in \S\ref{sectLsheaves} we are led to consider combinatorial data
$\M=(E_\cdot,f_{\cdot\cdot})$ consisting of a family of graded
representations $E_P$ of $L_P$, together with ``dual attaching'' morphisms
$f_{PQ}\colon H(\n_P^Q;E_Q)\to E_P[1]$ for all $P\subset Q$.
If $\M$ is required to satisfy a certain differential type condition
\eqref{eqnLsheafCondition} it is called an \emph{$\L$-module}.  From
an $\L$-module $\M$ one can construct in a natural way an object
$\Sheaf_{\Xhat}(\M)$ of $\Derived_\X(\Xhat)$ such that
$\H(\i_P^!\Sheaf_{\Xhat}(\M))= \EE_P$.

In fact we may define a category of $\L$-modules $\R(\L_W)$ for any locally
closed union of strata $W$ of $\Xhat$.  If $k\colon W\hookrightarrow Z$ is an
inclusion of such spaces, we have variants for $\L$-modules of the usual
functors $k_*$, $k_!$, $k^*$, and $k^!$, as well as truncation functors by
degree and by weight.  In \S\ref{sectRealizationLModules} we see there
exists a realization functor $\Sheaf_W\colon \R(\L_W) \to \Derived_\X(W)$ and
via $\Sheaf_W$ the functors $k_*$, etc.\ on $\L$-modules correspond to the
usual derived functors on $\Derived_\X(W)$.  The point is that these
derived functors preserve the extra structure given by an $\L$-module.

There are a number of different incarnations of $\Sheaf_W(\M)$ as a complex
of sheaves.  In the one we find simplest, a global section of
$\Sheaf_W(\M)$ corresponds to a family $\o=(\o_P)$, where $\o_P$ is a
special differential form \cite{refnGoreskyHarderMacPherson} on
$X_P\subset W$ with coefficients in $\EE_P$.  The differential on
$\Sheaf_W(\M)$ consists of the exterior derivative on each $\o_P$ together
with interaction terms based on the $f_{PQ}$.

If $\M$ is an $\L$-module on $\Xhat$, we define its cohomology
$H(\Xhat;\M)$ to be the hypercohomology of the associated complex of
sheaves $\Sheaf_{\Xhat}(\M)$.  The complexes of sheaves underlying many
cohomology theories on $\Xhat$ can be lifted to $\L$-modules.  For example,
let $E$ be an algebraic representation of $G$.  Then there exist bounded $\L$-modules
$\i_{G*}E$, $\IpC(E)$, and $\WnC(E)$ which correspond respectively to the
ordinary cohomology $H(X;\EE)=H(\Gamma;E)$, the intersection cohomology
$I_pH(\Xhat;\EE)$, and the weighted cohomology $W^\eta H(\Xhat;\EE)$ (see
\cite{refnGoreskyHarderMacPherson}).  The latter two $\L$-modules are
constructed in \S\S\ref{sectIntersectionCohomologyLmodule} and
\ref{sectWeightCohomologyLmodule}.

\subsection{Micro-support}
(\S\S\ref{sectMicroSupport}--\ref{sectAlternateMicroSupport})\ \ 
To an
$\L$-module $\M$ we associate its \emph{micro-support} $\mS(\M)$.
The micro-support is a subset of the set of all irreducible representations
of $L_P$ for all $P\in \Pl$.  As we will see below, the micro-support of
$\M$ yields a range of degrees outside of which $H(\Xhat;\M)$ must be zero.
(Micro-support here is vaguely analogous to the micro-support of a sheaf
introduced by Kashiwara and Schapira \cite{refnKashiwaraSchapira}, hence
the name.  There is no formal correspondence however.)

In order to define $\mS(\M)$, let $V$ be an irreducible $L_P$-module and
let $\xi_V$ denote the character by which the maximal $\QQ$-split torus
$S_P$ in the center of $L_P$ acts on $V$.  Also let $V|_{M_P}$ denote the
restriction of $V$ to $M_P=\lsp0 L_P$, the natural complement to $S_P$.
The closed strata $\Xhat_Q$ containing $X_P$ correspond to subsets $\D_P^Q$
of $\D_P$, the simple roots of $S_P$ acting on the Lie algebra of the
unipotent radical of $P$.  Let $\Xhat_Q\supset X_P$ be a closed stratum
such that
\begin{equation*}
\{\,\al\mid \langle\xi_V+\hsr,\al\spcheck\rangle < 0\,\} \subset \D_P^Q\subset  \{\,\al\mid
\langle\xi_V+\hsr,\al\spcheck \rangle\le  0\,\}
\end{equation*}
where as usual $\hsr$ denotes one half the sum of the positive roots.  Then
$V$ belongs to $\mS(\M)$ if and only if
\begin{gather}
(V|_{M_P})^* \cong \overline{V|_{M_P}}, \text{ and}
\label{eqnMicroSupportConjugateSelfContragredience} \\
H(\i_P^*\ihat_Q^!\M)_V\neq 0 \label{eqnMicroSupportVanishing}
\end{gather}
for some such $Q$.  Condition  \eqref{eqnMicroSupportConjugateSelfContragredience} here is reasonable since without it the $L^2$-cohomology $H(X_P;\VV)$ would vanish (Theorem ~\ref{ssectRaghunathanVanishing}).

\subsection{Vanishing theorem for $\L$-modules}
\label{ssectIntroVanishingTheoremLmodules}
(\S\S\ref{sectVanishingTheoremLModule}--\ref{sectProofGlobalVanishing})\ \
For an $L_P$-module $V$ in $\mS(\M)$, denote the least and greatest degrees for which
\eqref{eqnMicroSupportVanishing} can be nonzero by $c(V;\M)\le d(V;\M)$.
Set
\begin{align*}
c(\M) &= \inf_{V\in\mS(\M)} \tfrac12(\dim D_P - \dim D_P(V)) + c(V;\M)\ ,
\text{ and} \\
d(\M) &= \sup_{V\in\mS(\M)} \tfrac12(\dim D_P + \dim D_P(V)) + d(V;\M)\ ,
\end{align*}
where $D_P(V)$ is the symmetric space associated to the reductive
$\RR$-subgroup of $L_P$ whose roots (with respect to a fundamental Cartan
subalgebra) are orthogonal to the highest weight of $V$.  (This space
depends on a $\theta$-stable ordering; we choose one for which $\dim D_P(V)$ is
maximized.)  Then we have the

\begin{vanthm}
Let $\M$ be a bounded $\L$-module on $\Xhat$. Then  $H^i(\Xhat;\M)=0$ for $i
\notin [c(\M), d(\M)]$.
\end{vanthm}

The proof of the vanishing theorem proceeds by a combinatorial
generalization of the analytic arguments from \cite{refnSaperSternTwo}
together with a spectral sequence argument; this reduces the problem to
Raghunathan's vanishing theorem for ordinary $L^2$-cohomology, Theorem
~\ref{ssectRaghunathanVanishing}.  In \S\ref{sectLtwo} we recall ordinary
$L^2$-cohomology and a version in which restriction to boundary faces is
well defined.  Now we need to represent $H(\Xhat;\M)$ by a variant of
$L^2$-cohomology for our forms $(\o_P)\in \Sheaf_{\Xhat}(\M)$.  However
even though each $\Xhat_P$ is compact this is not entirely straightforward:
the $L^2$-norm on $\o_P$ with respect to the natural locally symmetric
metric on $\Xhat_P$ is not appropriate because, being a complete metric on
a noncompact space, it would impose $L^2$-growth conditions on $\o_P$ which
we do not want.  The solution in \S\ref{sectLTwoCohomologyLModules} is to
replace $\Xbar$ by a naturally diffeomorphic compact domain $\Xbar_t\subsetneq X$ as constructed in \cite{refnSaperTilings} and use the
induced metric; we then work on the corresponding $\Xhat_t$ instead of
$\Xhat$.

\subsection{Calculations of
micro-support}(\S\S\ref{sectMicroSupportWeightedCohomology}--\ref{sectPurityProofPartOne})\
\ For simplicity, in this introduction we only give the results for the
\emph{essential micro-support} $\emS(\M)$; this is the subset of
$\mS(\M)$ of $V$ for which \eqref{eqnMicroSupportVanishing} is nonvanishing
for all possible $Q$ (in a compatible fashion).  The essential
micro-support determines $\mS(\M)$ and may be used in calculating $c(\M)$
and $d(\M)$.

For $\M = \WnC(E)$, weighted cohomology, one can explicitly calculate $H(\i_P^*\ihat_Q^!\M)$
and finds the

\begin{SSWHthm}
Let $E$ be an irreducible algebraic $G$-module and let $\eta$ be a middle
weight profile.  The micro-support is nonempty if and only if $(E|_{\lsp0
G})^*\cong \overline{E|_{\lsp0 G}}$\textup; if this condition is satisfied
then $\emS(\WnC(E))=\{E\}$ and $c(E;\WnC(E))=d(E;\WnC(E))=0$.
\end{SSWHthm}

For $\M=\IpC(E)$, intersection cohomology, the situation is far more delicate.  We do not
have an explicit closed formula for $H(\i_P^*\ihat_Q^!\M)$ (or even for the
local cohomology $H(\i_P^*\M)$) and we are forced to use an inductive
argument.  However condition
\eqref{eqnMicroSupportConjugateSelfContragredience} is usually not
preserved under passage to a larger parabolic $\QQ$-subgroup.  Nonetheless
we have the following result:

\begin{SSIHcor}
Assume the irreducible components of the $\QQ$-root system of $G$ are of
type $A_n$, $B_n$, $C_n$, $BC_n$, or $G_2$.  Let $E$ be an irreducible
algebraic $G$-module and let $p$ be a middle perversity.  If $(E|_{\lsp0
G})^*\cong \overline{E|_{\lsp0 G}}$, then $\emS(\IpC(E))=\{E\}$ with
$c(E;\IpC(E))= d(E;\IpC(E))= 0$.
\end{SSIHcor}
\begin{rems*}
Unlike the situation for weighted cohomology, the micro-support can be
nonempty even if $(E|_{\lsp0 G})^*\not\cong \overline{E|_{\lsp0 G}}$; see
Theorem ~\ref{ssectIHMicroPurity} for details.

The hypothesis on the $\QQ$-root system ensures that the the corresponding
$\QQ$-split group for each almost $\QQ$-simple factor of $G$ has a
quasi-minuscule representation whose weights are linearly ordered (see the
proof of Lemma ~\ref{ssectQInductionLemma}).  We conjecture the results on
$\mS(\IpC(E))$ continue to hold without this hypothesis.  The hypothesis is
automatically satisfied if $D$ is Hermitian or has a Satake
compactification with all real boundary components equal-rank.
\end{rems*}

\subsection{Satake compactifications}
(\S\ref{sectSatakeCompactifications})\ \ We recall the theory of Satake
compactifications.  Briefly, given a faithful irreducible finite
dimensional representation $\sigma\colon G\to \GL(V)$, Satake
\cite{refnSatakeCompact}, \cite{refnSatakeQuotientCompact} constructed a
compactification $\lsb\RR \Dstar_\sigma= \coprod_R D_{R,h}$ of $D$ where
the strata $D_{R,h}$ are called the \emph{real boundary components}
and are indexed by their normalizers $R$.  The union  $\Dstar_\sigma$ of the so-called
\emph{rational boundary components}, under certain conditions on $\sigma$, yields a compactification
$\Xstar_\sigma = \G\back \Dstar_\sigma$ of $X=\G\back D$.  The strata $F_R
= \G_{L_{R,h}}\back D_{R,h}$ of $\Xstar_\sigma$ are indexed by
$R\in \Pl^\star_\sigma$, a subset of $\Pl$.  In the case that $D$ is Hermitian
symmetric, $\Xstar_\sigma$ for a certain $\sigma$ is the Baily-Borel-Satake
compactification.  It has a projective algebraic structure
\cite{refnBailyBorel} and its strata are indexed by $\G$-conjugacy classes
of maximal parabolic $\QQ$-subgroups.

In general Zucker showed \cite{refnZuckerSatakeCompactifications} that there is a
natural quotient map $\pi\colon \Xhat\to \Xstar_\sigma$.  For each stratum
$X_P$ of $\Xhat$ there exists a stratum $F_{P^\dag}$ of $\Xstar_\sigma$
such that $\pi|_{X_P}\colon X_P \to F_{P^\dag}$ is a flat bundle with
typical fiber $X_{P,\l}$.  In fact, there is a connected normal
$\QQ$-subgroup $L_{P,\l}\subset L_P$ so that $X_{P,\l}$ is an arithmetic
quotient of the symmetric space associated with $D_{P,\l}$, while
$F_{P^\dag}$ is an arithmetic quotient of $L_{P,h}= L_P/L_{P,\l}$.  The
inverse image $\pi^{-1}(F_R)$ of a stratum $F_R$ of $\Xstar_\sigma$ is the
union of those $X_P$ for which $P^\dag=R$; we denote it $X_R(L_{R,\l})$
since it is a partial compactification of $X_R$ in the vertical
(``$L_{R,\l}$'') directions.  The restriction $\pi|_{X_R(L_{R,\l})}\colon
X_R(L_{R,\l})\to F_R$ is again a flat bundle with typical fiber
$\Xhat_{R,\l}$.  Let $\ihat_{R,\l}\colon \Xhat_{R,\l}\hookrightarrow
X_R(L_{R,\l})$ denote the inclusion of a typical fiber.  In
\S\ref{ssectPullbackToFiber} we define functors $\ihat_{R,\l}^*$ and
$\ihat_{R,\l}^!$ on $\L$-modules; these commute with their analogues on the
derived category.

\subsection{Functoriality of micro-support}
(\S\S\ref{sectFunctorialityInverseImages}--\ref{sectRestrictMicroSupportToFiber})\
\ We need to understand how micro-support is affected when we apply
functorial operations such as $k^*$ or $k^!$ to $\M$.  For brevity we only
deal with $k^*$ in this introduction.  A few general results are presented
in \S\ref{sectFunctorialityInverseImages}.  For the case where $k$ is a closed
embedding, however, we only have a result for the \emph{weak
micro-support} $\mS_w(\M)$ which is defined similarly to $\mS(\M)$ but
omits condition \eqref{eqnMicroSupportConjugateSelfContragredience}.
The result briefly is that if $V$ is an irreducible $L_P$-module in $\mS_w(k^*\M)$ then there
exists an irreducible $L_{\tilde P}$-module $\tilde V\in \mS_w(\M)$ such
that $V$ is an irreducible constituent of $H^\l(\n_P^{\tilde P};\tilde V)$
with
\begin{equation}
(\xi_{V}+\hsr)|_{\sa_P^{\tilde P}} \in \lsp+\sa_P^{\tilde P*}
\label{eqnPositiveWeylHypotheses}
\end{equation}
where $\lsp+\sa_P^{\tilde P*}$ is the real convex cone generated by
$\D_P^{\tilde P}$.

The problem of $\mS_w(k^*\M)$ versus $\mS(k^*\M)$ is eliminated if $k$ is
the embedding of $X_R(L_{R,\l})$ for a Satake compactification in which all
rational boundary components are equal-rank; see
\S\ref{sectFunctorialityMicroSupportEqualRank}.  The point is that if $L_P
= \widetilde{L_{P,h}} L_{P,\l}$ is an almost direct product with the
symmetric space of $\widetilde{L_{P,h}}$ being equal-rank, then
$(V|_{\widetilde{L_{P,h}} })^* \cong \overline{V|_{\widetilde{L_{P,h}}}}$
is automatic.

We find that $c(k^*\M)$ and $d(k^*\M)$ may be estimated in terms of $c(\M)$
and $d(\M)$, together with the degrees $\l$ from above.  These are the
lengths $\l(w)$ of certain Weyl group elements $w$ and
\eqref{eqnPositiveWeylHypotheses} gives information on the geometry of the
corresponding Weyl chamber.  A basic lemma in \S\ref{sectBasicLemma}
translates this geometry into an estimate on $\l(w)$.  In the case that
$\#\D_P^{\tilde P} = 1$, the estimate is $\l(w)\le \frac12\dim\n_P^{\tilde
P}-\n_P(V)$.  See \S\ref{ssectRealSubmodules} for the definition of
$\n_P(V)$; in \S\ref{sectEqualRankBasicLemma} we estimate $\n_P(V)$ in more
geometric terms---this requires the stronger condition that all \emph{real}
boundary components are equal-rank.

The end result, after restricting to a typical fiber $\Xhat_{R,\l}\subset
X_R(L_{R,\l})$, is

\begin{SSFibercor}
Let $\Xstar_\sigma$ be a Satake compactification of $X$ and assume that all
the real boundary components $D_{R,h}$ of the associated Satake
compactification $\lsb\RR\Dstar_\sigma$ are equal-rank.  Let $\M$ be an
$\L$-module with $\emS(\M)=\{E\}$ for some irreducible algebraic $G$-module
$E$ and $c(E;\M)=d(E;\M)=0$.  Then for every proper stratum $F_R$ of
$\Xstar_\sigma$,
\begin{equation*}
d(\ihat_{R,\l}^*\M) \le \frac12\codim F_R - \dim \sa_{R}^{G}
\qquad\text{and}\qquad c(\ihat_{R,\l}^!\M) > \frac12\codim F_R + \dim
F_R\ .
\end{equation*}
\end{SSFibercor}

\subsection{The Conjecture of Rapoport and Goresky-MacPherson} (\S\ref{sectRapoportConjecture})\ \ 
Let $\M$ be an $\L$-module.  We wish to identify $\pi_*\Sheaf_{\Xhat}(\M)$
as an element of the derived category on $\Xstar_\sigma$.  If $x\in
F_R\subset \Xstar_\sigma$, then $\pi^{-1}(x) = \Xhat_{R,\l}$, a fiber of
$\pi|_{X_R(L_{R,\l})}$.  Thus the local cohomology at $x$ of
$\pi_*\Sheaf_{\Xhat}(\M)$ is simply $H(\ihat_{R,\l}^*\M)$, and the local
cohomology supported at $x$ is $H(\ihat_{R,\l}^!\M)$.  For $\M = \IpC(E)$,
we consequently see from Theorem ~\ref{ssectGlobalVanishing}, Corollary
~\ref{ssectIHMicroPurityCorollary}, and Corollary
~\ref{ssectRestrictMicroSupportToFiberCorollary} that
$\pi_*\Sheaf_{\Xhat}(\M)$ satisfies the local characterization of
intersection cohomology.  This proves the conjecture of Rapoport and
Goresky-MacPherson.  By replacing Corollary
~\ref{ssectIHMicroPurityCorollary} with Theorem
~\ref{ssectWeightCohomologyMicroSupport} we obtain the following
generalization of the main theorem of \cite{refnGoreskyHarderMacPherson}:

\begin{GoreskyHarderMacPhersonthm}
Under the assumptions of Theorem ~\textup{\ref{ssectRapoportConjecture}},
let $\eta$ be a middle weight profile and let $p$ be a middle perversity.
Then there is a natural quasi-isomorphism $\pi_*\WnC(\Xhat;\EE) \cong
\IpC(\Xstar;\EE)$.
\end{GoreskyHarderMacPhersonthm}

\subsection{Other applications of $\L$-modules}

The micro-support of $\i_{G*}E$ is easy to calculate and is treated
in \cite{refnSaperIHP}.  As an application of Vanishing Theorem ~\ref{ssectGlobalVanishing}, it is proved that if $E$ has regular highest weight, then $H^i(X;\EE)=0$ for $i<\frac12(\dim X - (\CCrank \lsp0 G - \rank K))$.
This answers a question posed by Tilouine in the Hermitian case (personal communication, Spring 2000).  Li and Schwermer \cite{refnLiSchwermer} have also proved this by entirely different methods.

In \cite{refnSaperLtwoLmoduleI} the theory of $\L$-modules is generalized
to allow $E_P$ to be a locally regular $L_P$-module (that is, any finite-dimensional subspace of $E_P$ is contained in an algebraic subrepresentation);
this allowed us to realize
$L^2$-cohomology as an $\L$-module and calculate its micro-support.

Finally note that our main results, Theorems \ref{ssectRapoportConjecture} and \ref{ssectGoreskyHarderMacPhersonTheorem}, imply the intersection cohomology of $\Xhat$ is isomorphic to the weighted cohomology of $\Xhat$ when $X$ is a Hermitian locally symmetric space.  In \cite{refnSaperMixed} we prove this holds for all locally symmetric spaces provided that types $D$, $E$, and $F$ do not occur in the $\QQ$-root system of $G$ as assumed in Corollary~\ref{ssectIHMicroPurityCorollary} above.

\subsection{Work to be done}
There are many other developments and applications of the theory of $\L$-modules
that would be interesting to pursue.  We mention for example the following
topics:
\begin{itemize}
\item $\QQ$-structures on an $\L$-module and its cohomology (see
  \S\ref{ssectOtherRealizationsLsheaves});

\item $\L$-modules on Satake
  compactifications 
(see \S\ref{ssectLModuleGeneralizations});

\item the removal of the
hypothesis on the $\QQ$-root system in Corollary ~\ref{ssectIHMicroPurityCorollary};

\item the action of Hecke correspondences on $\L$-modules and their cohomology;

\item
the localization of the homotopy category of $\L$-modules (see
\S\ref{ssectHomotopyLModules}) with respect to the
null system of $\L$-modules with empty micro-support;

\item a micro-support
spectral sequence for the cohomology of an $\L$-module;

\item the relationship of the structure of $H(\Gamma;E)$ coming from $\mS(\i_{G*}\EE)$ and the
structure coming from Eisenstein series as initiated by
Harder (see in particular the work of Franke \cite{refnFranke}).
\end{itemize}

A first step toward these last two topics is in \cite{refnSaperMixed}.

\subsection{Acknowledgments}
I would like to thank Mark Goresky for introducing me to this problem and
for discussing it with me and Mark McConnell many times during the academic
year 1988-89 while I visited Harvard University.  I would also like to
thank Michael Rapoport for reintroducing me to the problem in a letter to
Stern and myself in summer 1991, and encouraging us to publish our proof of
the $\QQ$-rank 1 case.  I am grateful for stimulating discussions with many
people which aided me in this research, in particular, A. Borel, R. Bryant,
R. Hain, M. Goresky, G. Harder, J.-P. Labesse, G. Lawler, M. Rapoport,
J. Rohlfs, J. Schwermer, J. Tilouine, and N. Wallach.  Finally I would like
to thank the referee for many valuable suggestions and corrections, and in
particular for suggesting the simplified proof of Proposition
~\ref{ssectWeightCohomologyLocalCohomologyWithSupports} that appears here.

\bigskip

\subsection{Notation}
\label{ssectNotation}
We now summarize the notation of the paper.

\subsubsection{}
\label{sssectBasics}
The cardinality of a set $T$ will be denoted
$\#T$\label{defCardinality}.  Let $\subset$ and $\subsetneq$ denote subsets and proper subsets respectively.
For a group $Z$ with $z\in Z$ and $S\subset Z$, let $\lsp z S = S ^{z^{-1}}$ denote
the conjugate $ z S z^{-1}$\label{defConjugation}.  If $Z$ is a topological space, the
topological closure of a subset $S\subset Z$ will be denoted by
$\cl{S}$\label{defClosure}. The notation $\id_Z$ indicates the identity map
of $Z$\label{defIdentity}.  If $\Cat$ is a category we will write $E\in \Cat$ to mean $E$ is an object of $\Cat$.  

\subsubsection{}
\label{sssectNotationCategory}
If $\Cat$ is an additive category, let
$\Complex(\Cat)$\label{defComplexes}  denote
the category of complexes of objects of $\Cat$ and let
$\Graded(\Cat)$\label{defGraded}
denote full subcategory of complexes with zero differential, that is,
graded objects of $\Cat$.  The full subcategory of $\Complex(\Cat)$ or
$\Graded(\Cat)$ consisting of bounded complexes will be denoted  $\Complex^b(\Cat)$ or
$\Graded^b(\Cat)$.

Let $\Cat$ be an additive category.  For $k\in \ZZ$, the \emph{shift}
$C[k]$\label{defShift} of an object $C\in\Complex(\Cat)$ is defined by
$C[k]^i=C^{i+k}$ and $d_{C[k]}=(-1)^k d_C$.  Given a functor $F\colon \Cat
\to \Complex(\Cat')$, where $\Cat$ and $\Cat'$ are additive categories, we implicitly extend it to a functor $F\colon \Graded(\Cat) \to
\Complex(\Cat')$ by setting $F(C)=\oplus_i F(C^i)[-i]$.  A common example is \eqref{eqnKostantDegeneration}.  Under this
convention we have the equality $F(C)[k]=F(C[k])$.  We further extend $F$ to a functor $\Complex(\Cat) \to \Complex(\Cat')$
by taking the associated total complex.

\subsubsection{}
\label{sssectSheaves}
For a topological space $W$, let $\Sh(W)$\label{defSheaves} denote the category of sheaves of
vector spaces over $\CC$ on $W$.  We write
$\Complex(W)$\label{defComplexesSheaves} for
$\Complex(\Sh(W))$.  Suppose $W$ is equipped with a stratification $\X$.
We will call a complex of sheaves $\Sheaf$ \emph{constructible}%
\footnote{In the terminology of \cite[V,
    3.3(ii)]{refnBorelIntersectionCohomology} this is $\X$-cohomologically
  locally constant, or $\X$-clc.}
if for all strata $X_P\in \X$, the restriction $H(\Sheaf)|_{X_P}$ of the
local cohomology sheaves is locally constant.  Let $\Complex_\X(W)$\label{defConstructibleComplexes} denote
the category of complexes of sheaves on $W$ which are constructible.
A morphism in $\Complex_\X(W)$ is a \emph{quasi-isomorphism} if it induces an isomorphism on all local cohomology sheaves; we denote this $\Sheaf \qiarrow \mathcal T$.
Let
$\Derived_\X(W)$\label{defConstructibleDerived} be the associated derived
category, obtained by formally
inverting quasi-isomorphisms in the corresponding homotopy category.  We
denote the usual derived functors on $\Derived_\X(W)$ without the prefix
``$R$'', for example, $i_*$ instead of $Ri_*$.

\subsubsection{}
\label{sssectDeRham}
All references to manifolds, smooth differential forms, fiber
bundles, locally constant sheaves, etc. should be taken in the sense of
$V$-manifolds \cite{refnSatakeVManifold} (also called a orbifolds
\cite{refnThurston}).

For $Y$ a manifold with (possibly empty) boundary or corners and $\EE$ a
locally constant sheaf on $Y$, let $\EE\mapsto
\A(Y;\EE)$\label{defDeRhamSheaf} denote the de
Rham resolution functor consisting of sheaves of smooth differential forms
on $Y$ with values in $\EE$.  The differential is exterior differentiation,
denoted $d=d_Y=d_{Y,\EE}$.  Note that if $\EE$ is graded we follow the sign
convention of \S\ref{sssectNotationCategory}.  The complex of global
sections is denoted $A(Y;\EE)$\label{defDeRhamComplex} and the subcomplex of forms with compact
support is denoted $A_c(Y;\EE)$\label{defCompactlySupportedForms}.

\subsubsection{}
\label{sssectCharacters}
For any algebraic group defined over $\RR$, we denote the Lie algebra of
its real points by the corresponding lowercase Gothic letter, for example, $\mathfrak
g = \operatorname{Lie} G(\RR)$\label{defLieAlgebra}.  The complexification
is denoted by a 
subscript $\CC$, for example $\mathfrak
g_{\CC}$\label{defComplexification}.  Now let $Z$ be an
algebraic group defined over $\QQ$.  The group of characters of $Z$ defined
over $\QQ$ is denoted $X(Z)$\label{defCharacters}.  The one-dimensional representation space
associated to a character $\chi$ is denoted $\CC_\chi$.  Let $\lsp 0Z =
\bigcap_{\chi\in X(Z)} \ker \chi^2$\label{defBorelPart} be the intersection of the kernels of
the squares of characters of $Z$.  The group $Z$ is an \emph{almost
direct product} of $\QQ$-subgroups $H$ and $K$ if the product morphism
$H\times K\to Z$ is surjective with finite kernel.

\subsubsection{}
\label{sssectRepresentations}
For a connected reductive $\QQ$-group $L$, let $\R(L)$\label{defRegRep} denote the category
of algebraic representations of $L$.  We write $\Complex(L)$\label{defComplexRegRep} and $\Graded(L)$\label{defGradedRegRep}
instead of $\Complex(\R(L))$ and $\Graded(\R(L))$.  If $M\subset L$ is a
reductive $\QQ$-subgroup and $E\in \R(L)$, let $\Res_M^L E$ or $E|_M$\label{defRestriction} denote the algebraic
representation of $M$ obtained from $E$ by restriction.

Let $\IrrRep(L)$\label{defIrrRep} denote the isomorphism classes of irreducible algebraic
$L$-modules.  For any $E\in \R(L)$ we have an \emph{isotypical
decomposition}
\begin{equation*}
E=\bigoplus_{V\in\IrrRep(L)} E_V
\end{equation*}
where $E_V=\Hom_{L}(V,E)\otimes V$\label{defIsotypicComponent}.  For any irreducible representation
$V\in\IrrRep(L)$, we let $\xi_V$\label{defCentralCharacter} denote the character by which $S_L$, the maximal $\QQ$-split torus in the center of $L$, acts
on $V$.  Since $S_L$ is $\QQ$-split, $\xi_V$ is defined over $\QQ$.
Although $\xi_V$ may not extend to a character of $L$, there exists a
natural number $k$ such that $\xi_V^k$ extends to a character of $L$
defined over $\QQ$.

\subsubsection{}
\label{sssectContragredient}
Let $Z$ be an algebraic group and let $\sigma\colon Z \to \GL(V)$ be an algebraic representation.  Recall the contragredient representation
$\sigma^*$ on $V^*$\label{defContragredient} is defined by $(\sigma^*(g)\psi)(v)=
\psi(\sigma(g^{-1})v)$ for $\psi\in V^*$.  Define the complex vector space $\overline V$ to
have vectors $\overline v$ indexed by $v\in V$, with operations $\overline
v + \overline w = \overline{v+w}$ and $\lambda \overline v =
\overline{\overline \lambda v}$.  If $Z$ is defined over $\RR$ there is a
complex conjugate representation $\overline\sigma$ on $\overline V$ defined
by $\overline\sigma(g)\overline v = \overline{\sigma(\overline g)v}$.

\subsubsection{}
\label{sssectUsualSuspects}
In this paper $G$ will be a connected reductive algebraic group defined
over $\QQ$ and $\G\subset G(\QQ)$ will be an arithmetic subgroup.  (We
briefly depart from this convention on $G$ in
\S\S\ref{ssectSatakeCompactificationsSymmetricSpaces},
\ref{ssectRealBoundaryComponents}.)  We let $D$ denote the associated
symmetric space $G(\RR)/KA_G$, where $K$ is a maximal compact subgroup of
$G(\RR)$ and $A_G$ is the connected component of the group of real points
of the maximal $\QQ$-split torus $S_G$ in the center of $G$.  Let
$X=\G\back D$ be the arithmetic quotient.  If $E$ is an algebraic
representation of $G$, the induced locally constant sheaf on $X$ will be
denoted by the corresponding blackboard bold letter, thus $\EE=D\times_{\G}
E$.

Since we do not assume $\G$ to be neat, $X$ (and the various other
arithmetic quotients that we will consider) need not be manifolds, but
rather $V$-manifolds.  Alternatively the reader may assume that $\G$ is neat; any
arithmetic group contains a neat subgroup with finite index.

\subsubsection{}
\label{sssectSplitComponents}
Let $P$ be a parabolic $\QQ$-subgroup of $G$, with unipotent radical $N_P$
and Levi quotient $L_P=P/N_P$. Set $M_P= \lsp 0 L_P$.  Let $S_P$ denote the
maximal $\QQ$-split torus in the center of $L_P$ and set $A_P =
S_P(\RR)^0$.  In terms of real points
we have $L_P(\RR) = M_P(\RR) \times A_P$.  Restriction of characters yields an injective morphism
$X(L_P)\to X(S_P)$ whose image has finite index.

If $P\subset Q$ are parabolic $\QQ$-subgroups, $P/N_Q$ is a parabolic
$\QQ$-subgroup of $L_Q$ with unipotent radical $N_P^Q=N_P/N_Q$.  One can
naturally identify $S_Q$ as a subgroup of $S_P$ and $X(L_Q)$ as a subgroup
of $X(L_P)$.  Namely consider the diagram
\begin{equation*}
\xymatrix{
{L_P} & {P/N_Q} \ar@{>>}[l]_-{j} \ar@{^{ (}->}[r]^-{i} & {L_Q\ ;}
}
\end{equation*}
since $P/N_Q$ is parabolic in $L_Q$, $i$ identifies its center with that of
$L_Q$, while $j$ identifies it with a subgroup of the center of $L_P$.
Furthermore
\begin{equation*}
\xymatrix{
{X(L_P)} \ar[r]^-{\sim} & {X(P/N_Q)}  & {X(L_Q)\ .}
\ar@{_{ (}->}[l]_-{i^*}
}
\end{equation*}

Set $S_P^Q = \bigl( \bigcap_{\chi\in X(L_Q)} \ker \chi\cap S_P\bigr)^0$.
Then $S_P$ is an almost direct product
\begin{equation*}
S_P = S_Q \cdot S_P^Q  \ .
\end{equation*}
If we write $A_P^Q= S_P^Q(\RR)^0$, then $A_P$ is an honest direct product
\begin{equation}
A_P = A_Q \times A_P^Q  \ .  \label{eqnUsualDecomp}
\end{equation}

\subsubsection{}
\label{sssectSplitComponentRoots}
Let $P$ be a parabolic $\QQ$-subgroup of $G$.  Let $\sa_P$ be the Lie
algebra of $A_P$ and let $\langle\ ,\ \rangle$ denote the natural pairing
$\sa_P^*\times\sa_P\to \RR$.  We may identify $X(S_P)\otimes_\ZZ \RR\cong
\sa_P^*$.  Consequently we will view elements of $X(S_P)$ both as
characters on $S_P$ and as linear functionals on $\sa_P$.

The group $L_P$ acts via conjugation on $N_P$ (and hence on $\n_{P\CC}$)
via any lift $\widetilde{L_P}\subset P$.  Although this action depends on
the lift, the characters $\al\in X(S_P)$ by which $S_P$ acts are
independent of the lift.  We denote these characters by $\Phi(\n_P,S_P)$ or
$\Phi(\n_P,\sa_P)$.  As sometimes done we call them ``roots'' even though in general
for $P$ nonminimal they are not the positive roots of a root system.  They
do, however, have the property that every such root is a unique nonnegative
integral linear combination of so-called simple roots; we let $\D_P$ denote
these simple roots.

In the case that $P=P_0$ is a minimal parabolic $\QQ$-subgroup, $\D_{P_0}$
is a basis of the $\QQ$-root system $\lsb\QQ\Phi$ of $G$.  In this case we
omit the subscript $P_0$ and simply write $\D$, $S$, $\sa$, etc.

\subsubsection{}
\label{sssectType}
If $P\subset Q$ are parabolic $\QQ$-subgroups, the \emph{type of $Q$
with respect to $P$} is the subset $\D_P^Q$ of $\D_P$ consisting of roots
restricting to 1 on $S_Q$.  Parabolic $\QQ$-subgroups $Q$ containing $P$
are in one-to-one correspondence with subsets of $\D_P$ via
$Q\leftrightarrow \D_P^Q$.

Let
\begin{equation*}
\sa_P = \sa_Q \oplus \sa_P^Q
\end{equation*}
be induced from \eqref{eqnUsualDecomp}.  By \S\ref{sssectSplitComponents}
we have a natural inclusion $\sa_Q^* \hookrightarrow \sa_P^*$.  Then
$\sa_Q$ is the subspace of $\sa_P$ annihilated by $\D_P^Q$, whereas
$\sa_P^Q$ is the subspace annihilated by $\D_Q$ together with a basis of
$X(S_G)$.  We also have a dual decomposition
\begin{equation*}
\sa_P^* = \sa_Q^* \oplus \sa_P^{Q*}\ .
\end{equation*}
As usual if $P=P_0$ is minimal we we omit the subscript $P_0$.
Given $\al\in \sa_P^*$ we write $\al = \al_Q + \al^Q$ for its corresponding
decomposition and likewise for elements of $\sa_P$.

\subsubsection{}
\label{sssectCoroots}
For each $\g\in \D$, let $\g\spcheck\in \sa$ be the corresponding coroot.
(We shall occasionally identify $\sa$ with $\sa^*$ via a Weyl group
invariant inner product $(\ ,\ )$, in which case $\g\spcheck =
\frac{2}{(\g,\g)}\g$.)  Let $\Dhat=\{\b_\al\}$ be the basis of $\sa^*$ dual to the coroots.

To define coroots where $P$ is an arbitrary parabolic $\QQ$-subgroup, we proceed as in \cite{refnArthurTraceFormula}.  Namely for $\al \in \D_P$, let $\g\in \D\setminus \D^P$ be the unique root restricting to $\al$ on $\sa_P$ and let $\al\spcheck= (\g\spcheck)_P$ be the projection of $\g\spcheck$ to $\sa_P$.
If $P\subset Q$, then $\sa_P^Q$ has a basis of coroots, $\{\al\spcheck\}_{\al\in\D_P^Q}$.  Let $\Dhat_P^Q=\{\b_\al^Q\}_{\al\in\D_P^Q}$ be the dual basis of $\sa_P^{Q*}$.
On the other hand, $\D_P^Q$ is a basis for $\sa_P^{Q*}$.  For $\b\in
\Dhat_P^Q$, define $\b\spcheck\in \sa_P^Q$ by $\langle
\al,\b\spcheck\rangle = \langle \b,\al\spcheck\rangle$ for all $\al \in
\D_P^Q$.  Then $\{\b\spcheck\}_{\b\in \Dhat_P^Q}$ is the basis of $\sa_P^Q$
dual to $\D_P^Q$.

\subsubsection{}
\label{sssectRootAndDominantCones}
Let $\lsp +\sa_P^*$ denote the real convex cone spanned by $\D_P$ and let
$\sa_P^{*+}$ denote the real convex cone spanned by $\Dhat_P$.  Then
$\sa_P^{*+}$ is the dominant cone
\begin{equation*}
\{\,\mu\in\sa_P^*\mid \langle\mu,\al\spcheck\rangle\ge 0
\text{ for all $\al\in\D_P$}\,\}
\end{equation*}
and $\sa_P^{*+} \subset \lsp +\sa_P^*$.  Similarly define $\lsp
+\sa_P^{Q*}$ and $\sa_P^{Q*+}$.

\subsubsection{}
\label{sssectRealParabolics}
Occasionally we will need to consider
parabolic $\RR$-subgroups.  The discussion in
\S\S\ref{sssectSplitComponents}--\ref{sssectRootAndDominantCones} applies
with $\QQ$ replaced by $\RR$.  We will add a left subscript $\RR$ to the
notation in this case, for example, $\lsb\RR S_P$, $\lsb\RR A_P$,
$\lsb\RR\D_P$, etc.

\subsubsection{}
\label{sssectComplementaryParabolic}
Let $P\subset Q$ be parabolic $\QQ$-subgroups.  Let $(P,Q)$ denote the
parabolic $\QQ$-subgroup ``complementary'' to $Q$ with respect to $P$ with
type $\D_P\setminus \D_P^Q$.  When $(P,Q)$ occurs as a subscript or a
superscript, we omit the parentheses, for example, $\D_P^{P,Q}$.  Thus we
have a disjoint decomposition
\begin{equation*}
\D_P = \D_P^{P,Q} \amalg \D_P^Q\ .
\end{equation*}
Restriction to $\sa_Q$ sends $\D_P^Q$ to zero and defines a bijection
$\D_P^{P,Q}\tildearrow \D_Q$.

Note that this decomposition of $\D_P$ induces a direct sum decomposition
\begin{equation*}
\sa_P = \sa_Q \oplus \sa_{P,Q}^G\ ,
\end{equation*}
however in general $\sa_{P,Q}^G \neq \sa_P^Q$.

\subsubsection{}
\label{sssectParabolicPoset}
Let
\begin{equation*}
\Pl=  \{\, \text{$\G$-conjugacy classes of parabolic $\QQ$-subgroups of
$G$} \,\}\ .
\end{equation*}
By abuse of notation we do not distinguish a parabolic $\QQ$-subgroup from its
$\G$-conjugacy class.  Consequently, for example, an algebraic representation $E$ of $L_P$ represents a family $(\sigma_\gamma\colon L_{\lsp\gamma P} \to \GL(E))_{\gamma\in\G})$ satisfying $\sigma_{\gamma'\gamma}(\gamma' r \gamma'^{-1}) = \sigma_\gamma(r)$ for all $\gamma',\gamma\in\G$ and $r\in L_{\lsp\gamma P}$.

Define a partial order on $\Pl$ by
\begin{equation*}
P \le Q \quad \Longleftrightarrow \quad P \subset Q^\gamma \text{ for
some $\gamma\in\G$;}
\end{equation*}
write $P<Q$ when $P\subsetneq Q^\gamma$ for any $\gamma\in\G$. When $P\le Q$ we will often replace $Q$ by a $\G$-conjugate
so as to arrange $P\subset Q$ as needed.  For $P$ and $Q$ in $\Pl$, let
$P\vee Q$ denote the least upper bound and let $P\cap Q$ denote
the greatest lower bound when it exists.

Given $P\le R$, let $[P,R]$ denote the interval $ \{\, Q\in \Pl \mid P\le Q
\le R\,\}$.  This is a distributive lattice.  It is also complemented---for every $Q\in [P,R]$ there
exists $Q'\in [P,R]$ such that $Q\cap Q' = P$ and $Q\vee Q' = R$, namely
$Q' = (P,Q)\cap R$.  We represent this situation by the tilted rectangle
\begin{equation}
\vcenter{\xymatrix @ur @M=1pt @R=1.5pc @C=1pc {
{Q'} \ar@{-}[r] \ar@{-}[d] & {R} \ar@{-}[d] \\
{P\rlap{\qquad .}} \ar@{-}[r] & {Q}
}} \label{eqnParallelogram}
\end{equation}

\subsubsection{}
\label{sssParabolicRoots}
Let $P$ be a parabolic $\QQ$-subgroup and let $\h$ be a Cartan subalgebra
of $\levi_P$.  Write $\h=\hb_P \oplus\sa_P$ where $\hb_P$ is a Cartan
subalgebra of $\m_P$.  (In \S\ref{sectEqualRankBasicLemma} we shall also
need $\h=\lsb\RR\hb_P\oplus\lsb\RR\sa_P$ for $P$ a parabolic $\RR$-subgroup,
where $\lsb\RR\hb_P$ is a Cartan subalgebra of $\lsb\RR\m_P$.)  The
corresponding roots will be denoted $\Phi^P=\Phi(\levi_{P\CC},\h_\CC)$ and a
positive system will be denoted $\Phi^{P+}=\Phi^+(\levi_{P\CC},\h_\CC)$.  (In the
case $P=G$ we shall often simply write $\Phi$ and $\Phi^+$.)

Via a lift, $\h$ acts on $\n_P$ via the adjoint action; the corresponding
roots are independent of the lift and will be denoted
$\Phi(\n_{P\CC},\h_{\CC})$.  Likewise via a lift, $\h$ may be considered
as a Cartan subalgebra of $\levi_Q$ for any parabolic $\QQ$-subgroup
$Q\supset P$; the decomposition $\h=\hb_Q \oplus \sa_Q$ and the roots
$\Phi(\levi_{Q\CC},\h_\CC)$ are independent of the lift.  We can always
choose a positive system $\Phi^+(\levi_{Q\CC},\h_\CC)$ containing
$\Phi^+(\levi_{P\CC},\h_\CC)$ and $\Phi(\n_{P\CC},\h_{\CC})$.

\subsubsection{}
\label{sssectHalfSumOfPositiveRoots}
Let $\hsr = \hsr^P + \hsr_P$ be one half the sum of roots in
$\Phi^+(\mathfrak g_\CC,\h_\CC)$, where $\hsr^P$ (resp.\ $\hsr_P$) is
one half the sum of the roots in $\Phi^+(\levi_{P\CC},\h_\CC)$
(resp.\ $\Phi(\n_{P\CC},\h_\CC)$).  We may identify $2\hsr_P$ with the
$\QQ$-rational character of $L_P$ by which $L_P$ (via any lift $\widetilde
L_P$) acts on $\bigwedge^{\dim \n_P} \n_P$.  Thus $\hsr_P$ is supported on
$\sa_P$ and $\hsr=\hsr^P+\hsr_P$ corresponds to the decomposition induced
by restriction to the summands in $\h=\hb_P \oplus \sa_P$.

\subsubsection{}
\label{sssectNilpotentCohomology}
Let $P$ be a parabolic $\QQ$-subgroup and let $E$ be an algebraic $G$-module.
The \emph{nilpotent Lie algebra cohomology functor} $E\mapsto H(\n_P;E)$
from $\R(G)$ to $\Graded(L_P)$ will play a crucial role in what follows.
As a vector space, $H(\n_P;E)$ is defined as the cohomology of the complex
$C(\n_P;E)=\Hom(\bigwedge\n_P,E)$ equipped with the usual differential (see
for example \cite[I, \S1.1]{refnBorelWallach}).
%
%
The coadjoint representation induces an action of $P$ on $C(\n_P;E)$ which
commutes with the differential, and hence there is an action of $P$ on
$H(\n_P;E)$.  The unipotent radical $N_P$ acts trivially \cite[I,
\S1.1(5)]{refnBorelWallach} so this descends to a representation of $L_P$.
Finally, a morphism of $G$-modules $f\colon E \to E'$ induces a morphism of
$P$-modules $C(\n_P;E) \to C(\n_P;E')$ which commutes with the
differentials and hence induces an $L_P$-module morphism $H(\n_P;E) \to
H(\n_P;E')$ which we denote $H(\n_P;f)$.

\subsubsection{}
\label{sssectWeyl}
Let $W$ denote the Weyl group of $\Phi$ and $\Phi^+$ a positive system.  For $w\in W$ let $\l(w)$ denote
the length of $w\in W$ with respect to the reflections in simple roots and
set $\Phi_w = \{\,\g\in \Phi^+\mid w^{-1}\g < 0\,\}$.  It is  well known
that $\l(w) = \#\Phi_w$.

If $P$ is parabolic $\QQ$-subgroup let $\Phi^+$ be a postive system containing  $\Phi^+ (\levi_{P\CC},\h_\CC)$ and $\Phi(\n_{P\CC},\h_\CC)$.
Let $W^P\subset W$ be the subgroup
corresponding to the  Weyl
group of $\Phi(\levi_{P\CC},\h_\CC)$ and let $W_P\subset W$ be the set of
minimal length coset representatives of $W^P\backslash W$.%
\footnote{Our use of subscript versus superscript here is opposite from the
usual convention, but this way is more compatible with the subscripts and
superscripts on $A$ and $\Delta$.}
Then $W^P$ is a normal subgroup of the automorphisms of $\Phi(\levi_{P\CC},\h_\CC)$ and we have $w\in W_P$ if and only if $\Phi_w \subset \Phi(\n_{P\CC},\h_\CC)$.
There is a factorization
\begin{equation*}
W=W^PW_P.
\end{equation*}
If the positive system is changed to $\Phitilde^+ = \widetilde w\Phi^+$ for $\widetilde w\in W^P$ then
\begin{equation}
  \label{eqnOrderingDependenceForW_P}
\widetilde W_P = \widetilde w \,W_P \widetilde w^{-1} \qquad \text{and}\qquad \Phitilde _w = \widetilde w \,\Phi_{\widetilde w^{-1}w \widetilde w}\ .
\end{equation}

The corresponding groups and coset representatives for the $\QQ$-root
system will be used in \S\ref{sectInductionLemmaProof}, and will be denoted
there $\lsb\QQ W$, etc.

\subsubsection{Kostant's Theorem}
\label{ssectKostantsTheorem}
Let $E$ be an algebraic $G$-module.  The complement $W_P$ plays an important
role in Kostant's theorem \cite{refnKostant}, which states that $H(\n_P;E)$
has a decomposition as an $L_P$-module indexed by $W_P$,
\begin{equation*}
H(\n_P;E) = \bigoplus_{w\in W_P} H(\n_P;E)_w \ ,
\end{equation*}
for which  $H(\n_P;E)_w$ is concentrated in degree $\l(w)$.

More precisely, fix a Cartan subalgebra $\h$ of $\levi_P$ (and hence of $\mathfrak g$) and a positive system $\Phi^+(\mathfrak g_{\CC},\h_{\CC})$ containing $\Phi^+(\levi_{P\CC},\h_{\CC})$ and $\Phi(\n_{\CC},\h_{\CC})$.  If $E$ is
irreducible with highest weight $\lambda$, then $H(\n_P;E)_w$ is the
irreducible $L_P$-module having highest weight $w(\lambda+\hsr)-\hsr$.

\subsubsection{}
\label{ssectKostantDegeneration}
Let $P\subset Q$ be parabolic $\QQ$-subgroups.  We have
\begin{equation*}
W = W^QW_Q =(W^PW_P^Q)W_Q
\end{equation*}
and one may check that
\begin{equation*}
W_P=W_P^QW_Q.
\end{equation*}
For $w\in W_P$ write $w=w^Qw_Q$ accordingly.  We have the important
equality
\begin{equation*}
\l(w)= \l(w^Q) + \l(w_Q)
\end{equation*}
and in fact
\begin{equation}
\Phi_w=\Phi_{w^Q} \textstyle \coprod w^Q \Phi_{w_Q}   \ .
\end{equation}
It immediately follows from Kostant's theorem that
\begin{equation*}
H(\n_P;E)_w \cong H(\n_P^Q;H(\n_Q;E)_{w_Q})_{w^Q}
\end{equation*}
and consequently \cite[(11.8)]{refnGoreskyHarderMacPherson},
\cite[\S4.10]{refnSchwermerGeneric}
\begin{equation}
H(\n_P;E) \cong H(\n_P^Q;H(\n_Q;E))\ .
\label{eqnKostantDegeneration}
\end{equation}


\bigskip
\bigskip

\section*{Index of Notation}
\setlength{\columnseprule}{.4pt}
\begin{multicols}{2}
  \setlength{\NotationWidth}{.8\linewidth}
  \begin{tabbing}
    \hspace{\NotationWidth}\=\hspace{35pt}\kill
\tabfill{$\#T$, $\subset$, $\subsetneq$, $\lsp z S$, $S ^z$, $\cl{S}$, $\id_Z$, $E\in\Cat$  } \>\`
\S\ref{sssectBasics}\\
\tabfill{$\Complex(\Cat)$, $\Complex^b(\Cat)$, $\Graded(\Cat)$, $\Graded^b(\Cat)$, $C[k]$  } \>\` \S\ref{sssectNotationCategory}\\
\tabfill{$\Sh(W)$, $\Complex(W)$, $\Complex_\X(W)$, $\Derived_\X(W)$,  $\Sheaf \qiarrow \mathcal T$} \>\` \S\ref{sssectSheaves}\\
\tabfill{$\A(Y;\EE)$, $A(Y;\EE)$, $A_c(Y;\EE)$} \>\` \S\ref{sssectDeRham}\\
\tabfill{$\mathfrak g$, $\operatorname{Lie} G(\RR)$, $\mathfrak g_{\CC}$, $X(Z)$, $\lsp 0Z$} \>\` \S\ref{sssectCharacters}\\
\tabfill{$\R(L)$, $\Complex(L)$, $\Graded(L)$, $E|_M$, $\Res_M^L E$, $\IrrRep(L)$, $E_V$, $\xi_V$} \>\` \S\ref{sssectRepresentations}\\
\tabfill{$V^*$, $\overline V$} \>\` \S\ref{sssectContragredient}\\
\tabfill{$G$, $\G$, $D$, $K$, $A_G$, $S_G$, $X$, $\EE$} \>\` \S\ref{sssectUsualSuspects}\\
\tabfill{$P$, $N_P$, $L_P$, $M_P$, $S_P$, $A_P$, $N_P^Q$, $S_P^Q$, $A_P^Q$} \>\` \S\ref{sssectSplitComponents}\\
\tabfill{$\sa_P$, $\langle\ ,\ \rangle$, $\widetilde{L_P}$, $\Phi(\n_P,S_P)$,
$\Phi(\n_P,\sa_P)$, $\D_P$, $\lsb\QQ\Phi$, $\D$, $S$, $\sa$} \>\` \S\ref{sssectSplitComponentRoots}\\
\tabfill{$\D_P^Q$, $\sa_P^Q$, $\D^Q$, $\sa^Q$, $\al_Q$, $\al^Q$} \>\` \S\ref{sssectType}\\
\tabfill{$\g\spcheck$, $\al\spcheck$, $\Dhat_P^Q$, $\b\spcheck$} \>\` \S\ref{sssectCoroots}\\
\tabfill{$\lsp +\sa_P^*$, $\sa_P^{*+}$, $\lsp +\sa_P^{Q*}$, $\sa_P^{Q*+}$} \>\` \S\ref{sssectRootAndDominantCones}\\
\tabfill{$\lsb\RR S_P$, $\lsb\RR A_P$, $\lsb\RR\D_P$, etc.} \>\` \S\ref{sssectRealParabolics}\\
\tabfill{$(P,Q)$,  $\D_P^{P,Q}$, $\sa_{P,Q}^G$ } \>\` \S\ref{sssectComplementaryParabolic}\\
\tabfill{  $\Pl$, $P \le Q$, $P\vee Q$,  $P\cap Q$,  $[P,R]$,
  $\xymatrix @ur @M=1pt @R=1.5pc @C=1pc { {Q'} \ar@{-}[r] \ar@{-}[d] & {R}
  \ar@{-}[d] \\ {P} \ar@{-}[r] & {Q} }$
 } \>\` \S\ref{sssectParabolicPoset}\\
\tabfill{$\h$, $\levi_P$, $\hb_P$, $\sa_P$, $\m_P$,
$\Phi(\levi_{P\CC},\h_\CC)$, $\Phi^+(\levi_{P\CC},\h_\CC)$, $\Phi^P$, $\Phi^{P+}$, $\Phi$, $\Phi^+$,  $\n_P$,  $\Phi(\n_{P\CC},\h_{\CC})$ } \>\` \S\ref{sssParabolicRoots}\\
\tabfill{$\hsr$, $\hsr^P$, $\hsr_P$} \>\` \S\ref{sssectHalfSumOfPositiveRoots}\\
\tabfill{$H(\n_P;E)$, $C(\n_P;E)$, $H(\n_P;f)$} \>\` \S\ref{sssectNilpotentCohomology}\\
\tabfill{$W$, $\l(w)$, $\Phi_w$, $W^P$, $W_P$,  $\lsb\QQ W$} \>\` \S\ref{sssectWeyl}\\
\tabfill{$H(\n_P;E)_w$ } \>\` \S\ref{ssectKostantsTheorem}\\
\tabfill{$W_P^Q$, $w^Q$, $w_Q$} \>\` \S\ref{ssectKostantDegeneration}
  \end{tabbing}
  \setlength{\NotationWidth}{.85\linewidth}
  \begin{tabbing}
    \hspace{\NotationWidth}\=\hspace{23.3pt}\kill
\tabfill{$L_{P,x}$, $z \geo y$, $z_x$ } \>\` \S\ref{ssectGeodesicAction}\\
\tabfill{$\mathscr A_P^G$, $e_P$, $\G_P$, $r_P$, $Y_P$} \>\` \S\ref{ssectGeodesicRetraction}\\
\tabfill{$\mathscr N_P(\RR)$, $D_P$, $\mu$, $n_x$, $\mathscr N_P(\RR)'$, $\G_{N_P}$, $\G_{L_P}$, $\mu'$,
$N_P(\RR)'$, $n_x'$, $p$, $X_P$} \>\` \S\ref{ssectNilmanifoldFibrations}\\
\tabfill{$\Abar_P^G$, $o_Q$, $D(P)$, $\Dbar$, $\overline{e_P}$, $\Xbar$, $Y_P\le
Y_Q$,  $\Ybar_Q$ } \>\` \S\ref{ssectBorelSerre}\\
\tabfill{$\bar{\mathscr A}_P^G$ } \>\` \S\ref{ssectLocalStructureXbar}\\
\tabfill{$\Xhat$, $p$, $\Xhat_Q$, $\i_P$, $\ihat_P$, $U_P$, $\j_P$, $\jhat_P$} \>\` \S\ref{ssectReductiveBorelSerre}\\
\tabfill{$\mathscr Z_P$ } \>\` \S\ref{ssectLocalStructureXhat}\\
\tabfill{$W$, $\Pl(W)$} \>\` \S\ref{ssectNotationLsheaves}\\
\tabfill{$\L_W$ } \>\` \S\ref{ssectLstructure}\\
\tabfill{$\L$, $\M$, $E_P$, $f_{PQ}$, $\i_P^!\M$, $\phi$, $\phi_{PQ}$, $\R(\L)$, $\R^b(\L)$} \>\` \S\ref{ssectLsheaves}\\
\tabfill{ $\i_P$, $\ihat_P$, $\j_P$, $\jhat_P$, $k^*\M$, $k_*\M$, $k^!\M$, $k_!\M$, $\i_P^*\M$} \>\` \S\ref{ssectFunctorsOnLsheaves}\\
\tabfill{$G_{\l}$, $\tilde G_h$, $G_{h}$,  $\G_{G_{\l}}$, $\G_{G_{h}}$,  $P_\l$,
$P_h$, $\tilde P_h$ $D_{h}$, $D_{\l}$,  $X_{h}$, $X_\l$, $X(G_\l)$, $\pi$,
$\ihat_{G,\l}$,  $L_{P,\l}$,  $\ihat_{G,\l}^*\M$, $\Res_{L_{P,\l}}^{L_P}$,
$\ihat_{G,\l}^!\M$, $\ihat_{R,\l}^*\M$,
$\ihat_{R,\l}^!\M$, $L_{R,\l}$, $X_R(L_{R,\l})$, $X_{R,h}$, $\Xhat_{R,\l}$,
$\widetilde{L_{R,h}}$} \>\` \S\ref{ssectPullbackToFiber}\\
\tabfill{$i_P^*\ihat_Q^!\M$ } \>\` \S\ref{ssectLocalCohomologyWithSupports}\\
\tabfill{$M(\phi)$,  $\M[1]$, $\al(\phi)$, $\b(\phi)$} \>\` \S\ref{ssectMappingCone}\\
\tabfill{$\K(\L)$, $\Derived(\L)$} \>\` \S\ref{ssectHomotopyLModules}\\
\tabfill{$(L_\sigma)_\sigma$, $\M_{\sigma\tau}$} \>\` \S\ref{ssectLModuleGeneralizations}\\
\tabfill{$\A_P$, $k_{PQ}$, $\A_W$, $\Sheaf_W$, $H(W;\M)$} \>\` \S\ref{ssectRealizationLModules}\\
\tabfill{$\A_W(\M)$, $\A_W(\phi)$} \>\` \S\ref{ssectRealizationLsheaves}\\
\tabfill{$q$, $\A\inv(Y_P;\EE)$, $A\inv(Y_P;\EE)$, $\Asp(\Xbar;\EE)$,
$\Asp(\Xhat;\EE)$, $A\sp(\Xbar;\EE)$, $A\sp(\Xhat;\EE)$} \>\` \S\ref{ssectSpecialDifferentialForms}\\
\tabfill{$A\inv(\mathscr N_P(\RR)';\EE)$, $\AA\inv(\mathscr N_P(\RR)';\EE)$,
$\o|_{Y_P}$, $\o|_{X_P}$, $\o|_{\Xhat_P}$} \>\` \S\ref{ssectRestrictionBoundaryStrata}\\
\tabfill{$\HH(\n_P;E)$ } \>\` \S\ref{ssectNilpotentCohomologySheaf}\\
\tabfill{$\coh_P$} \>\` \S\ref{ssectNilmanifoldCohomologyBundle}\\
\tabfill{$k_P$} \>\` \S\ref{ssectRestrictSpecialForms}\\
\tabfill{$\coh_{PQ}$, $k_{PQ}$} \>\` \S\ref{ssectProofOfTheoremRealization}\\
\tabfill{$\Acomb(\Xhat;\EE)$} \>\` \S\ref{ssectOtherRealizationsLsheaves}\\
\tabfill{$\t^{\leqslant n}C$, $\t^{> n}C$} \>\` \S\ref{ssectTruncation}\\
\tabfill{$\t_Q^{\leqslant n}\M$} \>\` \S\ref{ssectLocalDegreeTruncationLmodule}\\
\tabfill{$\t^{\leqslant n}\M$} \>\` \S\ref{ssectDegreeTruncationLmodules}\\
\tabfill{$p$, $p(k)$, $\IpC(E)$,  $\IpC(\Xhat;\EE)$ } \>\` \S\ref{ssectIntersectionCohomology}\\
\tabfill{$|\D_P|$, $|\D_P^Q|^\circ$, $c(|\D_P|)$, $p_w$, $I_{p_w}H(U)$, $I_{p_w}H(c(|\D_P|))$, $\ast$} \>\` \S\ref{ssectLocalIntersectionCohomology}\\
\tabfill{$X(S_P^G)$, $R_\al$, $\phi_P$, $\zeta\ge \eta$, $E_\chi$, $\t^{\geqslant
  \eta}E$, $\t^{\ngeqslant \eta}E$, $\eta$} \>\` \S\ref{ssectWeightTruncation}\\
\tabfill{$\t_Q^{\geqslant \eta}$, $\t^{\geqslant \eta}\M$ } \>\` \S\ref{ssectWeightTruncationLmodule}\\
\tabfill{$\WnC(E)$, $\WnC(\Xhat;\EE)$} \>\` \S\ref{ssectWeightedCohomologyLmodule}\\
\tabfill{$Q_V^W$, $Q^{\prime W}_V$, $Q_V$, $Q'_V$,} \>\` \S\ref{ssectMicroSupportNotation}\\
\tabfill{$\mS_w(\M)$, $\mS(\M)$, $\emS(\M)$} \>\` \S\ref{ssectMicroSupport}\\
\tabfill{$\t_P$, $\h^*_\RR$, $s^{P}$, $s^c$, $\t'_P$  } \>\` \S\ref{ssectDefineTau}\\
\tabfill{$\hb_{P,1}$ $\hb_{P,-1}$ } \>\` \S\ref{ssectConjugateSelfContragredient}\\
\tabfill{$\theta$} \>\` \S\ref{ssectFundamentalCSA}\\
\tabfill{ $V \preccurlyeq \tilde V$, $V \preccurlyeq_0 \tilde V$, $[\tilde V:V]$} \>\` \S\ref{ssectPartialOrdering}\\
\tabfill{$c_i(V;\M)$, $d_i(V;\M)$ } \>\` \S\ref{ssectAlternateMicroSupport}\\
\tabfill{$L_P(\mu)$, $\levi_P(\mu)$} \>\` \S\ref{ssectWeightCentralizer}\\
\tabfill{$D_P(\mu)$, $D_P(V)$} \>\` \S\ref{ssectVcentralizer}\\
\tabfill{$c(V;\M)$, $d(V;\M)$, $c(\M)$, $d(\M)$} \>\` \S\ref{ssectcMdMDefinition}\\
\tabfill{$c_i(\M)$, $d_i(\M)$, $\essc(V;\M)$, $\essd(V;\M)$} \>\` \S\ref{ssectEqualityDegreeRanges}\\
\tabfill{$\lVert\o\rVert$, $\lVert\o\rVert_Y$, $L_2(Y;\EE)$, $d$, $d_Y$,
$d_{Y,\EE}$, $\dom d$, $\dbar$, $H_{(2)}(Y;\EE)$, $\dbar^*$,
$\lVert\o\rVert_g$, $L_2(Y;\EE,g)$, $H_{(2)}(Y;\EE,g)$} \>\` \S\ref{ssectLtwoCohomology}\\
\tabfill{$\H_\epsilon(Y;\EE)$} \>\` \S\ref{ssectSpectralVanishingLtwoCohomology}\\
\tabfill{$\lvvv\o\rvvv$, $d_b$. $d_{Y,b}$, $d_{Y,b,\EE}$, $\Dom d_{b}$, $\dbar_{b}$,
$H_{(2),b}(Y;\EE)$} \>\` \S\ref{ssectLtwoBoundaryValues}\\
\tabfill{$\o|_F$, $\v$} \>\` \S\ref{ssectBoundaryRestriction}\\
\tabfill{$N_F$, $p_F$, $t_F$} \>\` \S\ref{ssectBoundaryAveraging}\\
\tabfill{$L_{2,b}(Y;\EE)$, $\dbar_b^*$, $\CoDifferential_F$, $\v_F^{\Ftilde}$} \>\` \S\ref{ssectLtwoBoundarySpace}\\
\tabfill{$\H_{\epsilon,b}(Y;\EE)$} \>\` \S\ref{ssectSpectralVanishingBoundaryValues}\\
\tabfill{$K$, $A_G$, $\theta$, $\k$, $\p$, $\sa_G$, $B$, $T_xD$, $(E,\sigma)$,
  $\EE$, $Z_P$, $\EEtilde_P$} \>\` \S\ref{ssectMetrics}\\
\tabfill{$\Xbar_t$, $s_t$, $\Dbar_t$, $e_{P,t}$} \>\` \S\ref{ssectTilings}\\
\tabfill{$A_P^G(s_P)$, $\langle A_P^G\rangle_t$, $\Omega_P$, $W_P$, $\Wbar_{P,t}$} \>\` \S\ref{ssectTilingInCylindricalSet}\\
\tabfill{$Y_{P,t}$, $X_{P,t}$, $\Xhat_t$, $\Xhat_{P,t}$, $\Xbar_{P,t}$, $\EE_{P,t}$} \>\` \S\ref{ssectMetrizeRetract}\\
\tabfill{$\A_{\Xhat_t}(\M)$, $d_P$, $\dd_{PQ}$
  } \>\` \S\ref{ssectLtwoMethods}\\
\tabfill{$\|\o\|_t$, $L_{2}(U;\M,t^{-\hsr})$,
  $d_{U;\M,t^{-\hsr}}$, $\dbar$, $\dbar_{\M}$, $\dbar_{U;\M,t^{-\hsr}}$,
  $H_{(2)}(U;\M,t^{-\hsr})$, $Y_{P,t}^Q$} \>\` \S\ref{ssectDDBoundIndependentT}\\
\tabfill{$\H_{\epsilon}^i(U;\M,t^{-\hsr})$} \>\` \S\ref{ssectSpectralVanishingLModules}\\
\tabfill{$\What_{P,t}$} \>\` \S\ref{ssectSpectralMayerVietoris}\\
\tabfill{$E'_Q$, $f'_{QR}$} \>\` \S\ref{ssectLModulesAPTimesXP}\\
\tabfill{$O_P$, $\nil$, $\nil_*\M$, $\nil_*\M\otimes \CC_{\hsr_P}$,  $\nil_{Q*}$,
  $H_Q$, $N_j^Q$, $Z_j^Q$, $\n_j^Q$  } \>\` \S\ref{ssectSpectralAverging}\\
\tabfill{$(\nil_*\M)_V$} \>\` \S\ref{ssectLPIsotypicalDecomposition}\\
\tabfill{$\Psi_t$} \>\` \S\ref{ssecttIndependence}\\
\tabfill{$\dim D(E)$,  $\sigma_1$, $\sigma_2$, $C(\lsp0\mathfrak g,\k; E\otimes F)$,
  $\epsilon(X_k)$, $\D_0$, $C$, $C_\k$ $\hb_{G,\k}$, $\hb_{G,\p}$,
  $\Phi_\k^+$, $\Phi^+(\p_\CC,\hb_{G,\k\CC})$, $G(\lambda)$, $\Phi^+(\p(\lambda)_\CC,\hb_{G,\k\CC})$} \>\` \S\ref{sectVanishingTheoremLtwo}\\
\tabfill{$\dbar_{R,b}$, $\i_R^{R'*}$, $A_P^\al$, $M$  } \>\` \S\ref{sectProofGlobalVanishing}\\
\tabfill{$T_V^\eta$} \>\` \S\ref{ssectLocalWeightedCohomologyWithAndWithoutProperSupports}\\
\tabfill{$\mu$, $\v$, $\langle \sa_P^*\rangle_R$, $\langle\langle
  \sa_P^*\rangle\rangle_R$, $\chi_R$ } \>\` \S\ref{ssectWeightCohomologyMicroSupport}\\
\tabfill{$m$, $m(k)$, $n$, $n(k)$} \>\` \S\ref{sectPurity}\\
\tabfill{$I_{p_w}H_Z(U)$} \>\` \S\ref{ssectPurityReductions}\\
\tabfill{$\n_\al$, $\Phi(\n_{\al\CC},\h_\CC)$, $\l_\al(w)$} \>\` \S\ref{ssectSimpleBasicLemma}\\
\tabfill{$\lsb \QQ W$ } \>\` \S\ref{sectInductionLemmaProof}\\
\tabfill{$\lsb\QQ\Phi^P$,  $\lsb\QQ W^P$, $\sb\QQ\Phi_P$ } \>\` \S\ref{ssectInductionLemmaProofNotation}\\
\tabfill{$\lsb\QQ W_P$, $\lsb\QQ W(P)$, $C$, $C^P$, $w\cdot C_P$} \>\` \S\ref{ssectSplitComponentChambers}\\
\tabfill{$\pr_P$,  $\lsb\QQ\Phi_w$} \>\` \S\ref{ssectAllOrNone}\\
\tabfill{$\nu\le \mu$, $\Phi(V)$ } \>\` \S\ref{ssectQInductionLemma}\\
\tabfill{$\IH(U)$, $\IH_Z$} \>\` \S\ref{ssectNotationPurityProofOne}\\
\tabfill{$dp(k)$, $\ModTwo(k)$, $dm$, $dn$, $\ModTwo$} \>\` \S\ref{ssectPerversityShift}\\
\tabfill{$S$, $\lsb\RR S$, $\lsb k\D$, $\D$, $\lsb\RR \D$, $\sigma$, $U$,  $\mu$, $\lsb k\mu$,
  $\kap(\Th)$,  $\o(\Th)$, $P^\dag$, $\restr$, $\Th\sphat$, $\Utilde$} \>\` \S\ref{ssectSatakeCompactificationsNotation}\\
\tabfill{$\lsb\RR A_G$, $h_0$, $T^*$, $S(U)$, $\tilde \sigma$, $\lsb\RR \Dstar_\sigma$} \>\` \S\ref{ssectSatakeCompactificationsSymmetricSpaces}\\
\tabfill{$D_{P,h}$, $L'_{P,\l}$, $L'_{P,h}$, $\widetilde{L'_{P,h}}$, $D_{R,h}\le D_{R',h}$} \>\` \S\ref{ssectRealBoundaryComponents}\\
\tabfill{$L_{R,\l}$, $L_{R,h}$, $\widetilde{L_{R,h}}$} \>\` \S\ref{ssectRationalBoundaryComponents}\\
\tabfill{$\Dstar_\sigma$, $\Xstar_\sigma$, $\Xstar$, $\Pl^*_\sigma$, $F_R$, $X_{R,h}$,
  $\G_{L_{R,h}}$, $F_R\le F_{R'}$ } \>\` \S\ref{ssectSatakeCompacticationsOfX}\\
\tabfill{$D_{P,\l}$, $K_{Q,\l}$, $q_P$, $q$, $\pi$} \>\` \S\ref{ssectSatakeAsQuotientOfRBS}\\
\tabfill{$X_S$, $Q_V^R$, $Q_V^{\prime R}$ } \>\` \S\ref{ssectFunctorialityInverseImagesNotation}\\
\tabfill{$V \preccurlyeq \tilde V$, $V \preccurlyeq_+ \tilde V$,  $V \preccurlyeq_- \tilde V$} \>\` \S\ref{ssectAdditionalPartialOrdering}\\
\tabfill{$k^*\mathscr S$, $k^!\mathscr S$, $k^{*/!}$, $N_V^{\tilde V}$ } \>\` \S\ref{ssectWeakMicroSupportClosedInverseImage}\\
\tabfill{$X_R(L_{R,\l})$, $\ihat_{R,\l}$, $V_h$, $V_\l$, $\Res_\l^*\mathscr
  S$, $\Res_\l^!\mathscr S$, $\Res_\l^* V$, $\Res_\l^! V$ } \>\` \S\ref{ssectWeakMicroSupportInclusionFiber}\\
\tabfill{$L_{P,\l}^{P',h}$} \>\` \S\ref{lemTwoBoundary}\\
\tabfill{$U(F_R)$} \>\` \S\ref{ssectMicroSupportInverseImages}\\
\tabfill{$\t'_P(F)$, $\n_P(\mu)$, $\n_\al(\mu)$, $\n_P(V)$} \>\` \S\ref{ssectRealSubmodules}\\
\tabfill{$D_{P,h}(V)$, $D_{P,\l}(V)$} \>\` \S\ref{ssectVcentralizerRationalFactors}\\
\tabfill{$\k_P$, $\p_P$, $\hb_{P,\k}$, $\hb_{P,\p}$, $\k_P(\mu)$, $\p_P(\mu)$,
  $\n_P(\r)$,  $\r_\g$  } \>\` \S\ref{ssectEqualRankLemma}\\
\tabfill{$D_{P_{i},\l}^{P_{i+1},h}$ } \>\` \S\ref{ssectRestrictMicroSupportToFiber}
\end{tabbing}
\end{multicols}

\newpage


\specialsection*{Part I. $\L$-modules and Micro-support}
In this part, after briefly recalling the reductive Borel-Serre compactification, we define the category of $\L$-modules.  We introduce the analogues of the Fary and Mayer-Vietoris spectral sequences for $\L$-modules which are crucial tools in our work.  We construct a realization functor from $\L$-modules to the derived category of constructible sheaves.  Finally we construct the principle examples for this paper, the intersection cohomology and the weighted cohomology $\L$-modules.

\section{Geodesic Action and Related Bundles}
\label{sectGeodesicActionBundles}
Recall that $D$ is the symmetric space $G(\RR)/KA_G$ associated to $G$ and
$X=\G\back D$ is an arithmetic quotient.  In this section we briefly recall
the geodesic retraction and nilmanifold fibration associated to a parabolic
$\QQ$-subgroup $P$ \cite[\S3]{refnBorelSerre},
\cite[\S3]{refnZuckerWarped}, \cite[ \S\S3,~ 4,
~7]{refnGoreskyHarderMacPherson}.  This section also serves to establish
our notation for the various spaces and bundles that will enter into the
construction of the reductive Borel-Serre compactification in
\S\ref{sectReductiveBorelSerre}.  One novelty of our presentation is that
we introduce a ``geodesic action'' for the entire Levi quotient $L_P(\RR)$
as opposed to merely $A_P$; this is convenient for describing the flat
connection on the nilmanifold fibration.

\subsection{Geodesic action}
\label{ssectGeodesicAction}
Given any basepoint $x\in D$, let $L_{P,x}$ denote the unique lift of $L_P$
to $P$ which is stable under the Cartan involution associated to $x$.
Since $P(\RR)$ acts transitively on $D$, we may express $y\in D$ as $px$
for $p\in P(\RR)$; since $P=N_P\rtimes L_{P,x}$ we may write $p = u r$ with
$u\in N_P(\RR)$ and $r\in L_{P,x}(\RR)$.  Define the \emph{geodesic
action} of $L_P(\RR)$ on $D$ by
\begin{equation}
z \geo y = u z_x r x, \label{eqnGeodesicAction}
\end{equation}
where $z_x\in L_{P,x}(\RR)$ denotes the lift of $z$.  Note that at the
basepoint $x$ the geodesic action of an element of $L_P$ agrees with the
ordinary action of its lift in $L_{P,x}$.  Clearly the restriction of this
action to $A_P$ agrees with the geodesic action as defined in
\cite{refnBorelSerre}.

\begin{lem*}
The geodesic action $(z,y)\mapsto z\geo y$ of $L_P(\RR)$ on $D$ is
well-defined \textup(that is, $z\geo y$ is independent of the choice of $p$
and $x$ above\textup).  It satisfies
\begin{equation}
q(z\geo y) = \lsp q z\geo qy \qquad (q\in P(\RR)), \label{eqnGeodesicAndUsual}
\end{equation}
where $z\mapsto \lsp q z$ denotes the action of $P(\RR)$ by conjugation on
$L_P(\RR)$.
\end{lem*}

\begin{proof}
Consider first the choice of $p$.  If $y =
p'x$ for $p'\in P(\RR)$, then $p'=pk$ for some $k\in K_x\cap P\subset
L_{P,x}$ and hence $p'=u r'$ with $r'=r k$.  Thus $u z_x r' x = u
z_x r x$ and so the right hand side of \eqref{eqnGeodesicAction} depends only
on $x$, $y$, and $z$.

Now consider another basepoint $x'$.  Write $x'=qx = vl x$ with $q\in P(\RR)$, $v\in
N_P(\RR)$, and $l\in L_{P,x}(\RR)$.  Then $L_{P,x'}=vL_{P,x}v^{-1}$.  We may
express $y= p'x'$ where $p'= pq^{-1}= u r l^{-1} v^{-1}= (uv^{-1})(vr
l^{-1} v^{-1}) = u' r'$.  Thus we compute $u' z_{x'}r'x'
=(uv^{-1})(vz_xv^{-1})(vr l^{-1} v^{-1})qx = u z_x r x$, which shows that
\eqref{eqnGeodesicAction} is independent of $x$.

To prove \eqref{eqnGeodesicAndUsual}, calculate $q(z\geo y) = quz_xrx = \lsp q u
\lsp q(z_x) \lsp qr qx$.  However $\lsp q(z_x)$ is the lift of $\lsp qz$ to
$L_{P,x'}$, where $x'=qx$.  Thus $q(z\geo y) =
\lsp q u (\lsp qz)_{x'} \lsp qr x' = \lsp qz \geo (\lsp q u  \lsp qr x') =
\lsp qz \geo qy$.
\end{proof}

\begin{cor*}
The geodesic action of $L_P(\RR)$ and the usual action of $N_P(\RR)$
commute.  The geodesic action of $A_P\subset L_P(\RR)$ commutes with the
usual action of $P(\RR)$.
\end{cor*}

\begin{proof}
By the lemma we need to show that $ \lsp q z=z$ for $q\in N_P(\RR)$ and
$z\in L_P(\RR)$ (and likewise for $q\in P(\RR)$ and $z\in A_P$).  Write
$z=z_x N_P(\RR)$ for some basepoint $x$.  Then $\lsp q z = qz_x
q^{-1}N_P(\RR) =z_x q^{z_x}q^{-1}N_P(\RR)$, so we need to show that
$q^{z_x}q^{-1}\in N_P(\RR)$.  This follows in the first case since $q\in
N_P(\RR)$ and $N_P(\RR)$ is normal.  In the second case, write $q=vl$ with
$v\in N_P(\RR)$ and $l\in L_{P,x}(\RR)$.  Then $q^{z_x}q^{-1} =v^{z_x}
l^{z_x} l^{-1} v^{-1} = v^{z_x}v^{-1}$, since $z_x\in A_{P,x}$ is in the
center of $L_P(\RR)$.  Now proceed as in the first case.
\end{proof}

\subsection{Comparison of geodesic actions}
\label{ssectCompareGeodesicActions}
Consider now two parabolic $\QQ$-subgroups $P\subset Q$; as in
\S\ref{sssectSplitComponents} the split component $A_Q$ may be viewed as a
subgroup of $A_P$.  Thus it has a geodesic action viewed as a subgroup of
either $L_Q(\RR)$ or $L_P(\RR)$; these two geodesic actions agree
\cite[\S3.10]{refnBorelSerre}.

For the extended geodesic action we have the

\begin{lem*}
The geodesic action of $N_P^Q(\RR) = N_P(\RR)/N_Q(\RR)
\subset L_Q(\RR)$ commutes with the geodesic action of $L_P(\RR)$.
\end{lem*}

\begin{proof}
Write $y\in D$ as $y=urx$, where $u\in N_P(\RR)$, $r\in L_{P,x}(\RR)$.
Since $N_P \subset Q = N_Q
\rtimes L_{Q,x}$ we can decompose $N_P= N_Q \rtimes N_{P,x}^Q$ where
$N_{P,x}^Q\subset L_{Q,x}$ is the lift of $N_P^Q$.  Thus we can
write $u =vw$, where $v\in N_Q(\RR)$ and $w\in
N_{P,x}^Q(\RR)$.  The geodesic action of $n\in N_P^Q(\RR)$ is now $n\geo y
= v n_x wrx$, which clearly commutes with the geodesic action of $z\in
L_P(\RR)$, $z\geo y = u z_x r x= vwz_x r x$.
\end{proof}

\subsection{Geodesic retraction}
\label{ssectGeodesicRetraction}
By Corollary ~\ref{ssectGeodesicAction} the geodesic action of $A_P$ and
the ordinary action of $\lsp0 P(\RR)$ commute.  The consequent action of
$A_P\times \lsp0 P(\RR)$ on $D$ is transitive and the stabilizer of any
point $x$ is $\{A_G\}\times\lsp0 P(\RR)\cap K_x$, where $K_x\subset
G(\RR)$ is the maximal compact subgroup fixing $x$.  Set%
\footnote{In \cite{refnBorelSerre} $e_P$ is denoted $e(P)$.}
$\mathscr A_P^G \equiv \lsp0 P(\RR)\back D$ and $e_P \equiv A_P\back D$.  The preceding
implies there is an isomorphism
\begin{equation}
D \cong \mathscr A_P^G \times e_P \label{eqnTrivializationPrincipalAPBundle}
\end{equation}
of $(A_P^G\times \lsp0 P(\RR))$-homogeneous spaces.  The quotient map
\begin{equation*}
D\longrightarrow e_P=A_P\back D
\end{equation*}
obtained by collapsing the orbits of the geodesic action of $A_P$ is a
trivial principal $A_P^G$-bundle with canonical trivializing sections given
by the orbits of $\lsp0 P(\RR)$.  Note that after choosing a basepoint
$x\in D$ and hence $\bar x\in \mathscr A_P^G$, there is an isomorphism of
$A_P^G$-homogeneous spaces $A_P^G\tildearrow \mathscr A_P^G$ given by
$a\mapsto a\geo \bar x$.  Thus $\mathscr A_P^G$ is an affine version of
$A_P^G$.

The action of $\G_P=\G\cap P$ on $D \cong \mathscr A_P^G \times e_P$
operates only on the second factor, so the principal $A_P^G$-bundle
structure persists after taking the quotient:
\begin{equation*}
r_P\colon  \G_P\back D \longrightarrow Y_P\equiv \G_P\back e_P.
\end{equation*}
The map $r_P$ is called \emph{geodesic retraction}.  Since $\G_P\subset
\lsp0 P(\RR)$, the quotients by $\G_P$ of the orbits of $\lsp0 P(\RR)$
still yield canonical trivializing sections of $r_P$.

If $P\subset Q$ this map descends to yield a canonically trivial
$A_P^Q$-bundle (also called geodesic retraction)
\begin{equation*}
r_P\colon  \G_P\back e_Q \longrightarrow Y_P.
\end{equation*}

\subsection{Nilmanifold fibrations}
\label{ssectNilmanifoldFibrations}
The geodesic action of $\lsp 0L_P(\RR)$ descends to $e_P$ and by Corollary
~\ref{ssectGeodesicAction} this action commutes with the $N_P(\RR)$-action.
The resulting action of $N_P(\RR)\times \lsp 0L_P(\RR)$ on $e_P$ is
transitive and the stabilizer of any point is contained in $\{1\}\times
\lsp 0L_P(\RR)$.  Set $\mathscr N_P(\RR)\equiv \lsp 0L_P(\RR)\back e_P$ and $D_P \equiv N_P(\RR)\back e_P$; the space $\mathscr N_P(\RR)$ is an affine version of
$N_P(\RR)$ and $D_P$ is the generalized symmetric space associated to
$L_P$.  Then there is an isomorphism
\begin{equation*}
e_P\cong \mathscr N_P(\RR)\times D_P
\end{equation*}
of $(N_P(\RR)\times \lsp 0L_P(\RR))$-homogeneous spaces. The principal
$N_P(\RR)$-bundle
\begin{equation*}
\mu\colon  e_P \longrightarrow D_P
\end{equation*}
obtained by collapsing the orbits of the action of $N_P(\RR)$ has canonical
trivializing sections given by the orbits of the geodesic action of $\lsp
0L_P(\RR)$.  Given a basepoint $x\in D$ and hence $\bar x\in\mathscr
N_P(\RR)$, there is an isomorphism $n_x\colon  N_P(\RR)\tildearrow\mathscr
N_P(\RR)$ of $N_P(\RR)$-homogeneous spaces defined by $n_x(u)= u \bar x$.

The action of $\lsp 0P(\RR)$ on $e_P$ descends to actions on both $\mathscr
N_P(\RR)$ and $D_P$, and the isomorphism $e_P\cong \mathscr N_P(\RR)\times
D_P$ is an isomorphism of $\lsp 0P(\RR)$-homogeneous spaces.  Set $\mathscr
N_P(\RR)' = \G_{N_P}\back \mathscr N_P(\RR)$ and note that the action of
$\G_P$ on $\mathscr N_P(\RR)$ descends to an action of
$\G_{L_P}=\G_P/\G_{N_P}$ on $\mathscr N_P(\RR)'$.  The previous
decomposition descends to an isomorphism $\G_{N_P}\back e_P \cong \mathscr
N_P(\RR)' \times D_P$ of $\G_{L_P}$-homogeneous spaces and thus there is a
canonically trivial $\mathscr N_P(\RR)'$-bundle
\begin{equation*}
\mu'\colon  \G_{N_P}\back e_P \to D_P\ .
\end{equation*}
Set%
\footnote{In \cite{refnGoreskyHarderMacPherson} the group $N_P(\RR)$ is
denoted $\mathpsscr U_P$ and the nilmanifold $N_P(\RR)'$ is denoted $N_P$.
What we denote here as $\mathscr N_P(\RR)$ and $\mathscr N_P(\RR)'$ are not
explicitly mentioned.}
$N_P(\RR)' = \G_{N_P}\back N_P(\RR)$.  Given a basepoint $x\in D$, the
isomorphism $n_x$ descends to $n_x'\colon  N_P(\RR)'\tildearrow \mathscr
N_P(\RR)'$.

Finally take the quotient by $\G_{L_P}$ of $\mu'$; we obtain a flat bundle
(with fibers noncanonically diffeomorphic to $\mathscr N_P(\RR)'$)
\begin{equation*}
\nil\colon  Y_P \longrightarrow X_P = \G_{L_P}\back D_P.
\end{equation*}
This is simply the flat bundle $ D_P \times_{\G_{L_P}} \mathscr N_P(\RR)'
\to X_P$ associated to the $\G_{L_P}$-space $\mathscr N_P(\RR)'$. It is
called the \emph{nilmanifold fibration}.

\subsection{}
\label{ssectActionOnNilmanifoldGivenBasepoint}
For later use, we calculate \cite[(7.8)]{refnGoreskyHarderMacPherson} the
action of $\G_{L_P}$ on $N_P(\RR)'$ that is needed in order to make $n_x'\colon 
N_P(\RR)'\tildearrow \mathscr N_P(\RR)'$ an isomorphism of
$\G_{L_P}$-spaces.  Let $\g\G_{N_P}\in \G_{L_P}$ and let $\G_{N_P}u\in
N_P(\RR)'$; write $\g=vl$ where $v\in N_P(\RR)$ and $l\in L_{P,x}(\RR)$.
Then
\begin{equation*}
\g\G_{N_P}\cdot n_x'(\G_{N_P}u) = \G_{N_P} \g. u.\bar x = \G_{N_P}
v.l.u.\bar x = \G_{N_P} v.\lsp lu. \bar x = n_x'(\G_{N_P}v.\lsp lu)\ .
\end{equation*}
Thus the induced action of $\g$ on $N_P(\RR)'$ is $\G_{N_P}u \mapsto
\G_{N_P}v.\lsp lu = \G_{N_P}\lsp \g u. v$.

\section{The Reductive Borel-Serre Compactification}
\label{sectReductiveBorelSerre}
Let $X=\G\back D$ be a locally symmetric space.  We describe the
Borel-Serre compactification $\Xbar$ \cite{refnBorelSerre} and its quotient
introduced by Zucker \cite{refnZuckerWarped}, the reductive Borel-Serre
compactification $\Xhat$ \cite[\S8]{refnGoreskyHarderMacPherson}.

\subsection{The Borel-Serre compactification}
\label{ssectBorelSerre}
Let $P$ be a parabolic $\QQ$-subgroup of $G$.  There is an isomorphism
$A_P^G \tildearrow (\RR^{> 0})^{\D_P}$ defined by $a\mapsto
(a^\al)_{\al\in\D_P}$ and thus $A_P^G$ may be embedded into a semigroup
$\Abar_P^G\tildearrow (\RR^{> 0}\cup\{\infty\})^{\D_P}$.

The orbits of the action of $A_P^G$ on $\Abar_P^G$ yield a stratification
indexed by the parabolic $\QQ$-subgroups $Q\supset P$.  Namely, let
$o_Q\in\Abar_P^G$ be defined by $\al(o_Q)=\infty$ for $\al\in\D_P^{P,Q}$,
and $\al(o_Q) = 1$ for $\al\in\D_P^Q$.  Then we have the stratification
\begin{equation}
\Abar_P^G = \coprod_{Q\supset P} A_P^G\cdot o_Q.
\label{eqnAbarStratification}
\end{equation}

Recall that the quotient map $D \to e_P$ by the geodesic action of $A_P^G$
is a principal $A_P^G$-bundle.  Let $D(P)= D \times_{A_P^G} \Abar_P^G$ be
the associated $\Abar_P^G$-bundle.  The orbits of $A_P^G\times \lsp0
P(\RR)$ yield a stratification
\begin{equation*}
D(P) = \coprod_{Q\supset P} e_Q\ .
\end{equation*}
If $P\subset Q$ there is a natural identification $D(Q)\subset D(P)$ as
an open union of strata.

Given the above identifications, set $\Dbar = \bigcup_P D(P)$.  This is a
manifold-with-corners whose stratification by faces is
\begin{equation*}
\Dbar =  \coprod_{P} e_P\ .
\end{equation*}
The construction of $\Dbar$ may be generalized to apply to $e_P$ (being a
space of type $S-\QQ$ \cite[\S\S2.3,~ 3.9]{refnBorelSerre}) and the resulting
$\overline{e_P}$ may be identified with the closure $\cl{e_P}$ of $e_P$ in
$\Dbar$.

The action of $G(\QQ)$ extends to $\Dbar$ and the quotient $\Xbar=\G\back
\Dbar$ is compact.  This is the \emph{Borel-Serre compactification}
of $X$.  On each stratum $e_P$ the equivalence relation induced by $\G$ is
that of $\G_P=\G\cap P$, and two strata $e_P$ and $e_{P'}$ become
identified in the quotient if and only if $P'=\lsp \g P$ for some $\g\in
\G$.  Thus the strata of $\Xbar$ are indexed by the $\G$-conjugacy classes
$\Pl$ of parabolic $\QQ$-subgroups and we have
\begin{equation*}
\Xbar= \coprod_{P\in \Pl} Y_P\ ,
\end{equation*}
where recall $Y_P=\G_P\back e_P$.  The
bijection $P\mapsto Y_P$ from $\Pl$ to the strata of $\Xbar$ is an
isomorphism of partially ordered sets, where $Y_P\le Y_Q$ if and only if
$Y_P\subset \cl{Y_Q}$.  Note that $\cl{Y_Q}$ may be identified with the
Borel-Serre compactification $\Ybar_Q$ of $Y_Q$.

\subsection{Local structure of $\Xbar$}
\label{ssectLocalStructureXbar}
Set $\bar{\mathscr A}_P^G = \mathscr A_P^G \times_{A_P^G} \Abar_P^G$.  The
factorization \eqref{eqnTrivializationPrincipalAPBundle} extends to $D(P)=
\bar{\mathscr A}_P^G \times e_P$.  It follows from reduction theory that
near a point of the stratum $e_P$ of $\Dbar$, the equivalence relation
induced by $\G$ is that of $\G_P$.  Thus $\Xbar$ is locally homeomorphic
along $Y_P$ to $\bar{\mathscr A}_P^G \times Y_P$ (where we embed $Y_P
\subset \bar{\mathscr A}_P^G \times Y_P$ as $\{o_P\}\times Y_P$).
Consequently if $\G$ is torsion-free then $\Xbar$ is a smooth
manifold-with-corners.  In general, since any arithmetic group has a
torsion-free subgroup of finite index there may be finite quotient
singularities and $\Xbar$ is a ``$V$-manifold-with-corners''.

The description of a neighborhood of the entire closed stratum $\Ybar_P$ is
much more subtle; see \cite[\S8]{refnSaperTilings}.
 
\subsection{Reductive Borel-Serre compactification}
\label{ssectReductiveBorelSerre}
Let $\Xhat$ be the stratified \emph{set}
\begin{equation*}
\Xhat= \coprod_{P\in \Pl} X_P
\end{equation*}
and define a stratified map $\nil\colon \Xbar \to \Xhat$ by setting
$\nil|_{Y_P}$ to be the flat nilmanifold fibration $Y_P\to X_P$ with fiber
$\mathscr N_P(\RR)'$.  We give $\Xhat$ the quotient topology induced by
$\nil$; this is the \emph{reductive Borel-Serre compactification} of
$X$.  (Note that $X=X_G=Y_G$.)  Again $P\mapsto X_P$ is an isomorphism of
$\Pl$ with the partially ordered set of strata of $\Xhat$.  The closure
$\cl{X_Q}$ is identified with the reductive Borel-Serre compactification
$\Xhat_Q$ of $X_Q$.

\subsection{Local structure of $\Xhat$}
\label{ssectLocalStructureXhat}
Note that the identification $A_P^G\cong \mathscr A_P^G$ from
\S\ref{ssectGeodesicRetraction} (which depends on the choice of a
basepoint) extends to an identification $\bar{\mathscr A}_P^G\cong
\Abar_P^G$.  The stratification \eqref{eqnAbarStratification} of
$\Abar_P^G$ thus induces one of $\bar{\mathscr A}_P^G$ and this
stratification is independent of the basepoint.

Set $\mathscr Z_P = (\bar{\mathscr A}_P^G \times \mathscr N_P(\RR)')/{\sim}$,
where $( a,\G_{N_P} n_1) \sim ( a,\G_{N_P} n_2)$ if and only if $n_2= u
n_1$ for $u\in N_Q(\RR)$ and $a$ belongs to the $Q$-stratum.  The action of
$\G_{L_P}$ on $\mathscr N_P(\RR)'$ induces an action on $\mathscr Z_P$.
Then $\Xhat$ is locally homeomorphic along $X_P$ to the bundle
$D_P\times_{\G_{L_P}} \mathscr Z_P$.

Note that if $\G$ is neat then the strata $X_P$ are smooth, however in
general they are only $V$-manifolds.  Also a description of a neighborhood
of the entire closed stratum $\Xhat_P$ may easily be deduced from
\cite[\S8]{refnSaperTilings}.

\section{$\L$-modules}
\label{sectLsheaves}

\subsection{Notation}
\label{ssectNotationLsheaves}
Let $W\subset \Xhat$ be a union of strata, that is, $W= \bigcup_{P\in
\Pl(W)} X_P$ for some subset $\Pl(W)\subset \Pl$.  We say $W$ is an
\emph{admissible} space if it is locally closed or, equivalently, if
$P$, $Q\in\Pl(W)$ imply that $[P,Q]\subset \Pl(W)$.  For $P\in\Pl(W)$ let $U_P = \coprod_{R\ge P} X_R$ be the open star neighborhood of a stratum $X_P$.  We have inclusions of admissible spaces
\begin{equation*}
\begin{CD}
X_P @>i_P>> W @<j_P<< (U_P\setminus X_P)\cap W \rlap{\ ,} \\
\Xhat_P\cap W @>\ihat_P>> W @<\jhat_P<<  W\setminus (\Xhat_P\cap W) \rlap{\ .}
\end{CD}
\end{equation*}
(To lighten notation we omit $W$ from the notation for these maps and let it be determined by the context.)

\subsection{}
\label{ssectLstructure}
Let $W$ be an admissible space.  Let $\L_W$ denote the system
consisting of
\begin{equation*}
\begin{cases}
\text{the reductive algebraic groups $L_P$}
&\text{for all $P\in\Pl(W)$, and}\\
\text{the functors }H(\n_P^Q;\cdot)\colon \R(L_Q) \longrightarrow
\Graded(L_P) &\text{for all $P\le Q\in \Pl(W)$.}
\end{cases}
\end{equation*}
Note that the functors $H(\n_P^Q;\cdot)$ are degree nondecreasing and there
are natural isomorphisms (in view of \eqref{eqnKostantDegeneration})
\begin{equation}
H(\n_P^R;\cdot) \cong H(\n_P^Q;H(\n_Q^R;\cdot))\qquad\text{for all
$P\le Q \le R$.} \label{eqnLinkComposition}
\end{equation}

\subsection{The category of $\L$-modules}
\label{ssectLsheaves}
Let $\L=\L_W$ for a fixed admissible space $W$.  An \emph{$\L$-module}
$\M=(E_\cdot,f_{\cdot\cdot})$ consists of
\begin{equation*}
\begin{cases}
\text{an object }E_P \text{ of }\Graded(L_P)
&\text{for all $P\in\Pl(W)$, and}\\
\text{morphisms }f_{PQ}\colon H(\n_P^Q;E_Q)\to E_P[1] &
\text{for all $P\le Q\in \Pl(W)$}
\end{cases}
\end{equation*}
satisfying the condition that
\begin{equation}
\sum_{Q\in [P,R]} f_{PQ}[1]\circ H(\n_P^Q;f_{QR}) = 0 \qquad\text{for all
$P\le R\in\Pl(W)$.} 
\label{eqnLsheafCondition}
\end{equation}
(Equation \eqref{eqnLinkComposition} shows this formula makes sense.)

Note that \eqref{eqnLsheafCondition} implies in particular  that $(E_P,f_{PP})$
is a complex which we denote $\i_P^!\M$.  (See
\S\ref{ssectFunctorsOnLsheaves} below for a general definition of $k^!$.)

An \emph{$\L$-morphism} $\M\to \M'$ is a family
$\phi=(\phi_{\cdot\cdot})$, where
\begin{equation*}
\phi_{PQ}\colon  H(\n_P^Q;E_Q)\to E'_P \qquad  \text{for all $P\le
Q\in\Pl(W)$},
\end{equation*}
such that the following condition is satisfied:
\begin{equation}
\sum_{Q\in [P, R]} \phi_{PQ}[1]\circ
H(\n_P^Q;f_{QR}) = \sum_{Q\in [P,R]} f'_{PQ}\circ H(\n_P^Q;\phi_{QR})
\qquad\text{for all $P\le R\in\Pl(W)$.}
\label{eqnMorphismLsheafCondition}
\end{equation}
The composition of two $\L$-morphisms $\phi\colon \M\to \M'$ and
$\phi'\colon \M'\to \M''$ is defined by
\begin{equation}
(\phi'\circ\phi)_{PR} = \sum_{Q\in [P,R]}
\phi'_{PQ} \circ H(\n_P^Q;\phi_{QR}) \qquad\text{for all $P\le R\in\Pl(W)$.}
\label{eqnCompositionLMorphism}
\end{equation}

The \emph{category of $\L$-modules} $\R(\L)$ has for objects the
$\L$-modules and for morphisms the $\L$-morphisms.  It is an additive category.  The full subcategory
$\R^b(\L)$ consists of $\L$-modules $\M$ for which all $E_P$ lie in
$\Graded^b(L_P)$; we call these \emph{bounded $\L$-modules}.

In \S\ref{sectRealizationLModules} we will define a functor $\A_W$ which associates to any $\L$-module $\M$ a constructible complex of sheaves.  The differential of $\A_W(\M)$ is constructed from the morphisms $(f_{\cdot\cdot})$.

Note that $\R(\L)$ is not defined as a
category of complexes.  However if $W=X_P$ consists of a single stratum, then
$\R(\L_{X_P})=\Complex(L_P)$.

\subsection{Standard functors on $\L$-modules}
\label{ssectFunctorsOnLsheaves}
Consider an inclusion $k\colon Z\hookrightarrow W$ of admissible spaces.  We will define functors $k^*$, $k_*$, $k^!$, and $k_!$ in some cases between the respective categories of $\L$-modules; these are motivated by the usual functors on the derived category of sheaves.  

The functor $k^!$ associates to an $\L_W$-module $\M$ an $\L_Z$-module with
the same data but with $P$ and $Q$ restricted to belong to $\Pl(Z)$.
Condition \eqref{eqnLsheafCondition} continues to hold since $Z$ is locally
closed.  One special case is the \emph{local
cohomology supported at $P$ complex}:
\begin{equation}\label{eqnLocalCohomologySupportedAtPFunctor}
  \i_P^!\M = (E_P,f_{PP})\ .
\end{equation}
The functor $k_*$ associates to
an $\L_Z$-module $\M$ an $\L_W$-module with the same data but extended by
$E_P=0$ and $f_{PQ}=0$ if any subscript is outside of $\Pl(Z)$.

We only define $k^*$ when $Z$ is open in $W$  or when $Z$ has a unique maximal face.  When $Z$ is open, $k^*=k^!$ has been
defined above.  When $Z$ has a unique maximal face $X_T$, define $k^*\M =
(E'_\cdot,f'_{\cdot\cdot})$ by
\begin{equation}\label{eqnPullBackFunctor}
\left\{
\begin{aligned}
E'_P&= \bigoplus_{R \in [P,(P,T)]} H(\n_P^R;E_R), \\
f'_{PQ}&= \sum_{\substack{R\cap T =P \\ S\cap T=Q \\ R\le S}} H(\n_P^R;f_{RS}),
\end{aligned}
\right.
\end{equation}
for all $P\le Q$ with $P$, $Q\in\Pl(Z)$.  (Recall that $(P,T)$ denotes the
parabolic $\QQ$-subgroup opposite to $T$ relative to $P$.)  The reader can
verify this is an $\L_Z$-module.  As a special cases, the \emph{local
cohomology at $P$ complex} is
\begin{equation}\label{eqnLocalCohomologyFunctor}
\i_P^*\M = \biggl( \bigoplus_{P\le R} H(\n_P^R;E_R), 
\sum_{P\le R\le S} H(\n_P^R;f_{RS}) \biggr)
\end{equation}
and the \emph{link cohomology at $P$ complex} is
\begin{equation}\label{eqnLinkCohomologyFunctor}
  \i_P^*j_{P*}j_P^*\M  =  \i_P^*\jhat_{P*}\jhat_P^*\M =\biggl( \bigoplus_{P < R} H(\n_P^R;E_R), 
\sum_{P<  R\le S} H(\n_P^R;f_{RS}) \biggr) \ .
\end{equation}
It should always be understood in these and similar formulas later that $R$ and
$S$ are restricted to belong to $\Pl(W)$.

We only define $k_!$ in the case we need, namely when
$Z$ is closed in $W$, in which case $k_!=k_*$ has been
defined above.

\begin{prop*}
Assume in the following that $j$ and $k$ are inclusions of admissible
spaces and that the indicated functors have been defined above.
\begin{enumerate}
\item There are identities\textup:
\begin{alignat*}{2}
j^!\circ k^! &= (k\circ j)^!\ , &\qquad\qquad k_*\circ j_* &= (k\circ j)_*
\ ,\\
j^*\circ k^* &= (k\circ j)^*\ , &\qquad\qquad k_!\circ j_! &= (k\circ j)_!
\ .
\end{alignat*}
\label{itemFunctorsComposition}
\item There are adjoint relations\textup:
\begin{equation}
\begin{aligned}
\Hom_{\R(\L_W)}(\M,k_*\M') &\cong
     \Hom_{\R(\L_Z)}(k^*\M,\M') , \\
\Hom_{\R(\L_W)}(k_!\M',\M) &\cong
     \Hom_{\R(\L_Z)}(\M',k^!\M).
\end{aligned}\label{eqnAdjointRelations}
\end{equation}
\label{itemFunctorsAdjoint}
\item Consider the commutative diagram of inclusions of admissible
spaces,
\begin{equation*}
\begin{CD}
Z @>k_2>> Z_2 \\
@V j_1 VV     @V j_2 VV \\
Z_1 @>k_1>> W \rlap{\ .}
\end{CD}
\end{equation*}
where $Z=Z_1\cap Z_2$.  Then
\begin{equation}
k_1^!\circ  j_{2*} = j_{1*}\circ k_2^!\ .
\end{equation}
Now assume that $Z_1$ has a unique maximal face $X_T$, that $Z$ has a
unique maximal face $X_{T_0}$, and furthermore that $Z_1$ and $Z_2$
intersect ``transversely'' in the sense that
\begin{equation}
U\cap Z_2 = U\cap W\cap\Xhat_{(T_0,T)},\label{eqnTransverse}
\end{equation}
where $U = \bigcup_{P\in\Pl(Z)}\bigcup_{Q\ge P} X_Q$ is the star
neighborhood of $Z$.  Then
\begin{equation*}
j_1^!\circ k_1^* = k_2^* \circ j_2^!.
\end{equation*}
\label{itemFunctorsSquare}
\end{enumerate}
\end{prop*}
\begin{proof}
Parts \itemref{itemFunctorsComposition} and \itemref{itemFunctorsAdjoint}
and the first part of \itemref{itemFunctorsSquare} are left to the reader.
For the rest of \itemref{itemFunctorsSquare}, note that $E'_P = \bigoplus_Q
H(\n_P^Q;E_Q)$ for both $j_1^! k_1^*\M$ and $k_2^* j_2^!\M$, where in the
first case the sum is over $Q\in [P,(P,T)]\cap\Pl(W)$, while in the second
case the sum is over $Q\in [P,(P,T_0)]\cap\Pl(Z_2)$.  The sums are equal by
\eqref{eqnTransverse} since one may check $(P,T)= (P,T_0)\cap (T_0,T)$.
\end{proof}

\subsection{Pullback of an $\L_W$-module to a fiber}
\label{ssectPullbackToFiber}
There is one other functor we will need to consider beginning in
\S\ref{sectFunctorialityInverseImages}.  Suppose $G$ has a connected normal
$\QQ$-subgroup which we denote $G_{\l}$.  Let $\tilde G_h$ be a
complementary connected normal $\QQ$-subgroup so that $G=\tilde G_h G_\l$
is an almost direct product.  Set $G_{h} = G/G_{\l}=\tilde
G_h/(G_\l\cap \tilde G_h)$ and let $\G_{G_{\l}}=\G\cap G_{\l}$ and
$\G_{G_{h}}=\G/ \G_{G_{\l}}$ be the induced arithmetic subgroups.
Parabolic $\QQ$-subgroups of $G$ correspond to pairs of a parabolic
$\QQ$-subgroup of $G_\l$ and a parabolic $\QQ$-subgroup of $G_h$: given
$P\subset G$ we associate $P_\l=P\cap G_\l$ and $P_h=P/P_\l$ and
conversely given $P_\l$ and $P_h$ we associate $P=\tilde P_h P_\l$
where $\tilde P_h\subset \tilde G_h$ is the inverse image of $P_h$ under the finite quotient map $\tilde G_h \to G_h$.

The geometric picture is as follows.  There is a factorization of symmetric
spaces $D = D_{h}\times D_{\l}$ and a flat bundle of arithmetic quotients
$X \to X_{h}$ with fibers $X_\l$.  (The action of $\G_{G_h}$ on $X_\l$ is induced
from the action of the finite quotient $\G/\G_{G_\l}(\G\cap \tilde G_h)$ so
this bundle becomes trivial over a finite cover of $X_h$.)  Let
$X(G_\l)\subset \Xhat$ be the partial compactification of $X$ obtained by
replacing the fibers by $\Xhat_\l$.  Specifically, note that the map
$P\mapsto P_h$ induces a surjection $\Pl(\Xhat) \to \Pl(\Xhat_h)$.  Then
set $X(G_\l)$ to be the union of strata $X_P$ for all $P\in\Pl(\Xhat)$
satisfying $P_h=G_{h}$.  The above flat bundle extends to a fibration
\begin{equation*}
\pi\colon  X(G_\l) \longrightarrow X_{h}
\end{equation*}
with fibers $\Xhat_{\l}$.  Let
$\ihat_{G,\l}\colon \Xhat_{\l}\hookrightarrow X(G_\l)$ be the inclusion of a
generic fiber.

The map $P_\l\mapsto P=\tilde G_h P_\l$ induces a surjection
$\Pl(\Xhat_\l)\to \Pl(X(G_\l))$.  (It may not be injective if the
finite group $\G/\G_{G_\l}(\G\cap \tilde G_h)$ is nontrivial.)  Define the
connected normal $\QQ$-subgroup $L_{P,\l} = P_\l/N_P\subset L_P$.
Let $\M=(E_\cdot,f_{\cdot\cdot})$ be an $\L_{X(G_\l)}$-module.  Define an
$\L_{\Xhat_{\l}}$-module $\ihat_{G,\l}^*\M =
(E'_\cdot,f'_{\cdot\cdot})$ by
\begin{equation*}
\left\{
\begin{aligned}
E'_{P_\l}&= \Res_{L_{P,\l}}^{L_P} E_P, \\
f'_{P_\l Q_\l}&= \Res_{L_{P,\l}}^{L_P} f_{PQ},
\end{aligned}
\right.
\end{equation*}
for all $P_\l\le Q_\l\in\Pl(\Xhat_\l)$.
This makes sense
since there is a natural isomorphism
\begin{equation*}
\Res_{L_{P,\l}}^{L_P} H(\n_P^Q;E_Q) \cong
H(\n_{P_\l}^{Q_\l};\Res_{L_{Q,\l}}^{L_Q} E_Q)\ .
\end{equation*}
Similarly define $\ihat_{G,\l}^!\M$ except that in this case
$E'_{P_\l} = \Res_{L_{P,\l}}^{L_P} E_P[-\dim D_h]$.

More generally for $R\in\Pl$ suppose  we are given a connected normal $\QQ$-subgroup
$L_{R,\l}$ of $L_R$ and hence as above an almost direct product $L_R=\widetilde{L_{R,h}} L_{R,\l}$.  We then have a fibration
$X_R(L_{R,\l})\to X_{R,h}$ and the inclusion of the generic fiber $\ihat_{R,\l}\colon \Xhat_{R,\l} \to X_R(L_{R,\l})$.  There is a surjection
$\Pl(\Xhat_{R,\l})\to \Pl(X_R(L_{R,\l}))$ given by $P_\l\mapsto P$
where $P/N_R = \widetilde{L_{R,h}} P_\l $.
Thus for $P\in\Pl(X_R(L_{R,\l}))$ we have a connected normal $\QQ$-subgroup $L_{P_\l}$ of $L_P$ and we can define as above the functors
\begin{equation*}
  \ihat_{R,\l}^*,\  \ihat_{R,\l}^! \colon \Rep(\L_{X_R(L_{R,\l})}) \longrightarrow \Rep(\L_{\Xhat_{R,\l}}) \ .
\end{equation*}
We note the following
identities for $P_\l\in\Pl(\Xhat_{R,\l})$:
\begin{equation}
\begin{alignedat}{2}
\i_{P_\l}^* \circ \ihat_{R,\l}^* &= \i_{P_\l}^* \circ \ihat_{P,\l}^*
\circ \i_P^*\ , &\qquad\qquad
\i_{P_\l}^* \circ \ihat_{R,\l}^! &= \i_{P_\l}^* \circ \ihat_{P,\l}^!
\circ \i_P^*\ , \\
\i_{P_\l}^! \circ \ihat_{R,\l}^* &= \i_{P_\l}^* \circ \ihat_{P,\l}^*
\circ \i_P^!\ , &\qquad\qquad
\i_{P_\l}^! \circ \ihat_{R,\l}^! &= \i_{P_\l}^* \circ \ihat_{P,\l}^!
\circ \i_P^!\ .
\end{alignedat}
\label{eqnPullbackToFiberFormulas}
\end{equation}
(The initial $\i_{P_\l}^*$ on the right-hand side of these equalities
would not be needed if $\G_{L_R}/\G_{L_{R,\l}}(\G_{L_R}\cap \widetilde {L_{R,h}})$ were trivial.)

\subsection{Local cohomology with supports}
\label{ssectLocalCohomologyWithSupports}
Let $\M$ be an $\L_W$-module, where $W$ is an admissible space.  There is a short exact sequence
\begin{equation}\label{eqnShortLocalCohomology}
0 \to \i_P^!\M \to \i_P^*\M \to  \i_P^*\j_{P*}\j_P^*\M \to 0
\end{equation}
relating the two types of local cohomology at $P\in \Pl(W)$ (see \eqref{eqnLocalCohomologySupportedAtPFunctor}, \eqref{eqnLocalCohomologyFunctor}, and \eqref{eqnLinkCohomologyFunctor}).
More generally, for $P\le Q\in \Pl(W)$ define the \emph{local
cohomology at $P$ supported on $Q$ complex} to be
\begin{equation}\label{eqnLocalCohomologyWithSupportsFunctor}
\i_P^*\ihat_Q^!\M = \biggl( \bigoplus_{R \in [P,Q]} H(\n_P^R;E_R), 
\sum_{R\le S \in [P,Q]} H(\n_P^R;f_{RS}) \biggr)\ ;
\end{equation}
this corresponds to $\i_P^!\M$ and $\i_P^*\M$ when $Q=P$ and $G$ respectively.
It is an essential part of the definition of the micro-support of $\M$ (see \S\ref{ssectMicroSupport}).
If $Q\le Q'\in\Pl(W)$ the short exact sequence \eqref{eqnShortLocalCohomology}
generalizes to
\begin{equation}
0\to \i_P^* \ihat_{Q}^!  \M \to \i_P^* \ihat_{Q'}^! \M \to \i_P^*
\jhat_{Q*}\jhat_Q^*\ihat_{Q'}^! \M \to 0 \ .
\label{eqnShortCompareLocalCohomologyWithSupports}
\end{equation}
The corresponding long exact sequence is
\begin{equation}
\cdots \longrightarrow
H^i(\i_P^* \ihat_Q^!  \M) \longrightarrow
H^i(\i_P^* \ihat_{Q'}^! \M) \longrightarrow
H^i(\i_P^* \jhat_{Q*}\jhat_Q^*\ihat_{Q'}^! \M) \longrightarrow \cdots\ .
\label{eqnLongCompareLocalCohomologyWithSupports}
\end{equation}

We wish to study $H(\i_P^* \jhat_{Q*}\jhat_Q^*\ihat_{Q'}^! \M)$.  Set
$P'= (P,Q)\cap Q'$ so that $P=P'\cap Q$ and $Q'=P'\vee Q$:
\begin{equation}
\vcenter{\xymatrix @ur @M=1pt @R=1.5pc @C=1pc {
{P'} \ar@{-}[r] \ar@{-}[d] & {Q'} \ar@{-}[d] \\
{P\rlap{\qquad .}} \ar@{-}[r] & {Q}
}}
\end{equation}
When $Q=Q'$ clearly $H^i(\i_P^* \jhat_{Q*}\jhat_Q^*\ihat_{Q'}^! \M)=0$.  In general we have the

\begin{lem}
\label{ssectRelativeLocalCohomologySupportsSS}
Given $P\le Q< Q'$, set $P'= (P,Q)\cap Q'$.  There are two spectral
sequences abutting to $H(\i_P^* \jhat_{Q*}\jhat_Q^*\ihat_{Q'}^!  \M)$.  The
Fary spectral sequence \cite[IV, \S12]{refnBredonSheafTheory}, \cite{refnFary} has
\begin{equation}
E_1^{-p,\cdot} =
\bigoplus_{\substack{P < \tilde P \le P' \\ \#\D_P^{\tilde P} = p}}
H(\n_P^{\tilde P};H(\i_{\tilde P}^*\ihat_{\tilde Q}^!\M))[-p] \ ,
\label{eqnFarySpectralSequence}
\end{equation}
where $\tilde Q = \tilde P\vee Q$.  The Mayer-Vietoris
spectral sequence has
\begin{equation}
E_1^{p,\cdot} =
\bigoplus_{\substack{P < \tilde P \le P' \\ \#\D_P^{\tilde P} = p+1}}
H(\n_P^{\tilde P};H(\i_{\tilde P}^*\ihat_{Q'}^!\M)) \ .
\label{eqnMayerVietorisSpectralSequence}
\end{equation}
\end{lem}

\begin{rem*}
  When applying the lemma it is helpful to note that since 
the functor $H(\n_P^{\tilde P};\cdot)$ is exact it commutes with
taking cohomology and thus
\begin{equation}\label{eqnPushForwardType}
H(\i_P^* \i_{\tilde P*}\i_{\tilde P}^*\ihat_{\tilde Q}^! \M) =
H(\n_P^{\tilde P}; H(\i_{\tilde P}^*\ihat_{\tilde Q}^! \M))\ .
\end{equation}
\end{rem*}

\begin{proof}
Any $R\in [P,Q']$ belongs to $[\tilde P,\tilde Q]$ for a unique $\tilde
P\in R\cap P'\in [P,P']$ and $\tilde Q = \tilde P\vee Q\in[Q,Q']$.  The situation is represented by
the following diagram, in which each rectangle has the analogous
meaning to \eqref{eqnParallelogram}:
\begin{equation}
\vcenter{\xymatrix @ur @M=1pt @R=1.5pc @C=1pc {
{P'} \ar@{-}[rr] \ar@{-}[d] & {} & {Q'} \ar@{-}[d] \\
{\tilde P} \ar@{-}[r] \ar@{-}[d] & {R} \ar@{-}[r] & {\tilde Q}
\ar@{-}[d] \\
{P\rlap{\qquad\qquad .}} \ar@{-}[rr] & {} & {Q}
}} \label{eqnParallelogramWithPtilde}
\end{equation}
(Of course this figure is not meant to suggest that the parabolic
$\QQ$-subgroups lying between $P$ and $P'$ are linearly ordered.)

Now $\i_P^* \jhat_{Q*}\jhat_Q^*\ihat_{Q'}^! \M$ is obtained by restricting the data of $\M$ to $(\Xhat_{Q'}\setminus \Xhat_Q)\cap W$, extending it by $0$ on $\Xhat_Q\cap W$, and finally applying $\i_P^*$.  Thus by \eqref{eqnParallelogramWithPtilde} we may decompose  $\i_P^* \jhat_{Q*}\jhat_Q^*\ihat_{Q'}^! \M$ (as an \emph{$L_P$-module}, not
as a complex) as
\begin{equation}
\label{eqnTildePDecomposition}
\bigoplus_{P<\tilde P\le P'}\i_P^* \i_{\tilde P*}\i_{\tilde P}^*
     \ihat_{\tilde Q}^! \M
=\bigoplus_{P<\tilde P\le P'} H(\n_P^{\tilde P};\i_{\tilde P}^*
\ihat_{\tilde Q}^! \M)
\ .
\end{equation}
 
For the Fary spectral sequence, the sum corresponding to \eqref{eqnTildePDecomposition} in which $\tilde P$ is restricted to those with
$\#\D_P^{\tilde P} \le p$ is a subcomplex and defines an increasing
filtration on $\i_P^* \jhat_{Q*}\jhat_Q^*\ihat_{Q'}^! \M$ as $p$ varies.
The associated graded complex is exactly equal to
\eqref{eqnTildePDecomposition}.  In view of \eqref{eqnPushForwardType}, the
$E_1$ term of the resulting spectral sequence has the desired form,
\eqref{eqnFarySpectralSequence}.

For the Mayer-Vietoris spectral sequence, consider the long exact sequence
\begin{equation*}
0\to \i_P^* \jhat_{Q*}\jhat_Q^*\ihat_{Q'}^! \M \to 
\bigoplus_{\#\D_P^{\tilde P}=1} \i_P^*\i_{{\tilde P}*}\i_{\tilde P}^*
         \ihat_{Q'}^!\M \to
\bigoplus_{\#\D_P^{\tilde P}=2} \i_P^*\i_{{\tilde P}*}\i_{\tilde P}^*
         \ihat_{Q'}^!\M \to
\cdots
\end{equation*}
in which $\tilde P$ satisfies $P< \tilde P\le P'$.  Thus $H(\i_P^*
\jhat_{Q*}\jhat_Q^*\ihat_{Q'}^! \M)$ is the cohomology of a double complex;
the spectral sequence \eqref{eqnMayerVietorisSpectralSequence} corresponds
to the ``first filtration'' $\#\D_P^{\tilde P} \ge p+1$.
\end{proof}

\subsection{}
\label{ssectMappingCone}
Define the \emph{shift} of an $\L$-module $\M$ by $\M[1]= (E_\cdot[1],-f_{\cdot\cdot})$.  The \emph{mapping cone} of an $\L$-morphism
$\phi\colon \M\to \M'$ is defined by
\begin{equation*}
M(\phi)= \left(E_\cdot[1]\oplus E'_\cdot\, , 
\left(\begin{matrix} -f_{\cdot\cdot} & 0 \\
            \phi_{\cdot\cdot} & f'_{\cdot\cdot} \end{matrix}\right)\right).
\end{equation*}
This is an $\L$-module; condition \eqref{eqnLsheafCondition} follows from
the corresponding condition for $\M$ and $\M'$ and from
\eqref{eqnMorphismLsheafCondition}.  The mapping cone defines a functor from the
category of diagrams $\phi\colon \M\to \M'$ to $\R(\L)$.  There are
natural $\L$-morphisms $\al(\phi)\colon \M'\to M(\phi)$ and $\b(\phi)\colon M(\phi)\to
\M[1]$ given by
\begin{alignat*}{2}
\al(\phi)_{PP} &= \begin{pmatrix} 0 \\ \id_{E'_P}\end{pmatrix} &&\qquad\text{for all
$P$ and} \\
\b(\phi)_{PP} &= \begin{pmatrix} \id_{E_P[1]} & 0\end{pmatrix} &&\qquad\text{for all
$P$}.
\end{alignat*}
(Set $\al(\phi)_{PQ}=0$ for $P\lneq Q$ and similarly for $\b(\phi)$.)  One can check that for $P\in\Pl(W)$,
\begin{equation}
  \label{eqnLocalCohomologyOfMappingCone}
  \i_P^* M(\M\xrightarrow{\phi} \M') = M(\i_P^*\M \xrightarrow{\i_P^*\phi} \i_P^*\M')\ .
\end{equation}

\subsection{}
\label{ssectHomotopyLModules}
Although it is not needed in this work, one can define $\L$-morphisms $\phi_1$, $\phi_2:\M\to \M'$ to be \textit{homotopic} if
there exists morphisms
$$
h_{PQ}: H(\n_P^Q;E_Q) \to E'_P[-1] \qquad
\text{for all $P\le Q\in\Pl(W)$}
$$
such that
\begin{multline}
\phi_{1PR} - \phi_{2PR} = \sum_{P\le Q\le R} f'_{PQ}[-1]\circ H(\n_P^Q;h_{QR}) \\
+ \sum_{P\le Q\le R} h_{PQ}[1]\circ H(\n_P^Q;f_{QR}) \qquad\text{for all
$P\le R\in\Pl(W)$.}
\label{eqnHomotopyLsheaf}
\end{multline}
One sets $\K(\L)$ to
be the category with $\L$-modules as objects and homotopy classes of
$\L$-morphisms as morphisms.  With the definition of mapping cone in \S\ref{ssectMappingCone}, one can show following \cite[\S\S1.4,~ 1.5]{refnKashiwaraSchapira} that
$\K(\L)$ has the structure of a triangulated category.
A {\it quasi-isomorphism}
of $\L$-modules is an $\L$-morphism that induces quasi-isomorphisms on
local cohomology complexes.  Since $\R(L_P)$ is a semisimple abelian
category for all $P$ one can show that any quasi-isomorphism of $\L$-modules
has a homotopy inverse.  Thus it is reasonable to notate $\K(\L)$ also as
$\Derived(\L)$.

\subsection{Generalizations}
\label{ssectLModuleGeneralizations}
The concept of $\L$-modules may be greatly generalized.  For example, in
\S\ref{ssectLstructure} one could start with a general family of reductive
$\QQ$-groups $(L_\sigma)_\sigma$ indexed by a partially ordered set and
replace $H(\n_P^Q;\cdot)$ by some ``link cohomology'' functors
$\M_{\sigma\tau}\colon \R(L_\tau) \longrightarrow \Graded(L_\sigma)$ satisfying
$\M_\sigma^\upsilon \cong \M_\sigma^\tau\circ \M_\tau^\upsilon$.  One
particularly simple situation will be introduced in
\S\ref{ssectLModulesAPTimesXP} in which all the groups $L_\sigma$ will be
identical.  One could also drop the condition that the groups are reductive
algebraic and replace $\R(L_\sigma)$ by a suitable category of
representations.  We will not consider such generalizations further
here.

\section{Realization of $\L$-modules}
\label{sectRealizationLModules}
Let $W$ be an admissible space and let $\X$ denote the stratification 
$\coprod_{P\in\Pl(W)} X_P$.  Recall that $\Complex_\X(W)$ denotes the
category of constructible complexes of sheaves on $W$ and that
$\Derived_\X(W)$ denotes the corresponding derived category.  Note that an algebraic $L_P$-module $E_P$ induces by restriction a $\G_{L_P}$-module and
hence a locally constant sheaf $\EE_P=E_P \times_{\G_P} D_P$ on $X_P$.
Thus we obtain a functor $\R(L_P)\to \Derived_\X(X_P)$.  In this section we
generalize this to obtain a realization functor $\Sheaf_W\colon
\R(\L_W)\to \Derived_{\X}(W)$ and to define the cohomology of an $\L_W$-module.

\begin{thm}
  \label{ssectRealizationLModules}
  Let $W$ be an admissible space.
\begin{enumerate}
\item
There exists a family of functors
\begin{equation*}
\A_P\colon \R(L_P) \to \Complex_{\X}(\Xhat_P) \qquad\text{for all
$P\in\Pl$}
\end{equation*}
and of natural morphisms  
\begin{equation}
k_{PQ}\colon \ihat_{Q*}\A_Q(\cdot)\longrightarrow \ihat_{P*}\A_P(H(\n_P^Q;\cdot))
\qquad\text{for all $P\le Q\in\Pl$}
\label{eqnRealizationPQ}
\end{equation}
with the following properties\textup: $\A_P(E_P)$ is quasi-isomorphic to $\i_{P*}\EE_P$
 in $\Derived_\X(\Xhat_P)$ and $\ihat_P^*(k_{PQ})$ is a
natural isomorphism in $\Derived_\X(\Xhat_P)$.
\label{itemExistenceOfBasicRealizationFunctors}
\item
A family of functors as in
\itemref{itemExistenceOfBasicRealizationFunctors}
determines functors
\begin{equation*}
\A_W\colon \R(\L_W)\to \Complex_{\X}(W)\qquad \text{and} \qquad \Sheaf_W\colon \R(\L_W)\to \Derived_{\X}(W)
\end{equation*}
that preserve boundedness.  If
$k\colon Z\hookrightarrow W$ is an inclusion of admissible spaces, then $k_*\circ
\Sheaf_Z= \Sheaf_W \circ k_*$ and similarly for the other functors defined
in \textup\S\textup\S\textup{\ref{ssectFunctorsOnLsheaves}} and
\textup{\ref{ssectPullbackToFiber}}.
\label{itemExistenceOfRealizationFunctors}
\item
If two families of functors as in
\itemref{itemExistenceOfBasicRealizationFunctors} are naturally
quasi-isomorphic compatibly with \eqref{eqnRealizationPQ}, then the
corresponding functors $\Sheaf_W$ are naturally isomorphic.
\label{itemUniquenessOfRealizationFunctors}
\end{enumerate}
\end{thm}

\begin{defn}
  The functor $\Sheaf_W$ is called a \emph{realization} of the category of $\L_W$-modules.   The \emph{cohomology $H(W;\M)$} of
an $\L_W$-module $\M$ with respect to the realization is defined to be the hypercohomology $H(W;\Sheaf_W(\M))$.  More
generally the cohomology $H(Y;\M)$ of any subset $Y\subset W$ with
coefficients in $\M$ is defined to be the hypercohomology of
$\Sheaf_W(\M)|_Y$.
\end{defn}
\begin{rem*}
  In this paper we use the incarnation $\A_W(\M)$
constructed below from special differential forms so cohomology is always with respect to this realization.  This is a complex of
fine sheaves so $H(W;\M)$ may be computed from the complex of
sections $\A_{\Xhat}(\M)(W)$.
\end{rem*}

\subsection{Proof of Theorem \ref{ssectRealizationLModules}\itemref{itemExistenceOfRealizationFunctors}\itemref{itemUniquenessOfRealizationFunctors}}
\label{ssectRealizationLsheaves}
We assume that a family of functors $\{\A_P\}_{P\in\Pl}$ exists as in
\itemref{itemExistenceOfBasicRealizationFunctors}.  Let
$\M=(E_\cdot,f_{\cdot\cdot})$ be an $\L_W$-module.  Define
\begin{equation}
\left\{
\begin{aligned}
\A_W(\M) &= \bigoplus_{P\in \Pl(W)} 
                            \ihat_{P*} \A_P(E_P), \\
d_{\A_W(\M)} &= \sum_{P\in \Pl(W)} d_P + \sum_{P\le Q\in \Pl(W)}
                            \A_P(f_{PQ}) \circ k_{PQ} \ ,
\end{aligned}
\right.
\label{eqnRealizationComplex}
\end{equation}
where $d_P$ is the differential of $\A_P(E_P)$.  (Note that in applying
$\A_P$ to $E_P$ we are using the conventions of
\S\ref{sssectNotationCategory}.)  This is a complex of sheaves.  For
 an $\L_W$-morphism $\phi=(\phi_{\cdot\cdot})$ define
\begin{equation*}
\A_W(\phi) = \sum_{P\le Q} \A_P(\phi_{PQ})\circ
k_{PQ} \ .
\end{equation*}
This is the desired functor for
\itemref{itemExistenceOfRealizationFunctors}; we let $\Sheaf_W$ be the
corresponding functor into $\Derived_\X(W)$.  It is easy to verify that
$\Sheaf_W$ commutes with the usual functors defined in
\S\S\ref{ssectFunctorsOnLsheaves} and \ref{ssectPullbackToFiber}
and that \itemref{itemUniquenessOfRealizationFunctors} holds.

\subsection{Special differential forms}
\label{ssectSpecialDifferentialForms}
We now prepare for the proof of Theorem
\ref{ssectRealizationLModules}\itemref{itemExistenceOfBasicRealizationFunctors}
which will be at the end of this section.

Let $E$ be an algebraic representation of $G$ and let $\EE$ denote the
corresponding locally constant sheaf on $\Xbar$ or on any stratum of
$\Xbar$.  For any stratum $Y_P$ of $\Xbar$, consider the diagram
\begin{equation*}
\Xbar = \G\back \overline{D} \overset{q}{\longleftarrow} \G_P\back D(P)
\overset{r_P}{\longrightarrow}  Y_P,
\end{equation*}
where $q\colon \G_P\back D(P)\hookrightarrow \G_P\back\overline{D} \to
\G\back\overline{D}$ is the inclusion followed by the quotient and $r_P$ is
the geodesic retraction.  Reduction theory implies that if $y\in Y_P$ there
exists a sufficiently small neighborhood $U\subset \Xbar$ of $y$ such
that $q\colon q^{-1}(U)\tildearrow U$ is a diffeomorphism.  An $\EE$-valued
differential form $\o$ on $\Xbar$ is \emph{locally lifted from the
boundary} \cite[\S8]{refnBorelStable} if for all strata $Y_P$ and all
$y\in Y_P$, there exists $U$ as above so that $q^*(\o|_U) = r_P^*\psi$
where $\psi$ is a smooth form on $U\cap Y_P$.  A form on a stratum
$Y_P=\G_P\back e_P$ is called \emph{$N_P(\RR)$-invariant} if its lift to
$e_P$ is invariant under the action of $N_P(\RR)$; let $\A\inv(Y_P;\EE)$
denote the sheaf of $N_P(\RR)$-invariant $\EE$-valued forms on $Y_P$ and
let $A\inv(Y_P;\EE)$ denote the global sections.

Let $\Asp(\Xbar;\EE)\subset \A(\Xbar;\EE)$ denote the subcomplex of forms
which are locally lifted from the boundary and whose restriction to each
stratum $Y_P$ is in $\A\inv(Y_P;\EE)$.  The sheaf complex of \emph{special
  differential forms} on $\Xhat$ is 
\cite[(13.2)]{refnGoreskyHarderMacPherson}
\begin{equation*}
\Asp(\Xhat;\EE) = \nil_{*} \Asp(\Xbar;\EE)\ ,
\end{equation*}
the pushforward under the quotient map $\nil\colon \Xbar\to \Xhat$.  Let
$A\sp(\Xbar;\EE)$ and $A\sp(\Xhat;\EE)$ denote the corresponding complexes
of global sections.

\subsection{Restriction to boundary strata}
\label{ssectRestrictionBoundaryStrata}
Recall from \S\ref{ssectNilmanifoldFibrations} that the nilmanifold
fibration $\nil\colon  Y_P \to X_P$ is the flat bundle associated to the
$\G_{L_P}$-space $\mathscr N_P(\RR)'$, and that given a basepoint $x\in D$
there is an isomorphism $n_x'\colon  N_P(\RR)'\tildearrow \mathscr N_P(\RR)'$.
The $\G_{L_P}$-action on $\mathscr N_P(\RR)'$ induces the structure of an
$\G_{L_P}$-module on $A\inv(\mathscr N_P(\RR)';\EE)$, the
$N_P(\RR)$-invariant forms; let $\AA\inv(\mathscr N_P(\RR)';\EE)$ denote%
\footnote{In \cite{refnGoreskyHarderMacPherson} corresponding complex of vector bundles is denoted
$\mathbf C^{\boldsymbol \cdot}(N_P;\mathbf E)$.}
the associated complex of locally constant sheaves on $X_P$.  The fiber
$\AA\inv(\mathscr N_P(\RR)';\EE)_x$ is the complex of $N_P(\RR)$-invariant
forms on $\nil^{-1}(x)$.

A special differential form $\o$ on $\Xhat$ determines an
$N_P(\RR)$-invariant form on $Y_P$ which we denote $\o|_{Y_P}$; this form
can be viewed as a differential form on $X_P$ with values in
$\AA\inv(\mathscr N_P(\RR)';\EE)$. In fact there
is a natural isomorphism of complexes \cite[\S\S12.6,~
13.4(2)]{refnGoreskyHarderMacPherson}
\begin{equation}
\label{eqnRestrictSpecialForms}
\i_P^*\Asp(\Xhat;\EE)  \longtildearrow
\A(X_P;\AA\inv(\mathscr N_P(\RR)';\EE))\ , \qquad \o\mapsto \o|_{X_P}\ ,
\end{equation}
where the complex on the right is the associated total complex.  This
operation extends to an isomorphism
\begin{equation*}
\ihat_P^*\Asp(\Xhat;\EE)\longtildearrow
\Asp(\Xhat_P;\AA\inv(\mathscr N_P(\RR)';\EE))
\end{equation*}
and hence there is a natural morphism
\begin{equation}\label{eqnRestrictclosedSpecialForms}
\Asp(\Xhat;\EE)\longrightarrow
\ihat_{P*}\Asp(\Xhat_P;\AA\inv(\mathscr N_P(\RR)';\EE)\ ,\quad \o\mapsto \o|_{\Xhat_P}\ .
\end{equation}

\subsection{}
\label{ssectResolutionBySpecialForms}
A theorem of Nomizu \cite{refnNomizu} and van Est \cite{refnvanEst} implies
that the natural inclusion of complexes of differential forms induces an
isomorphism
\begin{equation*}
H(A\inv(\mathscr N_P(\RR)';\EE))\cong H(\mathscr N_P(\RR)';\EE);
\end{equation*}
from this, the Poincar\'e lemma, and \eqref{eqnRestrictSpecialForms} one
sees that

\begin{lem*}
\textup(\cite[\S13]{refnGoreskyHarderMacPherson}\textup) There is a
natural isomorphism $\Asp(\Xhat;\EE)\tildearrow \i_{G*}\EE$ in
$\Derived_{\X}(\Xhat)$ for $E\in\R(G)$.
\end{lem*}

\subsection{}
\label{ssectNilpotentCohomologySheaf}
Recall that $C(\n_P;E)$ denotes the complex $\Hom(\bigwedge\n_P;E)$ with
differential as in \cite[I, \S1.1]{refnBorelWallach}.  Its cohomology
$H(\n_P;E)$ is a graded $L_P$-module; let $\HH(\n_P;E)$ denote the
corresponding graded locally constant sheaf on $X_P$.

\begin{lem}
\label{ssectNilmanifoldCohomologyBundle}
\textup{(\cite[Prop.~12.15]{refnGoreskyHarderMacPherson})} There is a
natural quasi-isomorphism on $X_P$\textup:
\begin{equation*}
\coh_P\colon \AA\inv(\mathscr N_P(\RR)';\EE) \qiarrow \HH(\n_P;E).
\end{equation*}
\end{lem}
\begin{proof}
Choose a basepoint $x\in D$ and let $L_{P,x}\subset P$ denote the lift of
$L_P$ stable under the Cartan
involution associated to $x$.  $L_P$ acts on $C(\n_P;E)$ through
the coadjoint action of $L_{P,x}$; we
denote by $C(\n_P;E)_x$ the resulting complex of $L_P$-modules and by
$\CC(\n_P;E)_x$ the corresponding complex
of locally constant sheaves.  

The quasi-isomorphism $\coh_P$ is given by a composition of quasi-isomorphisms
\begin{equation}
\AA\inv(\mathscr N_P(\RR)';\EE) \overset{\coh_x}{\longrightarrow}
\CC(\n_P;E)_{x} \longrightarrow
\HH(\n_P;E)\ .\label{eqnToCohomology}
\end{equation}
The first map is induced by an isomorphism
\cite[(12.13)]{refnGoreskyHarderMacPherson} of the underlying
$\G_{L_P}$-modules.  Namely, transfer a form from $\mathscr N_P(\RR)'$ to
$N_P(\RR)'$ via the isomorphism $n_x'$ from \S\ref{ssectNilmanifoldFibrations} and evaluate the form at the
identity; one may use \S\ref{ssectActionOnNilmanifoldGivenBasepoint} to
verify this respects the $\G_{L_P}$-action.  For the second map of
\eqref{eqnToCohomology}, note that if $E$ were irreducible, a theorem of
Kostant \cite{refnKostant} says that each irreducible component of
$H(\n_P;E)$ occurs in both $H(\n_P;E)$ and $C(\n_P;E)_x$ with multiplicity
one.  Thus for general algebraic $E$ there is a unique map $C(\n_P;E)_x\to
H(\n_P;E)$ in $\Complex(L_P)$ inducing the identity on cohomology.  This
induces the second map in \eqref{eqnToCohomology}.

If $x'$ is an another basepoint then $L_{P,x'}=v L_{P,x}v^{-1}$ for some
$v\in N_P(\RR)$.  Thus there is an isomorphism $\Ad^*(v)\colon  C(\n_P;E)_{x}
\tildearrow C(\n_P;E)_{x'}$ in $\Complex(L_P)$.  One may check from the
definition of $h_x$ that the left triangle of
\begin{equation*}
\xymatrix{
{} & {\CC(\n_P;E)_{x}} \ar[dd]_{\Ad^*(v)} \ar[rd] & {} \\
{\AA\inv(\mathscr N_P(\RR)';\EE)} \ar[ru]^{\coh_{x}} \ar[rd]_{\coh_{x'}} & {} &
{\HH(\n_P;E)} \\
{} & {\CC(\n_P;E)_{x'}} \ar[ru] & {}
}
\end{equation*}
commutes; since $\Ad^*(v)$ induces the identity on cohomology, the right
triangle commutes by the uniqueness noted above.  Thus $\coh_P$ is independent of
the basepoint $x$.

The verification of naturality for each map in \eqref{eqnToCohomology} is
left to the reader.
\end{proof}

\subsection{}
\label{ssectRestrictSpecialForms}
This lemma together with equation \eqref{eqnRestrictclosedSpecialForms}
yields the
\begin{cor*}
There is a natural morphism
\begin{equation*}
k_P\colon \Asp(\Xhat;\EE) \longrightarrow \ihat_{P*}\Asp(\Xhat_P;\HH(\n_P;E))
\end{equation*}
given by $\o\mapsto \Asp(\Xhat_P;\coh_P)(\o|_{\Xhat_P})$ such that
$\ihat_P^*(k_P)$ is a quasi-isomorphism.
\end{cor*}

\subsection{Proof of Theorem
 \ref{ssectRealizationLModules}\itemref{itemExistenceOfBasicRealizationFunctors}}
\label{ssectProofOfTheoremRealization}
Apply Lemma~\ref{ssectResolutionBySpecialForms},
Lemma~\ref{ssectNilmanifoldCohomologyBundle}, and Corollary
~\ref{ssectRestrictSpecialForms} with $X$ and $E$ replaced by $X_Q$ and
$E_Q$ respectively.  We obtain for all $P\le Q$ and all $E_Q\in\R(L_Q)$
natural quasi-isomorphisms
\begin{gather*}
\Asp(\Xhat_Q;\EE_Q)\tildearrow \i_{Q*}\EE_Q\ , \\
\coh_{PQ}\colon \AA\inv(\mathscr N_P^Q(\RR)';\EE_Q) \qiarrow
\HH(\n_P^Q;E_Q)\ ,
\end{gather*}
and a natural morphism
\begin{equation*}
k_{PQ}\colon \ihat_{Q*}\Asp(\Xhat_Q;\EE_Q) \longrightarrow
\ihat_{P*}\Asp(\Xhat_P;\HH(\n_P^Q;E_Q))\ ,
\end{equation*}
where $k_{PQ}(\o) = \Asp(\Xhat_P;\coh_{PQ})(\o|_{\Xhat_P})$ and
$\ihat_P^*(k_{PQ})$ is a quasi-isomorphism.  The proposition follows if we
set $\A_P(E_P) = \Asp(\Xhat_P;\EE_P)$.  \qed

\subsection{}
\label{ssectOtherRealizationsLsheaves}
There are other choices of $\{\A_P\}_{P\in\Pl}$ that yield different
incarnations of the same realization $\Sheaf_W$.
%
%
%
For example, if $E$ has a $\QQ$-structure one may set $\A_G(E)$
equal to the resolution $\i_{G*}\mathbf I(\EE)$ defined in
\cite[(27.9)]{refnGoreskyHarderMacPherson} and similarly for $\A_P$; in this way one obtains an
incarnation $\A_W$ defined over $\QQ$ for $\L_W$-modules defined over $\QQ$.
We will not go into further details here.

\section{Example: the Intersection Cohomology $\L$-module}
\label{sectIntersectionCohomologyLmodule}

\subsection{}
\label{ssectTruncation}
Recall that there are truncation functors on $\Complex(L_P)$ given by
\begin{align*}
\t^{\leqslant n}C\colon &\quad \dots\to C^{n-1} \to \ker d_C^n \to 0\to
\cdots\ ,\\
\t^{> n}C\colon &\quad \dots\to 0 \to C^{n+1}/\Im d_C^n \to C^{n+2}\to
\cdots\ .
\end{align*}
The natural morphisms $\t^{\leqslant n}C \to C$ and $C \to \t^{> n}C$ induce isomorphisms
%
%
\begin{equation*}
H^i(\t^{\leqslant n}C) \cong \begin{cases} H^i(C) & i\le n, \\
                                0      & i>n,
                         \end{cases}
\qquad\text{and}\qquad
H^i(\t^{> n}C) \cong \begin{cases} 0       & i\le n, \\
                        H^i(C)  & i>n
                  \end{cases}
\end{equation*}
and there are natural quasi-isomorphisms with mapping cones
\cite[\S1.7.5]{refnKashiwaraSchapira}
\begin{equation}
M(\t^{\leqslant n}C\!\to\! C) \qiarrow \t^{> n}C \qquad \text{and}
\qquad \t^{\leqslant n} C \qiarrow M(C\!\to\! \t^{> n}C)[-1]
\label{eqnTruncationIsCone}
\end{equation}
which suggest how to extend these functors to $\R(\L)$.

\subsection{}
\label{ssectLocalDegreeTruncationLmodule}
For $W$ admissible and $Q\in\Pl(W)$ consider the composition $\M\to \i_{Q*}\i_Q^*\M \to \i_{Q*}\t^{> n}\i_Q^*\M$ of natural maps  of $\L_W$-modules (the first being the adjunction morphism implied by \eqref{eqnAdjointRelations}).  Define a functor
$\t_Q^{\leqslant n}$ on $\R(\L_W)$ by
\begin{equation*}
\t_Q^{\leqslant n}\M = M(\M\to \i_{Q*}\t^{> n}\i_Q^*\M)[-1]\ ;
\end{equation*}
from \S\ref{ssectMappingCone} there is a natural morphism $\t_Q^{\leqslant n}\M \to \M$.  Note that
\begin{equation}
  \label{eqnLocalDegreeTruncationTrivial}
  \text{$\t_Q^{\leqslant n}\M = \M$ if $\t^{>n}\i_Q^*\M=0$}\ .
\end{equation}

By \eqref{eqnLocalCohomologyOfMappingCone}, $\i_P^* \t_Q^{\leqslant n}\M = M(\i_P^* \M \to \i_P^*\i_{Q*}\t^{> n}\i_Q^*\M)[-1]$,
and from \eqref{eqnTruncationIsCone} one sees there is a natural quasi-isomorphism
\begin{equation}
\label{eqnLocalCohomLocalTruncation}
\i_P^*\circ \t_Q^{\leqslant n} \cong 
                \begin{cases} \t^{\leqslant n} \circ \i_P^* & \text{if $P=Q$}\ , \\
                       \i_P^*                        & \text{if $P\nleq Q$}\ .
                \end{cases}
\end{equation}
If $P<Q$ a simple formula like \eqref{eqnLocalCohomLocalTruncation} does not
hold.  Instead there is a short exact sequence
\begin{equation*}
0 \to H(\n_P^Q;(\t^{> n}\i_Q^*\M)[-1]) \to \i_P^* \t_Q^{\leqslant n} \M
\to \i_P^*\M \to 0
\end{equation*}
and, since $H(\n_P^Q;\cdot)$ is exact and degree nondecreasing, the resulting long exact sequence implies there is a
quasi-isomorphism
\begin{equation}
\t^{\leqslant n} \circ \i_P^* \circ \t_Q^{\leqslant n} \qiarrow \t^{\leqslant n}
\circ \i_P^*\ .
\label{eqnTruncatedLocalCohomLocalTruncation}
\end{equation}
(This also holds if $P\nless Q$ but \eqref{eqnLocalCohomLocalTruncation} is stronger in that case.)
It is also clear that
\begin{equation}
\t_P^{\leqslant n}\circ \t_Q^{\leqslant n}=\t_Q^{\leqslant n}\circ
\t_P^{\leqslant n} \quad \text{if $P\nless Q$ and $Q\nless P$.}
\label{eqnCommuteLocalTruncation}
\end{equation}

\subsection{Degree truncation of $\L$-modules}
\label{ssectDegreeTruncationLmodules}
Let $W$ be admissible.  The functor $\t^{\leqslant n}$ on $\R(\L_W)$ is defined as
\begin{equation}\label{eqnDegreeTruncationLmodules}
\t^{\leqslant n} = \t_{Q_1}^{\leqslant n}\circ \t_{Q_2}^{\leqslant
n} \circ \dots \circ\t_{Q_N}^{\leqslant n}\ ,
\end{equation}
where we write $\Pl(W)=\{Q_1,\dots, Q_N\}$ so that $Q_i\le Q_j$ implies
$i\le j$; the functor is independent of the choice of such an ordering on
$\Pl(W)$ by \eqref{eqnCommuteLocalTruncation}.  There is a natural morphism
$\t^{\leqslant n}\M\to \M$.
It follows from repeated use of \eqref{eqnLocalCohomLocalTruncation} and \eqref{eqnTruncatedLocalCohomLocalTruncation} that there is a quasi-isomorphism
\begin{equation}\label{eqnTruncationAndPullbackCommute}
\i_Q^*\circ \t^{\leqslant n} \qiarrow \t^{\leqslant n} \circ \i_Q^*\qquad \text{ for $Q\in\Pl(W)$}\ .
\end{equation}

\begin{lem*}
Let $\Sheaf_{W}\colon \R(\L_W)\to
\Derived_{\X}(W)$ be a realization functor from Theorem
~\textup{\ref{ssectRealizationLModules}}.  Then $\Sheaf_{W}\circ
\t^{\leqslant n} \cong \t^{\leqslant n} \circ\Sheaf_{W}$.
\end{lem*}

\begin{proof}
It suffices to show that there is a natural isomorphism on each stratum $X_Q$.  For the left-hand side we calculate
\begin{equation}\label{eqnLHSTruncationLemma}
  \i_Q^*\Sheaf_{W}(\t^{\leqslant n}\M) = \Sheaf_{X_Q}(\i_Q^*\t^{\leqslant n}\M) \cong \Sheaf_{X_Q}(\t^{\leqslant n}\i_Q^*\M) = \Sheaf_{X_Q}(\t_Q^{\leqslant n}\i_Q^*\M)
\end{equation}
by Theorem ~\ref{ssectRealizationLModules}\ref{itemExistenceOfRealizationFunctors} and \eqref{eqnTruncationAndPullbackCommute} and for the right-hand side we similarly calculate
\begin{equation}\label{eqnRHSTruncationLemma}
  \i_Q^*\t^{\leqslant n}\Sheaf_{W}(\M) = \t^{\leqslant n}\i_Q^*\Sheaf_{W}(\M) = \t^{\leqslant n}\Sheaf_{X_Q}(\i_Q^*\M)\ .
\end{equation}
There is a quasi-isomorphism $\t^{\leqslant n}\Sheaf_{X_Q}(\i_Q^*\M)\qiarrow  \Sheaf_{X_Q}(\t_Q^{\leqslant n}\i_Q^*\M)$ by  \eqref{eqnTruncationIsCone}.
\end{proof}

\begin{rem*}
The isomorphism does not hold on the level of complexes.  The point is
that $\t^{\leqslant n}$ on $\L_W$-modules is defined ``externally'' via
mapping cones, whereas  $\t^{\leqslant n}$ on complexes of sheaves is defined
``internally'' as a sub-object.
\end{rem*}

\subsection{Intersection cohomology $\L$-module}
\label{ssectIntersectionCohomology}
Let $E$ be an algebraic $G$-module and let $p$ be an ordinary perversity, that
is, a function $p\colon \{2,\dots,\dim X\}\to \ZZ$ such that $p(2)=0$, and
$p(k+1) = p(k)$ or $p(k)+1$.  We will define the \emph{intersection cohomology
$\L_{\Xhat}$-module $\IpC(E)$ \textup(with perversity $p$ and coefficients
$E$\textup)} inductively.

Write the elements of $\Pl$ as $\{P_1=G,P_2,\dots,P_N\}$  where $(\codim_{\Xhat} X_{P_i})_{i=1}^N$ is nondecreasing.  If $P_k>P_j$ then  $\codim_{\Xhat} X_{P_k} < \codim_{\Xhat} X_{P_j}$ so $k<j$.  Thus the union $W_i = X_{P_1}\sqcup X_{P_2} \sqcup \dots \sqcup X_{P_i}$ is open and hence is admissible.
Now set $\M_{X_1} = E$ and, assuming that $\M_{W_{i}}$ has been defined, set
\begin{equation}
  \label{eqnIntersectionCohomologyLmodule}
  \M_{W_{i+1}} = \t^{\leqslant p(\codim_{\Xhat}X_{P_{i+1}})}\j_{P_{i+1}*} \M_{W_{i}}\ .
\end{equation}
Finally set $\IpC(E)=\M_{W_N}$.  There is a natural morphism of bounded $\L$-modules $\IpC(E) \to \i_{G*}E$.

\begin{lem*}
  The alternate inductive definition with $\widetilde\M_{X_1}=E$ and 
  \begin{equation*}
    \widetilde\M_{W_{i+1}} = \t_{P_{i+1}}^{\leqslant p(\codim_{\Xhat}X_{P_{i+1}})}\j_{P_{i+1}*} \widetilde\M_{W_{i}} \tag{$\widetilde{\ref{eqnIntersectionCohomologyLmodule}}$}
  \end{equation*}
yields $\widetilde \M_{W_N}=\IpC(E)$ as well.
\end{lem*}

\begin{proof}
  We inductively prove $\M_{W_i}=\widetilde \M_{W_i}$ and
{\newtagform{subscript}{(}{)$_i$}\usetagform{subscript}
\begin{equation}\label{eqnInductiveHypothesisIntersectionCohomologyLmodule}
  \text{$\t^{>  p(\codim_{\Xhat} X_{P_j})} i_{P_j}^*\M_{W_{i}}=0$ for $j\le i$}
\end{equation}
}
are true for all $i$. The case $i=0$ is trivial; we assume the result for $i$ and calculate
\begin{align*}
  \M_{W_{i+1}} &= \t_{P_{i+1}}^{\leqslant p(\codim_{\Xhat}X_{P_{i+1}})}\cdots\t_{P_1}^{\leqslant p(\codim_{\Xhat}X_{P_{i+1}})}\j_{P_{i+1}*} \M_{W_{i}} &&  \text{by \eqref{eqnIntersectionCohomologyLmodule} and \eqref{eqnDegreeTruncationLmodules}}\\
&\dots \\
  &= \t_{P_{i+1}}^{\leqslant p(\codim_{\Xhat}X_{P_{i+1}})}\t_{P_i}^{\leqslant p(\codim_{\Xhat}X_{P_{i+1}})}\j_{P_{i+1}*} \M_{W_{i}}  \\
  &= \t_{P_{i+1}}^{\leqslant p(\codim_{\Xhat}X_{P_{i+1}})}\j_{P_{i+1}*} \M_{W_{i}} = \widetilde \M_{W_{i+1}}\ .
\end{align*}
The last lines follow from repeated applications of
\eqref{eqnLocalDegreeTruncationTrivial}, \eqref{eqnInductiveHypothesisIntersectionCohomologyLmodule}$_i$,  and the fact that $p$ is nondecreasing.  Property \eqref{eqnInductiveHypothesisIntersectionCohomologyLmodule}$_{i+1}$ follows from \eqref{eqnTruncationAndPullbackCommute} and \eqref{eqnInductiveHypothesisIntersectionCohomologyLmodule}$_{i}$.
\end{proof}

Let $\EE= D\times_\G E$ be the locally constant sheaf on $X$ induced by $E$
and let $\IpC(\Xhat;\EE)$ be the corresponding intersection cohomology sheaf
\cite{refnGoreskyMacPhersonIHTwo} in $\Derived^b_{\X}(\Xhat)$.  The following
proposition is immediate from the definitions and the fact that
$\Sheaf_{\Xhat}$ commutes with the standard functors and with truncation
(Theorem
~\ref{ssectRealizationLModules}\itemref{itemExistenceOfRealizationFunctors}
and Lemma ~\ref{ssectDegreeTruncationLmodules}).

\begin{prop*}
There is a natural isomorphism $\Sheaf_{\Xhat}(\IpC(E)) \cong
\IpC(\Xhat;\EE)$.
\end{prop*}

\subsection{Local intersection cohomology}
\label{ssectLocalIntersectionCohomology}
We wish to give a formula for the local cohomology at $P$ of $\IpC(E)$.
There are two ingredients to the formula.  The first is the nilpotent
cohomology $H(\n_P;E)$ and its decomposition $ \bigoplus_{w\in W_P}
H(\n_P;E)_w $ from Kostant's theorem as in \S\ref{ssectKostantsTheorem}.
The second ingredient is a certain combinatorial invariant
$I_{p_w}H(c(|\D_P|))$, which we now describe.

Let $|\D_P|$ denote the geometric simplex with vertices indexed by the
elements of $\D_P$.  It has a natural stratification given by the
decomposition into open faces.  The strata correspond to parabolic
$\QQ$-subgroups $Q > P$ (that is, to nonempty subsets $\D_P^Q\subset
\D_P$) and are denoted $|\D_P^Q|^\circ$.  We give the (open) cone $c(|\D_P|)=|D_P| \times [0,1)/(x,0)\sim(x',0)$ the
induced stratification---the strata consist of the cones on
$|\D_P^Q|^\circ$ for $Q >P$ together with the cone point (indexed by
$P$).  Any constructible subset (that is, any union of open faces) of
$|\D_P|$ or $c(|\D_P|)$ is given the induced stratification.

Given an ordinary perversity $p$ and an element $w\in W$, define  $p_w\colon 
\Pl \to \ZZ$ by
\begin{equation*}
p_w(Q) = p(\dim \n_Q + \#\D_Q) - \l(w_Q),
\end{equation*}
where $w_Q$ denotes the projection of $w$ to the second factor of $W=W^Q
W_Q$.  This defines an integer-valued function on the strata of
$c(|\D_P|)$ or $|\D_P|$ (a perversity in the sense of
\cite{refnBeilinsonBernsteinDeligne}).  Let $I_{p_w}H(U)$ denote the
corresponding intersection cohomology group over $\ZZ$ for any open
constructible subset $U$ of $|\D_P|$ or $c(|\D_P|)$.  (Unlike the
intersection cohomology for an ordinary perversity $p$, the group here
depends on a choice of stratification; we always use the stratification
fixed above.)  Note that if $P\le Q$, then $w_Q = (w_P)_Q$ and so $p_w$ and
$p_{w_P}$ agree on the strata of $|\D_P|$.

\begin{prop*}
$H(\i_P^*\IpC(E)) \cong \bigoplus_{w\in W_P} H(\n_P;E)_w\otimes
I_{p_w}H(c(|\D_P|))$ where $L_P$ acts trivially on the second factor.
\end{prop*}
\begin{proof}
There is a Mayer-Vietoris spectral sequence that abuts to the link
cohomology $H(\i_P^*\j_{P*}\j_P^*\IpC(E))$ (Lemma ~
\ref{ssectRelativeLocalCohomologySupportsSS} with $Q'=G$ and $Q=P$) with
\begin{equation*}
E_1^{p,\cdot}= \bigoplus_{\substack{P<\tilde P \\ \#\D_P^{\tilde P}=p+1}} H( \n_P^{\tilde P}; H(\i_{\tilde P}^*\IpC(E)))\ .
\end{equation*}
By induction we may assume the proposition is true for ${\tilde P}>P$ and thus we
can compute using \S\ref{ssectKostantDegeneration} that
\begin{equation}
\begin{split}
E_1^{p,\cdot} &\cong \bigoplus_{\#\D_P^{\tilde P}=p+1}
          H( \n_P^{\tilde P};\bigoplus_{w_{\tilde P}\in W_{\tilde P}} H(\n_{\tilde P};E)_{w_{\tilde P}}\otimes 
          I_{p_{w_{\tilde P}}}H( c(|\D_{\tilde P}|) ) )\\
          &\cong \bigoplus_{\#\D_P^{\tilde P}=p+1} \bigoplus_{w_{\tilde P}\in 
          W_{\tilde P}} \bigoplus_{w^{\tilde P}_P\in W^{\tilde P}_P}
          H(\n_P^{\tilde P}; H(\n_{\tilde P};E)_{w_{\tilde P}})_{w^{\tilde P}_P} \otimes 
          I_{p_{w_{\tilde P}}}H( c(|\D_{\tilde P}|) ) \\
          &\cong \bigoplus_{w\in W_P} H(\n_P;E)_w \otimes
          \bigoplus_{\#\D_P^{\tilde P}=p+1} I_{p_w}H( c(|\D_{\tilde P}|) )\ .
\end{split}  \label{eqno}
\end{equation}

On the other hand, there is a topological Mayer-Vietoris spectral sequence
abutting to $I_{p_w}H(|\D_P|)$.  Let $(U_\al)_{\al\in\D_P}$ be the cover of $|\D_P|$ by the open stars of the vertices and let $\ast$ denote the \emph{join} of two simplices.  Then for
${\tilde P}> P$ the intersection
\begin{equation*}
  U_{\tilde P}=\bigcap_{\al\in \D_P^{\tilde P}} U_\al
  = \bigcup_{Q\in[P,\tilde P]}|\D_P^Q|^\circ
  = |\D_P^{\tilde P}|^\circ \ast |\D_P^{(\tilde P,P)}|
  = |\D_P^{\tilde P}|^\circ \times c(|\D_P^{(\tilde P,P)}|)
\end{equation*}
is the open star of $|\D_P^{\tilde P}|^\circ$.  We use the
homeomorphism $c(|\D_P^{(\tilde P,P)}|)\cong c(|\D_P^{\tilde P}|)$
to see that the Mayer-Vietoris spectral sequence is
\begin{equation}
\widetilde E_1^{p,\cdot} =\bigoplus_{ \#\D_P^{\tilde P}=p+1} I_{p_w}H(U_{\tilde P}) \cong 
\bigoplus_{ \#\D_P^{\tilde P}=p+1} I_{p_w}H(c(|\D_{\tilde P}|)) \Rightarrow  I_{p_w}H(|\D_P|) \ . \label{eqnp}
\end{equation}
Comparing \eqref{eqno} ~and \eqref{eqnp} we see that $E_1^{p,\cdot} \cong
\bigoplus_{w\in W_P} H(\n_P;E)_w \otimes \widetilde E_1^{p,\cdot}$.   This
isomorphism induces isomorphisms at all stages of the spectral sequences
and we find that
\begin{equation*}
H(\i_P^*\j_{P*}\j_P^*\IpC(E))\cong \bigoplus_{w\in W_P} 
H(\n_P;E)_w \otimes I_{p_w}H(|\D_P|)\ .
\end{equation*}
In order to obtain $H(\i_P^*\IpC(E))$, this link cohomology is to be
truncated in degrees greater than $p(\codim_{\Xhat} X_P) = p_w(P) + \l(w)$,
which corresponds exactly to replacing the second factor by
$I_{p_w}H(c(|\D_P|))$.
\end{proof}

\section{Example: the Weighted Cohomology $\L$-module}
\label{sectWeightCohomologyLmodule}

\subsection{}
\label{ssectWeightTruncation}
The reference for this subsection is \cite[\S9]{refnGoreskyHarderMacPherson} though we give a slightly different
description.

Let $P\in \Pl$ and let $X(S_P^G)$ denote the lattice of rational characters
on $S_P^G$.  For $\al\in\D_P$, let $R_\al\ge P$ be the maximal parabolic
$\QQ$-subgroup having type $\D_P\setminus \{\al\}$ with respect to $P$.  We
may canonically identify $X(S_{R_\al}^G)\cong \ZZ$ so that the unique
element of $\D_{R_\al}$ corresponds to a positive integer.  Consider the
map $X(S_P^G)\hookrightarrow \prod_{ \al\in\D_P} X(S_{R_\al}^G) \cong
\ZZ^{\D_P}$ defined by $\chi \mapsto (\chi|_{S_{R_\al}^G})_{\al\in\D_P}$.
The image is a sublattice of finite index; after tensoring with $\RR$ we
have an isomorphism $\phi_P\colon \sa_P^{G*} = X(S_P^G)\otimes_\ZZ \RR
\tildearrow \RR^{\D_P}$.  For $P=P_0$ minimal we omit the subscript.

Define a partial ordering on $\sa_P^{G*}$ by declaring $\zeta\ge \eta$ if
and only if $\phi_P(\zeta)\ge \phi_P(\eta)$ component-wise.  Thus if
$\lsp+\sa_P^{G*}$ denotes the real convex cone generated by all
$\al\in\D_P$, we have
\begin{equation}
\begin{aligned}
\zeta\ge \eta &\qquad\Longleftrightarrow\qquad \zeta-\eta\in
	\lsp+\sa_P^{G*} \\
 &\qquad\Longleftrightarrow\qquad \langle\zeta,\b\spcheck\rangle\ge
	\langle\eta,\b\spcheck\rangle \qquad \text{for all $\b \in \Dhat_P$.}
\end{aligned}
\label{eqnWeightOrdering}
\end{equation}

If $P\le Q$ and $\zeta \in \sa_P^{G*}$, let $\zeta_Q = \zeta|_{\sa_Q^G}$
denote the restriction.  We have a commutative diagram
\begin{equation*}
\begin{CD}
\sa_P^{G*}  @>\phi_P>> \RR^{\D_P} \\
@VVV     @VVV \\
\sa_Q^{G*} @>\phi_Q>> \RR^{\D_Q}
\end{CD}
\end{equation*}
where the left-hand vertical map is restriction and the right-hand vertical
map is the projection $\RR^{\D_P} = \RR^{\D_P^Q}\times \RR^{\D_P\setminus
\D_P^Q}\to \RR^{\D_P\setminus \D_P^Q} \cong \RR^{\D_Q}$.

For an algebraic $L_P$-module $E$ we write
\begin{equation*}
E= \bigoplus_{\chi\in X(S_P^G)} E_\chi
\end{equation*}
where $E_\chi$ is the submodule on which $S_P^G$ acts via $\chi$.  Given
$\eta \in \sa^{G*}$, there are \emph{weight truncation functors} on
$\R(L_P)$ defined by
\begin{equation*}
\t^{\geqslant \eta}E = \bigoplus_{\chi \ge \eta_P} E_\chi  \quad\text{and}\quad
\t^{\ngeqslant \eta}E = \bigoplus_{\chi \ngeq \eta_P} E_\chi \ ;
\end{equation*}
the element $\eta$ will be called a \emph{weight profile}.  (In
\cite{refnGoreskyHarderMacPherson} a weight profile is required to satisfy
$\phi(\eta)\in (\ZZ+\frac12)^{\D}$; this ensures that any truncation
functor $\t^{\geqslant \eta}$ arises from a unique weight profile $\eta$.
We will not assume this, but we will feel to replace $\eta$ by $\eta'$
provided $\t^{\geqslant \eta}$ is unchanged.)  Clearly we have a split
short exact sequence $0\to \t^{\geqslant \eta}E \to E \to \t^{\ngeqslant
\eta}E \to 0$.

\subsection{Weight truncation of $\L$-modules}
\label{ssectWeightTruncationLmodule}
Consider $P\le Q\in \Pl$ and let $E$ be an algebraic $L_Q$-module.  The
subtorus $S_Q^G\subset S_P^G$ acts on $H(\n_P^Q; E_\chi)$ via the
character $\chi$.  Thus
\begin{equation}
\t^{\geqslant \eta} H(\n_P^Q; E) \subset H(\n_P^Q; \t^{\geqslant \eta}
E)\ . \label{eqnWeightNonIncreasing}
\end{equation}

Given a weight profile $\eta$ we define local weight truncation functors
$\t_Q^{\geqslant \eta}$ on $\R(\L)$ via a mapping cone
\begin{equation}
\t_Q^{\geqslant \eta}\M = M(\M\to \i_{Q*}\t^{\ngeqslant
\eta}\i_Q^*\M)[-1]\ .
\end{equation}
These weight truncation functors satisfy properties analogous to those of the local degree truncation functors in
\S\ref{ssectLocalDegreeTruncationLmodule}.  The only change is in proving
the analogue of \eqref{eqnTruncatedLocalCohomLocalTruncation}.  Instead of
the fact that $H(\n_P^Q;\cdot)$ is degree nondecreasing, we use that it is
weight nonincreasing by \eqref{eqnWeightNonIncreasing}; we also use that
$\t^{\geqslant \eta}$ is an exact functor.

The weight truncation functor $\t^{\geqslant \eta}$ on $\R(\L)$ is
defined as a composition of these local truncation functors as in
\S\ref{ssectDegreeTruncationLmodules}.  There is a natural morphism
$\t^{\geqslant \eta}\M \to \M$ and a quasi-isomorphism
\begin{equation}
\i_Q^*\circ \t^{\geqslant \eta} \cong  \t^{\geqslant \eta}\circ \i_Q^* \ .
\label{eqnWeightTruncationPullbackToStratumCommute}
\end{equation}
As in Lemma~\ref{ssectDegreeTruncationLmodules} it follows that $\Sheaf_{\Xhat}\circ \t^{\geqslant \eta} \cong \t^{\geqslant \eta}\circ \Sheaf_{\Xhat}$.

\subsection{Weighted cohomology $\L$-module}
\label{ssectWeightedCohomologyLmodule}
Let $E$ be an algebraic $G$-module and let $\eta$ be a weight profile.  The
\emph{weighted cohomology $\L_{\Xhat}$-module $\WnC(E)$ \textup(with
weight profile $\eta$ and coefficients $E$\textup)} is defined as
\begin{equation*}
\WnC(E) = \t^{\geqslant \eta} \i_{G*}E \ .
\end{equation*}
There is a natural morphism of bounded $\L$-modules $\WnC(E) \to \i_{G*}E$.  Now let
$\WnC(\Xhat;\EE)$ be the corresponding weighted cohomology sheaf
\cite{refnGoreskyHarderMacPherson} in $\Derived_{\X}(\Xhat)$.  As in Proposition~\ref{ssectIntersectionCohomology} we can prove

\begin{prop*}
There is a natural isomorphism $\Sheaf_{\Xhat}(\WnC(E)) \cong \WnC(\Xhat;\EE)$.
\end{prop*}

\subsection{Local weighted cohomology}
\label{ssectLocalWeightedCohomology}
The formula for the local cohomology at $P$ of $\WnC(E)$ is considerably
simpler than that for $\IpC(E)$ since the former is defined by a single
truncation operation as opposed to the inductive truncation procedure of
the latter.  The following proposition is an immediate consequence of
\eqref{eqnWeightTruncationPullbackToStratumCommute} and the definitions.

\begin{prop*}
$H(\i_P^*\WnC(E)) \cong \t^{\geqslant \eta}
H(\i_P^*\i_{G*}E) \cong \t^{\geqslant \eta} H(\n_P;E)$.
\end{prop*}

\begin{rem}
In order to rephrase this in a form analogous to Proposition
~\ref{ssectLocalIntersectionCohomology}, let $E$ be irreducible and have
highest weight $\lambda$ with respect to a Cartan subalgebra of $\levi_P$
and a positive system $\Phi^+(\mathfrak g_\CC,\h_\CC)$ containing
$\Phi(\n_{P\CC},\h_{\CC})$.  Then by Kostant's theorem
\S\ref{ssectKostantsTheorem} we have
\begin{equation}
H(\i_P^*\WnC(E)) \cong \bigoplus_{w\in W_P} H(\n_P;E)_w\otimes
\t^{\geqslant \eta_P -(w(\lambda +\hsr)-\hsr)_P}\CC
\end{equation}
the second factor is either $0$ or $\CC$ with the trivial action of $L_P$.
As Rapoport pointed out to me, the second factor here depends on both $w$
and $E$, unlike $I_{p_w}H(c(|\D_P|))$ of Proposition
~\ref{ssectLocalIntersectionCohomology} which is independent of $E$.
\end{rem}

\section{Micro-support of an $\L$-module}
\label{sectMicroSupport}
In this section we define the micro-support of an $\L$-module using the local cohomology with supports functor \eqref{eqnLocalCohomologyWithSupportsFunctor}.  This is an
analogue of the notion for sheaves \cite{refnKashiwaraSchapira}; the
analogy is not precise, partially because it is most natural here to
restrict certain modules to be isomorphic to the complex conjugate of their
contragredient.

\subsection{Notation}
\label{ssectMicroSupportNotation}
Let $W$ be an admissible space possessing a unique maximal stratum $X_R$.
For $V\in\IrrRep(L_P)$ with $P\in\Pl(W)$ define parabolic $\QQ$-subgroups
$Q^{\prime W}_V\ge Q^W_V\ge P$ by
\begin{align*}
\D_{P}^{Q^W_V} &= \{\,\al\in \D_{P}^R \mid
\langle\xi_V+\hsr,\al\spcheck\rangle <0\,\}\ ,\\
\D_{P}^{Q^{\prime W}_V} &= \{\,\al\in \D_{P}^R \mid
\langle\xi_V+\hsr,\al\spcheck\rangle \le 0\,\}\ .
\end{align*}
Note that $Q_V^W$ and $Q^{\prime W}_V$ depend on $W$ through the parabolic
$\QQ$-subgroup $R$ indexing its maximal stratum.  When $W=\Xhat$ we
simply write $Q_V$ and $Q'_V$.

\begin{defn}
\label{ssectMicroSupport}
The \emph{weak micro-support} $\mS_w(\M)$ of an $\L_W$-module $\M$ is
the set of $V\in \IrrRep(L_P)$ for all $P\in\Pl(W)$ such that
\begin{enumerate}
\item $H(\i_{P}^* \ihat_{Q}^! \M)_V \neq 0$ for some
$Q\in [Q_V^W, Q^{\prime W}_V]$.
\label{itemMicroSupportQType}
\setcounter{saveenum}{\value{enumi}}
\end{enumerate}
The \emph{micro-support} $\mS(\M)$ consists of those $V\in
\mS_w(\M)$ which satisfy in addition
\begin{enumerate}
\setcounter{enumi}{\value{saveenum}}
\item $(V|_{M_P})^*\cong \overline {V|_{M_P}}$.
\label{itemMicroSupportConjugateSelfcontragedient}
\end{enumerate}
The \emph{essential micro-support} $\emS(\M)$ is the set of
those $V\in\mS(\M)$ such that
\begin{equation}
H(\i_{P}^* \ihat_{{Q_V^W}}^! \M)_V \longrightarrow H(\i_{P}^*
\ihat_{Q^{\prime W}_V}^! \M)_V
\label{eqnMicroNonVanishing}
\end{equation}
is nonzero.
\end{defn}

We will discuss
condition \itemref{itemMicroSupportConjugateSelfcontragedient} in \S\ref{sectConjugateSelfContragredientFundamentalWeyl}; we will present equivalent conditions to \itemref{itemMicroSupportQType} in \S\ref{sectAlternateMicroSupport} and clarify the
relationship between $\mS(\M)$ and $\emS(\M)$.

\section{Conjugate Self-contragredient Modules and Fundamental Weyl Elements}
\label{sectConjugateSelfContragredientFundamentalWeyl}
We recall work of Borel and Casselman
\cite{refnBorelCasselman} that gives a condition in terms of highest weights for Definition ~
\ref{ssectMicroSupport}\itemref{itemMicroSupportConjugateSelfcontragedient} to hold.  We also give a relation between this condition on an irreducible component $V\subset H(\n_P;E)$ and on $E$ itself.

\subsection{Notation}
Let $P$ be a parabolic $\RR$-subgroup and $\h$ a Cartan subalgebra of $\levi_P$.  For a parabolic $\RR$-subgroup $R\supset P$ we identify $\h$ with a Cartan subalgebra for $\mathfrak r$ by a lift.  Let $\Phi^P$ denote the root system $\Phi(\levi_{P\CC},\h_\CC)$ and let $\Phi$ denote the root system $\Phi(\mathfrak g_{\CC},\h_\CC)$.  Let $\Phi^{P+}\subset \Phi^P$ be a positive system.

\subsection{The involution $\boldsymbol\t_{\mathbf P}$}
\label{ssectDefineTau}
Let $V$ be an irreducible algebraic
$L_P$-module.  We say that $V$ is \emph{conjugate
self-contragredient} if
\begin{equation*}
V^*\cong \overline V
\end{equation*}
as representations of $L_P$.  Borel and Casselman
\cite[\S1]{refnBorelCasselman} construct an involution $\t_P\colon \h^*_\CC\to
\h^*_\CC$ that sends the highest weight of $V$ to that of $\overline {V^*}$.  Thus if $\mu$ is the highest weight of $V$, then
\begin{equation}
  \label{eqnTauPProperty}
\text{$V$  is conjugate
self-contragredient if and only if $\t_P(\mu) = \mu$} \ .
\end{equation}

To define $\t_P$, let $Z_P$ denote the center of $L_P$ and decompose $\mathfrak z_P = \mathfrak z_{P,d} + \mathfrak z_{P,u}$ according to the maximal $\RR$-split and $\RR$-compact subtori of $Z_P$.  Then let $\h^*_\RR$ be the real form of $\h^*_\CC$ spanned by the roots $\Phi^P$ together with $\mathfrak z_{P,d}^* + i\, \mathfrak z_{P,u}^*$; this ensures $\h^*_\RR$ contains the differentials of rational characters of $Z_P$.  Let $c$ denote the automorphism of  $\Phi^P$ defined by complex conjugation with respect to $\h^*$.  Let $s^{P}$ and $s^c$ be the unique
elements of the Weyl group $W^P$ of $\levi_P$ satisfying
\begin{equation}
  \label{eqnTauPConditions}
-s^{P}\Phi^{P+} =\Phi^{P+} \quad\text{and} \quad 
s^cc\,\Phi^{P+}=\Phi^{P+}\ .
\end{equation}
Note $-s^P$ transforms the highest weight of $V$ to that of $V^*$, while $s^c c$ transforms the highest weight of $V$ to
that $\overline V$.  Thus the desired
involution on $\h^*_\RR$ satisfying \eqref{eqnTauPProperty} is \cite[\S\S1.1, 1.2]{refnBorelCasselman} 
\begin{equation}
  \label{eqnTauPDefinition}
  \t_P = (- s^{P}) \circ (s^c c) \ .
  %
  %
  %
\end{equation}
We extend $\t_P$ $\CC$-linearly to $\h^*_\CC$ and set $\t'_P =-\t_P$.  Note that $\t_P$ leaves $\Phi^{P+}$ invariant.  It also preserves $\Phi^R$ for all $R\supset P$ and, in particular, $\Phi$.
%
%

\begin{lem}
  \label{lemTauPDependanceOnPositiveSystem}
  If $\Phitilde^{P+}=\widetilde w\,\Phi^{P+}$ for $\widetilde w\in W^P$, then the corresponding involution $\widetilde \t_P$ satisfies
  \begin{equation*}
    \widetilde \t_P = \widetilde w \t_P \widetilde w^{-1} \ .
  \end{equation*}
\end{lem}
\begin{proof}
  By definition $\widetilde\t_P = (-\tilde s^P)\circ (\tilde s^c c)$ where 
  $\tilde s^P$ and $\tilde s^c$ are the unique elements of $W^P$ satisfying
  $-\tilde s^{P}\Phitilde ^{P+} =\Phitilde^{P+}$ and $\tilde s^cc\,\Phitilde ^{P+}=\Phitilde ^{P+}$.  If we re-write these conditions in terms of $\Phi^{P+}$ we see that
  \begin{equation}
    \label{eqnTildeTauPConditions}
    -\widetilde w^{-1}\tilde s^{P}\widetilde w \,\Phi^{P+} =\Phi^{P+} \quad\text{and}\quad \widetilde w^{-1}\tilde s^cc \, \widetilde w \Phi^{P+}= \Phi^{P+}\ .
  \end{equation}
Since $W^P$ acts simply transitively on positive systems we deduce from the first parts of \eqref{eqnTauPConditions} and  \eqref{eqnTildeTauPConditions} that $\widetilde w^{-1}(-\tilde s^P)\widetilde w= -s^P$.  And since $W^P$ is a normal subgroup of the automorphisms of $\Phi^P$ we have $c\widetilde wc^{-1}\in W^P$.  Then we can conclude from the second parts of \eqref{eqnTauPConditions} and  \eqref{eqnTildeTauPConditions} that $ \widetilde w^{-1}(\tilde s^cc)\widetilde w=s^cc$ as desired.
\end{proof}

\subsection{Equivalent conditions to Definition
~\ref{ssectMicroSupport}\itemref{itemMicroSupportConjugateSelfcontragedient}}
\label{ssectConjugateSelfContragredient}

\begin{prop*}
  Let $P$ be an parabolic $\QQ$-subgroup and let $V\in \IrrRep(L_P)$ have highest weight $\mu$.  Decompose  $\h=\hb_P+\sa_P$ where $\hb_P$ is a Cartan subalgebra of $\m_P$.  Let  $\hb_{P\CC} = \hb_{P,1}+\hb_{P,-1}$ be the decomposition into the $\pm1$ eigenspaces of $\t_P$.  The following conditions are equivalent\textup:
\begin{enumerate}
\item $(V|_{M_P})^*\cong \overline{V|_{M_P}}$, that is, $V|_{M_P}$ is
conjugate self-contragredient, \label{itemVM_PConjugateSelfcontragradient}
\item $\t_P(\mu|_{\hb_P}) = \mu|_{\hb_P}$, \label{itemTauInvariant}
\item $\mu|_{\hb_{P,-1}}=0$. \label{itemEquivalentExpression}
\item $\t_P((\mu+\hsr)|_{\hb_P}) = (\mu+\hsr)|_{\hb_P}$, \label{itemTauPlusDeltaInvariant} 
\end{enumerate}
\end{prop*}

\begin{proof}
Since $\t_P\vert_{\sa_P}=-1$ and $\t_P(\hb_{P\CC})\subset \hb_{P\CC}$ the equivalence of \itemref{itemVM_PConjugateSelfcontragradient}--\itemref{itemEquivalentExpression} follows from  \eqref{eqnTauPProperty}.  Now recall from \S\ref{sssectHalfSumOfPositiveRoots} that $\hsr=\hsr^P+\hsr_P$.    Since  $2\hsr^P$ is the sum of
roots in $\Phi^{P+}$ and $\Phi^{P+}$ is $\t_P$-invariant we see $\t_P(\hsr^P) = \hsr^P$.  And since $2\hsr_P$ is the sum of roots in $\Phi(\n_{P\CC},\h_\CC)$ then $\hsr_P\in X(L_P)$ so $\hsr_P|_{\hb_P}=0$.  Thus $\t_P(\hsr|_{\hb_{P}}) = \hsr|_{\hb_{P}}$
which shows \itemref{itemTauPlusDeltaInvariant} is equivalent to the rest.
\end{proof}

%
%

\subsection{Fundamental Cartan subalgebras}
\label{ssectFundamentalCSA}
A convenient situation is when $\t_P=\theta$ for some Cartan involution
$\theta$ of $\levi_P$ and hence $\t'_P$ acts on roots as complex
conjugation.  For a fixed $P$, this can always be arranged by choosing $\h$
and $\Phi^{P+}$ appropriately.  Namely, $\t_P=\theta$ if
and only if $\h$ is a $\theta$-stable Cartan subalgebra of $\levi_P$ and
$\Phi^{P+}$ is $\theta$-stable.  A $\theta$-stable
positive system exists if and only if there exists a regular element in the
compact part of $\h$ and hence if and only if the compact part of $\h$ has
maximal dimension.  A Cartan subalgebra with this property for some Cartan
involution is called \emph{fundamental} (or \emph{maximally
compact}) for $\levi_P$; fundamental Cartan subalgebras always exist.

In the case that $\t_P=\theta$ the decomposition $\hb_{P\CC} =
\hb_{P,1}+\hb_{P,-1}$ is defined over $\RR$ and corresponds to the Cartan
decomposition into compact and $\RR$-split parts; the condition that
$V|_{M_P}$ is conjugate self-contragredient is equivalent to $\mu$
vanishing on the $\RR$-split part of $\hb_P$.

\subsection{Equal-rank groups}
\label{ssectEqualRankGroups}
A reductive algebraic $\RR$-group $Z$ will be called
\emph{equal-rank} if $\CCrank Z = \rank K_Z$, where $K_Z$ is a maximal
compact subgroup of $Z(\RR)$.  This holds if and only if a fundamental
Cartan subalgebra is compact and thus if and only if every algebraic
representation of $Z$ is conjugate self-contragredient.

\subsection{Fundamental parabolic subgroups}
\label{ssectFundamentalParabolic}
A Cartan subalgebra $\h$ for $\levi_P$ may be identified with a
Cartan subalgebra for $\mathfrak g$ via a lift; we will always assume this
has been done.  Suppose a parabolic $\RR$-subgroup $P$ has the property that a fundamental Cartan subalgebra $\h$ for $\levi_P$ is also fundamental
for $\mathfrak g$; in this case a parabolic $\RR$-subgroup $Q\supset P$ has the same property.   Such parabolic $\RR$-subgroups exist and a minimal one
is called \emph{fundamental}.  If $P_0$ is a fundamental parabolic
$\RR$-subgroup, then $L_{P_0}/\lsb\RR S_{P_0}$ is equal-rank; $P_0\neq G$
if and only if $G/\lsb\RR S_G$ is not equal-rank.  Various conditions
equivalent to $P$ containing a fundamental parabolic $\RR$-subgroup are
given in \cite[\S1.8]{refnBorelCasselman}; here are others we need.

\begin{lem*}
Let $P$ be a parabolic $\RR$-subgroup of $G$ and fix a fundamental $\theta$-stable Cartan subalgebra
$\h$ of $\levi_P$.  Let $\Phi^{P+}$ be a positive system
for $\levi_{P\CC}$ and extend it to a positive system $\Phi^+$ for $\mathfrak g_\CC$ containing $\Phi(\n_{P\CC},\h_\CC)$.  Let $\t_P$ and $\t_G$ be defined as in \textup{\S\ref{ssectDefineTau}} with respect
to $\Phi^{P+}$ and $\Phi^+$ respectively.  Then the following
conditions are equivalent\textup:
\begin{enumerate}
\item $P$ contains a fundamental parabolic $\RR$-subgroup\textup;
\label{itemFundamental}
\item There exists $w\in W_P$ such that $w\t_G
w^{-1}=\t_P$\textup;
\label{itemTauInvariance}
\item There exists $w\in W_P$ such that $\t_P'$ interchanges $\Phi_w$
and $\Phi(\n_{P\CC},\h_\CC)\setminus \Phi_w$\textup;
\label{itemPhiwSwap}
\item There exists $w\in W_P$ such that $w\Phi^+$ is
$\t_P$-stable.
\label{itemwPhiPositiveInvariance}
\end{enumerate}
The conditions in \itemref{itemTauInvariance},
\itemref{itemPhiwSwap} and
\itemref{itemwPhiPositiveInvariance} for a given $w$ are
equivalent.
\end{lem*}

\begin{proof}
By \eqref{eqnOrderingDependenceForW_P} and Lemma~ \ref{lemTauPDependanceOnPositiveSystem} it suffices to prove the result with $\Phi^+$ replaced by $\widetilde w\Phi^+$ for any $\widetilde w\in W^P$.  Thus by \S\ref{ssectFundamentalCSA} we may assume $\Phi^{P+}$ is $\theta$-stable and $\t_P=\theta$.

\itemref{itemFundamental}$\Longrightarrow$\itemref{itemTauInvariance}: Suppose
that $P$ contains a fundamental parabolic $\RR$-subgroup.  Then $\h$ is a
fundamental Cartan subalgebra for $\mathfrak g$ as well so there exists a
$\theta$-stable positive system $\Phitilde ^+$ for $\mathfrak g_\CC$ which
contains $\Phi^{P+}$.  Let $\widetilde \t_G$ denote the
operator of \S\ref{ssectDefineTau} with respect to $\Phitilde^+$; since
$\Phitilde^+$ is $\theta$-stable it satisfies $\widetilde \t_G=\theta=\t_P$,
while if $\Phitilde ^+=w\Phi^+$ for $w\in W$ then $\widetilde \t_G = w\t_G
w^{-1}$ by Lemma~ \ref{lemTauPDependanceOnPositiveSystem}.  Since both $\Phitilde^+$ and $\Phi^+$ contain
$\Phi^{P+}$, we must have $\Phi_w\subset \Phi(\n_{P\CC},\h_\CC)$ and hence $w\in W_P$ (see \S\ref{sssectWeyl}).

\itemref{itemTauInvariance}$\Longrightarrow$\itemref{itemPhiwSwap}: Suppose
$w\t_G w^{-1}=\t_P$ for some $w\in W_P$.  Then $\widetilde \t_G=\t_P$ where
$\widetilde \t_G$ is with respect to $\Phitilde^+=w\Phi^+$; this follows from Lemma ~\ref{lemTauPDependanceOnPositiveSystem} for $G$.  Decompose
$\Phi(\n_{P\CC},\h_\CC)=\Phi_w \cup (\Phi(\n_{P\CC},\h_\CC)\setminus
\Phi_w)$.  These subsets correspond to those roots in $\Phi(\n_{P\CC},\h_\CC)$ that are
negative (respectively positive) with respect to $\Phitilde ^+$.  Now  $\widetilde \t_G'$ maps  $\Phitilde ^+$ to $-\Phitilde ^+$ while $\t_P'$, being the identity on $\sa_P^*$, preserves $\Phi(\n_{P\CC},\h_\CC)$.  Thus $\t_P'=\widetilde \t_G'$ must interchange the two subsets.

\itemref{itemPhiwSwap}$\Longrightarrow$\itemref{itemwPhiPositiveInvariance}: If
$w$ is as in \itemref{itemPhiwSwap}, the positive system $w\Phi^+=
\Phi^{P+} \cup -\Phi_w \cup
(\Phi(\n_{P\CC},\h_\CC)\setminus \Phi_w)$ is $\t_P$-stable.

\itemref{itemwPhiPositiveInvariance}$\Longrightarrow$\itemref{itemFundamental}:
Since $w\Phi^+$ is a $\theta$-stable positive system for $\mathfrak g$, $\h$ is a fundamental Cartan subalgebra for $\mathfrak g$ by \S\ref{ssectFundamentalCSA}.
\end{proof}

\subsection{Fundamental Weyl elements}
\label{ssectFundamentalWeylElement}

\begin{defn*}
  In the setting of Lemma ~\ref{ssectFundamentalParabolic},  $w\in W_P$ satisfying $w\t_G w^{-1}=\t_P$ will be called
a \emph{fundamental Weyl element for $P$ in $G$}.  More generally,
if $P\le R$ an element $w \in W_P^R$ is fundamental for $P$ in $R$ if
$w\t_R w^{-1}=\t_P$ or, equivalently, it is fundamental for $P/N_R$ in $L_R$).
\end{defn*}

\begin{lem*}
  \leavevmode
\begin{enumerate}
\item A fundamental Weyl element $w\in W_P$ satisfies $\l(w)=\frac12\dim
\n_P$.
\label{itemFundamentalWeylHalf}
\item If $P\le R$ and $w=w^Rw_R\in W^R_P W_R=W_P$, then $w$ is
fundamental for $P$ in $G$ if and only if $w^R$ is fundamental for $P$ in
$R$ and $w_R$ is fundamental for $R$ in $G$.
\label{itemFundamentalWeylHereditary}
\end{enumerate}
\end{lem*}
\begin{proof}
Part \itemref{itemFundamentalWeylHalf} is clear from
Lemma~\ref{ssectFundamentalParabolic}\itemref{itemPhiwSwap}.  For part
\itemref{itemFundamentalWeylHereditary}, if $w\t_Gw^{-1}=\t_P$ then $w\Phi^+$ is $\t_P$-stable by  Lemma
~\ref{ssectFundamentalParabolic}\itemref{itemwPhiPositiveInvariance}.  Now note that $w\Phi^+\cap \Phi^R = w^R\Phi^{R+}$%
and $\t_P$ leaves $\Phi^R$ invariant.  Thus  $w\Phi^+$ being $\t_P$-stable implies $ w^R\Phi^{R+}$ is $\t_P$-stable. Then by  Lemma
~\ref{ssectFundamentalParabolic}\itemref{itemTauInvariance} applied to $P/N_R$ in $L_R$ it follows that $w^R\t_R(w^R)^{-1}=\t_P$.  Thus  $w\t_Gw^{-1}=w^R\t_R(w^R)^{-1}$ and hence $w_R\t_G w_R^{-1}=\t_R$.  The converse is clear.
\end{proof}

\begin{prop}
\label{ssectBorelCasselman}
Let $E$ be an irreducible algebraic $G$-module and let $P$ be a parabolic
$\QQ$-subgroup of $G$.  The following conditions are equivalent\textup:
\begin{enumerate}
\item $P$ contains a fundamental parabolic $\RR$-subgroup and
$(E|_{\lsp0G})^*\cong \overline{E|_{\lsp0G}}$\textup;
\label{itemFundamentalGTauInvariant}
\item There exists an irreducible algebraic $L_P$-submodule $V$ of $H(\n_P;E)$ satisfying
$(\xi_V+\hsr)|_{\sa^G_P}=0$ and $(V|_{M_P})^*\cong \overline{V|_{M_P}}$.
\label{itemSplitVanishingPTauInvariant}
\end{enumerate}
If either condition holds, $V= H^{\l(w)}(\n_P;E)_w$ for $w\in W_P$ a
fundamental Weyl element for $P$ in $G$ and
\itemref{itemSplitVanishingPTauInvariant} holds for any such $V$.
\end{prop}

\begin{proof}
  The equivalence of \itemref{itemFundamentalGTauInvariant} and
\itemref{itemSplitVanishingPTauInvariant} is proved in
\cite[3.6(iii)(iv)]{refnBorelCasselman} though with different notation.  For the final assertions adopt the setting of Lemma
~\ref{ssectFundamentalParabolic}.  Kostant's theorem says an irreducible  $L_P$-submodule $V$ of $H(\n_P;E)$ has highest weight $\mu=w(\lambda+\hsr)-\hsr$ where $\lambda$ is the highest weight of $E$ and $w\in W_P$.
Since $\t_P$ preserves the decomposition
$\hb_G = \hb_P + \sa^G_P$ and acts as $-1$ on the second factor, Proposition~ \ref{ssectConjugateSelfContragredient} implies that 
\itemref{itemSplitVanishingPTauInvariant} holds for $V$ if and only
if $w^{-1}\t_Pw((\lambda+\hsr)|_{\hb_G})=(\lambda+\hsr)|_{\hb_G}$.
Similarly \itemref{itemFundamentalGTauInvariant} implies  $\t_G((\lambda+\hsr)|_{\hb_G})=(\lambda+\hsr)|_{\hb_G}$ so $w^{-1}\t_Pw$ and $\t_G$ agree on the regular element $(\lambda+\hsr)|_{\hb_G}$.  However $w^{-1}\t_P w$ is the $\t_P$ involution for a different positive system by Lemma ~\ref{lemTauPDependanceOnPositiveSystem}.  Thus after composing with $-c$ both  $w^{-1}\t_Pw$ and $\t_G$ become  elements of $W$
and hence must be equal.  This proves that $w$ is a fundamental Weyl element for $P$ in $G$.
\end{proof}

%

\section{Alternate Forms of Definition
~\ref{ssectMicroSupport}\itemref{itemMicroSupportQType}}
\label{sectAlternateMicroSupport}
Let $W$ be an admissible space possessing a unique maximal stratum $X_R$.

\subsection{A partial ordering}
\label{ssectPartialOrdering}
Let $V\in \IrrRep(L_P)$ and $\tilde V\in \IrrRep(L_{\tilde P})$ for $P,\tilde P\in\Pl(W)$.  Set $V \preccurlyeq \tilde V$ if
the following two conditions are fulfilled:
\begin{enumerate}
\item \label{itemPartialOrderingPrerequisite} $P\le \tilde P$ and
\item \label{itemZeroKostantComponent} $V= H^{\l(w)}(\n_{P}^{\tilde P};
\tilde V)_w$ for an element $w\in W_P^{\tilde P}$.
\end{enumerate}
In this situation Kostant's theorem implies $(\xi_V+\hsr_P)\vert_{\sa_{\tilde P}}= \xi_{\tilde V}+\hsr_{\tilde P}$ since $w$ fixes $\sa_{\tilde P}$ and we define
\begin{equation*}
  [\tilde V:V]=\l(w)\ .
\end{equation*}
Set  $V \preccurlyeq_0 \tilde V$ if in addition to \itemref{itemPartialOrderingPrerequisite} and \itemref{itemZeroKostantComponent} we have
\begin{enumerate}[label=(\roman*)$_0$]
    \setcounter{enumi}{2}
\item \label{itemZero} $(\xi_V+\hsr_P)|_{\sa_P^{\tilde P}} = 0$.
\end{enumerate}
From \S\ref{ssectKostantDegeneration} it is easy to see that these are 
partial orders.  (In \S\ref{ssectAdditionalPartialOrdering} we consider other
partial orderings obtained by relaxing \itemref{itemZero}.)

\begin{prop}
\label{ssectAlternateMicroSupport}
Let $\M$ be an $\L_W$-module and for $P\in\Pl(W)$ let $V\in \IrrRep(L_P)$.  The following conditions on $V$ are equivalent\textup:
\begin{enumerate}
\item $H(\i_{P}^* \ihat_{Q}^! \M)_V \neq 0$ for some $Q\in [Q_V^W,
Q^{\prime W}_V]$.
\label{itemSSCondition}
\item $H(\i_P^* \i_{{\tilde P}*}\i_{\tilde P}^* \ihat_{\tilde Q}^!\M)_V
\neq 0$ for some $\tilde Q \in [Q_V^W, Q^{\prime W}_V]$ and $\tilde P =
(P,Q_V^W)\cap \tilde Q$.
\label{itemFaryTermCondition}
\item
\label{itemBasicSSCondition}
$H(\i_{\tilde P}^* \ihat_{{Q^W_{\tilde V}}}^!  \M)_{\tilde V} \neq 0$ for some irreducible $L_{\tilde P}$-module $\tilde V \succcurlyeq_0 V$.
\item
\label{itemessSSCondition}
There exists an irreducible $L_{\tilde{\tilde P}}$-module $\tilde{\tilde V}\succcurlyeq_0
V$ such that
\begin{equation*}
\Im\bigl ( H(\i_{\tilde{\tilde P}}^* \ihat_{{Q^W_{\tilde{\tilde V}}}}^!
\M)_{\tilde{\tilde V}} \longrightarrow H(\i_{\tilde{\tilde P}}^*
\ihat_{Q^{\prime W}_{\tilde{\tilde V}}}^!  \M)_{\tilde{\tilde V}} \bigr)
\neq 0\ .
\end{equation*}
\end{enumerate}
It can be arranged that the $\tilde P$ in
\itemref{itemFaryTermCondition} and
\itemref{itemBasicSSCondition} are the same and that
$\tilde{\tilde P}\ge \tilde P$.
\end{prop}

\begin{proof}
\itemref{itemSSCondition}$\Longleftrightarrow$\itemref{itemFaryTermCondition}:
Assume \itemref{itemSSCondition} holds and consider the long exact sequence
\begin{equation}
\cdots \longrightarrow
H^j(\i_P^* \ihat_{Q_V^W}^!  \M)_V \xrightarrow{\b}
H^j(\i_P^* \ihat_{Q}^! \M)_V \longrightarrow
H^{j}(\i_P^* \jhat_{Q_V^W*}\jhat_{Q_V^W}^*\ihat_{Q}^! \M)_V \longrightarrow
\cdots
\label{eqnCompareTypes}
\end{equation}
as in \eqref{eqnLongCompareLocalCohomologyWithSupports}.  If $\b$ is
nonzero, then $H^j(\i_{ P}^* \ihat_{{Q_V^W}}^!  \M)_{ V}\neq 0$ and
\itemref{itemFaryTermCondition} holds in degree $j$ with $\tilde Q = Q_V^W$ and
$\tilde P = P$.  Otherwise the last term of \eqref{eqnCompareTypes} is
nonzero which implies by the Fary spectral sequence of
Lemma~\ref{ssectRelativeLocalCohomologySupportsSS} and equation
\eqref{eqnPushForwardType} that \itemref{itemFaryTermCondition} holds in degree $j$
for some $\tilde P > P$ and $\tilde Q=\tilde P\vee Q_V^W\in (Q_V^W, Q]$.

Conversely suppose
\itemref{itemFaryTermCondition} holds for some $\tilde Q$ and $\tilde P$ and consider
the analogous long exact sequence
\begin{equation*}
\cdots \longrightarrow
H^j(\i_P^* \ihat_{Q_V^W}^!  \M)_V \xrightarrow{\b}
H^j(\i_P^* \ihat_{\tilde Q}^! \M)_V \longrightarrow
H^{j}(\i_P^* \jhat_{Q_V^W*}\jhat_{Q_V^W}^*\ihat_{\tilde Q}^! \M)_V \longrightarrow
\cdots \ .
\end{equation*}
in which $Q$ has been replaced by $\tilde Q$.  By
\itemref{itemFaryTermCondition} and \eqref{eqnPushForwardType} the Fary spectral sequence for  the
last term has $E_1$ nonzero at the top level $\tilde P=P'$.  Thus either the last
term is nonzero, in which case \itemref{itemSSCondition} holds for
$Q_V^W$ or $\tilde Q$, or else the last term is zero.  In the latter case $E_1$ must also be nonzero at a lower level $\tilde P<P'$ to cancel the top level.  Thus \itemref{itemFaryTermCondition} holds also for some smaller $\tilde Q$ and $\tilde P$ and we may repeat the argument.

\itemref{itemFaryTermCondition}$\Longleftrightarrow$\itemref{itemBasicSSCondition}:
Assume that \itemref{itemFaryTermCondition} holds in degree $j$.  Then by
\eqref{eqnPushForwardType} we see that $H^{j-[\tilde V:V]}(\i_{\tilde P}^*\ihat_{\tilde
Q}^!\M)_{\tilde V}\neq 0$ for some irreducible $L_{\tilde P}$-module
$\tilde V$ such that $V$ is an irreducible constituent of $H(\n_{P}^{\tilde
P}; \tilde V)$.  By the definition of $\tilde P$ and $\tilde Q\le Q_V^{\prime W}$ we see that
$(\xi_V+\hsr_P)|_{\sa_P^{\tilde P}}=0$ and so $V\preccurlyeq_0 \tilde V$.
Finally $(\xi_V+\hsr_P)|_{\sa_{\tilde P}} = \xi_{\tilde V}+\hsr_{\tilde P}$
implies that $\tilde Q = Q^W_{\tilde V}$ and thus
\itemref{itemBasicSSCondition} holds in degree $j-[\tilde V:V]$.

Conversely if \itemref{itemBasicSSCondition} holds for some $L_{\tilde P}$-module $\tilde V\succcurlyeq_0 V$ then $H(\n_P^{\tilde P};H(\i_{\tilde P}^*\ihat_{Q_{\tilde V}}^{W^!}\M))_V$ is nonzero and hence $H(\i_P^* \i_{{\tilde P}*}\i_{\tilde
P}^* \ihat_{Q^W_{\tilde V}}^!\M)_V \neq 0$ by \eqref{eqnPushForwardType}.  From $(\xi_V+\hsr_P)|_{\sa_P^{\tilde P}}=0$ and  $(\xi_V+\hsr_P)|_{\sa_{\tilde P}} = \xi_{\tilde V}+\hsr_{\tilde P}$ we can see $Q^W_{\tilde V} =
Q_V^W\vee \tilde P$, $Q^{\prime W}_{\tilde V} =
Q_V^{\prime W}\vee \tilde P$, and $\tilde P \le  Q_V^{\prime W}$.
Thus $Q^W_{\tilde V} \in [Q_V^W, Q^{\prime W}_V]$ and
\itemref{itemFaryTermCondition} holds.

\itemref{itemBasicSSCondition}$\Longleftrightarrow$\itemref{itemessSSCondition}:
Assume that \itemref{itemBasicSSCondition} holds.  In  the long
exact sequence
\begin{equation}
  \label{eqnCompareExtremeTypes}
\cdots \longrightarrow
H^{j-1}(\i_{\tilde P}^* \jhat_{Q^W_{\tilde V}*}\jhat_{Q^W_{\tilde
    V}}^*\ihat_{Q^{\prime W}_{\tilde V}}^! \M)_{\tilde V} \longrightarrow
H^j(\i_{\tilde P}^* \ihat_{Q^W_{\tilde V}}^!  \M)_{\tilde V} \xrightarrow{\b}
H^j(\i_{\tilde P}^* \ihat_{Q^{\prime W}_{\tilde V}}^! \M)_{\tilde V} \longrightarrow
\cdots
\end{equation}
either $\b$ is nonzero, in which case \itemref{itemessSSCondition} holds in degree $j$,
or the first term of the sequence is nonzero.  In the latter case the Fary spectral sequence and equation
\eqref{eqnPushForwardType} implies that
\itemref{itemFaryTermCondition} is true for some $\tilde Q \le Q_{\tilde V}^{\prime W}$ and $\tilde{\tilde P}>\tilde P$ (compare \itemref{itemSSCondition}$\Longrightarrow$\itemref{itemFaryTermCondition}).  Now apply \itemref{itemFaryTermCondition}$\Longrightarrow$\itemref{itemBasicSSCondition} to see that \itemref{itemBasicSSCondition}
holds for the new $\tilde Q$ and $\tilde{\tilde P}$.  Repeat the argument; since $\tilde{\tilde P}$ is increasing the process must end and \itemref{itemessSSCondition} holds.

As for degree, it goes down by $1$ with each iteration due to the $j-1$ in the first term of \eqref{eqnCompareExtremeTypes} resulting in a total decrease $\le \#\D_{\tilde P}^{\tilde{\tilde P}}=\dim \sa_{\tilde P}^{\tilde{\tilde P}}$.  In addition the degree goes down by $[\tilde{\tilde V}:\tilde V]$ with each use of \itemref{itemFaryTermCondition}$\Longrightarrow$\itemref{itemBasicSSCondition}.

Conversely if \itemref{itemessSSCondition} holds then $H(\i_{\tilde{\tilde P}}^* \ihat_{Q^W_{\tilde{\tilde V}}}^!\M)_{\tilde{\tilde V}}\neq0$ which establishes \itemref{itemBasicSSCondition}.
\end{proof}

Let $\M$ be an $\L_W$-module.  For $i=1,\dots,4$ let
\begin{equation}
  \label{eqnWeakAlternateMicrosupportDegreeRanges}
[c_i(\cdot;\M),d_i(\cdot;\M)]
\end{equation}
denote the smallest intervals containing the nonvanishing degrees of 
\itemref{itemSSCondition},\dots,\itemref{itemessSSCondition}
respectively.  In the parentheses we indicate the module and, if needed, 
the parabolic $\QQ$-subgroup to which we wish to restrict.  (It is possible that $c_i$ and $d_i$ can be $-\infty$ or $\infty$.)

\begin{corAlternateDegreeRanges}
  The intervals above satisfy
\begin{equation}
\label{eqnAlternateDegreeRanges}
\begin{aligned}
{[c_1(V;\M),d_1(V;\M)]} &\subset [c_2(V;\M), d_2(V;\M)] \ ,\\ 
[c_2(V,\tilde P;\M),d_2(V,\tilde P;\M)] &= [c_3(\tilde V;\M),d_3(\tilde V;\M)]
+[\tilde V:V] \ ,\\
[c_3(\tilde V;\M),d_3(\tilde V;\M)] &\subset [c_4(\tilde{\tilde
V};\M),d_4(\tilde{\tilde V};\M) + \dim \sa_{\tilde P}^{\tilde{\tilde P}}] +
[ \tilde{\tilde V}: \tilde V] \ .
\end{aligned}
\end{equation}
\end{corAlternateDegreeRanges}

 \begin{lem}
   \label{lemPartialOrdering}
If $V \preccurlyeq_0 \tilde V$ and $(V|_{M_P})^*\cong \overline{V|_{M_P}}$
then $ (\tilde V|_{M_{\tilde P}})^*\cong \overline{\tilde V|_{M_{\tilde
P}}}$, the Weyl element $w$ from \textup{\S\ref{ssectPartialOrdering}\itemref{itemZeroKostantComponent}}
is fundamental for $P$ in $\tilde P$, and $[\tilde V: V] = \tfrac12 \dim
\n_P^{\tilde P}$.
\end{lem}
\begin{proof}
  Use Proposition ~\ref{ssectBorelCasselman} and Lemma ~\ref{ssectFundamentalWeylElement}\itemref{itemFundamentalWeylHalf}.
\end{proof}

\begin{cor}
\label{ssectPartialOrderingOnMicroSupport}
If $V\in \mS(\M)$ there exists $\tilde V\in \emS(\M)$ with $V\preccurlyeq_0
\tilde V$.
If $\tilde V \in \emS(\M)$ and $V\in\IrrRep(L_P)$ for $P\in\Pl(W)$ such that $V|_{M_P}$ is conjugate
self-contragredient and $V\preccurlyeq_0 \tilde V$, then $V\in \mS(\M)$.
\end{cor}
\begin{proof}
  If $V\in\mS(\M)$ then $(V|_{M_P})^*\cong \overline{V|_{M_P}}$ and also $V$ satisfies \itemref{itemSSCondition}.  By Proposition~\ref{ssectAlternateMicroSupport} there exists  $\tilde V\succcurlyeq_0 V$ satisfying \itemref{itemessSSCondition} so $\tilde V\in \emS(\M)$ by Lemma~\ref{lemPartialOrdering}.
  If $\tilde V \in \emS(\M)$ then the implication \itemref{itemessSSCondition}$\Longrightarrow$\itemref{itemSSCondition} of Proposition~\ref{ssectAlternateMicroSupport} shows $V\in \mS(\M)$.
\end{proof}

Thus $\emS(\M)$ determines $\mS(\M)$; the converse is not true
however.

\specialsection*{Part II. A Global Vanishing Theorem for $\L$-modules}
In this part we present a vanishing theorem for the global cohomology of a bounded
$\L$-module $\M$.  The well-known local-global principle implies that if
the local cohomology of a sheaf vanishes, then the global hypercohomology
vanishes.  Our result implies that if the micro-support of a bounded $\L$-module
is empty, then the global cohomology vanishes.  More generally we give a
bound based on the micro-support for the degrees in which cohomology can be
nonzero.

\section{The Vanishing Theorem}
\label{sectVanishingTheoremLModule}

\subsection{}
\label{ssectWeightCentralizer}
Let $P\in\Pl$ and let $\h$ be a Cartan subalgebra of $\levi_P$.  Given
$\mu\in\h_\CC^*$ view $\mu$ as an element of $\levi_{P\CC}^*$ by extending
it to be zero on all root spaces.  Let $L_P(\mu)\subset L_P$ be the
stabilizer of $\mu$ under the coadjoint action.  This is a reductive
subgroup (defined over $\CC$) whose Lie algebra $\levi_P(\mu)$ is generated
by $\h$ and the root spaces for
\begin{equation}
\{\,\g\in\Phi(\levi_{P\CC},\h_{\CC}) \mid
         \langle\mu,\g\spcheck\rangle=0\,\}\ .
\label{eqnLPxiRoots}
\end{equation}
As extreme cases we have $L_P(0)=L_P$ and $L_P(\hsr)=H$, the Cartan
subgroup.

\subsection{}
\label{ssectVcentralizer}
Let $V$ be an irreducible algebraic $L_P$-module for which $(V|_{M_P})^*\cong
\overline{V|_{M_P}}$.  Let $\theta$ be a Cartan involution of $\levi_P$ and $\h=\hb_P+\sa_P$ a fundamental $\theta$-stable Cartan subalgebra of $\levi_P$.  Then for a $\theta$-stable positive system $\Phi^{P+}$ for
$\Phi(\levi_{P\CC},\h_{\CC})$ we obtain $\t_P=\theta$ as in \S\ref{ssectFundamentalCSA}.  Let $\mu\in \h^*_\CC$ be the corresponding
highest weight of $V$.  Then, by Proposition ~\ref{ssectConjugateSelfContragredient}, $\overline{\mu|_{\hb_{P}}} =
-\mu|_{\hb_{P}}$ and so $L_P(\mu)$ is defined over $\RR$.  Let
\begin{equation*}
D_P(\mu) = L_P(\mu)(\RR)/(K_P\cap L_P(\mu))A_P
\end{equation*}
denote the corresponding symmetric space (compare
\cite{refnBorelVanishingTheorem}).  The dimension of $D_P(\mu)$ can vary
depending upon the $\theta$-stable positive system $\Phi^{P+}$; we let
$D_P(V)$ denote any one of the $D_P(\mu)$ with maximal dimension.  Then
$\dim D_P(V)$ is well-defined and independent of the choice of $\h$ and
$\theta$.  It will be convenient beginning in
\S\ref{sectEqualRankBasicLemma} to allow $V$ in $D_P(V)$ to be merely
isotypical as opposed to irreducible.

\subsection{}
\label{ssectcMdMDefinition}
Let $\M$ be an $\L$-module on $\Xhat$ and $V\in\mS_w(\M)$ an $L_P$-module.  Let $[c(V;\M),d(V;\M)]$ be the smallest interval containing the nonvanishing degrees of $H(\i_{P}^* \ihat_{Q}^! \M)_V$ for any $Q\in[ Q_V,Q_V']$ as in Definition~ \ref{ssectMicroSupport}\itemref{itemMicroSupportQType}.  This is the same interval as $[c_1(V;\M),d_1(V;\M)]$ from \eqref{eqnWeakAlternateMicrosupportDegreeRanges}.

Define
\begin{align}
\label{eqncMDefn}
c(\M) &= \inf_{V\in\mS(\M)} \tfrac12(\dim D_P - \dim D_P(V)) + c(V;\M)\ ,
\text{ and} \\
\label{eqndMDefn}
d(\M) &= \sup_{V\in\mS(\M)} \tfrac12(\dim D_P + \dim D_P(V)) + d(V;\M)\ .
\end{align}
Since we are aiming toward a vanishing theorem for global cohomology the infimum and supremum here are taken over elements of $\mS(\M)$, not $\mS_w(\M)$.
Note that the shift $\tfrac12(\dim D_P \pm \dim D_P(V))$ corresponds to the
range of degrees in which $H_{(2)}(X_P;\VV)$ can be nonzero (see Theorem
~\ref{ssectRaghunathanVanishing}).

\begin{thm}
\label{ssectGlobalVanishing}
Let $\M$ be a bounded $\L$-module on $\Xhat$. Then  $H^i(\Xhat;\M)=0$ for $i
\notin [c(\M), d(\M)]$.  In particular if $\mS(\M)=\emptyset$ then $H^i(\Xhat;\M)=0$ for all $i$.
\end{thm}

The proof will appear in \S\ref{sectProofGlobalVanishing} after a number of
analytic preliminaries.  In the remainder of this section we will indicate
how $c(\M)$ and $d(\M)$ may be computed in terms of the essential
micro-support.

\begin{lem}
\label{ssectDPVInequality}
Let $V$ be an irreducible algebraic $L_P$-module for $P\in\Pl$ and assume
$V|_{M_P}$ is conjugate self-contragredient.  If $V\preccurlyeq_0
\tilde V$ then
\begin{equation}
\label{eqnComparedimDPV}
\dim D_P(V) \le \dim D_{\tilde P}(\tilde V) - \dim
\sa_P^{\tilde P}\ .
\end{equation}
Consequently
\begin{align}
\label{eqnLowerdimDPInequality}
\tfrac12\bigl(\dim D_P - \dim D_P(V)\bigr) &\ge \tfrac12\bigl(\dim
D_{\tilde P} - \dim D_{\tilde P}(\tilde V) \bigr) -\tfrac12\dim
\n_{P}^{{\tilde P}}\ , \\
\label{eqnUpperdimDPInequality}
\tfrac12\bigl(\dim D_P + \dim D_P(V)\bigr) &\le \tfrac12\bigl(\dim
D_{{\tilde P}} + \dim D_{{\tilde P}}({\tilde V}) \bigr)  -\tfrac12\dim
\n_{P}^{{\tilde P}} \\
&\qquad - \dim \sa_P^{{\tilde P}}\ . \notag
\end{align}
\end{lem}

\begin{proof}
  As in \S\ref{ssectVcentralizer} choose a $\theta$-stable Cartan subalgebra $\h=\hb_P+\sa_P$ for $\levi_P$ and a positive system $\Phi^{P+}$
such that $\t_P=\theta$.  Let
$\mu$ be the highest weight of $V$ and define
$D_P(\mu)$ to be the symmetric space of $L_P(\mu)(\RR)$.  Since $V|_{M_P}$ is conjugate self-contragradient and $V\preccurlyeq_0 \tilde V$, Lemma
~\ref{lemPartialOrdering} implies that $V=H^{\l(w)}(\n_P^{\tilde P};\tilde
V)_w$ where $w\in W_P^{\tilde P}$ is fundamental for $P$ in $\tilde P$.  Thus $w\t_{\tilde P}w^{-1} = \t_P$.  Now extend $\Phi^{P^+}$ to a positive system $\Phi^{\tilde P+}$ for $\levi_{\tilde P}$ by adjoining $\Phi(\n_{P\CC}^{\tilde P},\h_\CC)$.
By Lemma ~\ref{lemTauPDependanceOnPositiveSystem}, $\widetilde\t_{\tilde P}=
w\t_{\tilde P}w^{-1}$ is the involution of \S\ref{ssectDefineTau} defined
with respect to the positive system $w\Phi^{\tilde P+}$.  Thus $\widetilde\t_{\tilde P}=\t_P=\theta$.  If $\tilde\mu$ is the
highest weight of $\tilde V$ with respect to $\Phi^{\tilde P+}$, define
$D_{\tilde P}(w\tilde\mu)$ to be the symmetric space of $L_{\tilde
P}(w\tilde\mu)(\RR)$ as in \S\ref{ssectVcentralizer}.

By Kostant's theorem
\S\ref{ssectKostantsTheorem}, $\mu=w(\tilde\mu+\hsr)-\hsr = w\tilde\mu + (w\hsr-\hsr)$.  Now $w\hsr-\hsr$ is dominant for $\Phi^{P+}$ as is $w\tilde \mu$ since $\Phi_w\subset \Phi(\n_{P\CC}^{\tilde P},\h_{\CC})$.  Consequently if $\g$ is a root of $\levi_{P\CC}$ then
$\langle\mu,\g\spcheck\rangle=0$ implies $\langle
w\tilde\mu,\g\spcheck\rangle=0$.
Also note $\hb_{\tilde P}= \hb_P+\sa_P^{\tilde P}$. Thus
\begin{align*}
  \dim D_P(\mu) &= \#(\Phi^P\cap(\Span \mu)^\perp) + \dim \hb_P
  \\ &\le  \#(\Phi^{\tilde P}\cap(\Span w\tilde\mu)^\perp) + \dim\hb_{\tilde P} -\dim \sa_P^{\tilde P} =  \dim D_{\tilde P}(w\tilde\mu) - \dim \sa_P^{\tilde P}\ .
\end{align*}
Take the maximum of this inequality over all
possible $\Phi^{P+}$ to obtain \eqref{eqnComparedimDPV}.  One half of this is added and subtracted from $\frac12\dim D_P = \frac12(\dim D_{\tilde P} - \dim
\sa_P^{\tilde P} - \dim \n_P^{\tilde P})$ to finish the proof.
\end{proof}

\subsection{}
\label{ssectEqualityDegreeRanges}
Let $\M$ be an $\L$-module on $\Xhat$.  Recall the degree ranges $[c_i(\cdot;\M),d_i(\cdot;\M)]$ from \eqref{eqnWeakAlternateMicrosupportDegreeRanges} after Proposition
~\ref{ssectAlternateMicroSupport}.  Now define
\begin{align}
  \label{eqnD1Definition}
    d_1(\M) &= \sup_{V\in\mS(\M)} \tfrac12(\dim D_P + \dim D_P(V)) + d_1(V;\M)\ ,\\
  \label{eqnD2Definition}
    d_2(\M) &= \sup_{{\substack{V\in\mS(\M)\\ \tilde P\ge P}}} \tfrac12(\dim D_P + \dim D_P(V)) + \dim \sa_P^{\tilde P} + d_2(V,\tilde P;\M)\ ,\\
  \label{eqnD3Definition}
    d_3(\M) &= \sup_{\tilde V\in\mS(\M)} \tfrac12(\dim D_{\tilde P} + \dim D_{\tilde P}(\tilde V)) + d_3(\tilde V;\M)\ ,\\
  \label{eqnD4Definition}
    d_4(\M) &= \sup_{\tilde{\tilde V}\in\mS(\M)} \tfrac12(\dim D_{\tilde{\tilde P}} + \dim D_{\tilde {\tilde P}}(\tilde{\tilde V}) + d_4(\tilde{\tilde V};\M)\ .
\end{align}
analogously to \eqref{eqndMDefn}.  Similarly define $c_i(\M)$ 
analogously to \eqref{eqncMDefn} but without the $\dim \sa_P^{\tilde P}$ term for $i=2$.
Note $[c(\M),d(\M)]=[c_1(\M),d_1(\M)]$.  In fact we have

\begin{prop*}
The intervals $[c_i(\M),d_i(\M)]$ are equal for $i=1$, \dots, $4$.
\end{prop*}

\begin{proof}
We will prove all $d_i(\M)$ are equal; the proof of the other
equalities is similar.  For $V\in\mS(\M)$ let $V\preccurlyeq_0 \tilde V
\preccurlyeq_0 \tilde{\tilde V}$ be chosen as in Proposition
~\ref{ssectAlternateMicroSupport} such that $d_2(V;\M) = d_2(V,\tilde
P;\M)$.  Since $V|_{M_P}$ is conjugate self-contragredient the same is true for $\tilde V|_{M_{\tilde P}}$ and $\tilde{\tilde V}|_{M_{\tilde{\tilde P}}}$ by Lemma ~\ref{lemPartialOrdering} so $\tilde V$ and $\tilde{\tilde V}\in\mS(\M)$.  Furthermore the lemma implies $[\tilde V:V]= \frac12 \dim \n_P^{\tilde P}$ and $[\tilde{\tilde V}:\tilde V]= \frac12 \dim \n_{\tilde P}^{\tilde{\tilde P}}$.  Thus by \eqref{eqnAlternateDegreeRanges},
\begin{equation*}
\begin{split}
d_1(V;\M) \le  d_2(V,\tilde P;\M) &= d_3(\tilde V;\M) +
\tfrac12\dim\n_P^{\tilde P} \\
&\le d_4(\tilde{\tilde V};\M) + \dim \sa_{\tilde
P}^{\tilde{\tilde P}} + \tfrac12\dim \n_{P}^{\tilde{\tilde P}} \ .
\end{split}
\end{equation*}
These inequalities together with \eqref{eqnD1Definition}--\eqref{eqnD4Definition} and
\eqref{eqnUpperdimDPInequality} (applied both to $V\preccurlyeq_0 \tilde V$
and to $\tilde V\preccurlyeq_0 \tilde{\tilde V}$) yield $d_1(\M)\le d_2(\M)\le
d_3(\M) \le d_4(\M)$.  On the other hand, the inequality $d_4(\tilde{\tilde V};\M) \le d_1(\tilde{\tilde V};\M)$
implies that $d_4(\M)\le d_1(\M)$.
\end{proof}

\begin{cor*}
For and $L_P$-module $V\in\emS(\M)$ let $[\essc(V;\M), \essd(V;\M)]$ be the smallest interval containing the degrees in which
\begin{equation*}
\Im \bigl(H(\i_{P}^* \ihat_{{Q_V^W}}^! \M)_V \longrightarrow H(\i_{P}^*
\ihat_{Q^{\prime W}_V}^! \M)_V \bigr) \neq 0 \ .
\end{equation*}
Then
\begin{align*}
c(\M) &= \inf_{V\in\emS(\M)} \tfrac12(\dim D_P - \dim D_P(V)) + \essc(V;\M)\ ,
\text{ and} \\
d(\M) &= \sup_{V\in\emS(\M)} \tfrac12(\dim D_P + \dim D_P(V)) + \essd(V;\M)\ .
\end{align*}
\end{cor*}

\section{$L^2$-cohomology with Boundary Values}
\label{sectLtwo}
In this section we introduce a variant of $L^2$-cohomology which
admits restriction maps to boundary faces.

\subsection{$L^2$-cohomology}
\label{ssectLtwoCohomology}
We first recall ordinary $L^2$-cohomology; references are
\cite{refnCheeger} and \cite{refnZuckerWarped}.  Let $Y$ be a Riemannian
manifold (possibly with boundary or corners) and let $\EE$ be a metrized
locally constant sheaf.  (We allow $\EE$ to be graded but require it to be bounded.)  For a measurable $\EE$-valued form $\o$ on $Y$ let
\begin{equation*}
\lVert\o\rVert=\lVert\o\rVert_Y=\left(\int_Y \lvert\o\rvert^2
 \,dV\right)^{\frac12}
\end{equation*}
denote the usual $L^2$ norm (we will omit the subscript $Y$ if this would
not cause confusion).    Let $L_2(Y;\EE)$ be the Hilbert space of those
forms $\o$ for which $\lVert\o\rVert < \infty$; equivalently this is the
completion of $A_c(Y;\EE)$ with respect to the $L^2$ norm.  Let
$d=d_Y=d_{Y,\EE}$ denote the unbounded densely defined operator given by
exterior differentiation on the domain
\begin{equation*}
\dom d = \{\,\o\in A(Y;\EE)\mid
\lVert\o\rVert,\lVert d\o\rVert<\infty\,\}
\end{equation*}
and let $\dbar$ be its closure in the graph norm.  That is, $\o\in
L_{2}(Y;\EE)$ is in $\dom \dbar$ if and only if there exists a
sequence $\o_i\in\dom d$ such $\o_i$ and $d\o_i$ are Cauchy sequences and
$\o_i \to \o$ in $L^2$ norm.  The \emph{$L^2$-cohomology}
$H_{(2)}(Y;\EE)$ is defined to be the cohomology of the complex $\dom
\dbar$.  A smoothing argument \cite[\S15]{refndeRham},
\cite[\S8]{refnCheeger} shows that the inclusion $\dom d\subset \dom
\dbar$ induces an isomorphism $H_{(2)}(Y;\EE)\cong H(\dom d)$.
Consequently
\begin{equation}
\label{eqnLtwocohomologyIsOrdinaryForCompact}
H_{(2)}(Y;\EE)\cong H(Y;\EE) \qquad\text{if $Y$ is compact.}
\end{equation}

We give $H_{(2)}(Y;\EE)$ the quotient topology; this is Hausdorff if and
only if $\Range \dbar$ is closed.  If it is Hausdorff, then
$H_{(2)}(Y;\EE)$ is representable by harmonic forms, namely $\ker \dbar
\cap \ker \dbar^*$ where $\dbar^*$ is the adjoint of $\dbar$ in the sense
of unbounded operators on the Hilbert space $L_2(Y;\EE)$; conversely if it
is not Hausdorff the $L^2$-cohomology is infinite dimensional.

If $g>0$ is a positive function on $Y$ one may also consider the weighted
norm $\lVert\o\rVert_g^2 = \int_Y \lvert\o\rvert^2 g\,dV$, the
corresponding Hilbert space $L_2(Y;\EE,g)$, and the weighted $L^2$-cohomology
$H_{(2)}(Y;\EE,g)$.  Since this case may be absorbed into the previous by
modifying the metric on $\EE$, we do not mention it again in this section.

\subsection{Spectral vanishing criterion for $H_{(2)}(Y;\EE)$}
\label{ssectSpectralVanishingLtwoCohomology}
If $Y$ is complete as a metric space, a certain estimate implies the
vanishing of $L^2$-cohomology.  More precisely, define the space of
\emph{$\epsilon$-approximate harmonic forms}
\begin{equation*}
\H_\epsilon(Y;\EE)=
\{\,\o\in \dom \dbar^* \cap
A_c(Y;\EE) \mid \lVert d \o \rVert^2 + \lVert d^*\o\rVert^2 \le \epsilon
\lVert \o \rVert^2 \,\}\ .
\end{equation*}
It is useful to keep in mind that $\H_\epsilon(Y;\EE)=0$ is equivalent to the estimate
\begin{equation}
  \label{eqnSpectralVanishingEstimate}
  \epsilon \lVert \o \rVert^2 < \lVert d \o \rVert^2 + \lVert d^*\o\rVert^2\qquad \text{for all $\o\in \dom \dbar^* \cap
A_c(Y;\EE)$.}
\end{equation}
(Note that the intersection with $\dom \dbar^*$, which forces $\o$ to
satisfy Neumann boundary conditions, may be omitted if $\partial
Y=\emptyset$.)

\begin{prop*}
\textup(\cite[Prop.~1.2]{refnSaperSternTwo}, \cite[Thm. ~1.1.4]{refnHormander}\textup) Assume that $Y$ is
complete as a metric space.  Then $\H_{\epsilon}^i(Y;\EE) = 0$ for some
$\epsilon>0$ if and only if $H_{(2)}^i(Y;\EE)=0$ and $H_{(2)}^{i+1}(Y;\EE)$
is Hausdorff.
\end{prop*}

\subsection{$L^2$-cohomology with boundary values}
\label{ssectLtwoBoundaryValues}
For $\o\in A(Y;\EE)$ define the norm
\begin{equation*}
\lvvv\o\rvvv^2 = \|\o\|_Y^2 + \sum_F \| \o|_{F} \|_F^2,
\end{equation*}
where $F$ ranges over the closed proper boundary faces of $Y$.  We will
write, for example, ``$\lvvv\cdot\rvvv$-bounded'' if we mean to use this
norm.  Let $d_b= d_{Y,b}=d_{Y,b,\EE}$ denote the unbounded operator on
$L_{2}(Y;\EE)$ given by exterior differentiation on the domain
\begin{equation*}
\Dom d_{b} = \{\,\o\in A(Y;\EE) \mid \lvvv\o\rvvv, \lvvv
d\o\rvvv<\infty\,\}
\end{equation*}
and let $\dbar_{b}$ be its $\lvvv\cdot\rvvv$-closure in the graph norm.
That is, $\o\in L_{2}(Y;\EE)$ is in $\dom \dbar_b$ if and only if
there exists a sequence $\o_i\in\dom d_b$ such $\o_i$ and $d\o_i$ are
Cauchy sequences in the norm $\lvvv\cdot\rvvv$ and $\| \o_i - \o\|\to 0$
(usual $L^2$ norm).  Since $\left\|\cdot\right \| \le \lvvv\cdot\rvvv$, we
have an inclusion of complexes $\dom \dbar_b \subset \dom \dbar$, but it
is important to note that $\dbar_b$ is \emph{not} a closed operator on
$L_{2}(Y;\EE)$ in general.  (In \S\ref{ssectLtwoBoundarySpace} below
we will consider a larger Hilbert space in which $\dbar_b$ is closed.)  The
\emph{$L^2$-cohomology with boundary values} is defined to be
\begin{equation*}
H_{(2),b}(Y;\EE) = H(\Dom\dbar_b)\ .
\end{equation*}
(Actually in general one would want to consider various weight functions on
the boundary faces in defining $\lvvv\cdot\rvvv$, but this is not
necessary here; see Remark ~\ref{ssectGeneralBoundaryNorm} below.)

\subsection{}
\label{ssectBoundaryRestriction}
$L^2$-cohomology with boundary values is convenient since $\o\in \dom
\dbar_b$ admits a strong $L_2$ trace or restriction $\o|_F\in
L_2(F;\EE|_F)$ for each closed boundary face $F$.  Specifically

\begin{lem*}
Restriction to a closed bounded face $F$ induces a well-defined
$\lvvv\cdot\rvvv$-bounded map
\begin{equation*}
\Dom\dbar_{Y,b}\to \Dom\dbar_{F,b}
\end{equation*}
and consequently a map
\begin{equation*}
H_{(2),b}(Y;\EE) \longrightarrow H_{(2),b}(F;\EE|_F) \ .
\end{equation*}
\end{lem*}

\begin{proof}
Given $\o\in \dom \dbar_b$ let $\o_i$ be a sequence of smooth forms with
$\lvvv \o_i - \o\rvvv\to 0$ and $d\o_i$ a Cauchy sequence in the norm
$\lvvv\cdot\rvvv$.  Then $\o_i|_F $ is a Cauchy sequence in
$L_2(F;\EE|_F)$ and we set $\o|_F$ to be its limit.  Clearly
$\o|_F\in \dom \dbar_{F,b}$.

To see that $\o|_F$ is well-defined, it suffices to consider $\o = 0$
and $F$ a codimension 1 boundary face.  (For if $\o|_F$ is
well-defined and $F'\subset F$, then $\o|_F\in \dom\dbar_{F,b}$ and
$\o|_{F'} = (\o|_F)|_{F'}$.)  Let $\Int F$ denote the interior of $F$
and $\v$ the outward unit normal vector.  For any form $\al\in
A_c(\Int F;\EE|_F)$, there clearly exists a form $\tilde \al\in
A_c(Y;\EE)$ of one degree higher with Neumann boundary data $(\iota_{\v}\tilde
\al)|_F = \al$.  Then Stokes's theorem implies that
\begin{equation*}
(\o_i|_F,\al)_F = (d\o_i,\atilde) - (\o_i,\CoDifferential \atilde)
\end{equation*}
where $\CoDifferential$ is the formal adjoint to $d$.  Now if $\|\o_i\|_b
\to 0$ and $d\o_i$ is $\|\cdot\|_b$-Cauchy we have $\|d\o_i\| \to 0$ and
hence, by the above formula, $ (\o_i|_F,\al)_F\to 0$ for all $\al$.
Since such $\al$ are dense in $L_2(F;\EE|_F)$, we see that
$\o|_F = \lim \o_i|_F = 0$.
\end{proof}

\subsection{}
\label{ssectBoundaryAveraging}
On the other hand, we have an inclusion of complexes $\Dom\dbar_b
\hookrightarrow \Dom\dbar$ which induces a homomorphism
\begin{equation*}
H_{(2),b}(Y;\EE) \longrightarrow H_{(2)}(Y;\EE).
\end{equation*}
We will give a simple criterion for this to be an isomorphism which will
suffice for our purposes.

Consider a system of tubular neighborhoods
\begin{equation*}
(N_F, p_F, t_F)
\end{equation*}
of the closed boundary faces $F\subset Y$, where $N_F$ is an open
neighborhood of $F$, $p_F\colon N_F \to F$ is a smooth retraction, $t_F\colon N_F \to
[0,1)^{\codim F}$ is a normal variable with $F=t_F^{-1}(0)$, and
\begin{equation*}
(t_F,p_F)\colon  N_F \tildearrow [0,1)^{\codim F}\times F
\end{equation*}
is a smooth diffeomorphism.  We require that this data is compatible in the
sense that if $F$ and $F'$ intersect, then
\begin{equation}
\begin{aligned}
\relax&N_{F\cap F'}=N_F\cap N_{F'},\\
\relax&p_{F\cap F'} = p_Fp_{F'} = p_{F'}p_F \qquad\text{on $N_{F\cap F'}$},\\
\relax&t_{F\cap F'} = (t_F, t_{F'})\qquad\text{on $N_{F\cap F'}$, up to permutation
of coordinates.}
\end{aligned}
\label{eqnNeighborhoodCompatibility}
\end{equation}
It is not difficult to see that for any manifold with corners (always
assumed paracompact) such a system of tubular neighborhoods always exists.
Note also that such data induces canonical isomorphisms
\begin{equation*}
\EE|_{N_F} \cong p_F^* \EE|_F.
\end{equation*}

\begin{lem*}
Let $Y$ be a Riemannian
manifold with corners and $\EE$ a metrized locally constant sheaf on $Y$.
Assume that a system of tubular neighborhoods satisfying
\eqref{eqnNeighborhoodCompatibility} can be chosen so that the
Riemannian metric $ds^2$ on $Y$ and the metric $h$ on $\EE$ satisfy the
quasi-isometry%
\footnote{Two metrics $g$, $g'$ are said to be \emph{quasi-isometric}
if there exist constants $C$, $C'>0$ such that $Cg'\le g \le Cg'$; we
denote this by $g\sim g'$.}
relations
\begin{gather}
ds^2|_{N_F} \sim dt_F^2 + p_F^* ds^2_F
\label{eqnMetricalCollar} \\
\intertext{and}
h|_{N_F} \sim p_F^* (h|_F) \notag
\end{gather}
uniformly for all codimension $1$ closed boundary faces $F$.  \textup(This
holds in particular if $Y$ is compact.\textup)  Then the natural map induces
an isomorphism
\begin{equation*}
H_{(2),b}(Y;\EE) \longtildearrow H_{(2)}(Y;\EE).
\end{equation*}
\end{lem*}

\begin{proof}
Let $\eta(t)$ be a smooth function on $[0,1]$ which is identically 1 on
$[0,1/4]$ and has support in $[0,1/2)$.  Assume first that $Y$ has a smooth
boundary $F$ and set $t=t_F$ and $p=p_F$.
As in \cite{refnCheeger} and \cite{refnSaperIsol}, the homotopy operator
\begin{equation*}
H\o = 2 \eta(t) \int_0^{1/2}\!\left(\int_a^{t} \iota_{\partial/\partial
t}\o\right) da 
\end{equation*}
yields a homotopy formula 
\begin{equation}
dH + H d = I-P \label{eqnHomotopyFormula}
\end{equation}
on smooth forms, where
\begin{equation*}
P\o = (1-\eta(t))\o + 2 \eta(t)\, p^*\! \int_0^{1/2}\!\!(\o|_{a\times
F})\, da  -  
2 \eta'(t)\,dt\wedge \int_0^{1/2} \!\left(\int_a^t \iota_{\partial/\partial t}
\o \right) da
\end{equation*}
Given our hypotheses on the metrics, an easy estimate using the
Cauchy-Schwarz inequality demonstrates that $H$ and $P$ are bounded in
the usual $L^2$ norm.  Thus \eqref{eqnHomotopyFormula} holds on $\Dom \dbar$.
On the other hand,
\begin{gather*}
\|(H\o)|_F\| = 2 \biggl\| \int_0^{1/2}\! \left(\int_a^{0} \iota_{\partial/\partial
t}\o\right) da \biggr\| \le C \|\o\| \\
\intertext{and}
\|(P\o)|_F\| =  2 \biggl\| \int_0^{1/2}\!\!(\o|_{a\times F})\, da\biggr\| \le C
\|\o\|.
\end{gather*}
Thus $P$ defines a bounded map $\Dom\dbar\to \Dom\dbar_b$ and
\eqref{eqnHomotopyFormula} holds on $\Dom\dbar_b$; this proves the lemma
for the case of a single boundary face.  In general one performs this
construction for each codimension 1 closed boundary face and composes the
homotopies and projections as in for example \cite[\S15]{refndeRham}; the
hypotheses \eqref{eqnNeighborhoodCompatibility} ensure that the
$t_F$-invariance produced by $P_F$ is not disturbed by $P_{F'}$ or
$H_{F'}$.
\end{proof}

\begin{rem}
\label{ssectGeneralBoundaryNorm}
If $Y$ is not compact then analogous results may be proven even if
\eqref{eqnMetricalCollar} fails.  For example, if instead of
\eqref{eqnMetricalCollar} we have
\begin{equation*}
ds^2|_{N_F} \sim (p_F^*w_F)^2 dt_F^2 + p_F^* ds^2_F \ ,
\end{equation*}
where $w_F>0$ is a bounded function on $F$, then the lemma continues to
hold provided that $\lvvv\cdot\rvvv$ is defined using an $L^2$-norm for $F$
which is weighted by $w_F^2$  (and a product weight for faces of higher
codimension).  Even weight functions not independent of $t_F$ may be
handled---for a nontrivial example and application see
\cite[\S\S2--4]{refnSaperIsol}. (In particular Lemma ~3.2 of
\cite{refnSaperIsol} is the analogue of the lemma above.)
\end{rem}

\subsection{}
\label{ssectLtwoBoundarySpace}
Since $\dbar_b$ is not a closed operator on $L_2(Y;\EE)$, we need to be careful about the meaning of
the adjoint operator $\dbar_b^*$ before extending the spectral vanishing
criterion of Proposition ~\ref{ssectSpectralVanishingLtwoCohomology} to
$H_{(2),b}^i(Y;\EE)$.  Let $L_{2,b}(Y;\EE)$ denote the Hilbert space
completion of $A_c(Y;\EE)$ with respect to $\lvvv\cdot\rvvv$; by
associating to a form $\o\in A_c(Y;\EE)$ the collection $(\o|_F)_F$ of
restrictions to faces, we obtain an isometry $L_{2,b}(Y;\EE) \cong
\bigoplus_F L_2(F;\EE|_F)$.  Then $\dbar_b$ is a closed unbounded operator on $L_{2,b}(Y;\EE)$ and sends $(\o|_F)_F\mapsto (d\o|_F)_F$.  Define $\dbar_b^*$ to be the adjoint of this
operator in the sense of unbounded operators on $L_{2,b}(Y;\EE)$.

To understand $\dbar_b^*$ better, assume that $\psi=(\psi_F)_F\in \dom
\dbar_b^*$ with each $\psi_F$ a smooth form on $F$. Then from the relation
$(\dbar_b^*\psi,\o)_b= (\psi,\dbar_b \o)_b$ which is valid for all
$\o\in\dom \dbar_b$, but in particular for all $\o\in
A_c(Y;\EE)$, it is easy to see that
\begin{equation*}
(\dbar_b^* \psi)_F = \CoDifferential_F\psi_F + \sum_{ \Ftilde \supsetneq F}
(\iota_{\v_F^{\Ftilde}} \psi_{\Ftilde})|_F.
\end{equation*}
Here $\CoDifferential_F$ is the formal adjoint to $d$ on $F$ and the sum is
over boundary faces $\Ftilde$ for which $F$ is a codimension $1$ face with
unit normal vector $\v_F^{\Ftilde}$.

Note that $\psi$ belonging to $\dom \dbar_b^*$ does not necessarily imply
that the various $\psi_F$ satisfy Neumann boundary conditions.  Furthermore
even if $\psi$ is induced from a smooth form on $Y$, the same is not
necessarily true for $\dbar_b^*\psi$.

\subsection{Spectral vanishing criterion for $H_{(2),b}(Y;\EE)$}
\label{ssectSpectralVanishingBoundaryValues}
To state the vanishing criterion, define
\begin{multline*}
\H_{\epsilon,b}(Y;\EE)= \{\,\o\in \dom \dbar_b \cap\dom
\dbar_b^*,\ \supp \o \text{ compact} \mid  \\
\lvvv \dbar_b \o \rvvv^2 + \lvvv
\dbar_b^*\o\rvvv^2 \le \epsilon\lvvv \o \rvvv^2 \,\}.
\end{multline*}
(The condition that $\supp\o$ is compact for $\o\in\dom\dbar_b$ simply means $\o\in L_{2}(Y;\EE)$ has compact support.)
We give $\ker \dbar_b$ and $\range \dbar_b$ the
topology coming from the $\lvvv\cdot\rvvv$-norm and let
$H_{(2),b}(Y;\EE)$ have the quotient topology.

\begin{prop*}
Assume that $Y$ is complete as a metric space.  Then
$\H_{\epsilon,b}^i(Y;\EE) = 0$ for some $\epsilon>0$ if and only if
$H_{(2),b}^i(Y;\EE)=0$ and $H_{(2),b}^{i+1}(Y;\EE)=0$ is Hausdorff.
\end{prop*}

\begin{proof}
The desired equivalence would follow from \cite[Thm.~1.1.2, 1.1.4]{refnHormander} if it were not for the restriction that $\supp
\o$ be compact in the definition of $\H_{\epsilon,b}^i(Y;\EE)$.  To see
that this doesn't matter, first note that by Cohn-Vossen's generalization
\cite[\S2.4]{refnBallmann} of the Hopf-Rinow theorem, the completeness of $Y$
implies that the function $r(x)$ (distance in $Y$ from a fixed point) has
compact sublevel sets.  Then follow \cite{refnGaffney} to create a sequence of compactly supported functions $\eta_k$ such that $0\le \eta_k \le 1$, $|d\eta_k|\le
C$, and every compact $K\subset Y$ is contained in $\eta_k^{-1}(1)$ for
$k$ sufficiently large.  For $\o\in \dom \dbar_b \cap\dom \dbar_b^*$, we then obtain that $\eta_k\o\to \o$, $\dbar_b(\eta_k\o)\to \dbar_b\o$, and
$\dbar_b^*(\eta_k\o)\to \dbar_b^*\o$ in the $\lvvv\cdot\rvvv$ norm.  Thus
the existence of a nonzero $\o\in \dom \dbar_b \cap\dom \dbar_b^*$
satisfying $\bigl(\lvvv \dbar_b \o \rvvv^2 + \lvvv \dbar_b^*\o\rvvv^2\bigr)/ \lvvv \o \rvvv^2 \le \epsilon$ would imply that $\lim_{k\to\infty} \bigl(\lvvv \dbar_b (\eta_k \o) \rvvv^2 + \lvvv \dbar_b^*(\eta_k\o)\rvvv^2\bigr)/ \lvvv \eta_k\o \rvvv^2 \le
\epsilon$ and hence $\H_{\epsilon',b}^i(Y;\EE) \neq0$ for any $\epsilon'>\epsilon$.
\end{proof}

\begin{rem*}
The proposition would remain true if in the definition of
$\H_{\epsilon,b}(Y;\EE)$ we restrict to $\o$ smooth (as for ordinary
$L^2$-cohomology).  The proof requires a generalization of Friedrich's
regularization techniques to handle boundary values; since we will not need
this stronger result, we will not go into details.
\end{rem*}

%
%

\section{$L^2$-cohomology of $\L$-modules}
\label{sectLTwoCohomologyLModules}
In this section we will represent the cohomology $H(\Xhat;\M)$ of a bounded
$\L$-module $\M$ by a variant of $L^2$-cohomology.  To do this, we need to
consider $\Xbar_t \subsetneq X$, a diffeomorphic copy of $\Xbar$ which was
constructed in \cite{refnSaperTilings}, and replace $\Xhat$ by the
corresponding $\Xhat_t$.  The induced metric on a stratum $X_{P,t}$ of $\Xhat_t$ extends to a smooth metric on its Borel-Serre compactification unlike that for a stratum of $\Xhat$.  Thus we will be able to calculate $H(\Xhat;\M)$ using an $L^2$-complex with boundary values similarly to \S\ref{ssectLtwoBoundaryValues}.  The result is a spectral vanishing criterion for this 
$L^2$-cohomology analogous to those in Propositions
~\ref{ssectSpectralVanishingLtwoCohomology} and
\ref{ssectSpectralVanishingBoundaryValues}.

We also introduce a variant of $\L$-modules on $\Abar_P\times X_P$ and their $L^2$-cohomology.

\subsection{Locally symmetric metrics}
\label{ssectMetrics}
We will work with a fixed basepoint $x\in D$.  Let $KA_G\subset
G(\RR)$ be the stabilizer of the basepoint $x$ and let $\theta$ be the
Cartan involution associated to $K$ (see for example
\cite[V.24.6]{refnBorelLAG} or \cite[\S1.6]{refnBorelSerre} for the
nonconnected and nonsemisimple case) with Cartan decomposition
$\mathfrak g=\k+\p+\sa_G$.  Fix an $\Ad G(\RR)$-invariant,
$\theta$-invariant nondegenerate bilinear form $B$ on $\mathfrak g$
such that $(X,Y) \mapsto -B(X,\theta Y)$ is a positive definite inner
product on $\mathfrak g$.  (On the derived algebra one may take $B$ to
be the Killing form.)  The restriction of $B$ to $\p+\sa_G$ yields an
$\Ad (KA_G)$-invariant inner product on $T_x D$ which may be extended to
a $G(\RR)$-invariant metric on $D$; this descends to a locally
symmetric Riemannian metric on $X=\G\back D$.  In addition one obtains
from $B$ left-invariant metrics on all Lie subgroups of $G(\RR)$.

If $(E,\sigma)$ is an algebraic representation of $G$ let $\EE = D \times_{\G}
E$ be the corresponding flat bundle on $X$.  An \emph{admissible}
inner product on $E$ is a positive definite Hermitian bilinear form on $E$
which is $K$-invariant and for which $\p+\sa_G$ acts via self-adjoint
operators; equivalently it is a form for which $\sigma(g)\sigma(\theta
g)^*=\id_E$ for $g\in G(\RR)$.  Admissible metrics always exist (see
\cite[p.~ 375]{refnMatsushimaMurakami} in the case where $G(\RR)$ is
connected and semisimple; the generalization to our situation is clear).
Such an admissible inner product on $E$ induces a metric on $\EE$; in the
case that $E$ is isotypical it is given by the formula $| (gKA_G,v) | =
|\xi_E^k(g)|^{\frac1k} \cdot |\sigma(g^{-1})v|$, where $\xi_E$ is the
character by which $S_G$ acts on $E$ and $k\in\NN$ is such that $\xi_E^k$
extends to a character on $G$.

The above construction may be applied to all $L_P$ (using the basepoint
$A_P \geo N_Px\in D_P$, corresponding maximal compact subgroup
$K_P\subset L_P(\RR)$ and the bilinear form on $\levi_P$ induced from
$B$) in order to define locally symmetric metrics on the $X_P$ and metrics
on the flat bundles $\EE_P$ associated to algebraic representations $E_P$ of
$L_P$.

For a parabolic $P$ there is (analogously to
\S\ref{ssectNilmanifoldFibrations}) a flat nilmanifold fibration $\G_P\back
D \to Z_P \equiv \G_{L_P}\back (N_P(\RR)\back D)$.  The space
$N_P(\RR)\back D$ is naturally a homogeneous space for $L_P(\RR)$ with the
stabilizer of a point being a maximal compact subgroup times $A_G$.  The
metric on $D$ induces a $L_P(\RR)$-invariant metric on $N_P(\RR)\back D$
and hence a metric on $Z_P$.  There is a trivial principal $A_P^G$-bundle
$Z_P \to X_P$ given by geodesic retraction; canonical trivializing sections
are given by quotients of $L_P(\RR)$-orbits.  The metric on $Z_P$
decomposes under a canonical trivialization into the sum of an invariant
metric on $A_P^G$ and the locally symmetric metric on $X_P$.  Given an
algebraic representation $E_P$ of $L_P$, let $\EEtilde_P = (N_P(\RR)\back D)
\times_{\G_{L_P}} E_P$ be the corresponding flat bundle on $Z_P$;
alternatively $\EEtilde_P$ is the pullback of $\EE_P$ on $X_P$.  The
admissible inner product on $E_P$ induces a metric on $\EEtilde_P$ by the
formula $|(rK_PA_P,v)| = |\sigma(r^{-1})v|$ for all $r \in M_P(\RR)
A_P^G$; note that this is \emph{not} the pullback of the metric on
$\EE_P$ unless $A_P^G$ acts trivially.

\subsection{Diffeomorphisms of $\Xbar$ into $X$}
\label{ssectTilings}
By \S\ref{ssectGeodesicRetraction} we may identify $D$ with $A_P^G \times
e_P$ such that the $\lsp0P(\RR)$-orbit of the fixed basepoint $x$ is
identified with $\{1\}\times e_P$.  This induces identifications $\G_P\back
D \cong A_P^G \times Y_P$ and $Z_P \cong A_P^G\times X_P$.

By \cite[Thm.~ 6.1]{refnSaperTilings}%
\footnote{In fact \cite{refnSaperTilings} constructs a disjoint
decomposition of $D$ (indexed by parabolic $\QQ$-subgroups) called a
\emph{tiling} of which $\Dbar_t$ is the central tile corresponding to
$G$; see in particular Definition ~2.1, Theorem ~5.7, and Theorem ~6.1.
Note the differences in notation: in \cite{refnSaperTilings}, $X$ refers to
the symmetric space itself rather than its arithmetic quotient.  More
significantly, in \cite{refnSaperTilings} $X_P$ denotes the piece of the
tiling corresponding to $P$, as opposed to a stratum of $\Xhat$ as in the
current paper.}$^,$%
\footnote{The last assertion of \cite[Thm.~6.1]{refnSaperTilings}, that $s_t$ converges to the identity as $t\to\infty$, is false.  This assertion is not used in the current paper.  In future work we will indicate an alternate family $s'_t$ that satisfies this assertion.}
there exists a family of
$\G$-equivariant piecewise-analytic diffeomorphisms (depending on $t\in A^G$ sufficiently dominant) 
\begin{equation*}
s_t\colon  \Dbar\longtildearrow \Dbar_t\subsetneq D\ ,
\end{equation*}
such that $\Dbar_t$ is a subdomain
with corners.  For all%
\footnote{Actually equation \eqref{eqnBoundaryLocation} is only valid for
$P$ belonging to a fixed set of representatives of $\G$-conjugacy classes
of parabolic $\QQ$-subgroups; we will assume such representatives have been
chosen.  For other $P$ one has $e_{P,t}\subset \{t_Pb_P\}\times e_P$
where $b_P\in A_P^G$ is a certain parameter needed to make the construction
$\G$-invariant; this is because $\gamma\in \G$ may move the fixed basepoint
$x$.  See \cite[\S2]{refnSaperTilings}, particularly \S\S2.6, 2.7.}
$P\in\Pl$ each diffeomorphism $s_t$ satisfies:
\begin{equation}
e_P \xrightarrow{s_t|_{e_P}} e_{P,t}\equiv 
s_t(e_P) \text{ is $N_P(\RR)$-equivariant}
\label{eqnNPequivariance}
\end{equation}
and
\begin{equation}
e_{P,t} \subsetneq  \{t_P\}\times e_P \label{eqnBoundaryLocation}
\end{equation}
(where $t_P\in A_P^G$ is the image
of $t$ under the canonical projection $A^G\cong A_P^G\times A^P \to A_P^G$).
We obtain an induced diffeomorphism
\begin{equation*}
\Xbar \tildearrow \Xbar_t \equiv \G\back \Dbar_t \subsetneq X\ .
\end{equation*}

\subsection{}
\label{ssectTilingInCylindricalSet}
We can describe $\Xbar_t\subsetneq X$ more precisely by using a cylindrical
cover of $X$.  Given $s_P\in A_P^G$ and $t\in A^G$, set
\begin{align*}
A_P^G(s_P) &\equiv \{\, a\in A_P^G\mid a^\al > s_P^\al \text{ for all
$\al\in\D_P$}\,\}, \text{ and} \\
\langle A_P^G\rangle_t &\equiv \{\, a\in A_P^G \mid a^\b \le t^\b \text{ for all
$\b\in\Dhat_P$}\,\}.
\end{align*}
For $\Omega_P\subset Y_P$ open, relatively compact, and
$N_P(\RR)$-invariant, let $W_P \equiv A_P^G(s_P) \times \Omega_P\subsetneq
A_P^G\times Y_P \cong \G_P\back D$ be the associated \emph{cylindrical
set}.  We always assume that $s_P$ is sufficiently dominant (depending on
$\Omega_P$) so that $W_P$ may be identified with its image in $X=\G\back
D$.  For appropriate $s_P$ and $\Omega_P$ we obtain a \emph{cylindrical
cover} $\{W_P\}_{P\in\Pl}$ of $X$.  Given any cylindrical cover
$\{W_P\}$, the intersections $\Xbar_t \cap W_P$ form a cover of $\Xbar_t$;
for $t$ sufficiently dominant one can prove \cite[Thm.~ 5.7]{refnSaperTilings} that $\Xbar_t \cap W_P\neq\emptyset$ and that
\begin{equation}
\Xbar_t \cap W_P = (\langle A_P^G\rangle_t \times Y_P)\cap W_P\ .
\label{eqnAdaptedTiling}
\end{equation}
We set $\Wbar_{P,t} \equiv \Xbar_t \cap W_P$.

\subsection{}
\label{ssectMetrizeRetract}
Set $Y_{P,t} \equiv \G_P\back e_{P,t}$.
The diffeomorphism $\Xbar \tildearrow \Xbar_t$ induces diffeomorphisms
$Y_P \tildearrow Y_{P,t} \subsetneq \G_P\back D$ of the
boundary faces.   For each $P$ this diffeomorphism descends by
\eqref{eqnNPequivariance} to
\begin{equation}
X_P \longtildearrow X_{P,t} \subsetneq Z_P\ ,
\label{eqnStratumRetract}
\end{equation}
where $X_{P,t}$ is formed from $Y_{P,t}$ by collapsing the
$N_P(\RR)'$-fibers of the nilmanifold fibration (see \S\ref{ssectNilmanifoldFibrations}).  We set 
\begin{equation*}
\Xhat_t= \coprod_P X_{P,t}
\end{equation*}
and equip it with the topology so that \eqref{eqnStratumRetract} induces a
homeomorphism $\Xhat\tildearrow \Xhat_t$.  The closure of $X_{P,t}$ as a
subspace of $\Xhat_t$ will be denoted $\Xhat_{P,t}$.

The locally symmetric metric on $Z_P=A_P^G\times X_P$ induces on each stratum
$X_{P,t}$ of $\Xhat_t$ a metric which extends smoothly to its Borel-Serre
compactification $\Xbar_{P,t}$.  Given an algebraic representation $E_P$ of
$L_P$ there is a natural isomorphism of flat bundles
\begin{equation*}
\EE_P \longtildearrow \EE_{P,t} \equiv \EEtilde_P|_{X_{P,t}}
\end{equation*}
over the diffeomorphism \eqref{eqnStratumRetract}.  As in
\S\ref{ssectMetrics}, the choice of an admissible inner product on $E_P$
induces a metric on $\EEtilde_P$ and hence a metric on $\EE_{P,t}$ that
also extends smoothly over $\Xbar_{P,t}$.

\subsection{}
\label{ssectLtwoMethods}
Let $\M=(E_\cdot,f_{\cdot\cdot})$ be a bounded $\L$-module on $\Xhat$.  In
\S\ref{sectRealizationLModules} we constructed an incarnation
$\A_{\Xhat}(\M)$ of the realization using special differential forms; this
was a complex of fine sheaves on $\Xhat$.  The cohomology $H(\Xhat;\M)$ is
the cohomology of the global sections of $\A_{\Xhat}(\M)$.

We now consider an alternate incarnation and realization which is a complex
of sheaves on $\Xhat_t$ instead of $\Xhat$.  Namely we set $\A_P(E_P) =
\Asp(\Xhat_{P,t};\EE_{P,t})$ instead of $\Asp(\Xhat_P;\EE_P)$ and let
$\A_{\Xhat_t}$ be the corresponding realization functor provided by Theorem ~\ref{ssectRealizationLModules}.
If $U\subset
\Xhat_t$ is a domain define $H(U;\M)$ to be the cohomology of global
sections of $\A_{\Xhat_t}(\M)|_U$.  For $U=\Xhat_t$, clearly
\begin{equation}
  \label{eqnCohomologySameWithT}
H(\Xhat;\M) \cong H(\Xhat_t;\M)\ .
\end{equation}

Explicitly, by \S\S\ref{ssectRealizationLsheaves}, \ref{ssectProofOfTheoremRealization}
the cohomology $H(U;\M)$ is the
cohomology of the complex
\begin{equation}
\label{eqnCohomologyULModuleComplex}
\left\{
\begin{aligned}
A\sp(U;\M) & = \bigoplus_P
 A\sp(U\cap\Xhat_{P,t};\EE_{P,t})\ ,
 \\
d &= \sum_P d_P + \sum_{P \le Q} \dd_{PQ}\ ,
\end{aligned}
\right.
\end{equation}
where $d_P$ is the differential on $A\sp(U\cap\Xhat_{P,t};\EE_{P,t})$ and
$\dd_{PQ}$ is defined by
\begin{equation*}
\dd_{PQ}(\o_Q) \equiv A\sp(U\cap\Xhat_{P,t};f_{PQ}\circ
\coh_{PQ})(\o_Q|_{U\cap\Xhat_{P,t}}) \ .
\end{equation*}

\subsection{}
\label{ssectDDBoundIndependentT}
For $\o_Q\in  A\sp(U\cap\Xhat_{Q,t};\EE_{P,t})$ define a boundary value norm
\begin{equation}
\label{eqnComponentNorm}
\lvvv\o_Q\rvvv[Q]^2 \equiv \sum_{P\le Q} \|\o_Q|_{U\cap X_{P,t}}\|_{U\cap
X_{P,t}}^2\ ;
\end{equation}
this is analogous to \S\ref{ssectLtwoBoundaryValues} but we are using the restriction to $X_{P,t}$ from \eqref{eqnRestrictSpecialForms} instead of to $Y_{P,t}$.
Then for $\o=(\o_P)_P\in A\sp(U;\M)$ define the weighted norm
\begin{equation}
\|\o\|_t^2 \equiv \sum_P \lvvv t^{-\hsr_P} \o_P\rvvv[P]^2\ .
\label{eqnComplexNorm}
\end{equation}
We let $t^{-\hsr}$ denote the family of weights $(t^{-\hsr_P})_P$ and let
\begin{equation*}
L_{2}(U;\M,t^{-\hsr}) = 
   \bigoplus_P  L_{2,b}(U\cap\Xhat_{P,t};\EE_{P,t},t^{-\hsr_P})
\end{equation*}
denote the Hilbert space completions of $A\sp(U;\M)$ and
$A\sp(U\cap\Xhat_{P,t};\EE_{P,t})$ under the corresponding norms in
\eqref{eqnComplexNorm}.

For each $P$ the closed unbounded operator $\dbar_{U\cap
\Xhat_{P,t},b,t^{-\hsr_P}}$ is defined on
$L_{2,b}(U\cap\Xhat_{P,t};\EE_{P,t},t^{-\hsr_P})$ and has cohomology
$H_{(2),b}(U\cap\Xhat_{P,t};\EE_{P,t},t^{-\hsr_P})$.  On the other
hand, let $d_{U;\M,t^{-\hsr}}$ be the unbounded operator on
$L_{2}(U;\M,t^{-\hsr})$ corresponding to the subcomplex of
\eqref{eqnCohomologyULModuleComplex} consisting of $\o$ such that $\o$,
$d\o\in L_{2}(U;\M,t^{-\hsr})$ and denote its graph norm closure by
$\dbar= \dbar_{\M,t^{-\hsr}}= \dbar_{U;\M,t^{-\hsr}}$.  The cohomology of $\dom
\dbar_{U;\M,t^{-\hsr}}$ is by definition the \emph{$L^2$-cohomology of
$\M$ over $U$ with weights $t^{-\hsr}$} and is denoted
$H_{(2)}(U;\M,t^{-\hsr})$.  In order to relate these
$L^2$-cohomology groups to $H(U;\M)$ we need the

\begin{lem*}
The operators $\dd_{PQ}$ are bounded 
independently of $t$.
\end{lem*}
\begin{proof}
  The restriction map $\o_Q\mapsto \o_Q|_{U\cap\Xhat_{P,t}}$ is bounded since the  $\lvvv\cdot\rvvv[Q]$-norm incorporates the $L^2$-norm of boundary values.  The norm squared of $\coh_{PQ}$ (from Lemma ~\ref{ssectNilmanifoldCohomologyBundle} where $G$ is replaced by $L_Q$) is, up to a constant, the ratio of the weight function on $\Xhat_{P,t}$ to the weight function on $\Xhat_{Q,t}$.  This  ratio is $(t^{-\hsr_P})(t^{-\hsr_Q})^{-1}=(t^{\hsr_P^Q})^{-1}$ which is bounded as $t$ grows more dominant.
\end{proof}

\begin{cor}
There is a decomposition \textup(as vector spaces, not as complexes\textup)
\begin{equation*}
\dom \dbar_{U;\M,t^{-\hsr}}  =
\bigoplus_P \dom \dbar_{U\cap \Xhat_{P,t},b,t^{-\hsr_P}}
\end{equation*}
and a strict operator equality 
\begin{equation*}
\dbar_{U;\M,t^{-\hsr}} = \sum_P
\dbar_{U\cap \Xhat_{P,t},b,t^{-\hsr_P}} + \sum_{P \le Q} \dd_{PQ}\ .
\end{equation*}
\end{cor}

For $U=\Xhat_t$ we have the

\begin{cor}
\label{ssectCohomologyLModuleIsLtwoCohomology}
For a bounded $\L$-module $\M$, $H(\Xhat_t;\M)\cong H_{(2)}(\Xhat_t;\M,t^{-\hsr})$.
\end{cor}
\begin{proof}
Filter both \eqref{eqnCohomologyULModuleComplex} and $\dom
\dbar_{U;\M,t^{-\hsr}}$ by $\#\D^P\le p$ and compare the resulting
spectral sequences; they are isomorphic at $E_1$ by
\eqref{eqnLtwocohomologyIsOrdinaryForCompact} and Lemma
~\ref{ssectBoundaryAveraging}.
\end{proof}

\subsection{}
\label{ssectSpectralVanishingLModules}
Define
\begin{multline*}
\H_{\epsilon}(U;\M,t^{-\hsr})=
\{\,\o\in
\dom \dbar_{U;\M,t^{-\hsr}} \cap\dom \dbar_{U;\M,t^{-\hsr}}^*,\ \supp \o
\text{ compact}
\mid \\
\| \dbar \o \|_t^2 + \| \dbar^*\o\|_t^2 \leq
\epsilon\| \o \|_t^2 \,\}\ .
\end{multline*}
For $U=\Xhat_t$ we have as in Proposition
~\ref{ssectSpectralVanishingBoundaryValues}
\begin{prop*}
For a bounded $\L$-module $\M$, $\H_{\epsilon}^i(\Xhat_t;\M,t^{-\hsr})=0$ for some $\epsilon>0$ if  and
only if $H^i_{(2)}(\Xhat_t;\M,t^{-\hsr})=0$ and
$H^{i+1}_{(2)}(\Xhat_t;\M,t^{-\hsr})$ is Hausdorff.
\end{prop*}

\section{Vanishing Criterion for Cohomology of an $\L$-module}
\label{sectAnalyticLemmas}
Let $\M$ be a bounded $\L$-module.  The goal of this section is Proposition
~\ref{ssectParabolicAnalyticVanishingCriterion}, a vanishing criterion for
$H^i(\Xhat;\M)$ in terms of certain $L^2$-cohomology groups of $\langle
A_P^G\rangle_1 \times X_P$ for all $P\in\Pl$.  In view of the results of
\S\ref{sectLTwoCohomologyLModules} and particularly Proposition~\ref{ssectSpectralVanishingLModules} we achieve this by a
sequence of lemmas that yield sufficient conditions for the vanishing of
$\H_{\epsilon}^i(\Xhat_t;\M,t^{-\hsr})$.  As we proceed, the parameter $t$ will
have to be made successively more dominant.

\subsection{Spectral Mayer-Vietoris}
\label{ssectSpectralMayerVietoris}
Given a cylindrical cover $\{W_{P}\}$ assume $t$ is sufficiently dominant so
that \eqref{eqnAdaptedTiling} holds.  Each boundary stratum
$Y_{Q,t}\cap\Wbar_{P,t}$ of $\Wbar_{P,t}$ for $Q\ge P$ is part of a
nilmanifold fibration with fibers $N_Q(\RR)'$; let $\What_{P,t}\subset
\Xhat_t$ denote the result of collapsing these fibers.  We have
the following spectral analogue of the Mayer-Vietoris sequence:

\begin{lem*}
Given $\epsilon>0$ there
exists $\epsilon'>0$ and a cylindrical cover $\{W_{P}\}$ of $X$ such
that for all $t$ sufficiently dominant,
$\H_{\epsilon}^i(\What_{P,t};\M,t^{-\hsr})=0$ for all $P$
implies $\H_{\epsilon'}^i(\Xhat_t;\M,t^{-\hsr})=0$.
\end{lem*}
\begin{proof}
For every $c>0$ there exists
\cite[Prop.~ 2.1]{refnSaperSternTwo} a
cylindrical cover $\{W_{P}\}$ possessing a partition of unity
$\{\eta_P\}$ such that $|d\eta_P| < c$.  (One may eliminate the
neatness hypothesis.)  Examination of the construction shows that for $P\le
Q$ the function $\eta_P$ is $N_Q(\RR)$-invariant on $W_P\setminus
\bigcup_{P'\nleq Q} W_{P'}$.  But since $W_P\cap Y_{Q,t}$ is empty unless
$P\le Q$ by \eqref{eqnAdaptedTiling}, this implies that
$\eta_P|_{Y_{Q,t}}$ is $N_Q(\RR)$-invariant for \emph{any} $P$ and $Q$.
Thus we may pushdown to obtain a
partition of unity $\{\etahat_P\}$ for the covering $\{\What_{P,t}\}$ of
$\Xhat_t$.  Now for $\o\in \H_{\epsilon'}^i(\Xhat_t;\M,t^{-\hsr})$ (with
$\epsilon'>0$ to be chosen below) one may estimate (compare
\cite[Prop.~ 1.4]{refnSaperSternTwo}):
\begin{align*}
\| \dbar(\etahat_P\o) \|_t^2 + \|
\dbar^*(\etahat_P\o)\|_t^2
&\le 2 \bigl( \| \dbar\o \|_t^2 + \| \dbar^*\o\|_t^2 + 2c^2
 \| \o \|_t^2\bigr) \\
&< 2(\epsilon' +2c^2)   \| \o \|_t^2 \\
&\le  4(\#\Pl) (\epsilon' +2c^2)  \max_R \|  \etahat_R\o \|_t^2 \le
 \epsilon \max_R \|  \etahat_R\o \|_t^2,
\end{align*}
where we ensure the last line choosing $c>0$ and $\epsilon'>0$ small.  Thus if
$P$ realizes the maximum on the last line we have $\etahat_P\o\in
\H_{\epsilon}^i(\What_{P,t};\M,t^{-\hsr})$.  Since $\o\neq 0$ implies
$\etahat_P\o\neq0$ for this $P$, the lemma is proved.
\end{proof}

\subsection{$\L$-modules on $\langle A_P^G\rangle_t\times X_P$}
\label{ssectLModulesAPTimesXP}
The space $\langle A_P^G\rangle_t \times X_P \subsetneq Z_P$ is a manifold
with corners with strata indexed by the parabolic $\QQ$-subgroups $Q\ge P$.
Explicitly the $Q$-stratum is $\langle A_P^G\rangle_t^Q \times X_P$ where
(see \S\ref{sssectCoroots})
\begin{equation*}
\langle A_P^G\rangle_t^Q \equiv \{\, a\in A_P^G \mid a^\b \le t^\b \text{ for all
$\b\in\Dhat_P\setminus \Dhat_Q$, }a^\b = t^\b \text{ for all
$\b\in\Dhat_Q$}\,\}\ .
\end{equation*}
We may extend the theory of $\L$-modules extends to this (simpler) context:

\begin{defn*}
An \emph{$\L$-module
$\M'$ on $\langle A_P^G\rangle_t \times X_P$}
consists of a family $(E'_Q)_{Q\ge P}$ of graded algebraic representations of
the \emph{same} group $L_P$ and morphisms $f'_{QR}\colon
E'_R \to E'_Q[1]$.
\end{defn*}

Since the space has homotopically trivial
links the functor $H(\n_R^Q;\cdot)$ is replaced throughout in this context  by the identity
functor.  Furthermore special differential forms are no longer needed for
the realization.

The strata and the locally constant sheaves associated to $E'_Q$ are
metrized as in \S\ref{ssectMetrizeRetract} by restricting the locally
symmetric metrics on $Z_P$ and $\EEtilde'_Q$.  We will be using the
$L^2$-space $L_{2}(\langle A_P^G\rangle_t \times X_P; \M')$, and the corresponding $L^2$-cohomology.  There is also an analogue of
Proposition ~\ref{ssectSpectralVanishingLModules}.  Specifically, let
\begin{multline*}
\H_{\epsilon}( \langle A_P^G\rangle_t \times X_P; \M') \equiv
\{\,\o\in
\dom \dbar_{ \langle A_P^G\rangle_t \times X_P;\M'} \cap\dom \dbar_{ \langle A_P^G\rangle_t \times X_P;\M'}^*,\ \\
\supp \o
\text{ compact}
\mid 
\| \dbar \o \|_t^2 + \| \dbar^*\o\|_t^2 \leq
\epsilon\| \o \|_t^2 \,\}\ .
\end{multline*}

\begin{prop*}
For a bounded $\L$-module $\M'$,  $\H_{\epsilon}( \langle A_P^G\rangle_t \times X_P; \M')=0$ for some $\epsilon>0$ if  and
only if $H^i_{(2)}(\langle A_P^G\rangle_t \times X_P; \M')=0$ and
$H^{i+1}_{(2)}(\langle A_P^G\rangle_t \times X_P; \M')$ is Hausdorff.
\end{prop*}

\subsection{Spectral averaging}
\label{ssectSpectralAverging}
Assume that $W_{P}$ is a cylindrical set for a fixed $P$ and that $t$
is sufficiently dominant so that \eqref{eqnAdaptedTiling} holds.  Let $O_P$
be the image of $\Omega_P$ under the nilmanifold fibration $Y_P\to X_P$.
Each boundary stratum $X_{Q,t}\cap\What_{P,t}$ of $\What_{P,t}$ for $Q\ge
P$ is part of a nilmanifold fibration with fibers $N_P^Q(\RR)'$; let
\begin{equation*}
\What_{P,t} \xrightarrow{\nil} (\langle
A_P^G\rangle_t \cap A_P^G(s_P)) \times O_P 
\subset \langle A_P^G\rangle_t \times X_P
\end{equation*}
be the projection obtained by collapsing these fibers.

\begin{defn*}
Define the \emph{pushforward $\L$-module $\nil_*\M$ on
$\langle A_P^G\rangle_t \times X_P$} by $E'_Q = H(\n_P^Q;E_Q)$ and $f'_{QR}=
H(\n_P^Q;f_{QR})$.  Also define $\nil_*\M\otimes \CC_{\hsr_P} = (E'_\cdot\otimes
\CC_{\hsr_P},f'_{\cdot\cdot}\otimes \id_{\CC_{\hsr_P}})$.
\end{defn*}

\begin{lem*}
Given $\epsilon>0$ and a cylindrical set $W_P$
there exists $\epsilon'>0$ such that for all $t$ sufficiently dominant,
$\H_{\epsilon}^i( \langle A_P^G\rangle_t \times X_P; \nil_*\M\otimes\CC_{\hsr_P})=0$   implies
$\H_{\epsilon'}^i(\What_{P,t};\M,t^{-\hsr})=0$.
\end{lem*}
\begin{proof}
For $Q\ge P$ let $\nil_{Q*}$ be the operator on
$L_{2}(X_{Q,t}\cap \What_{P,t};\EE_{Q,t},t^{-\hsr_Q})$ which orthogonally
projects onto $N_P^Q(\RR)$-invariant and $\n_P^Q$-harmonic forms.
Zucker \cite[(4.24)]{refnZuckerWarped} shows there
is a homotopy formula 
\begin{equation}
\dbar_Q H_Q + H_Q \dbar_Q = I - \nil_{Q*}\label{eqnQHomotopy}
\end{equation}
where $H_Q$ and $\nil_{Q*}$ bounded in the $L^2$ norm.  These bounds are  indepedent of $t$ as $t$ grows more dominant though they are dependent on
$W_P$.  In fact since the
boundary faces of $\Xbar_{Q,t}\cap \What_{P,t}$ contain full
$N_P^Q(\RR)$-fibers, the result holds for the $\lvvv\cdot\rvvv[Q]$ norm
and we obtain a homotopy formula for $\dbar_{Q,b,t^{-\hsr_Q}}$.

For $Q\le R$ recall the operators $d_{QR}$ from \eqref{eqnCohomologyULModuleComplex}.  We now arrange that
\begin{equation}
\nil_{Q*} \dd_{QR} = \dd_{QR} \nil_{R*} \qquad \text{and}
             \qquad H_Q \dd_{QR} = - \dd_{QR} H_R \label{eqnCompatibility}
\end{equation}
(the sign is due to the fact that $\dd_{QR}$ is degree 1 in the
coefficients).  Zucker's
construction for $H_Q$ depends on the choice of a central series defined
over $\QQ$,
\begin{equation*}
1 = N_0^Q \subset N_1^Q \subset \dots \subset N_\v^Q = N_P^Q(\RR),
\label{eqnCentralSeries}
\end{equation*}
and unit $\QQ$-root vectors $Z_j^Q\in \n_j^Q$ which satisfy $\n_j^Q =
\n_{j-1}^Q + \RR\cdot Z_j^Q$.
Choose such a central series and root vectors for $N_P(\RR)$ (the case
$Q=G$) and use the induced data for each $N_P^Q(\RR)$.  (Some of the
vectors $Z_j^Q$ will now be zero but those stages in the series may be
omitted.)  One may now verify from the description in
\cite[\S4(c)]{refnZuckerWarped} that \eqref{eqnCompatibility} holds.

By summing \eqref{eqnQHomotopy} over $Q\ge P$ we obtain
\begin{equation}
\dbar_{\M}H + H \dbar_{\M} = I - \nil_*, \label{eqnFullHomotopy}
\end{equation}
where $H=\sum_{Q\ge P} H_Q$, $\nil_* = \sum_{Q\ge P} \nil_{Q*}$, and the
$\dd_{QR}$ terms vanish due to \eqref{eqnCompatibility}.
Since $p_{G*}\o=\o$ the map $p_*$ is injective.
Now for $\o\in \H_{\epsilon'}^i(\What_{P,t};\M,t^{-\hsr})$ calculate:
\begin{align*}
\| \dbar\o \|_t^2 + \| \dbar^*\o\|_t^2 &\le \epsilon'\| \o
\|_t^2 \\
&\le \epsilon' ( \|  \nil_*\o \|_t^2 + 
                 \| (I-\nil_*) \o \|_t^2 ) \\
&\le \epsilon'  \| \nil_*\o \|_t^2 + 
     2\epsilon'\| H\|_t^2\cdot (\| \dbar\o \|_t^2 +
                              \| \dbar^*\o\|_t^2 );
\end{align*}
the final line uses \eqref{eqnFullHomotopy} and  \cite[Prop.~ 1.6]{refnSaperSternTwo}.  Choose $\epsilon'>0$ so that
$2\epsilon'\| H\|_t^2 \le1/2$ and the last term
may be absorbed into the left-hand side.  Now we can estimate
\begin{equation*}
\| \dbar \nil_*\o \|_t^2 + \| \dbar^*\nil_*\o\|_t^2 \le
\| \nil_*\|_t^2 \cdot(\| \dbar\o \|_t^2 + \| \dbar^*\o\|_t^2) \le
2\epsilon'\| \nil_*\|_t^2 \cdot\| \nil_*\o \|_t^2 \ .
\end{equation*}
For  $\epsilon'>0$ so that
$2\epsilon'\| \nil_*\|_t^2\le \epsilon$ we see that $\nil_*\o \in \H_{\epsilon}^i(\What_{P,t};\M,t^{-\hsr})$.
Since $\| H\|_t$ and $\| \nil_*\|_t $ are bounded as $t$ grows more dominant,  $\epsilon'$ can be chosen independently of $t$.

To complete the proof, we consider the embedding
\[
i\colon \nil_*(\Dom \dbar_{\What_{P,t},\M,t^{-\hsr}})\hookrightarrow \Dom \dbar_{\langle A_P^G\rangle_t \times X_P,\nil_*\M\otimes \CC_{\hsr_P}}
\]
obtained by extending an $N_P^Q(\RR)$-invariant
compactly supported form by zero and applying Lemma
~\ref{ssectNilmanifoldCohomologyBundle} to the coefficients.  Here we also
to choose a flat trivializing section of the bundle associated to
$\CC_{\hsr_P}$.  To see that this is an isometry, note that the $L^2$-norm
integrand for the $Q$ summand of $\nil_*(\Dom \dbar_{\What_{P,t},\M,t^{-\hsr}})$
involves $t^{-2\hsr_Q}$ whereas for $\Dom \dbar_{\langle A_P^G\rangle_t \times X_P,\nil_*\M\otimes
\CC_{\hsr_P}}$ the coefficient system $\CC_{\hsr_P}$ contributes
$a^{-2\hsr_P}=a^{-2\hsr_Q}\cdot a^{-2\hsr_P^Q}=t^{-2\hsr_Q}\cdot
a^{-2\hsr_P^Q}$. The discrepancy is canceled by the map $\coh_{PQ}$ which
multiplies norm squared by $a^{2\hsr_P^Q}$ (see the proof of
Lemma~\ref{ssectDDBoundIndependentT}).  Thus $p_*\o\in \H_{\epsilon}^i( \langle A_P^G\rangle_t \times X_P; \nil_*\M\otimes\CC_{\hsr_P})$.
\end{proof}

\subsection{Isotypical reduction}
\label{ssectLPIsotypicalDecomposition}
We have a isotypical decomposition of $\L$-modules:
\begin{equation*}
\nil_*\M = \bigoplus_{V\in \IrrRep(L_P)} (\nil_*\M)_V\ .
\end{equation*}

\begin{lem*}
For any $\epsilon>0$ and any $t$, the vanishing
$\H_{\epsilon}^i ( \langle A_P^G\rangle_t \times X_P;
\nil_*\M\otimes \CC_{\hsr_P})=0$ is equivalent to the vanishing conditions 
\begin{equation*}
\H_{\epsilon}^i ( \langle A_P^G\rangle_t \times X_P;
(\nil_*\M)_V\otimes \CC_{\hsr_P})=0 \qquad\text{for all $V\in  \IrrRep(L_P)$}.
\end{equation*}
\end{lem*}
\begin{proof}
If $\o =\sum_V \o_V \in \H_{\epsilon}^i ( \langle A_P^G\rangle_t \times X_P;
\nil_*\M\otimes \CC_{\hsr_P})$ is nonzero, then for at least one $V$, $\o_V\neq
0$ and $\o_V \in \H_{\epsilon}^i ( \langle A_P^G\rangle_t \times X_P;
(\nil_*\M)_V\otimes\CC_{\hsr_P})$.  For otherwise we could sum the estimates
$\|\dbar \o_V\|^2 + \| \dbar^*\o_V\|^2 >\epsilon\| \o_V \|^2$ (compare \eqref{eqnSpectralVanishingEstimate}) and obtain a
contradiction.  The converse is obvious.
\end{proof}

\begin{lem}
\label{ssecttIndependence}
The vanishing condition
$\H_{\epsilon}^i ( \langle A_P^G\rangle_t \times X_P;
(\nil_*\M)_V\otimes\CC_{\hsr_P})=0$ is independent of $t$.
\end{lem}
\begin{proof}
Let $\Psi_t\colon  \langle
A_P^G\rangle_1 \times X_P\tildearrow \langle A_P^G\rangle_t \times X_P$ be the
diffeomorphism given by the left action of $t_P$.  Although $\Psi_t$ is an
isometry, the natural isomorphism of coefficient systems
$(\EE'_{Q,1})_V\otimes \CC_{\hsr_P}\tildearrow
\Psi_t^*((\EE'_{Q,t})_V\otimes \CC_{\hsr_P})$ multiplies the norm by
$t^{-(\xi_V+\hsr_P)}$,
where $\xi_V$ is the character by which $S_P$ acts on $V$.
Thus the norm of a form depends on
$t$ only via a scaling factor independent of $Q$.  This implies that
$\dbar^* \Psi_t^* = \Psi_t^*\dbar^*$ and that the scaling factor may be
canceled from all terms of the estimate
$\| \dbar\Psi_t^*\o\|_1^2 + \| \dbar^*\Psi_t^*\o\|_1^2 \le \epsilon \|
\Psi_t^*\o \|_1^2$ to yield $\| \dbar\o\|_t^2 + \| \dbar^*\o\|_t^2 \le \epsilon \|
\o \|_t^2$.
\end{proof}
\begin{rem*}
This lemma would not hold if we had not included the weight factors in
\eqref{eqnComplexNorm}.
\end{rem*}

\subsection{}
\label{ssectParabolicAnalyticVanishingCriterion}
We finally arrive at the following vanishing criterion for the cohomology
of an $\L$-module.
\begin{prop*}
Let $\M$ be a bounded $\L$-module on $\Xhat$. Assume
\begin{equation}
\begin{gathered}
H_{(2)}^i ( \langle A_P^G\rangle_1 \times X_P;
(\nil_*\M)_V\otimes\CC_{\hsr_P})=0\ ,\text{and}\\
H_{(2)}^{i+1} ( \langle A_P^G\rangle_1 \times X_P;
(\nil_*\M)_V\otimes\CC_{\hsr_P}) \text{ is Hausdorff.}
\end{gathered}
\label{eqnParabolicAnalyticVanishingCriterion}
\end{equation}
for all parabolic $\QQ$-subgroups $P$ and  $V\in \IrrRep(L_P)$.  Then
$H^i(\Xhat;\M)=0$.
\end{prop*}

\begin{proof}
In order to prove $H^i(\Xhat;\M)=0$ it suffices by \eqref{eqnCohomologySameWithT}, Corollary
~\ref{ssectCohomologyLModuleIsLtwoCohomology}, and Proposition
~\ref{ssectSpectralVanishingLModules} to show
$\H_{\epsilon}^i(\Xhat_t;\M,t^{-\hsr})=0$ for some $\epsilon>0$ and $t\in A^G$ sufficiently dominant (see \S\ref{ssectTilings}).
For this vanishing it suffices by Lemmas
~\ref{ssectSpectralMayerVietoris}--\ref{ssecttIndependence} to show there
exists $\epsilon>0$ such that
\begin{equation*}
\H_{\epsilon}^i ( \langle A_P^G\rangle_1
\times X_P; (\nil_*\M)_V\otimes\CC_{\hsr_P})=0
\end{equation*}
for all $P$ and $V$.  However by Proposition
~\ref{ssectLModulesAPTimesXP} this is implied by
\eqref{eqnParabolicAnalyticVanishingCriterion}.
\end{proof}

\section{Vanishing Theorem for $L^2$-cohomology}
\label{sectVanishingTheoremLtwo}
In the course of verifying \eqref{eqnParabolicAnalyticVanishingCriterion}
we will need to apply the following vanishing theorem to the ordinary
global $L^2$-cohomology $H_{(2)}(X_P;\VV)$.  The heart of the theorem
is a calculation originally due to Raghunathan \cite{refnRaghunathan},
\cite{refnRaghunathanCorrection}.  The theorem was proven in
\cite{refnSaperSternTwo} but since it was not stated explicitly in this
form we sketch the proof in this section.

\begin{thm}
\label{ssectRaghunathanVanishing}
Let $X = \G\back G(\RR)/KA_G$ be the locally symmetric space associated to
an arithmetic subgroup $\G$ of a reductive algebraic group $G$ defined over
$\QQ$.  Let $E$ be an irreducible algebraic representation of $G$ and let
$\EE$ be the associated locally constant sheaf on $X$.  Endow $X$ with a
locally symmetric Riemannian metric and let $\EE$ have a metric induced
from an admissible inner product on $E$.
\begin{enumerate}
\item If $(E|_{\lsp0 G})^* \not\cong \overline {E|_{\lsp0 G}}$, then
$H_{(2)}(X;\EE)=0$.
\label{itemRaghunathanVanishingNotTauInvariant}
\item If $(E|_{\lsp0 G})^* \cong \overline {E|_{\lsp0 G}}$, then
$H_{(2)}^i(X;\EE)=0$ for
\begin{equation*}
i \notin [(\dim D - \dim D(E))/2,(\dim D + \dim D(E))/2]\ ,
\end{equation*}
where $\dim D(E)$ is defined in \textup{\S\ref{ssectVcentralizer}}.
Furthermore, $H_{(2)}^i(X;\EE)$ is Hausdorff \textup(and hence
representable by harmonic forms\textup) for $i=(\dim D - \dim
D(E))/2$.
\label{itemRaghunathanVanishingDegree}
\end{enumerate}
\end{thm}

\begin{proof}[Sketch of Proof \textup{(compare
    \protect{\cite[\S\S8--10]{refnSaperSternTwo}})}] We will determine a
  sufficient condition on a degree $i$ so that the estimate
\begin{equation}
 \|d\o\|^2 + \|d^*\o\|^2 > \epsilon\|\o\|^2
\label{eqnBasicEstimate}
\end{equation}
holds for all $\o\in A^i_c(X;\EE)$ and some $\epsilon>0$.  Since $X$ is
complete and without boundary, Proposition
~\ref{ssectSpectralVanishingLtwoCohomology} shows that this estimate for a
given $i$ implies $H_{(2)}^i( X;\EE)=0$ and $H_{(2)}^{i+1}(X;\EE)$ is
Hausdorff.

Let $K\subset G(\RR)$ be a maximal compact subgroup with associated
Cartan involution $\theta$ and Cartan decomposition $\mathfrak
g=\k+\p+\sa_G$.  Let $\sigma_1$ denote the coadjoint action of $\k_\CC$ on
$\bigwedge \p_\CC^*$ and let $\sigma_2$ denote the action of $\mathfrak
g_\CC$ on $E$.  Since $X\cong \G\back (\lsp 0 G(\RR))^0/K^0$, the complex
$A_c(X;\EE)$ may as usual be isometrically identified with the complex
\cite[I, \S\S1.2,~ 5.1, VII, \S2.7]{refnBorelWallach}, \cite[
\S4]{refnMatsushimaMurakami}
\begin{equation*}
\left\{
\begin{gathered}
C(\lsp 0\mathfrak g,\k; E \otimes A_c^0( \G\back (\lsp0G(\RR))^0)) \cong
\Hom_\k(\bigwedge \p, E \otimes A_c^0( \G\back (\lsp0G(\RR))^0)) \ , \\
d = \smash[b]{\sum_k} \epsilon(X_k)(\sigma_2(X_k) + X_k)\ ,
\end{gathered}
\right.
\end{equation*}
where $\{X_k\}$ denotes an orthonormal basis of $\p$, acting on $A_c^0(
\G\back (\lsp0G(\RR))^0)$ via differentiation by the corresponding
left-invariant vector field, and $\epsilon(X_k)$ denotes exterior
multiplication by the corresponding element of the dual basis.  Note that
the operator $d$ extends to the larger space $\bigwedge \p_\CC^* \otimes E
\otimes A_c^0( \G\back (\lsp0G(\RR))^0)$ (without the $\k$-invariance
condition) though it is no longer a differential here.

Integration by parts and some multi-linear algebra allow one to extract
0-order terms and estimate that
\begin{equation}
\|d\o\|^2 + \|d^*\o\|^2 \ge (\D_0\o,\o) \ge 0 \label{eqnPestimate}
\end{equation}
\cite[Props.~9.2~ and~9.4]{refnSaperSternTwo} for a certain operator
$\D_0$.  Explicitly
\begin{equation}
\label{eqnAlgebraicLaplacianExpression}
\D_0 = \sigma_2(C) - \frac12\sigma_2(C_\k) + \frac12\sigma_1(C_\k) -
\frac12(\sigma_1\otimes\sigma_2)(C_\k) \ ,
\end{equation}
where $C$ (resp.\ $C_\k$) denotes the Casimir operator of $\lsp0\mathfrak g$
(resp.\ $\k$).  In fact $\D_0$ is induced from the algebraic Laplacian on
$C(\lsp0\mathfrak g,\k; E)= \Hom_{\k}(\bigwedge \p, E) \subset \bigwedge
\p^*_\CC \otimes E$.

We now need some more notation.  Let $\h=\hb_G + \sa_G =
\hb_{G,\k}+\hb_{G,\p} +\sa_G$ be a fundamental $\theta$-stable Cartan
subalgebra of $\mathfrak g$.  Fix a positive system
$\Phi_\k^+=\Phi^+(\k_\CC,\hb_{G,\k\CC})$ of roots for $\k_\CC$ and consider
any $\theta$-stable positive system $\Phi^+=\Phi^+(\mathfrak g_\CC,\h_\CC)$
compatible with $\Phi_\k^+$.  (A specific choice will be made below.)  Let
$\Phi^+(\p_\CC,\hb_{G,\k\CC})$ denote the positive nonzero weights of
$\hb_{G,\k\CC}$ in $\p_\CC$, counted with multiplicity.  (A weight is
``positive'' here if it is the restriction to $\hb_{G,\k}$ of a positive
root \cite[\S10.1]{refnSaperSternTwo}.)
Clearly
\begin{equation}
\label{eqnDimensionD}
\dim D = 2\cdot\#\Phi^+(\p_\CC,\hb_{G,\k\CC}) + \dim \hb_{G,\p}\ .
\end{equation}
Let $\lambda$ be the highest weight of $E$ and recall the reductive Lie
subgroup $G(\lambda)\subset G$ defined in \S\ref{ssectWeightCentralizer}
with positive roots $\Phi^+(\mathfrak g(\lambda)_\CC,\h_\CC)=
\{\,\g\in\Phi^+\mid \g\perp \lambda\,\}$; define
$\Phi^+(\p(\lambda)_\CC,\hb_{G,\k\CC})=
\{\,\al\in\Phi^+(\p_\CC,\hb_{G,\k\CC})\mid \al\perp \lambda\,\}$ similarly.
Note that if $\lambda|_{\hb_{G,\p}}=0$ (so that $G(\lambda)$ is defined
over $\RR$ as in \S\ref{ssectVcentralizer}) then
$\Phi^+(\p(\lambda)_\CC,\hb_{G,\k\CC})$ is indeed the set of positive
$\hb_{G,\k\CC}$-weights of the split part in the Cartan decomposition of
$\mathfrak g(\lambda)$.  In this case
\begin{equation}
\label{eqnDimensionDLambda}
\dim D(\lambda) = 2\cdot\#\Phi^+(\p(\lambda)_\CC,\hb_{G,\k\CC}) + \dim
\hb_{G,\p} \ .
\end{equation}

Equation \eqref{eqnAlgebraicLaplacianExpression} shows that for any fixed
irreducible submodule $R$ of $\sigma_1\otimes\sigma_2|_\k$, $\D_0$ acts as
a certain scalar on all forms in $R\otimes A_c^0( \G\back
(\lsp0G(\RR))^0))$.  We need to determine when this scalar can be zero.  It
suffices to assume $R\subset S_1\otimes S_2$, where $S_1$ is an
irreducible submodule of $\sigma_1$ with pure degree $i$ and highest weight
$\v_1$, and $S_2$ is an irreducible submodule of $\sigma_2|_\k$ with
highest weight $\v_2$.  From the usual basis of weight vectors for
$\bigwedge^i \p_\CC^*$ we see that $\v_1$ may be expressed as $\sum_{\al\in
A^+} \al - \sum_{\al\in A^-} \al$, where $A^+$, $A^-\subset
\Phi^+(\p_\CC,\hb_{G,\k\CC})$ and
\begin{equation}
\#A^+ + \#A^- \le i \le \#A^+ + \#A^- + \dim \hb_{G,\p}\ .
\label{eqnDegreeBounds}
\end{equation}
Choose $\Phi^+$ now so that $\v_2$ (extended by zero on $\hb_{G,\p}$) is
dominant.  A calculation \cite[Prop.~10.2]{refnSaperSternTwo} using the
well-known expression in terms of highest weights for the value of the
Casimir operator acting on an irreducible representation shows that $\D_0$
is zero on $R$ precisely when the following conditions hold:
\begin{enumerate}[label=(\arabic*)]
\item the highest weight of $R$ is $\v_1+\v_2$, \label{itemWeightDecompose}
\item $\lambda|_{\hb_G} = \v_2$ (where $\v_2$ has been extended to be zero
on $\hb_{G,\p}$), and
\label{itemWeightSplitVanish}
\item $(\Phi^+(\p_\CC,\hb_{G,\k\CC})\setminus A^+)$, $A^-\subset
\Phi^+(\p(\lambda)_\CC,\hb_{G,\k\CC})$.
\label{itemWeightOrthogonal}
\end{enumerate}

If $\D_0$ acts as a nonzero scalar for all $R$ with a fixed $i$ then
\eqref{eqnPestimate} will yield \eqref{eqnBasicEstimate} and so
$H_{(2)}^i(X;\EE)=0$.  So if $H^i_{(2)}(X;\EE)\neq 0$ conditions
\itemref{itemWeightDecompose}--\itemref{itemWeightOrthogonal} must hold for
some $R$.  Since we have a fundamental Cartan subalgebra
and $\theta$-stable positive system the involution $\t_G$  from \S\ref{ssectDefineTau} must equal $\theta$ by \S\ref{ssectFundamentalCSA}.  Thus  $\t_G(\lambda|_{\hb_G}) = \lambda|_{\hb_G}$  holds by condition \itemref{itemWeightSplitVanish}.   By  Proposition~ \ref{ssectConjugateSelfContragredient}, $(E|_{\lsp0 G})^* \cong \overline {E|_{\lsp0 G}}$  so part \itemref{itemRaghunathanVanishingNotTauInvariant} of the theorem follows.  Condition
\itemref{itemWeightOrthogonal} allows one to estimate the possible maximum
and minimums values of $\#A^\pm$; inserting these into
\eqref{eqnDegreeBounds} yields that
\begin{multline}
\#\Phi^+(\p_\CC,\hb_{G,\k\CC}) - \#\Phi^+(\p(\lambda)_\CC,\hb_{G,\k\CC})
\le i \le \\
\#\Phi^+(\p_\CC,\hb_{G,\k\CC})
+\#\Phi^+(\p(\lambda)_\CC,\hb_{G,\k\CC}) + \dim \hb_{G,\p}
\ .\label{eqnCohomologyDegreeBounds}
\end{multline}
This is equivalent to the degree range in part
\itemref{itemRaghunathanVanishingDegree} by applying  \eqref{eqnDimensionD}
and \eqref{eqnDimensionDLambda} and taking the maximum over $\Phi^+$.
\end{proof}

\begin{rem*}
  One also has vanishing of $L^2$-cohomology if there are elements of $\Gamma$ in the center of $G$ which act non-trivially on $V$ and trivially on functions on $X$.  This will not be used in the current paper; our results apply to all $\Gamma$.
\end{rem*}

\section{Proof of Theorem ~\ref{ssectGlobalVanishing}}
\label{sectProofGlobalVanishing}

By Proposition ~\ref{ssectParabolicAnalyticVanishingCriterion} we need to
show for all $P$ and $V$ that
\begin{equation}
H_{(2)}^i ( \langle A_P^G\rangle_1 \times X_P;
(\nil_*\M)_V\otimes\CC_{\hsr_P})=0 \qquad\text{for $i\notin [c(\M),d(\M)]$}
\label{eqnFinalGoal}
\end{equation}
and in addition
\begin{equation}
H_{(2)}^{c(\M)} ( \langle A_P^G\rangle_1 \times X_P;
(\nil_*\M)_V\otimes\CC_{\hsr_P}) \text{ is Hausdorff.}
\label{eqnFinalGoalAdded}
\end{equation}

We may write the complex $C$ which  computes
\eqref{eqnFinalGoal} as 
\begin{equation}
\left\{
\begin{aligned}
C &= \bigoplus_{R\ge P} \dom \dbar_{R,b}\otimes
\Hom_{L_P}(V,H(\n_P^R;E_R))\ , \\
d_C &= \sum_{R\ge P} \dbar_{R,b}\otimes
(-1)^{\text{degree}} + \sum_{R' \ge R \ge P} \i_R^{R'*}\otimes
\Hom_{L_P}(V,H(\n_P^{R'};f_{RR'}))\ ,
\end{aligned}
\right.
\label{eqnNPAveragedComplex}
\end{equation}
where $\dbar_{R,b} = \dbar_{\langle A_P^R\rangle_1\times X_P,b,\VV\otimes
\CC_{\hsr_P}}$ and $\i_R^{R'*}\colon \dom \dbar_{R',b} \to \dom
\dbar_{R,b}$ restricts a form to the smaller stratum.

Let $P'\ge P$ be the parabolic $\QQ$-subgroup with type $\D_P^{P'}=
\{\,\al\in \D_P \mid \langle\xi_V+\hsr,\al\spcheck\rangle = 0\,\}$ and set
$S=(P,P')$.  Then any $R$ with $P\le R$ belongs to $[\tilde P, \tilde S]$
for a unique $\tilde P\in [P, P']$ and $\tilde S =\tilde P\vee S$:
\begin{equation*}
\vcenter{\xymatrix @ur @M=1pt @R=1.5pc @C=1pc {
{P'} \ar@{-}[rr] \ar@{-}[d] & {} & {G} \ar@{-}[d] \\
{\tilde P} \ar@{-}[r] \ar@{-}[d] & {R} \ar@{-}[r] &
{\tilde S} \ar@{-}[d] \\
{P\rlap{\qquad\qquad .}} \ar@{-}[rr] & {} & {S}
}}
\end{equation*}
(Compare \eqref{eqnParallelogramWithPtilde}.)  For each such $\tilde P$ let
$C(\tilde P)$ be the subquotient complex obtained by restricting the sum in
\eqref{eqnNPAveragedComplex} to be over those $R\in [\tilde P,
\tilde S]$.  It suffices to prove the vanishing and Hausdorff assertions
of \eqref{eqnFinalGoal} and \eqref{eqnFinalGoalAdded} for each complex
$C(\tilde P)$.  (Compare the Fary spectral sequence from Lemma
~\ref{ssectRelativeLocalCohomologySupportsSS}.)

Filter the double complex $C(\tilde P)$ by the degree of the second factor
in \eqref{eqnNPAveragedComplex}; the associated spectral sequence has
\begin{equation}
E_1^{i-k,k} = \bigoplus_{ R\in [\tilde P,\tilde S] }
H^{k}_{(2)}(\langle A_P^R\rangle_1\times
X_P;\VV\otimes\CC_{\hsr_P})\otimes
\Hom_{L_P}^{i-k}(V,H(\n_P^R;E_R)). \label{eqnRawEoneTerm}
\end{equation}
Note that we may use ordinary $L^2$-cohomology here as opposed to
$L^2$-cohomology with boundary values in view of
Lemma~\ref{ssectBoundaryAveraging}.

For $\al\in \D_P$, let $A_P^\al$ denote the one-parameter subgroup
corresponding to $\al\spcheck$.  We have a factorization $\langle
A_P^R\rangle_1 = \prod_{\al\in \D_P^R} \langle A_P^{\al}\rangle_1$; for
each $\al$, Zucker's $L^2$-K\"unneth theorem
\cite[(2.34)]{refnZuckerWarped} may be applied to correspondingly factor
out the $L^2$-cohomology of $\langle A_P^{\al}\rangle_1$ in
\eqref{eqnRawEoneTerm} provided it is Hausdorff.  To calculate the result
write $\xi_V+\hsr = \sum_{\al\in\D_P}
\langle\xi_V+\hsr,\al\spcheck\rangle\b_\al$ and recall that
\begin{equation*}
H^{\cdot}_{(2)}(\langle
A_P^\al\rangle_1;\CC_{\langle\xi_V+\hsr,\al\spcheck\rangle\b_\al} )
\cong \begin{cases} 0 & \text{if $\langle\xi_V+\hsr,\al\spcheck\rangle>0$,} \\
		M[-1] & \text{if $\langle\xi_V+\hsr,\al\spcheck\rangle=0$,} \\
		\CC   & \text{if $\langle\xi_V+\hsr,\al\spcheck\rangle<0$,}
\end{cases}
\end{equation*}
where $M$ is an infinite dimensional non-Hausdorff space
\cite[Prop. ~3.2]{refnSaperCertain}, \cite
[(2.37)--(2.40)]{refnZuckerWarped}.  Thus let $Q_V\le S$ have type
$\D_P^{Q_V} =\{\,\al\in \D_P \mid \langle\xi_V+\hsr,\al\spcheck\rangle <
0\,\}$ as in \S\ref{ssectMicroSupport} and set $\tilde Q=\tilde P\vee Q_V$.
We find that
\begin{equation*}
E_1^{i-k,k} = 
H^{k}_{(2)}(\langle A_P^{\tilde P}\rangle_1\times
X_P;\CC_{\xi_V^{-1}}\otimes\VV)\otimes
\biggl(
\bigoplus_{ R\in [\tilde P, \tilde Q] }
\Hom_{L_P}^{i-k}(V,H(\n_P^R;E_R)) \biggr).
\end{equation*}

The advantage of this expression is that the differentials in the spectral
sequence act only on the final factor.  These differentials are induced by
$H(\n_P^R;f_{RR'})$ for $\tilde P\le R\lneq R'\le \tilde Q$ and so from the
definitions in \S\ref{ssectFunctorsOnLsheaves} we see that $H^i(C(\tilde
P))$ equals
\begin{equation}
\bigoplus_k
H^{k}_{(2)}(\langle A_P^{\tilde P}\rangle_1\times
X_P;\CC_{\xi_V^{-1}}\otimes\VV)\otimes
\Hom_{L_P}(V, H^{i-k}(\i_P^* \i_{{\tilde P}*}\i_{\tilde P}^* \ihat_{\tilde Q}^!
\M))\ . \label{eqnEinfinity}
\end{equation}
By the $L^2$-K\"unneth theorem this would be isomorphic to
\begin{equation}
\bigoplus_{k_1,k_2} H^{k_1}_{(2)}(\langle A_P^{\tilde P}\rangle_1)\otimes
H^{k_2}_{(2)}(X_P;\VV)\otimes \Hom_{L_P}(V, H^{i-k_1-k_2}(\i_P^* \i_{{\tilde
P}*}\i_{\tilde P}^* \ihat_{\tilde Q}^!  \M))\ ; \label{eqnDecomposedEinfinity}
\end{equation}
if at least one of the first two $L^2$-cohomology factors were
finite-dimensional in all degrees, however we cannot assume this.  Instead
note that by the technique in \cite[(2.29)--(2.34)]{refnZuckerWarped} one
can show that if each factor of \eqref{eqnDecomposedEinfinity} vanishes
outside a certain interval of degrees, then \eqref{eqnEinfinity} vanishes
outside the sum of the intervals.  The first factor is $0$ outside $[0,\dim
\sa_P^{\tilde P}]$, the vanishing interval of the second factor is given by
Theorem ~\ref{ssectRaghunathanVanishing}, and the third factor vanishes
outside $[c_2(V,\tilde P;\M),d_2(V,\tilde P;\M)]$ by definition (see
Proposition ~\ref{ssectAlternateMicroSupport} and \eqref{eqnWeakAlternateMicrosupportDegreeRanges}).  Thus by Proposition
~\ref{ssectEqualityDegreeRanges} we have vanishing for $i\notin
[c_2(\M),d_2(\M)]= [c(\M),d(\M)]$ which proves \eqref{eqnFinalGoal}.

Finally consider \eqref{eqnFinalGoalAdded}.  If $P=\tilde P$ this follows
from the final assertion of Theorem ~\ref{ssectRaghunathanVanishing}.  If
$P \neq \tilde P$ then $H^{k_1}_{(2)}(\langle A_P^{\tilde P}\rangle_1)$
vanishes in degree $0$ (and is non-Hausdorff for $1 \le k_1 \le \dim
\sa_P^{\tilde P}$ \cite[(4.51)]{refnZuckerWarped}) and thus
\eqref{eqnEinfinity} vanishes in degree $c(\M)$ as well.  \qed

\specialsection*{Part III. Micro-support Calculations}

\section{Micro-support of Weighted Cohomology}
\label{sectMicroSupportWeightedCohomology}

\subsection{}
\label{ssectLocalWeightedCohomologyWithAndWithoutProperSupports}
Fix a weight profile $\eta$.  For $V$ an irreducible $L_P$-module, let
$T_V^\eta\ge P$ have type%
\footnote{This type corresponds to what in
\cite{refnGoreskyMacPhersonTopologicalTraceFormula} would be denoted
$I_\eta(\xi_V)$.}
\begin{equation*}
\{\,\al\in\D_P\mid \langle\xi_V-\eta_P,\b_\al\spcheck\rangle<0\,\}\ .
\end{equation*}
relative to $P$.  Note that $T_V^\eta=P$ if and only if $\xi_V\in \eta_P +
\lsp+\sa_P^{G*}$ and $T_V^\eta =G$ if and only if $\xi_V \in \eta_P
-\Int(\lsp+\sa_P^{G*})$.

\begin{lem*}
Let $V$ be an irreducible $L_P$-module occurring in $H(\i_P^*\WnC(E))$.
Then $V$ has the form $H^{\l(w)}(\n_P;E)_w$ for some $w\in W_P$.
Furthermore, the natural morphism $\WnC(E) \to \i_{G*}E$ induces
\begin{equation*}
H(\i_P^*\WnC(E))_V \cong \begin{cases}
V[-\l(w)] & \text{if $T_V^\eta = P$,} \\
0 & \text{otherwise.}
\end{cases}
\end{equation*}
\end{lem*}

\begin{proof}
Apply Proposition ~\ref{ssectLocalWeightedCohomology}.
\end{proof}

\begin{prop}
\label{ssectWeightCohomologyLocalCohomologyWithSupports}
Let $V$ be an irreducible $L_P$-module occurring in
$H(\i_P^*\ihat_Q^!\WnC(E))$ for $P\le Q$.  Then $V$ has the form
$H^{\l(w)}(\n_P;E)_w$ for some $w\in W_P$.  Furthermore,
\begin{equation*}
H(\i_P^*\ihat_Q^!\WnC(E))_V \cong \begin{cases}
V[-\l(w)-\#\D_Q] & \text{if $Q = (P,T_V^\eta)$,} \\
0 & \text{otherwise.}
\end{cases}
\end{equation*}
\end{prop}

\begin{proof}
Assume first that $T_V^\eta=P$.  Then we claim that
$H(\i_P^*\ihat_Q^!\WnC(E))_V \cong H(\i_P^*\ihat_Q^!\i_{G*}E)_V$ from which
proposition follows immediately.  To prove the claim, consider the
commutative diagram
\begin{equation*}
\xymatrixnocompile@C-4mm@R-2mm{
{\cdots} \ar[r] &
{H^i(\i_P^*\ihat_Q^!\WnC(E))_V} \ar[r] \ar[d] &
{H^i(\i_P^*\WnC(E))_V} \ar[r] \ar[d] &
{H^i(\i_P^*\jhat_{Q*}\jhat_Q^*\WnC(E))_V} \ar[r] \ar[d] &
{\cdots} \\
{\cdots} \ar[r] &
{H^i(\i_P^*\ihat_Q^!\i_{G*}E)_V} \ar[r] &
{H^i(\i_P^*\i_{G*}E)_V} \ar[r] &
{H^i(\i_P^*\jhat_{Q*}\jhat_Q^*\i_{G*}E)_V} \ar[r] &
{\cdots}
}
\end{equation*}
in which both rows are long exact sequences (obtained from
\eqref{eqnLongCompareLocalCohomologyWithSupports} by setting $Q'=G$ and
then taking $V$-isotypical components).  The middle vertical arrow is an
isomorphism by Lemma
~\ref{ssectLocalWeightedCohomologyWithAndWithoutProperSupports}.  To show
the right-hand vertical arrow is an isomorphism, consider the corresponding
map on the Mayer-Vietoris spectral sequences from Lemma
~\ref{ssectRelativeLocalCohomologySupportsSS}.  For
$H(\i_P^*\jhat_{Q*}\jhat_Q^*\WnC(E))_V$ we have
\begin{equation}
\label{eqnWeightedMayerVietoris}
E_1^{p,\cdot} = \bigoplus_{\substack{P < \tilde P \le (P,Q) \\
\#\D_P^{\tilde P} = p+1}} H(\n_P^{\tilde P};H(\i_{\tilde
P}^*\WnC(E)))_V \ .
\end{equation}
The $\tilde P$-term can be rewritten $H(\n_P^{\tilde P};H(\i_{\tilde
P}^*\WnC(E))_{\tilde V})_V$ for an irreducible $L_{\tilde P}$-module
$\tilde V$ satisfying $\xi_{\tilde V}=\xi_V|_{\sa_{\tilde P}}$ (see
\S\ref{ssectKostantDegeneration}).  By Lemma
~\ref{ssectLocalWeightedCohomologyWithAndWithoutProperSupports} this is
isomorphic to the corresponding term in $E_1$ for
$H(\i_P^*\jhat_{Q*}\jhat_Q^*\i_{G*}E)_V$ provided that $T_{\tilde
V}^\eta=\tilde P$.  But this holds since $T_{\tilde V}^\eta = T_V^\eta \vee
\tilde P$.  The claim then follows from the $5$-lemma.

In the case that $T_V^\eta > P$, Lemma
~\ref{ssectLocalWeightedCohomologyWithAndWithoutProperSupports} shows that
$H(\i_P^*\WnC(E))_V = 0$ and hence $H(\i_P^*\ihat_Q^!\WnC(E))_V \cong
H(\i_P^*\jhat_{Q*}\jhat_Q^*\WnC(E))_V[-1]$.  It also shows that the sum in
\eqref{eqnWeightedMayerVietoris} yields one copy of $V[-\l(w)]$ for each
$\tilde P \in [T_V^\eta, (P,Q)]$.  The differential between two such copies
of $V[-\l(w)]$ for adjacent $\tilde P$ will be an isomorphism by comparison
with the spectral sequence for $H(\i_P^*\jhat_{Q*}\jhat_Q^*\i_{G*}E)_V$.
Consequently $H(\i_P^*\jhat_{Q*}\jhat_Q^*\WnC(E))_V = 0$ unless
$T_V^\eta=(P,Q)$, in which case we obtain $V[-\l(w)-\#\D_Q+1]$.  This
proves the proposition.
\end{proof}

\subsection{}
\label{ssectWeightCohomologyMicroSupport}
Let $\eta$ be one of the two \emph{middle weight profiles}:
\begin{equation*}
\begin{aligned}
\text{\emph{upper middle}: } \mu &=  - \hsr + \epsilon\hsr \ , \\
\text{\emph{lower middle}: } \v &=  - \hsr \ .
\end{aligned}
\end{equation*}
Here $\epsilon>0$ is chosen sufficiently small such that  $\epsilon\hsr_R <
\chi_R$ for all $R\in\Pl$ maximal, where $\chi_R$ is a positive generator
for $X(S_R^G)$.

\begin{thm*}
Let $E$ be an irreducible algebraic $G$-module and let $\eta$ be
a middle weight profile.  The weak micro-support $\mS_w(\WnC(E))$
consists of those irreducible $L_P$-modules $V$ such that
\begin{enumerate}
\item $V=H^{\l(w)}(\n_P;E)_w$ with $w\in W_P$, and
\label{itemWHMicroPurityNilpotent}
\setcounter{saveenum}{\value{enumi}}
\item $(\xi_V+\hsr)|_{\sa_P^G}=0$.
\label{itemWHMicroPurityVanishing}
\end{enumerate}
%
%
For such $V$ and $Q\in[Q_V,Q'_V]$,
\begin{equation*}
  H(i_P^*\ihat_Q^!\WnC(E))_V =
  \begin{cases}
      V[-\ell(w)-\#\D_P] &\qquad
      \text{if $Q=Q_V=P$ and $\eta=\mu$,} \\
      V[-\ell(w)] &\qquad
      \text{if $Q=Q'_V=G$ and $\eta=\v$,} \\
      0 &\qquad \text{otherwise.}
  \end{cases}
\end{equation*}
The micro-support $\mS(\WnC(E))$ consists of those $V\in \mS_w(\WnC(E))$
such that
\begin{enumerate}[resume*]
\item $V|_{M_P}$ is conjugate self-contragredient\textup;
  \label{itemWHMicroPurityTauInvariant}
\end{enumerate}
in this case $w$ from \itemref{itemWHMicroPurityNilpotent} is fundamental.
The micro-support is nonempty if and only if $E|_{\lsp0 G}$ is conjugate
self-contragredient\textup; in this case $\emS(\WnC(E))=\{E\}$.
\end{thm*}

\begin{proof}[Proof of Theorem
    ~\textup{\ref{ssectWeightCohomologyMicroSupport}}]
  By Proposition
~\ref{ssectWeightCohomologyLocalCohomologyWithSupports}, an irreducible
$L_P$-module $V$ is in $\mS_w(\WnC(E))$ if and only if
\itemref{itemWHMicroPurityNilpotent} and
\begin{enumerate}[label=(\roman*)$'$]
\setcounter{enumi}{\value{saveenum}}
\item $(P,T_V^\eta)\in[Q_V, Q'_V]$
\label{itemOppositeTInQRange}
\end{enumerate}
\noindent
hold.  Assume that $\eta=\v=-\hsr$.  Then
\itemref{itemOppositeTInQRange} is equivalent to
\begin{enumerate}[label=(\roman*)$''$]
\setcounter{enumi}{\value{saveenum}}
\item $\langle\xi_V+\hsr,\al\spcheck\rangle<0 \quad\Rightarrow\quad
\langle\xi_V+\hsr,\b_\al\spcheck\rangle\ge 0 \quad\Rightarrow\quad 
\langle\xi_V+\hsr,\al\spcheck\rangle \le 0 \quad$
\label{itemOppositeTInQRangeInequalities}
\end{enumerate}
\noindent
for all $\al\in\D_P$.  For $R\ge P$ define
\begin{equation*}
\begin{split}
\langle \sa_P^*\rangle_R &= \Int{\sa^{*+}_R} - \lsp+\sa_P^{R*} \\
&= \{\, \zeta\in  \sa_P^* \mid \langle\zeta,\g\spcheck\rangle>0 \text{ for
$\g\in\D_R$, and} \\
&\qquad\qquad\qquad\qquad  \langle\zeta,\b_\al^R{}\spcheck\rangle\le 0
\text{ for $\al\in\D_P^R$}\,\}.
\end{split}
\end{equation*}
By Langlands's ``geometric lemmas'' \cite[IV,~\S6.11]{refnBorelWallach}
these sets form a disjoint decomposition of $\sa_P^*$ so there exists a unique $R\ge P$ such that
\begin{equation}
\xi_V+\hsr_P\in -\langle \sa_P^*\rangle_R = -\Int{\sa^{*+}_R} + \lsp+\sa_P^{R*} . \label{eqnWeightInPartition}
\end{equation}
Now for any $\al\in\D_P\setminus \D_P^R$, write
$\al\spcheck = \al\spcheck_R + \al\spcheck{}^R= \g\spcheck +
\al\spcheck{}^R$ where $\g=\al_R\in \D_R$ and $\al\spcheck{}^R \in
-\sa_P^{R+}$ (since $(\al,\al')\le 0$ for all $\al'\in \D_P^R$).   Thus \eqref{eqnWeightInPartition} implies that
$\langle\xi_V+\hsr,\al\spcheck\rangle<0$.  By the first implication of
\itemref{itemOppositeTInQRangeInequalities}, $(\xi_V+\hsr)_R\in
\lsp+\sa^{*}_R$ which contradicts \eqref{eqnWeightInPartition}.
Consequently we have $R=G$, $\xi_V+\hsr_P\in \lsp+\sa^{*}_P$, and
$G=(P,T_V^\nu)$.  But then the second implication of
\itemref{itemOppositeTInQRangeInequalities} implies $\xi_V+\hsr_P\in -
\sa^{*+}_P$ and thus \itemref{itemWHMicroPurityVanishing} must hold.
Conversely \itemref{itemWHMicroPurityVanishing} implies
\itemref{itemOppositeTInQRangeInequalities} and $Q_V=P$ and $Q'_V=G$.  Now apply Proposition
~\ref{ssectWeightCohomologyLocalCohomologyWithSupports}.

The case of $\eta=\mu$ is similar except that now \itemref{itemOppositeTInQRange} is equivalent to
\begin{enumerate}[label=(\roman*)$'''$]
\setcounter{enumi}{\value{saveenum}}
\item $\langle\xi_V+\hsr,\al\spcheck\rangle<0 \quad\Rightarrow\quad
\langle\xi_V+\hsr,\b_\al\spcheck\rangle> 0 \quad\Rightarrow\quad 
\langle\xi_V+\hsr,\al\spcheck\rangle \le 0 \quad$
\label{itemOppositeTInQRangeInequalitiesUpperWeight}
\end{enumerate}
\noindent
and we consider the disjoint decomposition of $\sa_P^*$ given by
\begin{equation*}
\begin{split}
\langle\langle \sa_P^*\rangle\rangle_R &= \sa^{*+}_R -
\Int{\lsp+\sa_P^{R*}} \\
&= \{\, \zeta\in  \sa_P^* \mid \langle\zeta,\g\spcheck\rangle\ge0 \text{ for
$\g\in\D_R$, and} \\
&\qquad\qquad\qquad\qquad  \langle\zeta,\b_\al^R{}\spcheck\rangle< 0 \text{
for $\al\in\D_P^R$}\,\}.
\end{split}
\end{equation*}
If $\xi_V+\hsr_P \in - \langle\langle \sa_P^*\rangle\rangle_R$ for $R\ge P$ we show $R=G$, $\xi_V+\hsr_P \in \Int{\lsp+\sa_P^{*}}$, and $P=(P,T_V^\mu)$.  We find \itemref{itemOppositeTInQRangeInequalitiesUpperWeight} is equivalent
to \itemref{itemWHMicroPurityVanishing} and again use  Proposition
~\ref{ssectWeightCohomologyLocalCohomologyWithSupports}.

The assertion that $w$ is fundamental when $V\in\mS(\WnC(E))$ follows from
\itemref{itemWHMicroPurityVanishing} and
\itemref{itemWHMicroPurityTauInvariant} by Proposition
~\ref{ssectBorelCasselman}.  For the final assertion, note that we can
rephrase our result as $V\in\mS(\WnC(E))$ if and only if $V\preccurlyeq_0
E$ and $V|_{M_P}$ is conjugate self-contragredient.  By Lemma
~\ref{lemPartialOrdering} the existence of such a $V$ implies $E|_{\lsp0
G}$ is conjugate self-contragredient.  Also note that if $V\in
\emS(\WnC(E))$ then $Q_V=Q'_V$ by Proposition
~\ref{ssectWeightCohomologyLocalCohomologyWithSupports} and therefore
$\langle\xi_V+\hsr,\al\spcheck\rangle\neq 0$ for all $\al\in\D_P$.  In view
of \itemref{itemWHMicroPurityVanishing} we must have $P=G$ and thus $V=E$.
\end{proof}

\section{Micro-support of Intersection Cohomology}
\label{sectPurity}
In \S\S\ref{sectPurity}--\ref{sectPurityProofPartOne}, we consider a \emph{middle perversity}
$p$, either $m(k) = \left\lfloor\frac{k-2}2\right\rfloor$ or $n(k) =
\left\lfloor\frac{k-1}2\right\rfloor$.  We will calculate the micro-support of $\mS(\IpC(E))$.  Unlike in the case of weighted cohomology, we do not obtain a calculation for the weak micro-support.

\begin{thm}
\label{ssectIHMicroPurity}
Assume the irreducible
components of the $\QQ$-root system of $G$ are of type $A_n$, $B_n$,
$C_n$, $BC_n$, or $G_2$.  Let $E$ be an irreducible algebraic $G$-module
and let $p$ be a middle perversity.  The micro-support $\mS(\IpC(E))$
consists of those irreducible $L_P$-modules $V$ such that
\begin{enumerate}
\item $V=H^{\l(w)}(\n_P;E)_w$ with $w\in W_P$ a fundamental Weyl element
for $P$ in $G$\textup;
\label{itemIHMicroPurityNilpotent}
\item $\langle\xi_V+\hsr,\al\spcheck\rangle \begin{cases} \ge 0 &\qquad
\text{for all $\al\in\D_P$
                      \textup(if $p=m$\textup),}
\\
                      \le 0 &\qquad \text{for all $\al\in\D_P$
                      \textup(if $p=n$\textup)\textup;}
\end{cases}
$
\label{itemIHMicroPurityInequalities}
\item $V|_{M_P}$ is conjugate self-contragredient.
\label{itemIHMicroPurityTauInvariant}
\end{enumerate}
For such $V$ and $Q\in[Q_V,Q'_V]$,
\begin{equation*}
  H(i_P^*\ihat_Q^!\IpC(E))_V =
  \begin{cases}
      V[-\frac12\dim \n_P -\#\D_P] &\qquad
      \text{if $Q=Q_V=P$ and $p=m$,} \\
      V[-\frac12\dim\n_P] &\qquad
      \text{if $Q=Q'_V=G$ and $p=n$,} \\
      0 &\qquad \text{otherwise.}
  \end{cases}
\end{equation*}
The essential micro-support $\emS(\IpC(E))$ consists of the maximal elements
of $\mS(\IpC(E))$ under $\preccurlyeq_0$ \textup(see
\textup{\S\ref{ssectPartialOrdering}}\textup)\textup; this is equivalent to
requiring strict inequalities in
\itemref{itemIHMicroPurityInequalities}.
\end{thm}

\begin{cor}
\label{ssectIHMicroPurityCorollary}
In the setting of Theorem ~\textup{\ref{ssectIHMicroPurity}} assume that
$E|_{\lsp0 G}$ is conjugate self-contragredient.  Then
$\emS(\IpC(E))=\{E\}$ with $c(E;\IpC(E))= d(E;\IpC(E))= 0$ and
$(\xi_V+\hsr)|_{\sa_P^G}=0$ for all $V\in \mS(\IpC(E))$.
\end{cor}
\begin{proof}
Apply Proposition ~\ref{ssectBorelCasselman}.
\end{proof}

\subsection{}
\label{ssectPurityReductions}
In the remainder of this section we will reduce the theorem to a
combinatorial vanishing result which will be proven in the following
sections.

Recall the ``combinatorial'' intersection cohomology $I_{p_w}H(U)$ (for $U$
an open constructible subset of $|\D_P|$ or $c(|\D_P|)$) defined in
\S\ref{ssectLocalIntersectionCohomology} which depends on $w\in
W$ as well as $p$.  We also consider the version with supports
$I_{p_w}H_Z(U)$ defined for $Z$ constructible and relatively closed in
$U$. This may be interpreted in terms of relative intersection cohomology:
\begin{equation*}
I_{p_w}H_Z(U) = I_{p_w}H(U,U\setminus Z)\ .
\end{equation*}

\begin{prop}
\label{ssectLocalIHWithSupports}
Let  $P\le Q\in \Pl$.  Then
\begin{equation*}
H(\i_P^*\ihat_Q^!\IpC(E)) \cong \bigoplus_{w\in
W_P} H(\n_P;E)_w\otimes I_{p_w}H_{c(|\D_P^Q|)}(c(|\D_P|))\ .
\end{equation*}
\end{prop}
\begin{proof}
There is a long exact sequence (set $Q'=G$ in
\eqref{eqnLongCompareLocalCohomologyWithSupports})
\begin{equation*}
\dots \longrightarrow
H^j(\i_P^* \ihat_Q^! \IpC(E)) \longrightarrow H^j(\i_P^* \IpC(E)) \longrightarrow H^j(\i_P^* \jhat_{Q*}\jhat_Q^*
\IpC(E)) \longrightarrow \cdots\ .
\end{equation*}
On the other hand there is a topological long exact sequence
\begin{equation*}
\dots \rightarrow
I_{p_w}H^i_{c(|\D_P^Q|)}(c(|\D_P|)) \rightarrow I_{p_w}H^i(c(|\D_P|)) \rightarrow
I_{p_w}H^i( c(|\D_P|) \setminus c(|\D_P^Q|) ) \rightarrow \cdots\ .
\end{equation*}
We need to show that if the first term of the second sequence is tensored
with $H(\n_{P};E)_w$ and then direct summed over $w\in W_P$, then the
result is isomorphic to the first term of the first sequence.  This is true
for the middle term by Proposition ~\ref{ssectLocalIntersectionCohomology}.
It is also true for the last term: compare the Mayer-Vietoris sequences
abutting to each and apply Proposition
~\ref{ssectLocalIntersectionCohomology} again.  One may check that these
isomorphisms are compatible with the maps in the long exact sequences.
Application of the $5$-Lemma concludes the proof.
\end{proof}

\subsection{}
From the proposition it follows that if an irreducible $L_P$-module $V$
occurs in $\mS(\IpC(E))$, then $V=H^{\l(w)}(\n_P;E)_w$ for some $w\in W_P$,
and this $V$ will actually occur only if $V|_{M_P}$ is conjugate
self-contragredient and $I_{p_w}H_{c(|\D_P^Q|)}(c(|\D_P|) \neq 0$ for some
$Q\in[Q_V,Q'_V]$ (see \S\ref{ssectMicroSupportNotation}).
This condition on $Q$ is equivalent to
\begin{equation}
\begin{aligned}
\langle\xi_V+\hsr,\al\spcheck\rangle \le 0 \qquad &\text{for
$\al\in\Phi(\n_P^Q;\sa_P)$,} \\ 
\langle\xi_V+\hsr,\al\spcheck\rangle \ge 0 \qquad &\text{for
$\al\in\Phi(\n_P^{P,Q};\sa_P)$.}
\end{aligned} \label{eqnWeightInequality}
\end{equation}

\subsection{}
\label{ssectSimpleBasicLemma}
Thus we need to show $I_{p_w}H_{c(|\D_P^Q|)}(c(|\D_P|)$ vanishes except for
the cases indicated in Theorem ~\ref{ssectIHMicroPurity}; this is a
combinatorial statement involving the data $\l(w_R)$ for $R\ge P$.  However
our hypotheses give us information on $V$, or equivalently, its highest
weight.  Following an argument from \cite{refnSaperSternTwo}, the simple
basic lemma below will translate the weight information into combinatorial
information; it is a generalization of a result of Casselman
\cite[Prop.~2.6]{refnCasselman} in the real rank one case.  A stronger
result will be presented in \S\ref{sectBasicLemma}.

Let $\h=\hb_P+\sa_P$ be a Cartan subalgebra of $\levi_P$ and fix a positive
system $\Phi^+$ of $\Phi(\mathfrak g_\CC,\h_\CC)$ containing
$\Phi(\n_{P\CC},\h_\CC)$.  For $\al\in\Phi(\n_P,\sa_P)$, let
$\n_\al\subset \n_P$ be the corresponding root space and write
$\Phi(\n_{\al\CC},\h_\CC)= \{\,\g\in\Phi^+\mid \g|_{\sa_P}= \al\,\}$.  For
$w\in W$, the contribution to $\l(w)$ attributable to $\al$ is $\l_\al(w) =
\#(\Phi(\n_\al,\h_\CC)\cap \Phi_w)$.  Note that $\l_\al(w)=\l_\al(w_P)$ for
$\al\in\Phi(\n_P,\sa_P)$.  Also if $w\in W_P$ then $\l(w) =
\sum_{\al\in\Phi(\n_P,\sa_P)} \l_\al(w)$.

\begin{lem*}
Let $P$ be a parabolic $\QQ$-subgroup and let $w\in W_P$.  Let
$V=H^{\l(w)}(\n_P;E)_w$ for an irreducible $G$-module $E$ and assume
$V|_{M_P}$ is conjugate self-contragredient.  For any
$\al\in\Phi(\n_P,\sa_P)$ we have\textup:
\begin{enumerate} 
\item $\langle\xi_V+\hsr,\al\spcheck\rangle \le 
0 \implies 
\l_\al(w) \ge \frac12\dim \n_\al$.
\label{itemSimpleBasicLemmaLessThan}
\item $\langle\xi_V+\hsr,\al\spcheck\rangle = 0 \implies \l_\al(w) =
\frac12\dim \n_\al$. 
\label{itemSimpleBasicLemmaEqual}
\item $\langle\xi_V+\hsr,\al\spcheck\rangle \ge 0 \implies \l_\al(w) \le 
\frac12\dim \n_\al$.
\label{itemSimpleBasicLemmaGreaterThan}
\end{enumerate}
Furthermore if $\g\in \Phi(\n_\al,\h_\CC)$, then $\Phi_w$ contains
respectively \itemref{itemSimpleBasicLemmaLessThan} at least one,
\itemref{itemSimpleBasicLemmaEqual} exactly one,
\itemref{itemSimpleBasicLemmaGreaterThan} at most one of $\g$ and
$\t'_P\g$.
\end{lem*}

\begin{proof}
Let $\lambda\in \h_\CC$ be the highest weight of $E$; by Kostant's theorem
the highest weight of $V$ is then $w(\lambda+\hsr)-\hsr$ and $\xi_V+\hsr =
w(\lambda+\hsr)|_{\sa_P}$.

Let $\g\in \Phi(\n_\al,\h_\CC)$.  Since $w(\lambda+\hsr)$ is regular and
$\lambda$ is dominant, $\g\in\Phi_w$ if and only if
$\langle w(\lambda+\hsr),\g\spcheck\rangle \le 0$.  Decompose $\hb_P =
\hb_{P,1}+\hb_{P,-1}$ according to the 
eigenvalues of $\t_P$; by Proposition~ \ref{ssectConjugateSelfContragredient} the hypothesis that
$V|_{M_P}$ is conjugate self-contragredient is equivalent to
$w(\lambda+\hsr)|_{\hb_{P,-1}} = 0$.  On the other hand,
$\g\spcheck|_{\sa_P} = c \al$ where $c>0$.  Thus
\begin{equation*}
\langle w(\lambda+\hsr),\g\spcheck\rangle = \langle
w(\lambda+\hsr)|_{\hb_{P,1}},\g\spcheck\rangle  + c\langle w(\lambda+\hsr),
\al\spcheck\rangle\ .
\end{equation*}
Replacing $\g$ by $\t'_P\g$ in this equation negates the first term and
leaves the second term unchanged, so $\langle
w(\lambda+\hsr),\al\spcheck\rangle \le 0$ implies $\langle
w(\lambda+\hsr),\g\spcheck\rangle \le 0$ or $\langle
w(\lambda+\hsr),\t'_P\g\spcheck\rangle \le 0$.  Thus either $\g$ or
$\t'_P\g$ is in $\Phi_w$.  This proves
\itemref{itemSimpleBasicLemmaLessThan} and the proof of
\itemref{itemSimpleBasicLemmaGreaterThan} is similar;
\itemref{itemSimpleBasicLemmaEqual} follows from
\itemref{itemSimpleBasicLemmaLessThan} and
\itemref{itemSimpleBasicLemmaGreaterThan}.
\end{proof}

\begin{cor*}
Under the hypotheses of the lemma, if $Q\in[Q_V,Q'_V]$ then
lemma we have
\begin{equation}
\begin{aligned}
\l_\al(w) \ge \tfrac12\dim \n_\al \qquad &\text{for 
$\al\in\Phi(\n_P^Q;\sa_P)$,} \\
\l_\al(w) \le \tfrac12\dim \n_\al \qquad &\text{for 
$\al\in\Phi(\n_P^{P,Q};\sa_P)$.}
\end{aligned} \label{eqnLengthInequality}
\end{equation}
\end{cor*}

\begin{proof}
  Use \eqref{eqnWeightInequality}.
\end{proof}

\subsection{}
\label{ssectFundamentalCombinatorialPurity}
We begin with a special case of the vanishing theorem.  The proof will 
appear in \S\ref{sectFundamentalCombinatorialPurityProof}; the key point is
that in the case $Q=P$ or $G$ the hypotheses \eqref{eqnLengthInequality} are
preserved under passing by induction from $P$ to $R>P$.

\begin{prop*}
Let $p$ be a middle perversity and let $w\in W_P$.  Let $P\le Q$ be parabolic
subgroups so that \eqref{eqnLengthInequality} is satisfied.  Also assume
that $Q=P$ or $G$.  Then
\begin{equation*}
I_{p_w}H_{c(|\D_{P}^Q|)}(c(|\D_P|))= \begin{cases} 
  \ZZ[-\#\D_P] & \text{if $Q=P$, $p=m$, and
                                 $\l(w_P)=\tfrac12\dim \n_P$,} \\
  \ZZ          & \text{if $Q=G$, $p=n$, and
                                 $\l(w_P)=\tfrac12\dim \n_P$,} \\
  0            & \text{otherwise.}
\end{cases}
\end{equation*}
\end{prop*}

\subsection{}
\label{ssectInductionLemma}
For $P<Q<G$ the analogue of \eqref{eqnLengthInequality}
will no longer remain true when we pass by induction from $P$ to $R>P$.  Nor will
$H(\n_R;E)_{w_R}$ necessarily be conjugate self-contragredient (for any
$E$) so Lemma~\ref{ssectSimpleBasicLemma} will no longer apply directly.  The
following lemma yields an alternate hypothesis that is suitable for
induction and is sufficient.

\begin{lem*}
Assume the irreducible components of the $\QQ$-root system of $G$ are of
type $A_n$, $B_n$, $C_n$, $BC_n$, or $G_2$.  Let $P$ be a proper parabolic
subgroup and let $w\in W_P$.  Let $V=H^{\l(w)}(\n_P;E)_w$ for an irreducible
$G$-module $E$ and assume $V|_{M_P}$ is conjugate self-contragredient.
Assume $Q\ge P$ is a parabolic $\QQ$-subgroup such that
\eqref{eqnWeightInequality} is satisfied.  Then there exists $T\ge P$ so
that
\begin{equation}
1\le \#\D_T \le 2 \qquad\text{and} \qquad 0\le \#\D_T^{Q\vee T},\!
\ \#(\D_T\setminus \D_T^{Q\vee T}) \le 1 \label{eqnTrestriction}
\end{equation}
and that
\begin{equation}
\begin{aligned}
\l_\al(w) \ge \tfrac12\dim \n_\al \qquad &\text{for
$\al\in\Phi(\n_P;\sa_P)$ such that $\al|_{\sa_T}\in \Phi(\n_T^{Q\vee
T};\sa_T)$,} \\
\l_\al(w) \le \tfrac12\dim \n_\al \qquad &\text{for
$\al\in\Phi(\n_P;\sa_P)$ such that $\al|_{\sa_T}\in\Phi(\n_T^{T,Q\vee
T};\sa_T)$.}
\end{aligned} \label{eqnTLengthInequality}
\end{equation}
Furthermore, if $P<Q<G$ we can assume that
\begin{equation}
\l(w_T) \neq \tfrac12\dim \n_T \text{ when  $\#\D_T=1$.} \label{eqnTNonHalf}
\end{equation}
\end{lem*}

The lemma will be proved in \S\ref{sectInductionLemmaProof}.  When $P<Q<G$
this lemma provides the hypotheses for the following theorem:

\begin{thm}
\label{ssectCombinatorialPurity}
Let $p$ be a middle perversity and let $w\in W$.  Let $P\le Q$ be parabolic
$\QQ$-subgroups and assume there exists $T\ge P$ satisfying
\eqref{eqnTrestriction}, \eqref{eqnTLengthInequality}, and
\eqref{eqnTNonHalf}.  Then
\begin{equation*}
I_{p_w}H_{c(|\D_{P}^Q|)}(c(|\D_P|))= 0.
\end{equation*}
\end{thm}

The proof will appear in \S\ref{sectPurityProofPartOne}.  

\subsection{Proof of Theorem ~\ref{ssectIHMicroPurity}}
Note that the nonvanishing cases of Proposition
~\ref{ssectFundamentalCombinatorialPurity} are relevant only when $Q_V=P$ (for
$p=m$) and $Q_V'=G$ (for $p=n$).  Together with Theorem
~\ref{ssectCombinatorialPurity}, this implies that an irreducible $L_P$-module
$V$ belongs to $\mS(\IpC(E))$ if and only if
\begin{enumerate}
\renewcommand{\theenumi}{\roman{enumi}}
\renewcommand{\labelenumi}{(\theenumi$'$)}
\item $V=H^{\l(w)}(\n_P;E)_w$ with $w\in W_P$ satisfying $\l(w) =
\frac12 \dim \n_P$;
\label{itemIHMicroPurityNilpotentPrime}
\end{enumerate}
together with \itemref{itemIHMicroPurityTauInvariant} and
\itemref{itemIHMicroPurityInequalities} of Theorem
~\ref{ssectIHMicroPurity} hold.  The last assertion of
Lemma~\ref{ssectSimpleBasicLemma} implies in this case that $\t_P'$
interchanges $\Phi_w$ and $\Phi(\n_{P\CC},\h_\CC)\setminus \Phi_w$; in
other words, $w$ is a fundamental Weyl element for $P$ in $G$ by Lemma
~\ref{ssectFundamentalParabolic}\itemref{itemTauInvariance}\itemref{itemPhiwSwap} and \S\ref{ssectFundamentalWeylElement}.  Thus
\itemref{itemIHMicroPurityNilpotent}$\Longleftrightarrow$(\ref{itemIHMicroPurityNilpotentPrime}$'$)
given \itemref{itemIHMicroPurityTauInvariant} and
\itemref{itemIHMicroPurityInequalities}.  Such a $V$ will lie in
$\emS(\IpC(E))$ furthermore if and only if $Q_V=Q_V'$, that is, if and only
if the inequalities in \itemref{itemIHMicroPurityInequalities} are all
strict.

Now if $V \preccurlyeq_0 \tilde V$ then $(\xi_V+\hsr_P)|_{\sa_P^{\tilde P}} = 0$.  Thus $V\in \emS(\IpC(E))$ implies $V=\tilde V$ due to strictness in \itemref{itemIHMicroPurityInequalities}.  Conversely if $V \in \mS(\IpC(E))$  is maximal writh respect to $\preccurlyeq_0$ then Corollary  ~\ref{ssectPartialOrderingOnMicroSupport} implies $V\in \emS(\IpC(E))$.
\qed

\section{Proof of Proposition ~\ref{ssectFundamentalCombinatorialPurity}}
\label{sectFundamentalCombinatorialPurityProof}

First assume that $\l(w_P) = \frac12\dim\n_P$.  This implies
$\l(w_R)=\frac12\dim \n_R$ for all $R\ge P$ by
\eqref{eqnLengthInequality} and so
\begin{equation*}
p_w(R)= \begin{cases} \left\lfloor\frac{\#\D_R-2}2\right\rfloor & \text{if
               $p=m$,}\\[2ex]
               \left\lfloor\frac{\#\D_R-1}2\right\rfloor & \text{if $p=n$.}
        \end{cases}
\end{equation*}
For $p=m$ we see that $p_w(R) < \#\D_R-1$ for $R<G$ and thus the link 
cohomology for the $R$-stratum (which inductively is $\ZZ[-\#\D_R 
+1]$) will always be truncated.  For $p=n$ we see that $p_w(R) \ge 0$ 
and thus the link cohomology (which inductively is $\ZZ$) will never 
be truncated.  This proves the proposition in this case.

Now assume that $\l(w_P) \neq \frac12\dim\n_P$ and consider the case
$Q=P$.  Let $S> P$ be the unique parabolic $\QQ$-subgroup which is minimal
with the property that $\l(w_S)=\frac12\dim \n_S$.  Consider the Fary
spectral sequence
\begin{equation*}
E_1^{-p,i+p} = \bigoplus_{\#\D_P^R = p}
I_{p_w}H_{|\D_P^R|^\circ}^i(|\D_P|)  \cong
\bigoplus_{\#\D_P^R = p}
I_{p_w}H_{c(\emptyset)}^i(c(|\D_R|)) \Rightarrow 
I_{p_w}H^i(|\D_P|).
\end{equation*}
The hypothesis \eqref{eqnLengthInequality} remains valid if we replace $P$
by $R$ (with $Q=R$) so by induction all the $E_1$ terms vanish except in the
case that $p=m$.  In this case, only the terms for $R\in [S, G]$ are
nonvanishing and
the spectral sequence is isomorphic to the Fary spectral sequence 
abutting to 
$I_{p_w}H^i(c(|\D_S|))$.  This vanishes by the argument in 
the previous paragraph unless $S=G$.  Thus
\begin{equation*}
I_{p_w}H(|\D_P|) = \begin{cases} 0   & \text{if $S < G$ or $p=n$,} \\
                          \ZZ & \text{if $S = G$ and $p=m$.}
                   \end{cases}
\end{equation*}
However in this second case, $p_w(P)
\ge \left\lfloor\frac{1 +\#\D_P-2}2\right\rfloor \ge 0$ and so a 
degree 0 link cohomology class will not be truncated at the 
$P$-stratum.  Thus $I_{p_w}H_{c(\emptyset)}(c(|\D_P|)) = 0$.

The argument for the case $Q=G$ is similar except that we 
use the Mayer-Vietoris spectral sequence converging to $I_{p_w}H^i(|\D_P|)$
from \eqref{eqnp}:
\begin{equation*}
E_1^{p,i-p} =\bigoplus_{\#\D_P^R=p+1} I_{p_w}H^{i-p}(U_R) \cong 
\bigoplus_{\#\D_P^R=p+1} I_{p_w}H^{i-p}_{c(|\D_R|)}(c(|\D_R|))\ . \quad
\qed
\end{equation*}

\section{Proof of Lemma ~\ref{ssectInductionLemma}}
\label{sectInductionLemmaProof}

We begin with some generalities on chambers in $\sa_P$ and an analogue of
the lemma for the $\QQ$-Weyl group $\lsb \QQ W$ of $G$ (Lemma
~\ref{ssectQInductionLemma}).  This is a purely combinatorial result
regarding the $\QQ$-root system and its Weyl group.  We then reduce Lemma
~\ref{ssectInductionLemma} to this result.  Thus $w$ in this section will
refer to an element of $\lsb \QQ W$ until \S\ref{ssectInductionLemmaProof}.

\subsection{}
\label{ssectInductionLemmaProofNotation}
Let $\lsb \QQ\Phi\subset \sa^*$ denote the $\QQ$-root system of $G$ with
simple roots $\D$.  Let $P$ be a parabolic $\QQ$-subgroup and let
$\lsb\QQ\Phi^P$ denote the subroot system with basis $\D^P$ in $\sa^{P*}$.
The Weyl group of $\lsb\QQ\Phi^P$ is the parabolic subgroup $\lsb\QQ
W^P\subset \lsb\QQ W$ generated by the simple reflections
$\{s_\al\}_{\al\in \D^P}$.  Let $\sb\QQ\Phi_P$ denote the restriction of
the elements $\lsb\QQ\Phi\setminus \lsb\QQ\Phi^P$ to $\sa_P$.  An element of
$\lsb\QQ\Phi_P$ may be expressed as a $\ZZ$-linear combination of the
elements of $\D_P$ with coefficients all $\ge0$ or $\le0$, however in
general $\lsb\QQ\Phi_P$ is not a root system.

\subsection{Chambers in $\sa_P$}
\label{ssectSplitComponentChambers}
The ``chambers'' of $\sa_P$
are the connected components of $\sa_P\setminus \bigcup_{\al\in\lsb\QQ\Phi_P}
\ker \al$.  Even though $\lsb\QQ\Phi_P$ is not a root system there is a
combinatorial object that parametrizes these chambers.  Let $\lsb\QQ
W_P$ denote the minimal length representatives of the cosets in
$\lsb\QQ W^P\backslash \lsb\QQ W$ and define%
\begin{equation}
\lsb\QQ W(P) = \{\, w\in \lsb\QQ W\mid w^{-1} \D^P \subset \D\,\}
\subset \lsb\QQ W_P
\label{eqnSplitWeyl}
\end{equation}
(cf. \cite[I.1.7]{refnMoeglinWaldspurger}%
\footnote{Our notation differs in several ways from
\cite{refnMoeglinWaldspurger}.  For one thing, the meaning of
$\lsb\QQ W_P$ and $\lsb\QQ W^P$ are reversed.  Furthermore, since we
consider right cosets of $\lsb\QQ W^P$ rather than left cosets, we
have $w^{-1}$ where \cite{refnMoeglinWaldspurger} would have $w$.}%
).  This is not a subgroup in general.

There is a nice geometric way to think of $\lsb\QQ W_P$ and $\lsb\QQ W(P)$.
Namely let $C$ denote the strictly dominant cone of $\lsb\QQ\Phi$ in $\sa$
and let $C^P$ denote the strictly dominant cone of $\lsb\QQ\Phi^P$ in
$\sa^P$.  Then $w\in \lsb\QQ W_P$ if and only if the chamber $wC$ projects
into $C^P$ under the projection $\sa\to \sa^P$.  Furthermore such an
element belongs to $\lsb\QQ W(P)$ if and only if the closure of $wC$ is
adjacent to $\sa_P$ in the sense that $w\overline C \cap \sa_P$ is a
boundary stratum of $w\overline C$ with codimension $\# \D^P$.  In this
case we let $w\cdot C_P \equiv \Int_{\sa_P}(w\overline{C}\cap \sa_P)$
denote the corresponding chamber of $\sa_P$.  (Here $\Int_{\sa_P}$ denotes
the interior of a subset of $\sa_P$.)

The correspondence
\begin{equation}
w \quad \longleftrightarrow \quad  w\cdot C_P \label{eqnUnipotentWeylChamber}
\end{equation}
is a bijection of $\lsb\QQ W(P)$ with the chambers of $\sa_P$
\cite[I.1.10]{refnMoeglinWaldspurger}.

\subsection{An Alternative characterization of $\lsb\QQ
W(P)$}
\label{ssectAllOrNone}
Let $\pr_P$ denote the map $\lsb\QQ\Phi\setminus \lsb\QQ\Phi^P \to
\lsb\QQ\Phi_P$ given by restriction to $\sa_P$; we write $\pr_P^{-1}(\al)$ for $\pr_P^{-1}(\{\al\})$.  For $w\in \lsb\QQ W$ let
$\lsb\QQ\Phi_w = \{\g \in\lsb\QQ\Phi^+| w^{-1}\g<0\}$.

\begin{lem*}
  \begin{equation*}
\lsb\QQ W(P) = \{\, w\in \lsb\QQ W_P \mid \lsb\QQ\Phi_w\cap\pr_P^{-1}(\al)
= \pr_P^{-1}(\al) \text{ or } \emptyset \text{ for all
}\al\in\lsb\QQ\Phi_P^+\,\}.
  \end{equation*}
The cases $\pr_P^{-1}(\al)$ and $\emptyset$ corresponding respectively to
$\langle \al, w\cdot C_P\rangle<0$ and $\langle \al, w\cdot C_P
\rangle>0$.
\end{lem*}

\begin{proof}
If $w\in \lsb\QQ W(P)$ and $\g\in \pr_P^{-1}(\al)$ we have
$\g\in\lsb\QQ\Phi_w\Leftrightarrow \langle \g, wC \rangle <0 \Leftrightarrow
\langle \g, w\overline{C}\rangle \le 0 \Leftrightarrow \langle \al, w\cdot
C_P \rangle <0$.  Conversely,
if $w\in \lsb\QQ W_P$ and $\lsb\QQ\Phi_w\cap\pr_P^{-1}(\al) =
\pr_P^{-1}(\al)$ or $\emptyset$ for all $\al\in\lsb\QQ\Phi_P^+$, consider
$\g\in \D^P$.  If $w^{-1}\g=\d_1+\d_2$ ($\d_1$ and $\d_2$
positive roots) is not simple, then $w\d_1>0$ (say) and $w\d_2<0$.
Thus $-w\d_2\in\lsb\QQ\Phi_w$ and $-w\d_2+\g \notin \lsb\QQ\Phi_w$,
while $\pr_P(-w\d_2) = \pr_P(-w\d_2+\g)$.  This contradiction
implies $w^{-1}\D^P\subset \D$ and thus $w\in \lsb\QQ W(P)$.
\end{proof}

\begin{lem}
\label{ssectQInductionLemma}
Assume the irreducible components of the $\QQ$-root system of $G$ are of
type $A_n$, $B_n$, $C_n$, $BC_n$, or $G_2$.  Let $P$ be a proper parabolic
$\QQ$-subgroup and let $w\in \lsb\QQ W(P)$.  Then there exists $T\ge P$
such that one of the following hold\textup:
\begin{description}
\item[Case 1] $\#\D_T= 1$.  In this case
$\lsb\QQ\Phi_w\cap\pr_T^{-1}(\lsb\QQ\Phi_T^+) = \pr_P^{-1}(\lsb\QQ\Phi_T^+)$ or $\emptyset$.
\item[Case 2] $\#\D_T= 2$.  In this case $\D_T=\{\b_-,\b_+\}$ and
\begin{equation*}
\begin{aligned}
\lsb\QQ\Phi_w\cap\pr_T^{-1}(\NN \b_-) &= \pr_T^{-1}(\NN \b_-),\\
\lsb\QQ\Phi_w\cap\pr_T^{-1}(\NN \b_+) &= \emptyset.
\end{aligned}
\end{equation*}
\end{description}
\end{lem}
\begin{rem*}
The lemma does not quite imply that
$w_T \in \lsb\QQ W(T)$ since in Case 2 nothing is asserted regarding
$\lsb\QQ\Phi_w\cap\pr_T^{-1}(\b)$ when $\b\in\lsb\QQ\Phi_T^+$ is not a multiple of a
simple root.
\end{rem*}

\begin{proof}
It suffices to assume that $\lsb\QQ\Phi$ is irreducible and
even reduced (the case $BC_n$ follows from the case $B_n$).  Let
$\widetilde G$ be the simply connected semisimple $\QQ$-split group with
root system $\lsb\QQ\Phi$.  There always exists an irreducible
representation $V$ of $\widetilde G$ whose nonzero weights form a single
Weyl orbit \cite[\S2]{refnLakshmibaiMusiliSeshadriIII},
\cite{refnKempf}.  We partially order the weights $\Phi(V)$ as usual:
$\mu\prec \v$ if $\v-\mu$ is a sum of nonnegative integral multiples of
elements of $\D$.  Then the root systems $A_n$, $B_n$, $C_n$, and $G_2$ are
precisely those for which we can arrange that $\Phi(V)$ is totally ordered.
(In the notation of the appendices of \cite{refnBourbakiLie} we choose $V$
to have highest weight $\omega_1$.)  Write $\Phi(V) = \{\e_1\succ
e_2\succ\dots\succ \e_N\}$.  Since the representation is faithful every
root $\g\in \lsb\QQ\Phi^+$ may be expressed (nonuniquely) as $\e_i-\e_j$,
where $1\le i< j\le N$.  Moreover, if $\al\in \D$ is in the support of $\g$
and we represent $\al$ as $\e_k-\e_{k+1}$ for a fixed $k$, we may choose $i$
and $j$ so that $i\le k <j$.  This can be proven in general but is easily
checked for the cases at hand.  Note that to simplify the exposition that
follows we will treat $\e_i-\e_j$ as if it were always a root; any
assertion regarding $\e_i-\e_j$ should be ignored if it is not a root.  For
$w\in \lsb\QQ W$ we will also denote by $w$ the corresponding permutation
of $\{1,2,\dots,N\}$; thus $w(\e_i-\e_j)=\e_{w(i)}-\e_{w(j)}$.

First assume that $w(1)> 1$.  Then
\begin{equation}
\begin{aligned}
\e_i - \e_{w(1)}	&\in \lsb\QQ\Phi_w	&& \text{for } i<w(1), \\
\e_{w(1)} - \e_k	&\notin \lsb\QQ\Phi_w	&& \text{for } w(1)<k.
\end{aligned} \label{eqnNegativeColumnPositiveRow}
\end{equation}
Since $w\in \lsb\QQ W_P$ it follows that $\e_{w(1)-1}-\e_{w(1)}\notin
\D^P$.  If $\e_{w(1)}-\e_{w(1)+1}$, \dots, $\e_{N-1}-\e_N\in
\D^P$, then $\pr_P(\e_i-\e_j)=\pr_P(\e_i- \e_{w(1)})$ for all $i<w(1)\le j$
and so by Lemma~\ref{ssectAllOrNone} and \eqref{eqnNegativeColumnPositiveRow} such
$\e_i-\e_j$ must belong to $\lsb\QQ\Phi_w$.  We can thus let $T$ have type
$\D\setminus\{\e_{w(1)-1}-\e_{w(1)}\}$ and Case 1 is satisfied.  Otherwise
let $b\ge w(1)$ be the least index such that $\e_b-\e_{b+1}\notin \D^P$.
Then Lemma~\ref{ssectAllOrNone} and \eqref{eqnNegativeColumnPositiveRow} imply that
\begin{equation*}
\begin{aligned}
\e_i-\e_j	&\in \lsb\QQ\Phi_w	&& \text{for } i< w(1)\le j\le b, \\
\e_j-\e_k	&\notin \lsb\QQ\Phi_w	&& \text{for } w(1)\le j \le b< k.
\end{aligned}
\end{equation*}
Thus we can let $T$ have type
$\D\setminus\{\e_{w(1)-1}-\e_{w(1)},\e_b-\e_{b+1}\}$ and Case 2 is satisfied.

Now assume that $w(1)=1$.  If $w(i)=i$ for all $i$ then
$\lsb\QQ\Phi_w=\emptyset$ and letting $T$ have type $\D\setminus
\{\e_i-\e_{i+1}\}$ for any $\e_i-\e_{i+1}\notin \D^P$ will do.  Otherwise
let $a>1$ be the least integer such that $w(a)>a$.  Since
$\e_a-\e_{w(a)}\in \lsb\QQ\Phi_w$ whereas $\e_{a-1}-\e_{w(a)}\notin
\lsb\QQ\Phi_w$ we must have $\e_{a-1}-\e_a\notin \D^P$ by
Lemma~\ref{ssectAllOrNone}.  Furthermore $\e_i-\e_j\notin \lsb\QQ\Phi_w$ for
$i<a\le j$ so we can let $T$ have type $\D\setminus\{\e_{a-1}-\e_a\}$.
\end{proof}

\subsection{Proof of Lemma ~\ref{ssectInductionLemma}}
\label{ssectInductionLemmaProof}
We return to the setting of Lemma ~\ref{ssectInductionLemma}; thus $w$ is
now an element of $W$.  Let $\lambda\in\h_\CC$ be the highest weight of
$E$; the conjugate self-contradient hypothesis on $V|_{M_P}$ and the assumption that $Q$ satisfies \eqref{eqnWeightInequality} imply that $\lambda$ satisfies
\begin{equation}
\begin{aligned}
\langle w(\lambda+\hsr),\al\spcheck \rangle  \le 0 \qquad &\text{for 
$\al\in\Phi(\n_P^Q;\sa_P)$,} \\
\langle w(\lambda+\hsr),\al\spcheck \rangle  \ge 0 \qquad &\text{for 
$\al\in\Phi(\n_P^{P,Q};\sa_P)$.}
\end{aligned} \label{eqnWeightInequalityAgain}
\end{equation}
and
\begin{equation*}
\lambda\in
	B\equiv \{\,\lambda \text{ dominant}\mid
		w(\lambda+\hsr)|_{\hb_{P,-1}}=0\,\}
\end{equation*}
(see Proposition~\ref{ssectConjugateSelfContragredient}).  The transformation
$\lambda\mapsto c\lambda+(c-1)\hsr$ for $c\ge 1$ preserves these conditions
(since it corresponds to rescaling $w(\lambda+\hsr)$ by $c$) so we may
assume that $\lambda$ is strictly dominant and hence in the interior of
$B$.  By perturbing $\lambda$ within $B$ slightly we can arrange that all
the inequalities in \eqref{eqnWeightInequalityAgain} are strict
inequalities and that in addition that $\langle w(\lambda+\hsr),\al\spcheck \rangle \neq 0$ for
all $\al\in \Phi(\n_P,\sa_P)$.

Since we have arranged that $w(\lambda+\hsr)|_{\sa_P}$ lies in a chamber of
$\sa_P$, the discussion of \S\ref{ssectSplitComponentChambers} and in
particular \eqref{eqnUnipotentWeylChamber} shows that it determines an element
$\lsb\QQ w\in \lsb\QQ W(P)$.  By \eqref{eqnWeightInequalityAgain} and
Lemma~\ref{ssectAllOrNone} we have
\begin{equation}
\begin{aligned}
\lsb\QQ\Phi_{\lsb\QQ w}\cap\pr_P^{-1}(\al) &= \pr_P^{-1}(\al)
&& \Longleftrightarrow &\langle w(\lambda+\hsr),\al\spcheck \rangle  &< 0
&& \Longleftarrow  &&\al\in\Phi(\n_P^Q;\sa_P), \\ 
\lsb\QQ\Phi_{\lsb\QQ w}\cap\pr_P^{-1}(\al) &= \emptyset
&& \Longleftrightarrow & \langle w(\lambda+\hsr),\al\spcheck \rangle  &> 0
&& \Longleftarrow  &&\al\in\Phi(\n_P^{P,Q};\sa_P). \\ 
\end{aligned} \label{eqnAllOrNoneVersusWeightInequality}
\end{equation}
for all $\al\in\Phi(\n_P,\sa_P)$.  Apply Lemma ~\ref{ssectQInductionLemma}.
If $\#\D_T=1$, re-notate the unique element $\b$ of $\D_T$ by $\b_-$ or $\b_+$
depending on whether $\lsb\QQ\Phi_{\lsb\QQ w}\cap\pr_T^{-1}(\b) =
\pr_T^{-1}(\b)$ or $\emptyset$.   In general let $\al_\pm\in \D_P$
restrict to $\b_\pm$.  The lemma and \eqref{eqnAllOrNoneVersusWeightInequality}
imply that $\al_-\in \D_P^Q$ and $\al_+\in \D_P\setminus \D_P^Q$ (which
establishes \eqref{eqnTrestriction}) and furthermore that for any
$\al\in\Phi(\n_P,\sa_P)$,
\begin{equation}
\begin{aligned}
\langle w(\lambda+\hsr),\al\spcheck \rangle  < 0 \qquad &\text{if } \al|_{\sa_T}\in \NN\b_-, \\
\langle w(\lambda+\hsr),\al\spcheck \rangle  > 0 \qquad &\text{if } \al|_{\sa_T}\in \NN\b_+.
\end{aligned} \label{eqnTWeightInequality}
\end{equation}
Lemma ~\ref{ssectSimpleBasicLemma} concludes the proof of \eqref{eqnTLengthInequality}.

The final assertion may be proved by induction.  Suppose that $P<Q<G$,
$\D_T=1$ and $\l(w_T)=\tfrac12\dim \n_T$.  Then $\l_\al(w)=\tfrac12 \dim
\n_\al$ for all $\al\in\Phi(\n_T,\sa_P)$.  If $P=Q\cap T$, then
$\D_P^T=\D_P\setminus \D_P^Q$ and so $\l_\al(w) \le \tfrac12 \dim \n_\al$
for all $\al\in\Phi(\n_P^T,\sa_P)$ by Lemma ~\ref{ssectSimpleBasicLemma}.
Then $\widetilde T$ with $\D_P^{\widetilde T}=\D_P^T\setminus \{\al\}$ for
any $\al\in\D_P^T$ has the desired properties.  An analogous argument
applies if $Q\cap T= T$.  Otherwise we can apply the lemma by induction to
the parabolic $\QQ$-subgroups $P/N_T <(Q\cap T)/N_T < T/N_T$ with $w$ and
$\lambda$ replaced by $w^T$ and $w_T(\lambda+\hsr)-\hsr$ respectively to
obtain a parabolic $\QQ$-subgroup $\widetilde T/N_T$.  Then $\widetilde T$
has the desired properties.  \qed

\section{Proof of Theorem 
~\ref{ssectCombinatorialPurity}}
\label{sectPurityProofPartOne}

If $T=P$ and $Q=P$ or $G$ we can apply Proposition
~\ref{ssectFundamentalCombinatorialPurity} to obtain the desired
vanishing; the hypothesis \eqref{eqnLengthInequality} follows from
\eqref{eqnTLengthInequality} and the nonvanishing cases are excluded by
\eqref{eqnTNonHalf}.  So in this section we can assume that if $T=P$, then
$P<Q<G$; this is necessary for Lemma ~\ref{ssectACLemmaOne} below.

\renewcommand{\IH}{IH}
\subsection{Notation}
\label{ssectNotationPurityProofOne}
In order to simplify the
notation, we will simply write $\IH(U)$ for $I_{p_w}H(U)$.  Furthermore, if
$Z\subset |\D_P|$ is a locally closed constructible subset, we write
\begin{equation*}
\IH_Z = \IH_Z(U)
\end{equation*}
where $U\subset |\D_P|$ is any open constructible subset containing $Z$
as a relatively closed subset.  This is well-defined by excision.

\subsection{}
Consider the commutative diagram
\begin{equation*}
\xymatrixnocompile@C-3mm@R-2mm{
{\cdots} \ar[r] & {\IH_{c(|\D_P^Q|)}^i(c(|\D_P|))} \ar[r]
& {\IH^i(c(|\D_P|))} \ar[r]^-{r} \ar[d]_{s} &
{\IH^i(c(|\D_P|)\setminus c(|\D_P^Q|))} \ar[r] \ar[d]_{\cong} & {\cdots} \\
{\cdots} \ar[r] &
\composite{i{\IH_{c(|\D_P^Q|)}^i(c(|\D_P|))}*{\IH^i_{|\D_P^Q|}}} \ar[r] &
\composite{i{\IH^i(c(|\D_P|))}*{\IH^i(|\D_P|)}} \ar[r] &
\composite{i{\IH^i(c(|\D_P|)\setminus c(|\D_P^Q|))}*{\IH^i(|\D_P|\setminus |\D_P^Q|)}} \ar[r] & {\cdots}
}
\end{equation*}
in which both rows are long exact sequences.  The vanishing assertion
of Theorem ~\ref{ssectCombinatorialPurity} is equivalent to $r$ being an
isomorphism.  However by the local characterization of intersection
cohomology, the map $s$ is an isomorphism in degrees $i \le p_w(P)$ and the
inclusion of $0$ in degrees $i > p_w(P)$.  Thus we need to show that
\begin{subequations}
\begin{align}
&\IH^{i-1}(|\D_P|\setminus |\D_P^Q|) = 0  &&\text{for $i >
p_w(P) + 1$,} \label{eqnPurityVanishingA} \\
&\Im \left(\IH^{p_w(P)}(|\D_P|\setminus |\D_P^Q|)
      \longrightarrow \IH_{|\D_P^Q|}^{p_w(P)+1}\right) = 0
      &&\text{for $i =
p_w(P) + 1$,} \label{eqnPurityVanishingB} \\
&\IH_{|\D_P^Q|}^i = 0 &&\text{for $i < p_w(P)+1$.}
\label{eqnPurityVanishingC}
\end{align}
\end{subequations}

To prove \eqref{eqnPurityVanishingA} we will use the long exact sequence
\begin{equation}
\dots
\longrightarrow \IH^{i-1}_{|\D_P^{Q\vee T}|\setminus |\D_P^Q|}
\longrightarrow \IH^{i-1}(|\D_P|\setminus |\D_P^Q|)
\longrightarrow \IH^{i-1}(|\D_P|\setminus |\D_P^{Q\vee T}|)
\longrightarrow \dots \label{eqng}
\end{equation}
and show the outside terms vanish in the appropriate degrees.
Equation \eqref{eqnPurityVanishingC} will follow similarly using
\begin{equation}
\dots
\longrightarrow \IH^i_{ |\D_P^Q\setminus\D_P^{Q\cap T}|}
\longrightarrow \IH^i_{|\D_P^Q|}
\longrightarrow \IH^i_{|\D_P^Q|\setminus |\D_P^Q\setminus\D_P^{Q\cap T}|}
\longrightarrow \dots \ . \label{eqnDualSequence}
\end{equation}
For \eqref{eqnPurityVanishingB} we will use these arguments to  show that
either $\IH^{p_w(P)}(|\D_P|\setminus |\D_P^Q|)=0$ or
$\IH_{|\D_P^Q|}^{p_w(P)+1} = 0$.

\begin{lem}
\label{ssectALemmaTwo}
$\IH^{i-1}_{|\D_P^{Q\vee T}|\setminus |\D_P^Q|}=0$ for all $i$.
\end{lem}

\begin{proof}
For $R> P$, let $U_R$ be the star neighborhood of the open face
$|\D_P^R|^\circ$.  Decompose
\begin{equation*}
|\D_P^{Q\vee T}|\setminus |\D_P^Q| =
\coprod_{\emptyset \neq \D_P^R \subset \D_P^T\setminus\D_P^{Q \cap T}}
           U_R\cap|\D_P^{Q\vee R}|.
\end{equation*}
The corresponding Fary spectral sequence abutting to
$\IH^{i-1}_{|\D_P^{Q\vee T}|\setminus |\D_P^Q|}$ has (compare
Lemma~\ref{ssectRelativeLocalCohomologySupportsSS})
\begin{equation*}
E_1^{-p,i-1+p} =
\bigoplus_{\substack{\emptyset \neq \D_P^R \subset 
                   \D_P^T\setminus\D_P^{Q \cap T}
\\ \#\D_P^{R} = p}}
\negmedspace \IH^{i-1}_{U_R\cap|\D_P^{Q\vee R}|} \cong
\bigoplus_{\substack{\emptyset \neq \D_P^R \subset
\D_P^T\setminus\D_P^{Q \cap T} \\
\#\D_P^{R} = p}}
\negmedspace\IH^{i-1}_{c(|\D_{R}^{Q\vee R}|)}(c(|\D_{R}|)) \ .
\end{equation*}
Theorem ~\ref{ssectCombinatorialPurity} may be applied by induction (with $P\le
Q$ replaced by $R\le Q\vee R$ and $T$ remaining the same) to prove that
this vanishes.
\end{proof}

Similarly we have the

\begin{lem}
\label{ssectCLemmaTwo}
$\IH^i_{|\D_P^Q|\setminus |\D_P^Q\setminus\D_P^{Q\cap T}|} = 0$ for 
all $i$.
\end{lem}

\begin{proof}
Cover $|\D_P^Q|\setminus |\D_P^Q\setminus\D_P^{Q\cap T}|$ by 
$\{U_\al\}_{\al\in \D_P^{Q\cap T}}$, where $U_\al$ is the open star of
the vertex $\al$.  The corresponding Mayer-Vietoris 
spectral sequence abutting to $\IH^i_{|\D_P^Q|\setminus 
|\D_P^Q\setminus\D_P^{Q\cap T}|}$ has
(compare Lemma~\ref{ssectRelativeLocalCohomologySupportsSS})
\begin{equation*}
E_1^{p,i-p} =
\bigoplus_{\substack{\emptyset \neq \D_P^R \subset \D_P^{Q\cap T}
\\ \#\D_P^{R} = p+1}}
IH^{i-p}_{U_R\cap|\D_P^Q|} \cong
\bigoplus_{\substack{\emptyset \neq \D_P^R \subset \D_P^{Q\cap T}
\\ \#\D_P^{R} = p+1}}
\IH^{i-p}_{c(|\D_{R}^{Q}|)}(c(|\D_{R}|)) \ .
\end{equation*}
Theorem ~\ref{ssectCombinatorialPurity} may be applied by induction (with $P\le
Q$ replaced by $R\le Q$ and $T$ remaining the same) to prove that
this vanishes.
\end{proof}

\subsection{}
\label{ssectPerversityShift}
In view of the preceding two lemmas, Lemma ~\ref{ssectACLemmaOne} below
will conclude the proof of the theorem.  First we need a useful formula.
For a perversity $p$ let
\begin{equation*}
dp(k) = p(k+1)-p(k)\qquad \text{and} \qquad \ModTwo(k) = 
        \begin{cases} 1 &
	\text{if $k$ odd,} \\
	0 & \text{if $k$ even,}
	\end{cases}
\end{equation*}
Thus for the middle perversities we have $dm=\ModTwo$ and $dn=1-\ModTwo$.

\begin{lem*}
For $p$ a middle perversity,
$w\in W$, and $P\le R$,
\begin{equation*}
\begin{split}
p_w(P) &= p_w(R) +
\left( 
	\tfrac12(\dim\n_P^R+\#\D_P^R) - 
		\l(w_P^R)
\right)  \\
&\qquad +
\ModTwo(\dim\n_P^R+\#\D_P^R)\left(\tfrac12-dp(\dim\n_P+\#\D_P)\right).
\end{split}
\end{equation*}
\end{lem*}

\begin{proof}
This is a simple verification using
$\l(w_P)=\l(w_P^R)+\l(w_R)$ and the definition of a middle perversity.
\end{proof}

\begin{lem}
\label{ssectACLemmaOne}
In addition to the hypotheses of Theorem
~\textup{\ref{ssectCombinatorialPurity}}, assume that if $T=P$, then
$P<Q<G$.  Then the following vanishing results hold\textup:
\begin{enumerate}
\renewcommand{\theenumi}{\alph{enumi}}
\renewcommand{\labelenumi}{(\theenumi)}
\item
$\IH^{i-1}(|\D_P|\setminus |\D_P^{Q\vee T}|)=0 \text{ for }
	{\begin{cases} i > p_w(P)+1\text{, or} \\
		i= p_w(P)+1\text{ and} \\
		\qquad dp(\dim\n_P+\#\D_P)=0.
	\end{cases}} $
\label{itemACLemmaA}
\smallskip
\addtocounter{enumi}{1}
\item
$\IH^i_{ |\D_P^Q\setminus\D_P^{Q\cap T}|}=0    \text{ for }
	{\begin{cases} i= p_w(P)+1\text{ and }dp(\dim\n_P+\#\D_P)=1\text{, or}\\
		i < p_w(P)+1.
	\end{cases}} $
\label{itemACLemmaC}
\end{enumerate}
\end{lem}

\begin{proof}
Consider part \itemref{itemACLemmaA}.  Either $\D_P \setminus \D_P^{Q\vee T}$ is empty
(in which case the lemma is trivial) or by \eqref{eqnTrestriction} it has
just one element.  Let $R=(P,Q\vee T)$ be the corresponding parabolic
subgroup with $\#\D_P^R=1$.  We know that $R\neq G$ since otherwise
$Q\vee T=P$ and hence $Q=T=P$.  Thus $\IH^{i-1}(|\D_P|\setminus
|\D_P^{Q\vee T}|) \cong
\IH^{i-1}(c(|\D_R|))$ will be zero due to truncation at the cone point if
$i-1 > p_w(R)$.  By Lemma ~\ref{ssectPerversityShift} this can be re-expressed
as
\begin{equation}
\begin{split}
&\bigl[ i-(p_w(P)+1) \bigr] +
\bigl[ \tfrac12(\dim\n_P^R+1) - \l(w_P^R) \bigr] \\
&\qquad +
\bigl[ \ModTwo(\dim\n_P^R+1)\left(\tfrac12-dp(\dim\n_P+\#\D_P)\right) \bigr]
> 0.
\end{split} \label{eqnNonvanishingRange}
\end{equation}
Since $\l(w_P^R)\le \tfrac12\dim\n_P^R$ by \eqref{eqnTLengthInequality}, the
second bracketed term of \eqref{eqnNonvanishingRange} is at least $\frac12$.
We are considering $i\ge p_w(P)+1$ so the first term is nonnegative
(and at least $1$ if $i> p_w(P)+1$), while the third term is at least
$-\tfrac12$ (and nonnegative if $dp(\dim\n_P+\#\D_P)=0$).  This proves
part \itemref{itemACLemmaA}.

Part \itemref{itemACLemmaC} is similar.  Let $R$ have type $\D_P^R =
\D_P^Q\setminus\D_P^{Q\cap T}= \D_P^{Q\vee T}\setminus \D_P^T$.  Again
$R\neq G$ since otherwise $Q=G$ and $T=P$.  Then $\IH^i_ {
|\D_P^Q\setminus\D_P^{Q\cap T}|}\cong \IH^i_{c(\emptyset)}(c(|\D_{R}|))$
will be zero if $i-1 \le p_w(R)$.  The application of Lemma
~\ref{ssectPerversityShift} and the inequality $\l(w_P^R)\ge
\tfrac12\dim\n_P^R$ from \eqref{eqnTLengthInequality} concludes the proof.
\end{proof}

\specialsection*{Part IV. Satake Compactifications and Functoriality of
Micro-support}

In this part we prove our main result, Theorem~
\ref{ssectRestrictMicroSupportToFiber}, which describes the micro-support of the restriction of an $\L$-module $\M$ to the fibers of $\pi\colon \Xhat\to \Xstar_\sigma$; here $\Xstar_\sigma$ is an equal-rank Satake compactification.  As an application we deduce the Rapoport and Goresky-MacPherson conjecture.

We begin with \S\ref{sectSatakeCompactifications} which introduces Satake compactifications
$\Xstar_\sigma$ of $X$ and recalls Zucker's result
\cite{refnZuckerSatakeCompactifications} that there is a natural quotient
map $\pi\colon \Xhat\to \Xstar_\sigma$.  Then in \S\S\ref{sectFunctorialityInverseImages} and \ref{sectFunctorialityMicroSupportEqualRank}
we discuss the functorial behavior of
the micro-support of an $\L_W$-module under $k^*$ and $k^!$, where
$k\colon Z\hookrightarrow W$ is an inclusion of admissible spaces.  The resulting degree bounds for these micro-supports involve the lengths of certain Weyl group elements $\l(w)=[\tilde V:V]$ which are estimated by the ``Basic Lemma'' in \S\S\ref{sectBasicLemma} and \ref{sectEqualRankBasicLemma}.

\section{Satake Compactifications}
\label{sectSatakeCompactifications}
We will briefly outline the theory of Satake compactifications
$\Xstar_\sigma$ and give Zucker's realization
\cite{refnZuckerSatakeCompactifications} of $\Xstar_\sigma$ as a quotient
space of $\Xhat$.  References are \cite{refnSatakeCompact},
\cite{refnSatakeQuotientCompact}, \cite{refnBorelEnsemblesFondamentaux},
\cite{refnZuckerSatakeCompactifications}, and
\cite{refnCasselmanGeometricRationality}.

\subsection{Notation}
\label{ssectSatakeCompactificationsNotation}
Let $S\subset \lsb\RR S$ be maximal tori of $G$ split over $\QQ$ and
$\RR$ respectively and choose compatible orderings on the roots.  For
$k=\QQ$ or $\RR$ let $\lsb k\D$ denote the simple restricted $k$-roots; as
usual we omit the left subscript when $k=\QQ$.  We will be working with an
algebraic representation $\sigma\colon G \to GL(U)$.  Let $\mu$ denote the highest
weight of $\sigma$ and let $\lsb k\mu=\mu|_{\lsb k\sa}$ be the restricted
highest $k$-weight.  A subset $\Th\subset \lsb k\D$ is said to be
\emph{$\sigma$-connected} if $\Th\cup \{\lsb k\mu\}$ is connected.
Here a finite subset $\Th$ of an inner product space is said to be
\emph{connected} if the associated graph (with vertices $\Th$ and an
edge between $\al$ and $\b$ whenever $(\al,\b)\neq 0$) is connected.  For
any subset $\Th\subset\lsb k\D$, let $\kap(\Th)$ be the largest
$\sigma$-connected subset of $\Th$.  If $\Th$ is a $\sigma$-connected
subset of $\lsb k\D$, let $\o(\Th)$ be the largest subset of $\lsb k\D$
such that $\kap(\o(\Th))=\Th$.
As in
\cite{refnCasselmanGeometricRationality} we extend this notation to
arbitrary subsets $\Th\subset \lsb k\D$ by setting
$\o(\Th)=\o(\kap(\Th))$;
it is easy to see that $\o(\Th)\supset \Th$.

If $P$ is a parabolic $\RR$-subgroup, let $P^\dag$ be the parabolic
$\RR$-subgroup containing $P$ with $\RR$-type $\o(\DRR^P)$.  This is the
unique largest parabolic $\RR$-subgroup containing $P$ with
$\kap(\DRR^{P^\dag})=\kap(\DRR^P)$.  In general a parabolic $\RR$-subgroup
will be called \emph{$\sigma$-saturated} if it has $\RR$-type
$\o(\Th)$ for some $\Th\subset \DRR$.

Let $\restr\colon \DRR\to \D\cup\{0\}$
denote%
\footnote{Our use in this section of this standard notation should not be
confused with our equally standard use elsewhere of $\hsr$ to denote one half
the sum of the positive roots.}
the restriction map on simple roots.  For a subset $\Th\subset \DRR$ set
$\Th\sphat= \restr(\Th)\setminus\{0\}$ and for a subset $\U\subset \D$
set $\Utilde= \restr^{-1}(\U\cup\{0\})$.  Subsets of $\DRR$ of the form
$\Utilde$ are called \emph{$\QQ$-rational} and are precisely the
$\RR$-types of the parabolic $\QQ$-subgroups.  Clearly $(\Utilde)\sphat=\U$
and $\widetilde{(\Th\sphat)}$ is the smallest $\QQ$-rational set containing
$\Th$.

We will frequently make use of the fact that a subset $\U\subset \D$ is $\sigma$-connected if and only if there exists a $\sigma$-connected $\Psi\subset \DRR$ with $\Psi\sphat = \U$ \cite[(2.4)]{refnZuckerSatakeCompactifications}.

\subsection{Satake compactifications of $D$}
\label{ssectSatakeCompactificationsSymmetricSpaces}
In this subsection and the next we depart from our usual context and assume
that $G$ is merely defined over $\RR$ with $K\subset G(\RR)$ a maximal
compact subgroup and that $D$ denotes the associated symmetric space
$G(\RR)/K\lsb\RR A_G$.

Let $\sigma\colon G \to \GL(U)$ be an irreducible algebraic representation of $G$
which is nontrivial on every simple factor of $G(\RR)$ and fix an
admissible metric $h_0$ on $U$ with respect to $K$ as in
\S\ref{ssectMetrics}.  Let $T^*$ denote the adjoint with respect to $h_0$
of an endomorphism $T$ of $U$ and let $S(U)$ denote the real vector space
of self-adjoint endomorphisms.  The representation $\sigma$ induces a
natural representation $\tilde \sigma\colon G(\RR) \to \GL(\PP S(U))$ given by
\begin{equation*}
\tilde\sigma(g)[T]= [\sigma(g)\circ T \circ\sigma(g)^*]\ .
\end{equation*}
The class of the identity
endomorphism is fixed by $K\lsb\RR A_G$ under $\tilde\sigma$ and thus the
map $g \mapsto [\sigma(g)\sigma(g)^*]$ descends a $G(\RR)$-equivariant map
\begin{equation*}
D \longrightarrow \PP S(U)\ .
\end{equation*}
This map is an embedding; the closure of the image is denoted $\lsb\RR
\Dstar_\sigma$, the \emph{Satake compactification%
\footnote{Here we have followed most closely Zucker
\cite{refnZuckerSatakeCompactifications}; Borel
\cite{refnBorelEnsemblesFondamentaux} uses an action on the right for $U$
and $S(U)$ and hence replaces $\mu$ by the lowest weight of $U$.  Satake's
construction \cite{refnSatakeCompact} is identical except that he
represents $S(U)$ by Hermitian matrices.  Casselman
\cite{refnCasselmanGeometricRationality} dispenses with $U$ altogether and
works directly with an irreducible representation $V$ containing a
$K$-fixed vector.  In \cite[\S3]{refnSaperGeometricRationality} it is
shown that given $U$ there exists $V$ (in fact contained in $S(U)_\CC$)
such that Casselman's construction agrees with Satake's construction.  The
highest weight of this $V$ is $2\lsb\RR\mu$ and hence the notion of
$\sigma$-connected defined here agrees with that in
\cite{refnCasselmanGeometricRationality} for subsets of $\D$ or $\DRR$.}
of $D$ associated to $\sigma$}.  The Satake compactification is
independent of the choice of $K$ and $h_0$.  If $D$ is Hermitian symmetric then the
closure of the Harish-Chandra embedding of $D$ as a bounded domain is
called the \emph{natural compactification} and it is topologically
equivalent to the Satake compactification for a certain $\sigma$.

\subsection{Real boundary components}
\label{ssectRealBoundaryComponents}
A \emph{real Satake boundary component} $D_{P,h}\subset \lsb\RR
\Dstar$ is the set of points fixed under the action of $N_P(\RR)$ for some
parabolic $\RR$-subgroup $P$.  The action of $P(\RR)$ descends to an action
of $L_P(\RR)$ on $D_{P,h}$.  The centralizer of $D_{P,h}$ under this action
is the group of real points of a normal reductive $\RR$-subgroup
$L'_{P,\l}\subset L_P$.  The action then descends further to the real
points of $L'_{P,h} = L_P/L'_{P,\l}$, and thus $D_{P,h}$ is realized as
$L'_{P,h}(\RR)/(K_P\cap L'_{P,h}(\RR))$, the symmetric space associated to
$L'_{P,h}(\RR)$.  (We use the subscripts $h$ and $\l$ to mimic the notation
frequently used in the Hermitian case.)

We can describe these groups concretely in terms of roots.  Set
$\Th=\kap(\DRR^P)$; then the orthogonal decomposition $\DRR^P = \Th \coprod
(\DRR^P\setminus\Th)$ yields a decomposition of the Levi factor as an
almost direct product,
\begin{equation}
L_P = \widetilde{L'_{P,h}} L'_{P,\l}. \label{eqnLevisFactors}
\end{equation}
Here $\widetilde{L'_{P,h}}$ is the minimal normal connected semisimple
$\RR$-subgroup with $\Th$ as simple $\RR$-roots; it projects onto $L'_{P,h}$ with finite kernel.  The group $L'_{P,\l}$ is the almost direct product of the
minimal normal connected semisimple $\RR$-subgroup with
$\DRR^P\setminus\Th$ as simple $\RR$-roots together with the center of
$L_P$ and all almost simple factors of $L_P$ with compact groups of real
points \cite[\S5.11]{refnBorelTits}.

The parabolic $\RR$-subgroup $P$ giving rise to a boundary component
$D_{P,h}$ as above is not unique.  To index $D_{P,h}$ in a unique fashion
we use its normalizer $\{\,g\in G(\RR) \mid gD_{P,h}=D_{P,h}\,\}$; this is
the group of real points of a $\sigma$-saturated parabolic $\RR$-subgroup,
namely $P^\dag$.  Recall that $P^\dag$ is the unique parabolic
$\RR$-subgroup containing $P$ with $\RR$-type $\o(\DRR^P)$.  For later use
we note that the centralizer $\{\,g\in G(\RR) \mid gy=y \text{ for all
$y\in D_{P,h}$}\,\}$ of $D_{P,h}$ is the group of real points of the
inverse image of $L'_{P^\dag,\l}$ under $P^\dag\to L_{P^\dag}$.  On the
other hand, if $R$ is a $\sigma$-saturated parabolic $\RR$-subgroup,
$R(\RR)$ is the normalizer of a unique boundary component $D_{R,h}$, namely
the points fixed by $N_R(\RR)$.  We have
\begin{equation*}
\lsb\RR \Dstar_\sigma = \coprod_{\substack{\text{$R$ $\sigma$-saturated}\\
				\text{parabolic $\RR$-subgroup}}} D_{R,h}\ .
\end{equation*}
Define a partial order on the boundary components by setting $D_{R,h}\le
D_{R',h}$ if and only if $R\cap R'$ is a parabolic $\RR$-subgroup and $\kap(\DRR^R)\subset \kap(\DRR^{R'})$ \cite[\S2]{refnSaperGeometricRationality}.  This corresponds
to the condition that $D_{R,h}\subset \cl{D_{R',h}}$.

\subsection{Geometrically rational Satake compactifications of $D$}
\label{ssectRationalBoundaryComponents}
Now consider again a connected reductive algebraic $\QQ$-group $G$ with
$D=G(\RR)/KA_G$.  We assume that $\lsb\RR S_G = S_G$.  (This is not
necessary but simplifies the exposition; the assumption is satisfied if $D$
is Hermitian or $\lsp0 G$ is equal-rank.)  Given an appropriate representation $\sigma$ of $G$, the preceding two
subsections may be applied to yield a Satake compactification $\lsb\RR
\Dstar_\sigma$.

A real boundary component $D_{R,h}$ of $\lsb\RR \Dstar_\sigma$ is called
\emph{rational} \cite[\S\S3.5,~3.6]{refnBailyBorel} if
\begin{enumerate}
\item its normalizer $R$ is defined over $\QQ$, and
\label{itemRationalNormalizer}
\item
\label{itemRationalCentralizer}
the group $L'_{R,\l}$ contains a normal subgroup $L_{R,\l}$ of $L_R$
defined over $\QQ$ such that $L'_{R,\l}(\RR)/L_{R,\l}(\RR)$ is compact.%
\footnote{In \cite{refnBorelEnsemblesFondamentaux} it is assumed that
$L'_{R,\l}$ is defined over $\QQ$, but as noted in
\cite{refnBailyBorel} this is too restrictive.}
\end{enumerate}
Here Condition ~\itemref{itemRationalNormalizer} is equivalent to
$\G_{N_R}\back N_R(\RR)$ being compact while Condition
~\itemref{itemRationalCentralizer} ensures that $D_{R,h}$ may be realized
as the symmetric space associated to the reductive algebraic $\QQ$-group
$L_{R,h}\equiv L_R/L_{R,\l}$.

The real boundary components that lie in the closure of a Siegel set are, by  \cite[Lemma
~8.1]{refnCasselmanGeometricRationality}, precisely those whose normalizers $R$ have type
\begin{equation}
\label{eqnSiegelBoundaryComponents}
\DRR^R = \o(\Utilde) \qquad \text{for some $\U\subset \D$.}
\end{equation}
Following Casselman \cite{refnCasselmanGeometricRationality} we say the Satake compactification $\lsb\RR
\Dstar_\sigma$ is \emph{geometrically rational} if every real
boundary component in the closure of a Siegel set (that is, satisfying \eqref{eqnSiegelBoundaryComponents}) is rational in the above sense.%
\footnote{The assumption that Condition
~\itemref{itemRationalNormalizer} holds for the real boundary components with normalizers
satisfying \eqref{eqnSiegelBoundaryComponents} is the assertion that
$\o(\Utilde)$ is $\QQ$-rational for any subset
$\U\subset \D$.  A different assumption is made in
\cite[(3.3)]{refnZuckerSatakeCompactifications} (Assumption 1) and used in  \cite[Prop.~3.3(ii)]{refnZuckerSatakeCompactifications} to prove $\QQ$-rationality.  However, as pointed out in \cite[\S9]{refnCasselmanGeometricRationality}, this proposition is
incorrect.}
Casselman \cite{refnCasselmanGeometricRationality}
gives a necessary and sufficient criterion for geometric rationality in
terms of $\lsb\RR\mu$ and the Tits index of $G$.  
Examples of geometrically rational Satake compactifications include those
where $\sigma$ is defined over $\QQ$ \cite[Thm.~ 8]{refnSaperGeometricRationality} and the Baily-Borel case where $D$ is a Hermitian symmetric
space and $\lsb\RR \Dstar_\sigma$ is the natural compactification
\cite{refnBailyBorel}.  Geometric rationality also holds if all real boundary components $D_{R,h}$ of $\lsb\RR \Dstar_\sigma$ are equal-rank (and an extra condition in $\QQ$-rank $1$ and $2$) \cite[Theorems
~20, 22]{refnSaperGeometricRationality}.

By \cite[Prop.~ 3.3(i)]{refnZuckerSatakeCompactifications} any real boundary component
satisfying Condition ~\itemref{itemRationalNormalizer} satisfies
\eqref{eqnSiegelBoundaryComponents}.  Thus for geometrically rational Satake compactifications, Condition ~\itemref{itemRationalNormalizer} in the definition of a rational boundary component implies Condition ~\itemref{itemRationalCentralizer}.  It follows that the rational boundary components of a geometrically rational Satake
compactification are indexed by the $\sigma$-saturated parabolic
$\QQ$-subgroups.  We now show these are precisely the
parabolic $\QQ$-subgroups with $\QQ$-type $\o(\U)$ for $\U\subset \D$.

\begin{lem*}
  Let $\Th\subset \DRR$ and assume $\Th\subset \Xi\subset \o(\Th)$ for some $\QQ$-rational $\Xi$.  Then
  \begin{equation}
    \kap(\Th\sphat) = \kap(\Th)\sphat \ .
  \end{equation}
\end{lem*}
\begin{proof}
Since $\kap(\Th)\sphat$ is $\sigma$-connected we have $\kap(\Th)\sphat \subset \kap(\Th\sphat)$.  On the other hand, $\kap(\Th\sphat)$ is $\sigma$-connected so there exists a $\sigma$-connected $\Psi\subset \DRR$ such that $\Psi\sphat =  \kap(\Th\sphat)$.  Thus $\Psi\subset \widetilde{\kap(\Th\sphat)} \subset \widetilde{\Th\sphat} \subset \widetilde{\Xi\sphat}$.  Since $\Xi$ is $\QQ$-rational it follows that $\Psi\subset \Xi \subset \o(\Th)$ and since $\Psi$ is $\sigma$-connected we obtain $\Psi\subset \kap(\Th)$.  Thus $\kap(\Th\sphat) = \Psi\sphat \subset \kap(\Th)\sphat$.
\end{proof}

\begin{prop*}
If $\lsb\RR \Dstar_\sigma$ is geometrically rational,
\begin{alignat}{2}
\label{eqnQQsaturatedImpliesRRsaturated}
\widetilde{\o(\U)}&= \o(\Utilde) &\qquad &\text{for $\U\subset \D$,} \\
\label{eqnRationalRRsaturatedImpliesQQsaturated}
\o(\Th)\sphat&= \o(\Th\sphat) &\qquad &\text{for $\Th\subset \DRR$ with
$\o(\Th)$ $\QQ$-rational.}
\end{alignat}
\end{prop*}
\begin{proof}
To prove \eqref{eqnQQsaturatedImpliesRRsaturated} we first claim that
$\kap(\widetilde{\o(\U)})= \kap(\Utilde)$.  Since $\widetilde{\o(\U)}$
contains $\kap(\Utilde)$ it follows that $\kap(\widetilde{\o(\U)})\supset
\kap(\Utilde)$.  On the other hand, $\kap(\widetilde{\o(\U)})\sphat$ is
$\sigma$-connected and contained in $\o(\U)$ and therefore
$\kap(\widetilde{\o(\U)})\sphat \subset \kap(\U) \subset \U$; it
follows that $\kap(\widetilde{\o(\U)}) \subset \kap(\Utilde)$ which
finishes the proof of the claim.  The claim implies that
$\widetilde{\o(\U)}\subset \o(\Utilde)$.  For the opposite inclusion note both sets are $\QQ$-rational (for $\o(\Utilde)$ use geometric rationality) and so it suffices to
prove that $\o(\U)\supset \o(\Utilde)\sphat$.  However the lemma (applied twice with $\Th=\Xi$ being $\Utilde$ and $\o(\Utilde)$) yields $\kap(\o(\U)) = \kap(\Utilde)\sphat = \kap(\o(\Utilde)\sphat)$ which implies the desired inclusion.  For
\eqref{eqnRationalRRsaturatedImpliesQQsaturated} suppose $\o(\Th)$ is
$\QQ$-rational for $\Th\subset \DRR$ and set $\U=\kap(\Th)\sphat$.  By
\cite[Prop.~3.3(i)]{refnZuckerSatakeCompactifications},
$\kap(\Th)=\kap(\Utilde)$ and hence $\o(\Th)=\o(\Utilde)=
\widetilde{\o(\U)}$ (use \eqref{eqnQQsaturatedImpliesRRsaturated} for the
last equality).  Thus $\o(\Th)\sphat= \o(\kap(\Th)\sphat)$; the lemma (with $\Xi=\o(\Th)$) completes the proof.
\end{proof}

\begin{cor*}
  Let $\lsb\RR \Dstar_\sigma$ be a geometrically rational Satake compactification of $D$.  Let $P$ be a parabolic $\QQ$-subgroup and let $R=P^\dag$.
  \begin{enumerate}
  \item \label{itemSaturated} $P$ is $\sigma$-saturated if and only if $\D^P=\o(\U)$ for some $\U\subset \D$.
  \item \label{itemSaturationRational} $R$ is a parabolic $\QQ$-subgroup and $\D^{R}=\o(\D^P)$.
  \item \label{itemKapOfSaturationRational}$\kap(\D^{R}) = \kap(\D^P)$.
  \end{enumerate}
\end{cor*}

\begin{proof}
  Since $P$ is defined over $\QQ$, $\lsb\RR\D^P =\widetilde{\D^P}$ is $\QQ$-rational.  For \itemref{itemSaturated}, if $P$ is $\sigma$-saturated then $\lsb\RR\D^P = \o(\Theta)$ for some $\Theta\subset \lsb\RR\D$ and so $\D^P=(\lsb\RR\D^P)\sphat = \o(\Theta)\sphat= \o(\Theta\sphat)$ by \eqref{eqnRationalRRsaturatedImpliesQQsaturated}.  Conversely if $\D^P=\o(\U)$ then $\sb\RR\D^P = \widetilde{\D^P} = \widetilde {\o(\U)} = \o(\Utilde)$ by \eqref{eqnQQsaturatedImpliesRRsaturated}.  For \itemref{itemSaturationRational}, $\lsb\RR\D^R = \o(\lsb\RR\D^P)=\o(\widetilde{\D^P}) =\widetilde{\o(\D^P)}$ by  \eqref{eqnQQsaturatedImpliesRRsaturated} which shows $R$ is defined over $\QQ$ and also $\D^R=(\lsb\RR\D^R)\sphat
  = \o(\D^P)$.  For \itemref{itemKapOfSaturationRational} apply $\kap$ to \itemref{itemSaturationRational}.
\end{proof}

For a rational boundary component $D_{R,h}$ of a geometrically rational Satake compactification, equation \eqref{eqnLevisFactors} may be replaced by an almost direct product of
$\QQ$-subgroups
\begin{equation}
L_R = \widetilde{L_{R,h}} L_{R,\l} \label{eqnSaturatedLevisRationalFactors}
\end{equation}
and we set $L_{R,h}=L_R/L_{R,\l}$.  If $P\subset R$ is a parabolic $\QQ$-subgroup with $P^\dag= R$ the corollary above implies that $P/N_R = \widetilde{L_{R,h}}P_\l$ where $P_\l$ is a parabolic $\QQ-$subgroup of $L_{R,\l}$.  Thus there is an almost direct product of $\QQ$-subgroups
\begin{equation}
  L_P = \widetilde{L_{P,h}} L_{P,\l} \label{eqnLevisRationalFactors}
\end{equation}
where $\widetilde{L_{P,h}} =  \widetilde{L_{R,h}}$ and $L_{P,\l}= P_\l/N_P^R$.  

\subsection{Satake compactifications of $X$}
\label{ssectSatakeCompacticationsOfX}
Assume $\lsb\RR \Dstar_\sigma$ is a geometrically rational Satake
compactification of $D$.  Let $\Dstar_\sigma$ be the union of $D$ and the
rational boundary components.  By the preceding discussion we have
\begin{equation*}
\Dstar_\sigma = \coprod_{\substack{\text{$R$ $\sigma$-saturated}\\
				\text{parabolic $\QQ$-subgroup}}} D_{R,h}
\end{equation*}
and $G(\QQ)$ acts on $\Dstar_\sigma$.  There is a topology on
$\Dstar_\sigma$ (the ``Satake topology'') for which $\Xstar_\sigma \equiv
\G\back \Dstar_\sigma$ is a compact Hausdorff space containing $X$ as a
dense open subset; $\Xstar_\sigma$ is the \emph{Satake compactification
of $X$ associated to $\sigma$}.  For example, if $D$ is a Hermitian symmetric space and
$\lsb\RR \Dstar_\sigma$ is the natural compactification, then $\Xstar_\sigma$
is the Baily-Borel-Satake compactification \cite{refnBailyBorel}.  Below in \S\ref{ssectSatakeAsQuotientOfRBS} we indicate an alternate construction of the topology on $\Xstar_\sigma$.

Let $\Pl^*_\sigma$ denote the $\G$-conjugacy classes of $\sigma$-saturated
parabolic $\QQ$-subgroups.  For $R\in \Pl^*_\sigma$, set $F_R =
\G_{L_{R,h}}\back D_{R,h}$, where we view $D_{R,h}$ as the symmetric space
associated to $L_{R,h}$ as above and $\G_{L_{R,h}} = \G_{L_R}/(\G_{L_R}\cap
L_{R,\l})$. In the notation of \S\ref{ssectPullbackToFiber}, $F_R = X_{R,h}$.  Then $\Xstar_\sigma$ has a stratification
\begin{equation*}
\Xstar_\sigma = \coprod_{R\in \Pl^*_\sigma} F_R  \ .
\end{equation*}
We define a partial order on the strata of $\Xstar_\sigma$ by setting $F_R\le
F_{R'}$ if and only if $R^\g \cap R'$ is a parabolic $\QQ$-subgroup for
some $\g\in \G$ and $\kap(\D^R)\subset \kap(\D^{R'})$.  This corresponds
to the condition that $F_R\subset \cl{F_{R'}}$.

\subsection{$\Xstar_\sigma$ as a quotient of $\Xhat$}
\label{ssectSatakeAsQuotientOfRBS}
Continue to assume that $\lsb\RR \Dstar_\sigma$ is a geometrically
rational Satake compactification of $D$.  Let $P$ be a parabolic
$\QQ$-subgroup.  By Corollary ~\ref{ssectRationalBoundaryComponents}, $P^\dag$ is a $\sigma$-saturated parabolic
$\QQ$-subgroup and thus defines a stratum $D_{P^\dag,h}$ of $\Dstar_\sigma$.
Let $D_{P,\l} = L_{P,\l}(\RR)/K_{P,\l}A_P$ be the symmetric space of
the second factor of \eqref{eqnLevisRationalFactors} (where $K_{P,\l} =
K_P\cap L_{P,\l}(\RR)$).  The factorization \eqref{eqnLevisRationalFactors}
induces a decomposition
\begin{equation}
D_{P} = L_P(\RR)/(K_PA_P) = D_{P,h}\times D_{P,\l}
\label{eqnEhatsFactors}
\end{equation}
and we let
\begin{equation*}
q_P\colon   D_{P} \to D_{P^\dag,h}
\end{equation*}
be the projection onto the first factor.  (Recall that
$D_{P,h}=D_{P^\dag,h}$.)
Together these maps yield a surjection
\begin{equation*}
q\colon  \widehat D \to \Dstar_\sigma \ .
\end{equation*}
Let $\pi\colon \Xhat\to \Xstar_\sigma$ be the $\Gamma$-equivariant map induced by $q$.  Zucker's
main result in \cite{refnZuckerSatakeCompactifications} is that $\pi$ is a
quotient map.  (Note however that the quotient topology on $\Dstar_\sigma$
induced by $q$ may be finer than the Satake topology on $\Dstar_\sigma$
\cite[(3.10)]{refnZuckerSatakeCompactifications}.)  In the
notation of \S\ref{ssectPullbackToFiber} we have $\pi^{-1}(F_R)=
X_R(L_{R,\l})$ for any $\sigma$-saturated parabolic $\QQ$-subgroup $R$.

\section{Functoriality of Weak Micro-support under Inverse Images}
\label{sectFunctorialityInverseImages}
Let $k\colon Z\hookrightarrow W$ be an inclusion map of admissible spaces
possessing unique maximal strata and let $\M$ be an $\L_W$-module.  Let
$X_S$ denote the unique maximal stratum of $W$.  We wish to describe the
weak micro-support of the inverse images $k^*\M$ and $k^!\M$ in terms of the weak micro-support of $\M$.  We will consider three cases: when $Z$ is open, when $Z$ is closed, and when $Z$ is a fiber as in \S\ref{ssectPullbackToFiber}.  The results in the last two cases do not apply for the full micro-support; we will add hypotheses to treat this in the next section.

\subsection{Notation}
\label{ssectFunctorialityInverseImagesNotation}
Given an admissible space $W$ with unique maximal stratum $X_S$ and an $\L_W$-module $\M$, let  $\mS_w(\M)$ denote the weak micro-support of $\M$ as defined in \S\ref{sectMicroSupport}.   Specifically for $V\in \IrrRep(L_P)$ with $P\in\Pl(W)$ let $Q^W_V$ be the parabolic $\QQ$-subgroup containing $P$
with type $\{\,\al\in \D_P^S \mid
\langle\xi_V+\hsr,\al\spcheck\rangle <0\,\}$ and let $Q^{\prime W}_V$ be similarly defined with the non-strict inequality.  Then $\mS_w(\M)$ is the set of all such $V$ with $H(\i_{P}^* \ihat_{Q}^! \M)_V \neq 0$ for some $Q\in [Q_V^W, Q^{\prime W}_V]$.
For $V\in \mS_w(\M)$ let $[c(V;\M),d(V;\M)]$ be the smallest interval containing all degrees $j$ such that  $H^j(\i_{P}^* \ihat_{Q}^! \M)_V \neq 0$ for some $Q$ as in \S\ref{ssectcMdMDefinition}.  For $R\in\Pl(W)$ we will write $Q_V^R$ instead of $Q_V^{\Xhat_R\cap W}$ and
similarly $Q_V^{\prime R}$.

We will use the functors on $\L$-modules induced by the inclusions $i_P$, $\ihat_P$, $\j_P$, and $\jhat_P$ into $W$ as in \S\ref{ssectFunctorsOnLsheaves}.  Since $W$ changes within a composition of such functors, each instance of these functors may refer to a different $W$ which can be determined from the context.

Lastly we use the notation $(P,Q)$ for the parabolic $\QQ$-subgroup complementary to $Q$ with respect to $P$ as in \S\S\ref{sssectComplementaryParabolic}, as well as the notation of \S\S\ref{ssectPullbackToFiber} and
\ref{sectSatakeCompactifications}.

\subsection{The case of an open embedding}
\label{ssectMicroSupportOpenInverseImage}
The following proposition is clear:

\begin{prop*} Let $k\colon Z\hookrightarrow W$ be the inclusion of a open
admissible subspace and let $\M$ be an $\L_W$-module.  Then
\begin{equation*}
\mS_w(k^*\M)=\{\,V\in\mS_w(\M)\mid V\in\IrrRep(L_P)\text{ for some }P\in\Pl(Z)\,\}
\end{equation*}
and for $V\in \mS_w(k^*\M)$,
\[
  [c(V;k^*\M),d(V;k^*\M)] = [c(V;\M),d(V;\M)]\ .
  \]
The same result holds for $\mS(\cdot)$ as well as for $k^!$ \textup(which equals $k^*$ in this case\textup).
\end{prop*}

\subsection{More partial orderings}
\label{ssectAdditionalPartialOrdering}
Before treating the case of a closed embedding, we need to generalize the
partial orderings of \S\ref{ssectPartialOrdering} and present several lemmas.  Let $V\in \IrrRep(L_P)$ and $\tilde V\in \IrrRep(L_{\tilde P})$ for $P,\tilde P\in\Pl(W)$.  Recall we write $V \preccurlyeq \tilde V$
if \begin{enumerate*} \item \label{itemParabolic} $P\le \tilde P$ and
\item \label{itemKostantComponent} $V= H^{\l(w)}(\n_{P}^{\tilde P};
\tilde V)_w$ for some  $w\in W_P^{\tilde P}$.
\end{enumerate*}
In this setting we define $[\tilde V:V]=\l(w)$.

Let $\lsp+\sa_P^{\tilde P*}$ denote the convex cone
generated by all $\al\in\D_P^{\tilde P}$.  Write $V \preccurlyeq_+ \tilde V$ (resp.\ $V\preccurlyeq_- \tilde V$) if  $V \preccurlyeq \tilde V$ holds and in addition
\begin{enumerate}[label=(\roman*)$_\pm$]
  \setcounter{enumi}{2}
  \item \label{itemRootCone} $(\xi_V+\hsr)|_{\sa_P^{\tilde P}} \in
    \lsp+\sa_P^{\tilde P*}$ \quad (resp.\ $(\xi_V+\hsr)|_{\sa_P^{\tilde P}} \in
    -\lsp+\sa_P^{\tilde P*}$)
\end{enumerate}
is satisfied.  If $V \preccurlyeq_+ \tilde V$ and $V \preccurlyeq_-
\tilde V$ hold then $V\preccurlyeq_0 \tilde V$ .  Write $V\preccurlyeq_\pm
\tilde V$ if either $V \preccurlyeq_+ \tilde V$ or $V \preccurlyeq_-
\tilde V$ are true.  In all of these cases write $V\prec \tilde V$, etc., if $P<\tilde P$.

\begin{lem}
\label{ssectPartialOrderingLemma}
All of $\preccurlyeq$, $\preccurlyeq_0$, $\preccurlyeq_+$, and $\preccurlyeq_-$ are partial orderings.  They are generated by the relations where $\#\D_P^{\tilde
P}=1$.
\end{lem}

\begin{proof}
Consider $\preccurlyeq_+$.  Transitivity is clear for
\itemref{itemParabolic} and for \itemref{itemKostantComponent} it follows
from \eqref{eqnLinkComposition}.  As for \itemref{itemRootCone}$_+$, suppose
$V\preccurlyeq_+ \tilde V \preccurlyeq_+ \tilde{\tilde V}$.  Then
\begin{equation*}
(\xi_V+\hsr)|_{\sa_P^{\tilde{\tilde P}}} = 
(\xi_V+\hsr)|_{\sa_P^{\tilde P}} +
(\xi_{\tilde V}+\hsr)|_{\sa_{\tilde P}^{\tilde{\tilde P}}}
\end{equation*}
so we need to check that the elements of both $\D_P^{\tilde P}$ and
$\D_{\tilde P}^{\tilde{\tilde P}}$ are in $\lsp+\sa_P^{\tilde{\tilde
P}*}$ when viewed as linear functions on $\sa_P^{\tilde{\tilde P}} =
\sa_P^{\tilde P} + \sa_{\tilde P}^{\tilde{\tilde P}}$.  But this is
one of Langlands's ``geometric lemmas'' \cite[IV,
\S6.5(2)]{refnBorelWallach}.

For the last assertion, suppose that $V\preccurlyeq_+ \tilde{\tilde V}$.
If $\langle \xi_V+\hsr,\al\spcheck\rangle< 0$ for all $\al \in
\D_P^{\tilde{\tilde P}}$, then $(\xi_V+\hsr)|_{\sa_P^{\tilde{\tilde P}}}$
belongs to the interior of $-\lsp+\sa_P^{\tilde{\tilde P}*}$ \cite[IV,
\S6.2]{refnBorelWallach}, which is a contradiction.  So there exists
$\al_0\in \D_P^{\tilde{\tilde P}}$ such that
$\langle\xi_V+\hsr,\al_0\spcheck\rangle\ge 0$.  Let $\tilde P\ge P$ have
type $\D_P^{\tilde P}=\{\al_0\}$ and decompose $w=w^{\tilde P}w_{\tilde
P}\in W_P^{\tilde P}W_{\tilde P}^{\tilde{\tilde P}}$.  Then it is easy to
check that $\tilde V= H^{\l(w_{\tilde P})}(\n_{\tilde P}^{\tilde{\tilde
P}}; \tilde{\tilde V})_{w_{\tilde P}}$ satisfies $V\preccurlyeq_+ \tilde
V\preccurlyeq_+ \tilde{\tilde V}$ and we can use induction.
\end{proof}

\begin{lem}
  \label{ssectBehaviorOfQsubVUnderPartialOrder}
  If $V\preccurlyeq_+ \tilde V$ then $Q_V^W\vee \tilde P \ge Q_{\tilde V}^W$ and $Q_V^{\prime W}\vee \tilde P \ge Q_{\tilde V}^{\prime W}$.  If $V\preccurlyeq_- \tilde V$ then $Q_V^W\vee \tilde P \le Q_{\tilde V}^W$ and $Q_V^{\prime W}\vee \tilde P \le Q_{\tilde V}^{\prime W}$.
  \end{lem}

\begin{proof}
  By Lemma~\ref{ssectPartialOrderingLemma} it suffices to consider the case   where $\D_P^{\tilde P}=\{\al_0\}$.  Suppose $V\preccurlyeq_+ \tilde V$ and thus  $\langle \xi_V+\hsr,\al_0\spcheck\rangle\ge0$.  Let $\tilde \al\in \D_{\tilde P}^S$ and let $\al\in \D_P^S\setminus \D_P^{\tilde P}$ be the unique element for which
$\al|_{\sa_{\tilde P}}=\tilde \al$.  Since $\xi_{\tilde V} + \hsr_{\tilde P}=
(\xi_V+\hsr)|_{\sa_{\tilde P}}$ we have
\begin{equation*}
\langle \xi_{\tilde V} + \hsr,\tilde \al{}\spcheck\rangle =
\langle \xi_V+\hsr,\al\spcheck\rangle - 
\langle \b_{\al_0}^{\tilde P}, \al\spcheck\rangle
\langle \xi_V+\hsr,\al_0\spcheck\rangle \ .
\end{equation*}
Since $\langle \b_{\al_0}^{\tilde P}, \al\spcheck\rangle = c\langle \al_0,
\al\spcheck\rangle\le 0$ (where $c>0$), the inequality $\langle\xi_{\tilde V} +
\hsr,\tilde \al{}\spcheck\rangle<0$ implies
$\langle\xi_V+\hsr,\al\spcheck\rangle<0$ and thus $Q_V^W\vee \tilde P \ge
Q_{\tilde V}^W$.  The other assertions are proven similarly.
\end{proof}

\begin{lem}\label{lemStarShriekMicroSupport}
  Consider $P\le Q\le R\le S\in\Pl(W)$ and an irreducible $L_P$-module $V$.  Set $T=(Q,R)\cap S$\textup:
  \begin{equation*}
\xymatrix @ur @M=1pt @R=1.5pc @C=1pc {
{} & {T}
\ar@{-}[r] \ar@{-}[d] &  {S} \ar@{-}[d]  \\
{P} \ar@{-}[r] & {Q} \ar@{-}[r]
& {R}
}\ .
  \end{equation*}
  Then
  \begin{equation}
    \label{eqnStratumStarMicroSupport}
    H(\i_P^* \ihat_Q^!(\ihat_R^*\M))_V =H(\i_P^* \ihat_T^!\M)_V
  \end{equation}
  and
  \begin{equation}
    \label{eqnStratumShriekMicroSupport}
    H(\i_P^* \ihat_Q^!(\ihat_R^!\M))_V = H(\i_P^* \ihat_Q^!\M)_V\ .
  \end{equation}
  Furthermore\textup:
  \begin{enumerate}
  \item\label{itemStratumStarMicroSupport}
    $V\in\mS_w(\ihat_R^*\M)$ if and only if \eqref{eqnStratumStarMicroSupport} is $\neq0$ for some $Q\in[Q_V^R, Q_V^{\prime R}]$.  In this case $[c(V;\ihat_R^*\M),d(V;\ihat_R^*\M)]$ is the smallest interval containing all nonzero degrees of \eqref{eqnStratumStarMicroSupport} for all such $Q$.
  \item\label{itemStratumShriekMicroSupport}
    $V\in \mS_w(\ihat_R^!\M)$ if and only if \eqref{eqnStratumShriekMicroSupport} is $\neq0$
    for some $Q\in[Q_V^R, Q_V^{\prime R}]$.  In this case $[c(V;\ihat_R^!\M),d(V;\ihat_R^!\M)]$ is the smallest interval containing all nonzero degrees of \eqref{eqnStratumShriekMicroSupport} for all such $Q$.
  \end{enumerate}
\end{lem}

%

\begin{proof}
The equalities in \eqref{eqnStratumStarMicroSupport} and \eqref{eqnStratumShriekMicroSupport} follow from Proposition
~\ref{ssectFunctorsOnLsheaves}\itemref{itemFunctorsComposition}\itemref{itemFunctorsSquare}.
\end{proof}

Now fix a Satake compactification $\Xstar_\sigma$ with strata $(F_R)_{R\in \Pl^*_\sigma}$.

\begin{lem}
  \label{lemTwoBoundary}
  Say $P\le \tilde P$ and set $R=P^\dag$ and $\tilde R={\tilde P}^\dag$.  Let $V$ and $\tilde V$ irreducible $L_P$- and $L_{\tilde P}$-modules respectively and assume $V\preccurlyeq \tilde V$.
  \begin{enumerate}
    \item \label{itemTwoBoundarySets} $\tilde R \le \tilde P \vee R$.
    \item \label{itemTwoBoundaryOrder} $F_R\le F_{\tilde R}$.
  \end{enumerate}
  For the remainder of the lemma assume as well that $\tilde P\cap R = P$.
  \begin{enumerate}[resume*]
    \item \label{itemTwoBoundaryRoot} $\D^{\tilde P}\setminus \D^P \subset \kap(\D^{\tilde P})$.
    \item \label{itemTwoBoundaryLeviFactors} $N_P^{\tilde P} \subset \widetilde{L_{\tilde P,h}}$ and the decomposition \eqref{eqnLevisRationalFactors} for $L_P$ can be refined to
    \[
    L_P = \widetilde{L_{P,h}} L_{P,\l}^{\tilde P,h} L_{\tilde P,\l} \ .
    \]
  \item \label{itemTwoBoundaryOrthogonal} If $\tilde \g\in \D^{\tilde R}\setminus \D^{\tilde P}$ then $\tilde \g|_{\sa_P^{\tilde P}}=0$.
  \item \label{itemQTransfer} $(Q_V^R\vee \tilde P)\cap\tilde R = Q_{\tilde V}^{\tilde R}$ and $(Q_V^{\prime R}\vee \tilde P)\cap\tilde R = Q_{\tilde V}^{\prime \tilde R}$.
  \end{enumerate}
\end{lem}

\begin{proof}
  First some preliminaries.  By Corollary ~\ref{ssectRationalBoundaryComponents}\itemref{itemSaturationRational} $R$ is a parabolic $\QQ$-subgroup with $\D^R=\o(\kap(\D^P))$.  Thus we have an orthogonal decomposition $\D^R = \kap\sqcup \zeta$ where $\kap=\kap(\D^P)$ and $\zeta$ is the largest subset of $\D$ such that $\zeta\perp (\kap\cup\{\lsb\QQ\mu\})$ (here $\lsb\QQ\mu$ is the restricted highest $\QQ$-weight of $\sigma$).  We have an analogous decomposition $\D^{\tilde R} = \tilde \kap \sqcup \tilde\zeta$.  The following relations between these decompositions will be used frequently:
  \begin{equation*}
    \kap\subset \tilde\kap \qquad \text{and} \qquad  \tilde \zeta \subset \zeta\ .
  \end{equation*}
The first holds since $P\le\tilde P$, the second  since $\tilde\zeta\perp (\tilde\kap\cup\{\lsb\QQ\mu\})$ implies $\tilde\zeta\perp (\kap\cup\{\lsb\QQ\mu\})$.

Part \itemref{itemTwoBoundarySets} holds since $\D^{\tilde R} =  \tilde \kap \sqcup \tilde\zeta \subset \tilde \kap \cup \zeta \subset \D^{\tilde P}\cup \D^R$.  Part \itemref{itemTwoBoundaryOrder} follows from the definition of the partial order in \S\ref{ssectSatakeCompacticationsOfX} since $R\cap\tilde R\supset P$ is a parabolic $\QQ$-subgroup and $\kap\subset \tilde \kap$.

For part \itemref{itemTwoBoundaryRoot}, consider $\g\in\D^{\tilde P}\setminus \D^P$.  Either $\g \in \tilde\kap$ or $\g \in\tilde\zeta\cap\D^{\tilde P} \subset \tilde \zeta\subset \zeta$.  But in the second case $\g\in\D^R$ which, since $\g\notin \D^P$, contradicts $\tilde P\cap R=P$.  Thus $\g\in\tilde \kap$ which proves \itemref{itemTwoBoundaryRoot}.  To prove  \itemref{itemTwoBoundaryLeviFactors}, note by \eqref{eqnLevisRationalFactors} we have an almost direct product $L_{\tilde P}= \widetilde{L_{\tilde P,h}}L_{\tilde P,\l}$ which corresponds to $\D^{\tilde P} = \tilde \kap\sqcup (\tilde \zeta\cap\D^{\tilde P})$.  Part \itemref{itemTwoBoundaryRoot} shows that $L_P$, viewed as the Levi quotient of the parabolic subgroup $P/N_{\tilde P}$ of $L_{\tilde P}$, is obtained by removing simple roots from $\tilde\kap$, the first factor.  This yields an orthogonal decomposition
\begin{equation}
  \label{eqnTwoBoundaryRootDecomposition}
  \D^P = \kap \sqcup (\tilde \kap \cap \zeta\cap \D^P) \sqcup (\tilde \zeta \cap \D^{\tilde P})
\end{equation}
and thus the almost direct product of \itemref{itemTwoBoundaryLeviFactors}.

Now let $\tilde \g\in \D^{\tilde R}\setminus \D^{\tilde P}$.  To prove \itemref{itemTwoBoundaryOrthogonal} it suffices to show that $(\tilde \g,\g_P )=0$ for all $\g\in \D^{\tilde P}\setminus \D^P$ and for this it suffices to show $\tilde \g \perp \g$ and $\tilde \g \perp \g^P$.  The first assertion follows from  $\tilde \g \in \tilde \zeta$ and $\g\in \tilde \kap$ (by \itemref{itemTwoBoundaryRoot}).  For the second, $\g\in\tilde \kap$ is orthogonal to the third factor of \eqref{eqnTwoBoundaryRootDecomposition} so $\g^P$ is in the span of the first two factors.  But $\tilde \g\in\tilde\zeta$  is orthogonal to the first two factors which finishes the proof.

%
%

To prove \itemref{itemQTransfer}, let $\tilde \al\in\D_{\tilde P}^{\tilde R}$ and let $\al$ be the unique element of
$\D_P^{\tilde R}\setminus\D_P^{\tilde P}$ for which $\al|_{\sa_{\tilde P}}=\tilde \al$.  We know $\alpha\in\D_P^R$ by \itemref{itemTwoBoundarySets}.  We calculate
\[
\langle\xi_{\tilde V} + \hsr,\tilde \al\spcheck\rangle =  \langle\xi_V+\hsr,\tilde\al\spcheck\rangle =
\langle\xi_V+\hsr,\al\spcheck\rangle
\]
by using in turn $V\preccurlyeq \tilde V$ and \itemref{itemTwoBoundaryOrthogonal}.
It follows that $\tilde\al\in\D_{\tilde P}^{Q_{\tilde V}^{\tilde R}}$
if and only if $\al\in \D_P^{Q_V^R}$ which is equivalent to $\tilde\al\in \D_{\tilde P}^{Q_V^R\vee
  \tilde P}$.  Since  $\D_{\tilde P}^{(Q_V^R\vee \tilde P)\cap \tilde R} = \D_{\tilde P}^{(Q_V^R\vee \tilde P)}
\cap \D_{\tilde P}^{\tilde R}$, this proves the first equality of \itemref{itemQTransfer}; the other follows similarly.
\end{proof}

\subsection{The case of a closed embedding:\ weak form}
\label{ssectWeakMicroSupportClosedInverseImage}

If $k\colon Z\hookrightarrow W$ is the inclusion of an admissible subspace
possessing a unique maximal stratum and  $\mathscr S\subset \coprod_{P\in\Pl(W)} \IrrRep(L_P)$, define
\begin{equation}
\begin{aligned}
k^*\mathscr S &= \{\,V \mid V\in \IrrRep(L_P) \text{ for }P\in\Pl(Z) \text{ and }
V\preccurlyeq_+ \tilde V \text{ for }\tilde V\in\mathscr S\,\} \\
k^!\mathscr S &= \{\,V \mid V\in \IrrRep(L_P) \text{ for }P\in\Pl(Z) \text{ and }
V\preccurlyeq_- \tilde V \text{ for }\tilde V\in\mathscr S\,\}\ .
\end{aligned}
\label{eqnRestrictMicroSupportClosedStratum}
\end{equation}

Let $\pi\colon \Xhat \to \Xstar_\sigma$ be the projection defined in \S\textup{\ref{ssectSatakeAsQuotientOfRBS}}.  For a stratum $F_R$ of $\Xstar_\sigma$ let $U(F_R)= \coprod_{F_{R'}\ge F_R} F_{R'}$ denote its open star neighborhood.

\begin{prop*}
  Assume $W=\Xhat_S\cap \pi^{-1}(U(F_R))$ where $F_R$ is a stratum of $\Xstar_\sigma$ and $S\ge R$.  Let $\ihat_R\colon \Xhat_R\cap W\hookrightarrow W$ be the inclusion of the closed admissible subspace and let $\M$ be an $\L_W$-module.  Then
\begin{equation*}
\mS_w(\ihat_R^*\M) \subset \ihat_R^*\mS_w(\M) \qquad
\text{and}\qquad
\mS_w(\ihat_R^!\M) \subset \ihat_R^!\mS_w(\M)\ .
\end{equation*}
In fact if $V\in\IrrRep(L_P)$ is in $\mS_w(\ihat_R^{*}\M)$ \textup(resp.\ $\mS_w(\ihat_R^{!}\M)$\textup) then there exists a sequence
\begin{align*}
  V&\equiv V_0\prec_+ V_1\prec_+ \dots\prec_+ V_N \equiv \tilde V \\
  \text{\textup(resp.\ } V&\equiv V_0\prec_- V_1\prec_- \dots\prec_- V_N \equiv \tilde V
  \text{\ \textup)}
\end{align*}
with an $L_{\tilde P}$-module $\tilde V \in \mS_w(\M)$.  Here $V_i\in \IrrRep(L_{P_i})$ belongs to $\mS_w(\ihat^{*}_{R_i}\ihat_{S_i}^!\M)$
for $R_i=P_i^\dag$ and some $S_i\ge R_i$.
%
%
Furthermore
\begin{align}
  \label{eqnPversusR}
  P_{i}\cap R_{i-1} = P_{i-1}  &\qquad \text{for $0 < i \le N$,}\\
  \label{eqnStrictDominantCone}
(\xi_{V_{i-1}}+\hsr)|_{\sa_{P_{i-1}}^{P_{i}}}\in \pm \sa_{P_{i-1}}^{P_{i}*+} &\qquad
  \text{for $0< i \le N$,}\\
    \label{eqnNumberOfInductionSteps}
  N_V^{\tilde V} \equiv N \in [0,\dim \sa_P^{\tilde P}]
\end{align}
and
\begin{align}
  \label{eqnWeakMSInverseImageDegree}
    [c(V;\ihat_R^*\M),d(V;\ihat_R^*\M)] &\subset [c(\tilde V;\M), d(\tilde V;\M)] + [\tilde V:V]\ ,\\
\label{eqnWeakMSInverseImageProperSupportDegree}
  [c(V;\ihat_R^!\M),d(V;\ihat_R^!\M)] &\subset
     [c(\tilde V;\M), d(\tilde V;\M)] + [\tilde V:V] + N_V^{\tilde V} \ .
\end{align}
\end{prop*}

\begin{rems*}
\leavevmode
\begin{enumerate}
  \item\label{itemStrictnessForStar} In the $\ihat_R^*$ case \eqref{eqnStrictDominantCone} holds with $\Int \sa_{P_{i-1}}^{P_{i}*+}$.  We do not use this.
  \item\label{itemWeakMSClosedEmbeddingFalseForMS} The proposition is false in general if $\mS_w$ is replaced by $\mS$.  However we will see in \S\ref{sectFunctorialityMicroSupportEqualRank} a situation in which this replacement is valid.
  \item Properties \eqref{eqnPversusR} and \eqref{eqnStrictDominantCone} on the relationship between $P_{i-1}$ and $P_{i}$ may not extend to a relationship between $P$ and $\tilde P$.  Specifically although it is true that $\tilde P \cap R\ge P$, equality may fail.  Also while $(\xi_{V}+\hsr)|_{\sa_{P}^{\tilde P}}\in \pm \lsp+ \sa_{P}^{\tilde P*}$, the codominant cone, it may not be in the dominant cone.  This is the importance of the sequence $(V_i)_i$.
    %
    %
    %
  \end{enumerate}
\end{rems*}

\begin{proof}
  First note the following facts which follow from Lemma~ \ref{lemTwoBoundary}\itemref{itemTwoBoundarySets}\itemref{itemTwoBoundaryOrder}:
  \begin{enumerate}
  \item $W=\Xhat_S\cap \pi^{-1}(U(F_R))$ is admissible with unique maximal stratum $X_S$;
  \item If $P\in \Pl(\Xhat_R\cap W)$ then $P^\dag=R$;
  \item If $\tilde P\in \Pl(W)$ with $\tilde P\ge P$ then $\tilde P^\dag \in \Pl(W)$.
  \end{enumerate}
As a consequence the parabolic subgroups that will arise in the constructions below all belong to $\Pl(W)$.

  Let $V$ be an irreducible $L_P$-module for $P\in \Pl(\Xhat_R\cap W)$.
  We prove the case where $V\in \mS_w(\ihat_R^*\M)$; the changes necessary for $V\in \mS_w(\ihat_R^!\M)$ will be indicated at the end.  We use induction on $\#\D_P$ so we may assume
\begin{equation}\label{eqnInductionStarMicroSupport}
  \parbox[t]{.8\linewidth}{if $P_1>P$ the result is true for an $L_{P_1}$-module $V_1\in \mS_w(\ihat_{R_1}^*\M_1)$ where $R_1=P_1^\dag$, $S_1\ge R_1$,  
and $\M_1$ is a $\L_{\Xhat_{S_1}\cap \pi^{-1}(U(F_{R_1}))}$-module.}    
\end{equation}

Let $Q\in[Q_V^R, Q_V^{\prime R}]$, set $U=(P,R)\cap S$, and let $T=(Q,R)\cap S = Q\vee U$.  
Let $P'\ge P$ have type $\D_P^{P'}= \{\,\al\in \D_P^S\setminus \D_P^R
\mid \langle\xi_V+\hsr,\al\spcheck\rangle > 0 \,\}$ and let $Q' = Q \vee P'$.  Then let $U' = (P,P')\cap U$ and $T' = (Q,Q')\cap T=Q\vee U'$.  The situation is pictured in Figure ~\ref{figureThreeDDiagram}.
\begin{figure}
\begin{equation*}
\xymatrix @ur @M=1pt @R=1.5pc @C=1pc {
{} & {U} \ar@{-}[rr] \ar@{-}[dd] & {}  & {T}
\ar@{-}[rr] \ar@{-}[dl] \ar@{-}[dd] & {}  & {S} \ar@{-}[dd] \\
{U'} \ar@{-}'[r][rr] \ar@{-}[ru] \ar@{-}[dd] & *+{\hole}  & {T'} \ar@{-}'[r][rr]
\ar@{-}'[d][dd] & *+{\hole} & {R\vee U'} \ar@{-}[ru] \ar@{-}'[d][dd] \\
{} & {P'} \ar@{-}[rr] & *+{\hole}  & {Q'}
\ar@{-}[rr] \ar@{-}[dl] & *+{\hole}  & {R\vee P'} \\
{P} \ar@{-}[rr] \ar@{-}[ru] & {}  & {Q} \ar@{-}[rr]
& {} & {R} \ar@{-}[ru]
}
\end{equation*}
\caption{}\label{figureThreeDDiagram}
\end{figure}
Note that $\D_P^{U'}= \{\,\al\in \D_P^S\setminus \D_P^R
\mid \langle\xi_V+\hsr,\al\spcheck\rangle \le 0 \,\}$ so
\begin{equation}
  \label{eqnUprimeConsequencesStar}
  Q_V^W \le Q_V^R \vee U' \le T' \quad\text{ and }\quad Q_V^{\prime W} = Q_V^{\prime R} \vee U'\ .
\end{equation}

Consider the short exact sequence
\begin{equation*}
0\to \i_P^* \ihat_{T'}^!  \M \to \i_P^* \ihat_{T}^! \M \to
\i_P^* \jhat_{T'*}(\jhat_{T'}^*\ihat_{T}^!  \M) \to 0
\end{equation*}
from \eqref{eqnShortCompareLocalCohomologyWithSupports} and the corresponding long exact sequence:
\begin{equation}
\dots\to
H^j(\i_P^* \ihat_{T'}^!  \M)_V \to
H^j(\i_P^* \ihat_{T}^! \M)_V \to
H^j(\i_P^* \jhat_{T'*}(\jhat_{T'}^*\ihat_{T}^!  \M))_V \to \cdots\ .
\label{eqnExactTripleTransfer}
\end{equation}
By Lemma ~\ref{lemStarShriekMicroSupport}\itemref{itemStratumStarMicroSupport} the middle term is
nonzero for some $j\in [c(V;\ihat_{R}^* \M),d(V;\ihat_{R}^*\M)]$ so either the first or last term is as well.  If
the first term of \eqref{eqnExactTripleTransfer} is nonzero, then we
have that
$V\in\mS_w(\M)$ (since $T'=Q\vee U'\in [Q_V^W ,Q_V^{\prime W}]$ by \eqref{eqnUprimeConsequencesStar}).  If
instead the last term of \eqref{eqnExactTripleTransfer} is nonzero, the Fary spectral sequence \eqref{eqnFarySpectralSequence} implies there exists $P_1\in (P, P']$ such that $H^\l(\n_P^{P_1};H^{j-\l}(\i_{P_1}^*\ihat_{T' \vee P_1}^!\M))_V \neq 0$ for some $\l$.  Thus
\begin{equation}
H^{j-[V_1:V]}(\i_{P_1}^*\ihat_{T_1}^!\M)_{V_1}\neq 0
\label{eqnPPrimeStarTShriekNonvanishing}
\end{equation}
for some $V_1\in \IrrRep(L_{P_1})$ and $T_1 = T' \vee P_1$.   Equations \eqref{eqnPversusR} and  \eqref{eqnStrictDominantCone} hold for $i=1$, the latter due to the definition of $P'$.  In particular we have $V\prec_+ V_1$.

Now set $R_1= P_1^\dag$ and $S_1 = T_1 \vee R_1$.
%
%
%
By Lemma~\ref{lemTwoBoundary}\itemref{itemTwoBoundarySets}, $R_1\in [P_1,P_1\vee R]$.  Set $Q_1= T_1\cap R_1$:
\begin{equation*}
\xymatrix @ur @M=1pt @R=1.5pc @C=1pc {
{} & {T_1}
\ar@{-}[r] \ar@{-}[d] &  {S_1} \ar@{-}[d] \\
{P_1} \ar@{-}[r] & {Q_1} \ar@{-}[r]
& {R_1}
}\ .
\end{equation*}
Here \eqref{eqnStratumStarMicroSupport} of Lemma ~\ref{lemStarShriekMicroSupport} applied to $\ihat_{S_1}^!\M$ shows \eqref{eqnPPrimeStarTShriekNonvanishing} is equivalent to
\begin{equation}
H^{j-[V_1:V]}(\i_{P_1}^* \ihat_{Q_1}^!(\ihat_{R_1}^*\ihat_{S_1}^!\M))_{V_1}\neq 0\ .
\label{eqnVPrimeInTerm}
\end{equation}

Now $Q_1=  (Q\vee U'\vee P_1)\cap R_1$ and $ U'\cap R_1\le P < P_1$.  Then, by  Lemma~ \ref{lemTwoBoundary}\itemref{itemQTransfer},
\begin{align*}
  Q_{V_1}^{R_1} &= (( Q_V^R\vee P_1)\cap R_1) \vee  (U'\cap R_1) =
  (Q_V^R\vee U'\vee P_1)\cap R_1 \le Q_1 \\
  &\le (Q_V^{\prime R}\vee U'\vee P_1)\cap R_1
  = ((Q_V^{\prime R}\vee P_1)\cap R_1) \vee (U'\cap R_1)
  = Q_{V_1}^{\prime R_1}
\end{align*}
which implies $Q_1\in[Q_{V_1}^{R_1},Q_{V_1}^{\prime R_1}]$.  As a result \eqref{eqnVPrimeInTerm} yields that $V_1\in\mS_w(\ihat_{R_1}^*\ihat_{S_1}^!\M)$.  For $j\in\{ c(V;\ihat_R^*\M), d(V;\ihat_R^*\M) \}$
Lemma ~\ref{lemStarShriekMicroSupport}\itemref{itemStratumStarMicroSupport} and \eqref{eqnVPrimeInTerm} imply
\begin{equation}
  \label{eqnWeakMSVOneInverseImageDegree}
  [c(V;\ihat_R^*\M),d(V;\ihat_R^*\M)] \subset [c(V_1;\ihat_{R_1}^*\ihat_{S_1}^!\M), d(V_1;\ihat_{R_1}^*\ihat_{S_1}^!\M)] + [V_1:V]\ .
\end{equation}

Since $P_1>P$ we can apply induction (see \eqref{eqnInductionStarMicroSupport}) and obtain  irreducible $L_{P_i}$-modules $V_1 \preccurlyeq_+ V_2 \preccurlyeq_+ \dots \preccurlyeq_+ V_N$  with $\tilde V\equiv V_N\in \mS_w(\ihat_{S_1}^!\M)$ and $\tilde P \equiv P_N \le S_1$.
%
%
%
Thus $S_1 \ge T_1\vee \tilde{P} \ge Q_{V}^W \vee \tilde{P} \ge Q_{\tilde V}^W$ by \eqref{eqnUprimeConsequencesStar} and Corollary
~\ref{ssectBehaviorOfQsubVUnderPartialOrder} so $Q_{\tilde V}^{S_1} \ge Q_{\tilde V}^W$.  Since clearly $Q_{\tilde V}^{\prime S_1} \le Q_{\tilde V}^{\prime W}$ we conclude $\tilde V\in \mS_w(\M)$.
%
%
Induction also completes the proof of equations \eqref{eqnPversusR} and \eqref{eqnStrictDominantCone} and shows that 
\begin{equation}
  \label{eqnWeakMSTildeVInverseImageDegree}
 [c(V_1;\ihat_{R_1}^*\ihat_{S_1}^!\M), d(V_1;\ihat_{R_1}^*\ihat_{S_1}^!\M)] \subset  [c(\tilde V;\M), d(\tilde V;\M)]+ [\tilde V:V_1]\ .
\end{equation}
The composition of \eqref{eqnWeakMSVOneInverseImageDegree} and \eqref{eqnWeakMSTildeVInverseImageDegree} yield \eqref{eqnWeakMSInverseImageDegree}.

The proof for $\mS_w(\ihat_R^!\M)$ is similar except that we let
$P'\ge P$ have type $\D_P^{P'}= \{\,\al\in \D_P^S\setminus \D_P^R \mid
\langle\xi_V+\hsr,\al\spcheck\rangle \le 0 \,\}$ and use the long exact sequence
\begin{equation}
\dots\to
H^j(\i_P^* \ihat_{Q}^!  \M)_V \to
H^j(\i_P^* \ihat_{Q'}^! \M)_V \to
H^j(\i_P^* \jhat_{Q*}(\jhat_{Q}^*\ihat_{Q'}^!  \M))_V \to \cdots\ .
\end{equation}
It is now the first term that is nonzero for some $j\in [c(V;\ihat_{R}^! \M),d(V;\ihat_{R}^!\M)]$ and hence either the middle
or last term (in degree $j-1$) must be nonzero as well.  In the latter case we can use the Fary spectral sequence to find $P_1\in (P,P']$ and $V_1\in\IrrRep(L_{P_1})$ satisfying $V\prec_- V_1$ and
\begin{equation}
  \label{eqnPPrimeStarTShriekNonvanishingProperSupportFinal}
H^{j-1-\l}(\i_{P_1}^* \ihat_{Q_1}^!\M)_{V_1}  = H^{j-1-\l}(\i_{P_1}^* \ihat_{Q_1}^!(\ihat_{S_1}^!\M))_{V_1} \neq 0
\end{equation}
where $Q_1 = Q\vee P_1$ and $S_1=Q_1 \vee R_1$.  Set $T_1 = Q_1\cap R_1$ and  $Y_1=(P_1,T_1)\cap Q_1$:
\begin{equation*}
\xymatrix @ur @M=1pt @R=1.5pc @C=1pc {
*+!R!<-.5em,0pt>{Y_1}  \ar@{-}[r] \ar@{-}[d] & {Q_1}
\ar@{-}[r] \ar@{-}[d] &  {S_1} \ar@{-}[r] \ar@{-}[d] & {G} \\
{P_1} \ar@{-}[r] & {T_1} \ar@{-}[r]
& {R_1}
}
\end{equation*}
We can use Lemma~ \ref{lemTwoBoundary}\itemref{itemQTransfer},
Corollary~ \ref{ssectBehaviorOfQsubVUnderPartialOrder}, and Lemma ~\ref{lemStarShriekMicroSupport}\itemref{itemStratumShriekMicroSupport} to prove $V_1\in\mS_w(\ihat_{S_1}^!\M)$ and the degree estimate.

The proof concludes by induction; the degree shift of $-1$ in \eqref{eqnPPrimeStarTShriekNonvanishingProperSupportFinal} for each inductive step induces the summand $N_V^{\tilde V}$ in \eqref{eqnWeakMSInverseImageProperSupportDegree}.
\end{proof}

\begin{cor}
   Let $\ihat_R\colon \Xhat_R\cap W\hookrightarrow W$ be the inclusion of the closed admissible subspace and let $\M$ be an $\L_W$-module.  Then
\begin{equation*}
\mS_w(\ihat_R^*\M) \subset \ihat_R^*\mS_w(\M) \qquad
\text{and}\qquad
\mS_w(\ihat_R^!\M) \subset \ihat_R^!\mS_w(\M)\ .
\end{equation*}
For $V\in \mS_w(\ihat_R^*\M)$ there exists $\tilde V\in \mS_w(\M)$ so that $V\preccurlyeq_+ \tilde V$ and \eqref{eqnWeakMSInverseImageDegree} holds.
For $V\in \mS_w(\ihat_R^!\M)$ there exists $\tilde V\in \mS_w(\M)$ so that $V\preccurlyeq_- \tilde V$ and \eqref{eqnWeakMSInverseImageProperSupportDegree} holds.
\end{cor}

\begin{proof}[Proof of Corollary]
  We will show that for any $P\le R$ there exists a Satake compactification $\Xstar_\sigma$ such that $R=P^\dag$.  The corollary then follows from the above proposition.
  
  Let $\U_1,\dots,\U_k$ be all the connected components of $\D^R$ which are subsets of $\D^P$.  Since $P\le R$ each $\U_i$ is also a connected component of $\D^P$.  For $1\le i\le k$ choose a simple $\QQ$-root $\al_i\in\U_i$ and set $\Psi = \{\al_1,\dots,\al_k\}\cup(\D\setminus\D^R)$.  Let $\sigma$ be a irreducible represention of $G$, projectively rational over $\QQ$, whose restricted highest weight is a multiple of $\sum_{\al\in\Psi} \b_\al$.   The Satake compactification $\lsb\RR D_\sigma$ is geometrically rational by \cite[Cor.~7]{refnSaperGeometricRationality}.  Now $\kap(\D^R) = \U_1\cup\dots\U_k$ and $\kap(\D^R\cup\{\al\})= \kap(\D^R)\cup\{\al\}$ for any $\al\in\D\setminus \D^R$.  Thus $R$ is $\sigma$-saturated.  Since $\kap(\D^P)=\kap(\D^R)$ it follows that $P^\dag=R$.
\end{proof}

\subsection{The case of an inclusion of a fiber}
\label{ssectWeakMicroSupportInclusionFiber}
As in  \S\ref{ssectPullbackToFiber} suppose we are given a connected normal $\QQ$-subgroup $L_{R,\l}\subset L_R$
for a given $R\in\Pl$ and set $L_{R,h}=L_R/L_{R,\l}$ with associated
symmetric space $D_{R,h}$; we have an almost
direct product decomposition $L_R= \widetilde{L_{R,h}}L_{R,\l}$ where $\widetilde{L_{R,h}}$ projects onto $L_{R,h}$.  Let $X_R(L_{R,\l})\subset
\Xhat_R$ be the partial compactification of $X_R$ in the
``$L_{R,\l}$-directions''.  Let
$\ihat_{R,\l}\colon \Xhat_{R,\l}\hookrightarrow X_R(L_{R,\l})$ be the inclusion
of a generic fiber of the flat bundle $\pi\colon  X_R(L_{R,\l})\to X_{R,h}$.

In \S\ref{ssectPullbackToFiber} we defined functors $\ihat_{R,\l}^*$ (resp.\  $\ihat_{R,\l}^!$) from $\L_{ X_R(L_{R,\l})}$-modules to $\L_{\Xhat_{R,\l}}$-modules by applying $\Res_{L_{P,\l}}^{L_P}$ (resp.\ $\Res_{L_{P,\l}}^{L_P}[-\dim D_{R,h}]$) to the data.  Note that any $V\in \IrrRep(L_P)$ with $P\in \Pl(X_R(L_{R,\l}))$ may be written as
$V_h\otimes V_\l$, where $V_h\in\IrrRep(\widetilde{L_{R,h}})$ and
$V_\l\in\IrrRep(L_{P,\l})$; the restriction $\Res_{L_{P,\l}}^{L_P} V$ is a
sum of copies of $V_\l$.  For a subset $\mathscr S \subset
\coprod_{P\in\Pl(X_R(L_{R,\l}))}\IrrRep(L_P)$, define
\begin{multline}
  \Res_\l^*\mathscr S = \Res_\l^!\mathscr S = \\
\biggl\{\,V_\l \biggm|
\parbox[c]{.55\linewidth}{$V_\l\in \IrrRep(L_{P,\l})$ is the isotypical type of
$\Res_{L_{P,\l}}^{L_P} V$ for an irreducible $L_P$-module $V \in \mathscr S$ }\,\biggr\}\ .
\label{eqnRestrictMicroSupportRellFiber}
\end{multline}
For $V\in\IrrRep(L_P)$ we let $\Res_\l^* V=\Res_\l^! V$ denote the unique element of $\Res_\l^*\{V\}$.

\begin{prop*}
Let $\M$ be an $\L_{X_R(L_{R,\l})}$-module.  Then $\mS_w(\ihat_{R,\l}^*\M) = \Res_\l^* \mS_w(\M)$ and  $\mS_w(\ihat_{R,\l}^!\M) = \Res_\l^! \mS_w(\M)$.  For $V_\l\in \mS_w(\ihat_{R,\l}^{*/!}\M)$ we have
\begin{align}
[c(V_\l;\ihat_{R,\l}^*\M),d(V_\l;\ihat_{R,\l}^*\M)] &= \bigcup_{
\Res_\l^*V=V_\l} [c(V;\M),d(V;\M)] \text{ and} \\
[c(V_\l;\ihat_{R,\l}^!\M),d(V_\l;\ihat_{R,\l}^!\M)] &= \bigcup_{
\Res_\l^!V=V_\l} [c(V;\M),d(V;\M)] + \dim D_{R,h}\ .
\end{align}
\end{prop*}

\begin{proof}
For $P_\l\le Q_\l\in \Pl(\Xhat_{R,\l})$ we compute via
\eqref{eqnPullbackToFiberFormulas} that
\begin{equation}
\i_{P_\l}^*\ihat_{Q_\l}^!(\ihat_{R,\l}^*\M) =
\i_{P_\l}^*\ihat_{P,\l}^*(\i_P^*\ihat_Q^!\M)\ ,
\label{eqnMicroSupportOnFiberCalculation}
\end{equation}
where $P\le Q\in \Pl(X_R(L_{R,\l}))$ are defined as in
\S\ref{ssectPullbackToFiber}.  Thus an irreducible module $V_\l\in
\IrrRep(L_{P_\l})$ appears in the cohomology of
\eqref{eqnMicroSupportOnFiberCalculation} if and only if
$H(\i_P^*\ihat_Q^!\M)_V\neq 0$ for some $V=V_h\otimes V_\l\in
\IrrRep(L_P)$.  Furthermore since $\widetilde{L_{P,h}}=\widetilde{L_{R,h}}$, $ Q_\l \in [Q_{V_\l}^{\Xhat_{R,\l}},
Q_{V_\l}^{\prime \Xhat_{R,\l}}]$ if and only if $Q\in [Q_V^R, Q_V^{\prime R}]$.
\end{proof}

\section{Functoriality of Micro-support Assuming Equal-Rank}
\label{sectFunctorialityMicroSupportEqualRank}
In this section we present versions of Propositions
~\ref{ssectWeakMicroSupportClosedInverseImage} and \ref{ssectWeakMicroSupportInclusionFiber} that treat the full micro-support $\mS$ rather than just the weak micro-support $\mS_w$.  In order to preserve the condition of conjugate self-contragredience, we will assume certain strata of $\Xstar_\sigma$ are ``equal-rank''.

A symmetric space is called \emph{equal-rank} if it can be expressed as
$G(\RR)/K$, where $G$ is a reductive $\RR$-group, $K$ is a maximal compact
subgroup of $G(\RR)$, and $\CCrank G=\rank K$; this last condition holds if and only if every algebraic representation of $G$ is conjugate self-contradient (see
\S\ref{ssectEqualRankGroups}).   If $G$ is defined over $\QQ$ and $D$
is defined as usual to be $G(\RR)/KA_G$, then $D$ is equal-rank if and only
if $\lsp0G$ is equal-rank, that is, if and only if $V|_{\lsp0G}$
is conjugate self-contragredient for any algebraic $G$-module $V$.  A locally
symmetric space is called \emph{equal-rank} if it is an arithmetic quotient
of an equal-rank symmetric space.

We will use the notation of \S\S\ref{ssectPullbackToFiber},
\ref{sectSatakeCompactifications}, and \ref{sectFunctorialityInverseImages}.  Let $\Xstar_\sigma = \G\back \Dstar_\sigma$ be a Satake compactification with strata $(F_R)_{R\in \Pl^*_\sigma}$.

\begin{lem}
  \label{lemTransferLemma}
  Say $P\le \tilde P$ and set $R=P^\dag$ and $\tilde R={\tilde P}^\dag$.  Assume $\tilde P\cap R = P$.  Let $V$ and $\tilde V$ irreducible $L_P$- and $L_{\tilde P}$-modules respectively and assume $V\preccurlyeq \tilde V$.
  \begin{enumerate}
    \item \label{itemRestrictionVToPtildeEll} $\Res_{L_{\tilde P,\l}}^{L_P} V$ and  $\Res_{L_{\tilde P,\l}}^{L_{\tilde P}} \tilde V$ are isotypical of the same type.
  \item \label{itemTwoBoundaryEqualRankTransfer} If $V|_{M_P}$ is conjugate self-contragredient and
$F_{\tilde R}$ is equal-rank, then $\tilde V|_{M_{\tilde P}}$ is conjugate
self-contragredient.
  \end{enumerate}
\end{lem}

\begin{proof}
Write $\tilde V = \tilde V_h\otimes \tilde V_\l$, a tensor
product of irreducible modules according to the almost direct product
decomposition $L_{\tilde P}=\widetilde{L_{\tilde P,h}}L_{\tilde P,\l}$.  Since $N_P^{\tilde P} \subset \widetilde{L_{\tilde P,h}}$ by Lemma~\ref{lemTwoBoundary}\itemref{itemTwoBoundaryLeviFactors} we can write $H(\n_P^{\tilde P};\tilde V)$ as $H(\n_P^{\tilde P};\tilde V_h)\otimes \tilde V_\l$ according to the almost direct product $L_P = (\widetilde{L_{P,h}}L_{P,\l}^{\tilde P,h})\cdot L_{\tilde P,\l}$.  
Since $V\preccurlyeq \tilde V$, this shows that $\Res_{L_{\tilde P,\l}}^{L_P} V \subset \Res_{L_{\tilde P,\l}}^{L_P} H^\l(\n_P^{\tilde P};\tilde V)$ is isotypic of type $\tilde V_\l$.  This proves \itemref{itemRestrictionVToPtildeEll}.

For \itemref{itemTwoBoundaryEqualRankTransfer} $V|_{M_P}$ and thus $V|_{M_{\tilde P,\l}}$ are conjugate self-contragredient.  Thus by \itemref{itemRestrictionVToPtildeEll} $\tilde V_\l|_{M_{\tilde P,\l}}$ is conjugate
self-contragredient.  In addition $\tilde V_h$ is conjugate self-contragredient since
$\widetilde{L_{\tilde P,h}} = \widetilde{L_{\tilde R,h}}$ has no $\QQ$-split
center and $F_{\tilde R}$ is equal-rank.
\end{proof}

\subsection{Micro-support of the inverse image under a closed embedding}
  \label{ssectMicroSupportInverseImages}
Fix a Satake compactification $\Xstar_\sigma$ and let $\pi\colon \Xhat \to \Xstar_\sigma$ be the projection  defined in \S\textup{\ref{ssectSatakeAsQuotientOfRBS}}.  For a stratum $F_R$ of $\Xstar_\sigma$ let $U(F_R)= \coprod_{F_{R'}\ge F_R} F_{R'}$ be its open star neighborhood.
\begin{prop*}
  Let $F_R$ be a stratum of $\Xstar_\sigma$.  If all strata  $F_{R'}\ge F_R$ are equal-rank then all conclusions of Proposition ~\textup{\ref{ssectWeakMicroSupportClosedInverseImage}} hold true with $\mS_w$ replaced by $\mS$.  In particular
\begin{equation}
  \label{eqnConjugateSelfContragredient}
V_i|_{M_{P_i}} \text{ is conjugate self-contragredient} \qquad \text{for
  $0\le i\le N$.}
\end{equation}
\end{prop*}

\begin{proof}
If $V$ belongs to $\mS(\ihat_R^*\M)$ or $\mS(\ihat_R^!\M)$ then $V|_{M_{P_0}}$ is conjugate self-contragredient so if all $F_{R'}\ge F_R$ are equal-rank then Lemma  ~\ref{lemTransferLemma}\itemref{itemTwoBoundaryEqualRankTransfer} implies \eqref{eqnConjugateSelfContragredient} holds.
\end{proof}

\subsection{Micro-support of the inverse image under the inclusion of a fiber}
\label{ssectMicroSupportInclusionFiber}
We now generalize Proposition ~\ref{ssectWeakMicroSupportInclusionFiber} from weak micro-support to micro-support.

\begin{prop*}
Let $\M$ be an $\L_{X_R(L_{R,\l})}$-module and assume that $F_R$ is
equal-rank.  Then $\mS(\ihat_{R,\l}^*\M) = \Res_\l^* \mS(\M)$ and  $\mS(\ihat_{R,\l}^!\M) = \Res_\l^! \mS(\M)$.  For $V_\l\in \mS(\ihat_{R,\l}^{*/!}\M)$ we have
\begin{align}
[c(V_\l;\ihat_{R,\l}^*\M),d(V_\l;\ihat_{R,\l}^*\M)] &= \bigcup_{
\Res_\l^*V=V_\l} [c(V;\M),d(V;\M)] \text{ and}
\label{eqnDegreeFiberRestrictionEstimates} \\
[c(V_\l;\ihat_{R,\l}^!\M),d(V_\l;\ihat_{R,\l}^!\M)] &= \bigcup_{
\Res_\l^!V=V_\l} [c(V;\M),d(V;\M)] + \dim D_{R,h}\ .
\label{eqnDegreeFiberRestrictionEstimatesShriek}
\end{align}
\end{prop*}

\begin{proof}
As in \S\ref{ssectPullbackToFiber} let $P_\l\in\Pl(\Xhat_{R,\l})$ and let $P\in \Pl(X_R(L_{R,\l}))$ be such that $P/N_R=\widetilde{L_{R,h}}P_\l$; we have $L_P = \widetilde{L_{R,h}}L_{P,\l}$.
By Proposition ~\ref{ssectWeakMicroSupportInclusionFiber} it suffices to show that if $V_\l\in\IrrRep(L_{P_\l})$ and  $V=V_h\otimes V_\l \in\IrrRep(L_P)$ then $V_\l |_{M_{P,\l}}$ is conjugate self-contragredient if and only if $V|_{M_P}$ is  conjugate self-contragredient.  But since $F_R$ is equal-rank, $V_h$ is automatically conjugate self-contragredient so this holds.
\end{proof}

\section{The Basic Lemma}
\label{sectBasicLemma}
We now address the issue of estimating the term $[\tilde V:V]=\l(w)$ in
\eqref{eqnWeakMSInverseImageDegree} and
\eqref{eqnWeakMSInverseImageProperSupportDegree} under the equal-rank assumption of
Proposition \ref{ssectMicroSupportInverseImages}.  Our tool is a basic
lemma which, under the assumption that $V|_{M_P}$ is conjugate
self-contragredient, establishes the fundamental relationship between the
geometry of $\xi_V$ \emph{vis a vis} an $S_P^{\tilde P}$-root $\al$
and the contribution to the length of $w$ attributable to $\al$.  This
basic lemma is important in other contexts; a simple form already occurred
as Lemma~\ref{ssectSimpleBasicLemma}.

For simplicity we will assume that $\tilde V=E$ is an irreducible
$G$-module.  We use the notation of \S\ref{ssectSimpleBasicLemma}.

\subsection{}
\label{ssectRealSubmodules}
Let $\h$ be a fundamental Cartan subalgebra of $\levi_P$ and let
$\Phi^+(\levi_{P\CC},\h_\CC)$ be a $\theta$-stable positive system for
$\levi_P$.  Let $V$ be an irreducible $L_P$-module with highest weight
$\mu\in\h_\CC^*$ and assume $V|_{M_P}$ is conjugate self-contragredient.  We
have defined a reductive $\RR$-subgroup $L_P(\mu)\subset L_P$ in
\S\S\ref{ssectWeightCentralizer} and \ref{ssectVcentralizer} with roots
\begin{equation}
\{\,\g\in\Phi(\levi_{P\CC},\h_{\CC}) \mid
\langle\mu,\g\spcheck\rangle=0\,\}\ .
\end{equation}

Via a lift of $L_P$ to $P$, we obtain an adjoint representation of
$\levi_P$ (and hence $\levi_P(\mu)$) on $\n_{P\CC}$.  There is a unique
decomposition
\begin{equation*}
\n_{P\CC} = \bigoplus_{\substack{\text{ $\levi_P(\mu)$-irreducible}\\
                         \text{submodules} \\
                         F\subset \n_{P\CC}}}       F,
\end{equation*}
which coarsens the decomposition into root spaces.  Set $\t'_P = -\t_P$; note $\t'_P$ permutes the roots of $\n_{P\CC}$ and, by Proposition ~\ref{ssectConjugateSelfContragredient}\itemref{itemTauInvariant}, the roots of  $\levi_{P\CC}(\mu)$. It follows that if
$F$ is an irreducible component as above with weights
$\Phi(F,\h_\CC)$, then $\t'_P(\Phi(F,\h_\CC))$ is the set of weights of
another irreducible component, which we denote $\t'_P(F)$.   In this way we obtain an action of $\t'_P$
on the $\levi_P(\mu)$-irreducible components of $\n_{P\CC}$,
generalizing the action of $\t'_P$ on the roots of $\n_{P\CC}$.  Let%
\footnote{The quantities $\levi_P(\mu)$ and $\n_P(\mu)$ are
defined quite differently despite the similarity of notation.  We trust
that with this warning there should be no confusion.}
\begin{equation*}
\n_P(\mu) = \bigoplus_{\t'_P(F)=F} F ;
\end{equation*}
since the action of $\levi_P(\mu)$ on $\n_{P\CC}$ preserves $\n_{\al\CC}$
for any $\al\in\Phi(\n_P,\sa_P)$, we may likewise define $\n_\al(\mu)$.
Although $\n_P(\mu)$ and $\n_\al(\mu)$ depend on the choice of a lift of
$L_P$ to $G$, their dimensions are independent of this choice.  The
dimension of $\n_P(\mu)$ can vary depending on the choice of
$\Phi^+(\levi_{P\CC},\h_\CC)$ however; let $\n_P(V)$ denote any one of the
$\n_P(\mu)$ with maximal dimension.

\begin{BasicLemma}
\label{ssectBasicLemma}
Let $P$ be a parabolic $\QQ$-subgroup and let $w\in W_P$.  Let
$V=H^{\l(w)}(\n_P;E)_w$ have highest weight $\mu$ with respect to a
$\theta$-stable positive system.  Assume $V|_{M_P}$ is conjugate
self-contragredient.  For any $\al\in\Phi(\n_P,\sa_P)$ we have\textup:
\begin{enumerate} 
\item $\langle\xi_V+\hsr,\al\spcheck\rangle \le 
0 \implies 
\l_\al(w) \ge \frac12(\dim \n_\al+\dim \n_\al(\mu))$.
\label{itemBasicLemmaLessThan}
\item 
$\langle\xi_V+\hsr,\al\spcheck\rangle = 0 \implies \l_\al(w) = \frac12\dim
\n_\al$ and $\n_\al(\mu)=0$.
\label{itemBasicLemmaEqual}
\item $\langle\xi_V+\hsr,\al\spcheck\rangle\ge 0 \implies \l_\al(w) \le 
\frac12(\dim \n_\al-\dim \n_\al(\mu))$.
\label{itemBasicLemmaGreaterThan}
\end{enumerate}
Furthermore if $F\subset \n_{\al\CC}$ is an irreducible
$\levi_P(\mu)$-submodule, then $\Phi_w$ contains respectively
\itemref{itemBasicLemmaLessThan} at least one,
\itemref{itemBasicLemmaEqual} exactly one,
\itemref{itemBasicLemmaGreaterThan} at most one of
$\Phi(F,\h_\CC)$ and $\Phi(\t'_P(F),\h_\CC)$.
\end{BasicLemma}

\begin{proof}
Let $W^{L_P(\mu)}$ denote the Weyl group of $L_P(\mu)$ and let
$\lambda\in\h_\CC$ be the highest weight of $E$ relative to $\Phi^+ =
\Phi^+(\levi_{P\CC},\h_\CC) \cup \Phi(\n_{P\CC},\h_\CC)$.  We begin by
establishing
\begin{equation}
s\in W^{L_P(\mu)} \quad \Rightarrow \quad s\Phi_w =\Phi_w \text{ and }
sw\lambda = w\lambda \label{eqnCartier}
\end{equation}
following \cite{refnCartier}.  Express the highest weight $\mu = w(\lambda +
\hsr)-\hsr$ of $V$ as $w\lambda -{\sum_{\g\in\Phi_w}\g}$.  Since $s$ may be
written as a product of reflections in simple coroots orthogonal to $\mu$ we
see that $\mu=s\mu$ and hence
\begin{equation*}
w\lambda -{\sum_{\g\in\Phi_w}\g} =\mu = s\mu =s w\lambda -
{\sum_{\g\in s\Phi_w} \g}.
\end{equation*}
Note that $s\Phi_w\subset \Phi^+$; by canceling the terms of the right
hand sum that lie in $\Phi_w$ with the corresponding terms of the left hand
sum one may compute that
\begin{equation}
\lambda = w^{-1}sw\lambda - {\sum_{\g\in (w^{-1}\Psi)\cap \Phi^+} \g},
\label{eqnLambdaFormula}
\end{equation}
where $\Psi = s\Phi_w\cup -(\Phi^+\setminus s\Phi_w)$.  Since $\lambda$ is
dominant, we may replace $w^{-1}sw\lambda$ in \eqref{eqnLambdaFormula} by
$\lambda - \sum_i m_i
\al_i$ for $\al_i$ simple and $m_i\ge 0$.  In addition since the
roots of $(w^{-1}\Psi)\cap \Phi^+$ being subtracted in \eqref{eqnLambdaFormula}
are positive, it follows that all $m_i=0$ and that
$(w^{-1}\Psi)\cap \Phi^+ = \emptyset$; these two facts are equivalent to
the desired \eqref{eqnCartier}.

Now consider an $\levi_P(\mu)$-irreducible submodule $F\subset
\n_{\al\CC}$.  Since by \eqref{eqnCartier}
$\Phi_w$ is stable under $W^{L_P(\mu)}$, either all extremal weights of
$\Phi(F,\h_\CC)$ belong to $\Phi_w$ or none do.  Since furthermore membership
in $\Phi_w$ is given by a linear inequality on positive
roots and $\Phi(F,\h_\CC)$ is the convex hull of the extremal weights, this
implies that either all of
$\Phi(F,\h_\CC)$ is contained in $\Phi_w$ or none of it is.
An application of Lemma ~\ref{ssectSimpleBasicLemma} to any $\g\in
\Phi(F,\h_\CC)$ concludes the proof.
\end{proof}

\begin{cor}
\label{ssectBasicLemmaCor}
Let $P$ be a parabolic $\QQ$-subgroup and let $V$ be an irreducible
constituent of $H^\l(\n_P;E)$ for which $V|_{M_P}$ is conjugate
self-contragredient.  Then
\begin{align*}
(\xi_V+\hsr)|_{\sa_P} \in -\sa^{*+}_P &\implies \l \ge \frac12(\dim
\n_P+\dim\n_P(V)) \\
\intertext{and}
(\xi_V+\hsr)|_{\sa_P} \in \sa^{*+}_P &\implies \l \le \frac12(\dim
\n_P-\dim\n_P(V)). 
\end{align*}
\end{cor}

\begin{rem}
All results in this section and their proofs generalize immediately to the
situation in which $E$ is replaced by a ``virtual'' irreducible $G$-module,
whose highest weight $\lambda$ is allowed to be any dominant element of
$\h_\CC^*$, not necessarily integral.  They also generalize to parabolic
$k$-subgroups where $k$ is any subfield of $\RR$.
\end{rem}

\section{The Basic Lemma in the Presence of Equal-Rank Real Boundary
Components}
\label{sectEqualRankBasicLemma}
We need to replace the term $\dim \n_P(V)$ appearing in Corollary
~\ref{ssectBasicLemmaCor} by a more geometric expression.  This can be
done in the presence of a Satake compactification with equal-rank
\emph{real} boundary components.

\subsection{}
\label{ssectVcentralizerRationalFactors}
Let $V$ be an irreducible $L_P$-module and assume that $V|_{M_P}$ is
conjugate self-contragredient.  Let $D_P(V)$ be defined as in
\S\ref{ssectVcentralizer}.  Let $D\subset \lsb\RR\Dstar_\sigma$ be a
Satake compactification and write $V=V_\l\otimes V_h$ according to the almost direct product $L_P= \widetilde{L_{P,h}} L_{P,\l}$ \eqref{eqnLevisRationalFactors}.  This decomposition
induces a decomposition of symmetric spaces
\begin{equation}
D_P(V)= D_{P,h}(V)\times D_{P,\l}(V)\ .
\end{equation}
Here $D_{P,h}(V)$ is defined using $\Res_{\widetilde{L_{P,h}}}^{L_P} V$ so $D_{P,h}(V)  = D_{P,h}(V_h)$ (see the comment in \S\ref{ssectVcentralizer} regarding isotypical representations).  Similarly $D_{P,\l}(V) = D_{P,\l}(V_\l)$.

\begin{lem}
\label{ssectEqualRankLemma}
Let $P$ be a parabolic $\QQ$-subgroup and let $D\subset \lsb\RR
\Dstar_\sigma$ be a Satake compactification for which we assume all
real boundary components $D_{R',h}\ge D_{P^\dag,h}$ are equal-rank.
Let $V$ be an irreducible constituent of $H(\n_P;E)$ for which
$V|_{M_P}$ is conjugate self-contragredient.  Then
\begin{equation*}
\dim \n_P(V) \ge \dim \sa_P^G + \dim D_{P,\l}(V) \ .
\end{equation*}
\end{lem}
\begin{rem*}
A similar result was proved by Borel in
\cite[\S5.5]{refnBorelVanishingTheorem} under the weaker assumptions that
$D$ and $D_{P^\dag,h}$ are equal-rank and a certain condition (B) holds.
The condition (B) needed to be verified for each desired case and Borel
expressed the hope that it could be replaced by a more conceptual argument.
The proof offered below, while not transparent, is at least free from
case-by-case analysis.

Note also that condition (B) does not hold%
\footnote{In \cite[\S6.6]{refnBorelVanishingTheorem} it is mistakenly
asserted that condition (B) holds in this case.}
in one case of interest where the lemma above \emph{does} apply:
$G=\SO(p,q)$ with $p+q$ odd with the Satake compactification in which
$\sigma$ is connected only to the short simple $\RR$-root.  In this case
there are imaginary noncompact roots in $\levi_P$ (relative to a
fundamental Cartan subalgebra) which are strongly orthogonal to all real
roots in $\n_P$; the hypothesis that the intermediate real boundary
components are equal-rank and its use below in Sublemma
~\ref{ssectEqualRankSubLemma} was needed precisely to handle this
situation.
\end{rem*}

\begin{proof}[Proof of Lemma ~\textup{\ref{ssectEqualRankLemma}}]
We first establish notation.  Let $\theta$ be a Cartan involution of
$L_P(\RR)$ with Cartan decomposition $\levi_P = \k_P + \p_P + \sa_P$.  Let
$\h= \hb_{P,\k}+\hb_{P,\p} +\sa_P$ be a fundamental $\theta$-stable Cartan
subalgebra of $\levi_P$ and choose a $\theta$-stable positive system
$\Phi^+(\levi_{P\CC},\h_\CC)$ for $\levi_{P\CC}$; thus $\t_P=\theta$.  Let
$\mu$ be the highest weight of $V$ and let $L_P(\mu)$ be defined as in
\S\ref{ssectWeightCentralizer} with Cartan decomposition $\k_P(\mu) +
\p_P(\mu) +\sa_P$.  Let $\Phi^+(\p_P(\mu)_\CC,\hb_{P,\k\CC})$ denote the
positive nonzero weights of $\hb_{P,\k\CC}$ in $\p_P(\mu)_\CC$, counted with
multiplicity.  We extend the preceding notation to $L_{P,h}$ and $L_{P,\l}$
by adding the appropriate subscript.  Finally note that any reductive
subalgebra $\r\subset \levi_{P\CC}$ acts on $\n_{P\CC}$ by the adjoint
action of a lift of $\levi_{P\CC}$ and define $\n_P(\r)$ to be the sum of
the irreducible $\r$-submodules whose set of weights are $\t'_P$-stable;
thus $\n_P(\levi_P(\mu))=\n_P(\mu)$ as in \S\ref{ssectRealSubmodules}.  It
will suffice to prove $\dim \n_P(\mu) \ge \dim \sa_P^G + \dim D_{P,\l}(\mu)$.

Clearly
\begin{equation}
\label{eqnLinearDimension}
\dim \sa_P^G + \dim D_{P,\l}(\mu) =
2\cdot\#\Phi^+(\p_{P,\l}(\mu)_\CC,\hb_{P,\k\CC}) + \dim \hb_{P,\p} + \dim
\sa_P^G \ .
\end{equation}
We begin by bounding the last two terms.  Consider $\n_P(\h_\CC)$ as
defined above; the set $\Phi(\n_P(\h_\CC),\h_\CC)$ consists of the positive
roots $\g$ satisfying $\t'_P \g = \g$, that is, the real roots.  We claim
that
\begin{equation}
\label{eqnRealRoots}
\dim \n_P(\h_\CC) \ge \dim \hb_{P,\p} + \dim \sa_P^G\ .
\end{equation}
In fact there exists a set of strongly orthogonal roots
$\{\d_1,\dots,\d_r\}$ satisfying $\t'_P \d_i = \d_i$ and which span
$\hb_{P,\p}^*+\sa_P^G{}^*$.  To see this, proceed as in
\cite{refnBorelVanishingTheorem} to write the centralizer in
$\lsp0\mathfrak g$ of $\hb_{P,\k}$ as $\hb_{P,\k} + \mathfrak z_1$.  Then
$\mathfrak z_1$ has both a split Cartan subalgebra $\hb_{P,\p}+\sa_P^G$
and, since by hypothesis $\lsp0G =L_{G,h}$ is equal-rank, a compact Cartan
subalgebra.  The existence of $\{\d_1,\dots,\d_r\}$ then follows from
\cite{refnKostantConjugacyCartan}, \cite[Prop.~ 11]{refnSugiura},
\cite{refnSugiuraCorrection}.

We now attempt to bound the rest of \eqref{eqnLinearDimension}.
For $\g_\k\in \Phi^+(\p_{P,\l}(\mu)_\CC,\hb_{P,\l,\k\CC})$ choose
$\g\in\Phi^+(\levi_{P,\l}(\mu)_\CC,\hb_{P,\l\CC})$ such that
$\g|_{\hb_{P,\l,\k\CC}} = \g_\k$.  The root $\t_P\g$ is also in
$\Phi^+(\levi_{P,\l}(\mu)_\CC,\hb_{P,\l\CC})$ since
$\t_P(\mu|_{\hb_P})=\mu|_{\hb_P}$.  Let $\r_\g$ denote the reductive
subalgebra $\h_\CC\subset \r_\g \subset \levi_P(\mu)_\CC$ with root
system $\Phi\cap \Span\{\g,\t_P\g\}$.  Clearly
\begin{equation*}
\Phi(\n_P(\h_\CC),\h_\CC) \subset \Phi(\n_P(\r_\g),\h_\CC) \subset
\Phi(\n_P(\levi_P(\mu)_\CC),\h_\CC)= \Phi(\n_P(\mu),\h_\CC).
\end{equation*}
We will see below in Sublemma ~\ref{ssectEqualRankSubLemma}%
\itemref{itemEqualRankTauNoninvariant}\itemref{itemEqualRankTauInvariant}
that there exists a pair of roots $\eta\neq\t'_P\eta \in
\Phi(\n_P(\r_\g),\h_\CC) \setminus \Phi(\n_P(\h_\CC),\h_\CC) \subset
\Phi(\n_P(\mu),\h_\CC)$.  This would imply
\begin{equation}
\dim \n_P(\mu) - \dim \n_P(\h_\CC) \ge
2\cdot\#\Phi^+(\p_{P,\l}(\mu)_\CC,\hb_{P,\k\CC})
\end{equation}
(and hence prove the lemma in view of \eqref{eqnLinearDimension} and
\eqref{eqnRealRoots}) except that we have to deal with the possibility that
the sets $\Phi(\n_P(\r_\g),\h_\CC) \setminus \Phi(\n_P(\h_\CC),\h_\CC)$ for
the different $\g_\k$ may not be disjoint.

Since $|\g|=|\t_P\g|$ and no multiple of $\g-\t_P\g$ can be a root (recall
that $\levi_{P\CC}$ has no roots negated by $\t_P$), one easily checks that
$\{\g,\t_P\g\}$ is a base for the root system of $\r_\g$ and the type is
$A_1$ (if $\g=\t_P\g$), $A_1\times A_1$, or $A_2$.  Now consider $\g_\k\neq
\g_\k'\in \Phi^+(\p_{P,\l}(\mu)_\CC,\hb_{P,\l,\k\CC})$ and choose
corresponding roots $\g\neq \g'
\in\Phi^+(\levi_{P,\l}(\mu)_\CC,\hb_{P,\l\CC})$.  Suppose that
$\Phi(\n_P(\r_\g),\h_\CC)\setminus \Phi(\n_P(\h_\CC),\h_\CC)$ and
$\Phi(\n_P(\r_{\g'}),\h_\CC)\setminus \Phi(\n_P(\h_\CC),\h_\CC)$ have
nonempty intersection.  Then there exists a root $\eta\in
\Phi(\n_{P\CC},\h_\CC)$ for which $\t'_P\eta-\eta$ is a nonzero
$\ZZ$-linear combination of both $\{\g,\t_P\g\}$ and $\{\g',\t_P\g'\}$.
Since $\t'_P\eta-\eta$ is $\t_P$-invariant, this means that $\g+\t_P\g$ and
$\g'+\t_P\g'$ are scalar multiples of each other.  But $(1/2)(\g+\t_P\g)$
belongs to either the reduced root system $\Phi(\levi_{P\CC},\h_\CC)$ (if
$\g=\t_P\g$) or the reduced root system $\Phi(\k_{P\CC},\hb_{P,\k\CC})$ (if
$\g\neq \t_P\g$), and similarly for $(1/2)(\g'+\t_P\g')$.  Thus $\g+\t_P\g$
and $\g'+\t_P\g'$ can be scalar multiples of each other only if $\g$ (say)
satisfies $\g\neq \t_P\g$ and $\g'$ satisfies $\g'=\t_P\g'$.  It follows
from our comments earlier then that $\g'=\g+\t_P\g$ and $\r_\g$ has type
$A_2$.

So $\Phi(\n_P(\r_\g),\h_\CC)\setminus \Phi(\n_P(\h_\CC),\h_\CC)$ can
intersect with at most one other such set and it will suffice to show that
in this case $\Phi(\n_P(\r_\g),\h_\CC)\setminus \Phi(\n_P(\h_\CC),\h_\CC)$
contains an additional pair of roots $\eta'\neq \t'_P\eta'$.  But if such a
pair of roots did not exist, the only nontrivial $\r_\g$-submodule of
$\n_P(\r_\g)$ would have weights (after interchanging $\eta$ and
$\t'_P\eta$ if necessary) $\{\eta, \eta+\g=\t'_P\eta-\t_P\g, \t'_P\eta\}$
and hence $\{\eta,\g,\t_P\g\}$ would be the base of a root subsystem of
type $A_3$.
Sublemma~\ref{ssectEqualRankSubLemma}\itemref{itemEqualRankAThree}
following shows this situation cannot occur.
\end{proof}

\begin{sublem}
\label{ssectEqualRankSubLemma}
As in the proof of Lemma ~\textup{\ref{ssectEqualRankLemma}} let
$\g_\k\in \Phi^+(\p_{P,\l\CC},\hb_{P,\l,\k\CC})$ and choose $\g\in
\Phi^+(\levi_{P,\l\CC},\hb_{P,\l\CC})$ such that
$\g|_{\hb_{P,\l,\k\CC}}=\g_\k$.
\begin{enumerate}
\item If $\t_P\g \neq \g$, then there exists a root $\d=\t'_P\d
\in \Phi(\n_{P\CC},\h_\CC)$ such that {\rm(}after interchanging $\g$ and
$\t_P\g$ if necessary\/{\rm)} $\eta= \d-\g$ and
$\t'_P\eta=\d+\t_P\g$ are also roots.
\label{itemEqualRankTauNoninvariant}
\item If in \itemref{itemEqualRankTauNoninvariant} the roots
$\{\eta,\g,\t_P\g\}$ are a base of a root subsystem of type $A_3$, then
there also exists a root $\dtilde\neq \d$ with $\dtilde=\t'_P\dtilde$ and
one of $\dtilde\pm \g$ a root.
\label{itemEqualRankAThree}
\item If $\t_P\g=\g$, then there exists a root $\eta \in
\Phi(\n_{P\CC},\h_\CC)$ such that $\t'_P\eta=\eta+k\g$ for some
$k\in\ZZ\setminus\{0\}$.
\label{itemEqualRankTauInvariant}
\end{enumerate}
\end{sublem}
\begin{proof}
Let $\{\d_1,\dots,\d_r\}$ be a set of strongly orthogonal roots spanning
$\hb_{P,\p}^*+\sa_P^{G\,*}$ and satisfying $\t'_P \d_i = \d_i$ as in the
proof of Lemma ~\ref{ssectEqualRankLemma}.  If $\t_P\g\neq \g$ then
$\g|_{\hb_{P,\p}+\sa_P^G}\neq 0$ and so there exists $\d=\d_i$ such that
$\langle\g,\d\spcheck\rangle\neq 0$.  By interchanging $\g$ and $\t_P\g$ if
necessary, we may assume $\langle\g,\d\spcheck\rangle >0$.  Then $\eta = \d
-\g$ is a root which proves \itemref{itemEqualRankTauNoninvariant}.  Part
\itemref{itemEqualRankAThree} will follow if we show there exists
$\dtilde=\d_j\neq \d_i$ with $\langle\g,\dtilde\spcheck\rangle\neq 0$ as
well.  But if such a $\dtilde$ did not exist, $\g|_{\hb_{P,\p}+\sa_P^G} = c
\d$ for some scalar $c$ and it is easy to check that the conditions on $c$
implied by the equalities $|\d-\g|=|\g|=|\d|$ and $2(\g,\t_P\g)/|\g|^2 =
-1$ for a Weyl group invariant inner product on a root system of type $A_3$
are inconsistent.

If $\t_P\g = \g$ we claim that there is a nontrivial $\g$-root string
through some root in $\RR\g +\hb_{P,\p}^*+\sa_P^{G\,*}$.  On such a root
string the operator $\t'_P$ acts as reflection in $\g$ and hence the root
string is $\t'_P$-stable.  Thus part ~\itemref{itemEqualRankTauInvariant}
follows from the claim if we set $\eta$ to be any root in the root string
other than at the center.

The claim is obvious if $\d_i+\g$ or $\d_i-\g$ is a root for some $i$ we
are done, so we assume that $\g$ is strongly orthogonal to
$\{\d_1,\dots,\d_r\}$.
Recall \cite[Chapter~ VII, ~\S7]{refnKnapp} that a parabolic $\RR$-subgroup
$R$ is called \emph{cuspidal}%
\footnote{This usage is  different from that in
\cite{refnLanglandsFE} or \cite{refnBorelCasselman}; $R$ being cuspidal
in those references corresponds here to $R$ being defined over $\QQ$.}
if $L_R/Z(L_R)$ is equal-rank, where $Z(L_R)$ denotes the center.  Let
$R\subset P$ be a cuspidal parabolic $\RR$-subgroup with $\lsb\RR\hb_{R}
= \hb_{P,\k}$ and split component $\lsb\RR \sa_{R} = \hb_{P,\p}+\sa_P$; to
construct $R$ apply \cite[Prop.~7.87]{refnKnapp}.  Since
$\t_{R}=\t_P=\theta$ on $\h$, $\g$ remains an imaginary noncompact root for
$\levi_{R\CC}$.  Thus there is a Cayley transformation \cite[Chapter~ VI,
~\S7]{refnKnapp} which transforms $\g$ into a real root.  That is, there
exists a inner automorphism $\mathbf c_{\g}$ of $\levi_{R\CC}$ such that
$\h'\equiv \levi_R \cap \mathbf c_{\g}(\h_\CC)$ is a $\theta$-stable Cartan
subalgebra and $\mathbf c_{\g}\circ(s_{\g}\circ\theta) = \theta \circ
\mathbf c_{\g}$ on $\h$.  Let $R'\subset R$ be a cuspidal parabolic
$\RR$-subgroup with $\lsb\RR\hb_{R'}= \mathbf
c_{\g}(\ker\g\cap\lsb\RR\hb_R)$ and $\lsb\RR\sa_{R'}=\RR \mathbf
c_\g(\g\spcheck)+ \lsb\RR\sa_R$.  The roots $\{\d_1,\dots,\d_r,\mathbf
c_\g(\g)\}$ form a strongly orthogonal basis of $\lsb\RR\sa_{R'}^{G\,*}$
and to prove our claim we need to demonstrate that there is a nontrivial
$\mathbf c_\g(\g)$-root string through some root in
$\lsb\RR\sa_{R'}^{G\,*}$.

From our hypotheses on the Satake compactification we now construct another
independent set $\{\d'_1,\dots,\d'_{r+1}\}$ of roots spanning
$\lsb\RR\sa_{R'}^{G\,*}$.  Define a sequence of parabolic $\RR$-subgroups
\begin{equation*}
R'=R'_0\subset R'_1 \subset \dots \subset R'_{r+1} = G
\end{equation*}
inductively by choosing $R'_i$ to have type $\DRR_{R'}^{R'_i}=\{\al_1,\dots
,\al_i\}$ where $\al_i\in\DRR_{R'}^{R',R_{i-1}^{\prime\dag}}$.  Thus the
corresponding real Satake boundary component strictly increases at each stage,
\begin{equation*}
D_{R^{\prime\dag},h}= D_{R_0^{\prime\dag},h} \lneq  D_{R_1^{\prime\dag},h}
\lneq \dots \lneq D_{R_{r+1}^{\prime\dag},h} = D \ .
\end{equation*}
On the other hand, since $L_{P,h}$ is equal-rank and $\g$ was a root of
$\levi_{P,\l}$, we have $L_{R',h} = L_{P,h}$ and hence
$D_{P^\dag,h}=D_{R^{\prime\dag},h}$.  We now claim that each parabolic
subgroup $R'_i$ is cuspidal.  This follows since $L_{R'_i,\l}/Z(L_{R'_i})$
is an almost direct factor of $L_{R'_{i-1},\l}/Z(L_{R'_{i-1}})$, which is
equal-rank by induction, and $L_{R'_i,h}$ is equal-rank since the boundary
component $D_{R'_i{}^\dag,h}$ is equal-rank by hypothesis.  Thus for each
$i=1,\dots,r+1$ there exists a root $\d'_i$ which spans
$\lsb\RR\sa_{R'_{i-1}}^{R'_i\,*}$ \cite{refnKostantConjugacyCartan},
\cite{refnSugiura}, \cite{refnSugiuraCorrection}.

Since $\mathbf c_\g(\g)\in \lsb\RR\sa_{R'}^{G\,*}$ is nonzero there exists
$1\le i \le r+1$ such that $\langle\mathbf
c_\g(\g),\d'_i{}\spcheck\rangle\neq 0$.  The desired nontrivial $\mathbf
c_\g(\g)$-root string would now exist provided $\mathbf c_\g(\g)\neq \pm
\d'_i$.  To show $\mathbf c_\g(\g)$ cannot equal $\pm\d'_i$, write
$\DRR_{R'}^R=\{\al_k\}$ where $1\le k \le r+1$.  Thus $\mathbf c_\g(\g)$ is
a multiple of $\al_k$.  On the other hand, $\d'_i$ is proportional to the
projection of $\al_i$ to $\lsb\RR\sa_{R_{i-1}'}^{R_i'}$, that is,
$\d'_i\sim\al_i + \sum_{j < i} c_j \al_j$ where all $c_j=0$ if and only if
$\al_i\in \lsb\RR\sa_{R_{i-1}'}^{R_i'}$.  So if $\mathbf c_\g(\g)$ were to
equal $\pm \d'_i$ then $k=i$ and $\al_i\in \lsb\RR\sa_{R_{i-1}'}^{R_i'}$.
This implies that the simple $\RR$-roots underlying $\al_i$ as opposed to
$\{\al_1,\dots,\al_{i-1}\}$ belong to different components of $L_{R'_i,h}$.
This implies that the condition $\al_i\in\DRR_{R'}^{R',R_{i-1}'{}^\dag}$
actually means $\al_i\in\DRR_{R'}^{R',R'{}^\dag}$.  But since $\g$ is a root
of $\levi_{P,\l}$, $\mathbf c_\g(\g)$ is a root of $\levi_{R,\l}$ and hence
a linear combination of roots in $\DRR_{R'}^{R{}^\dag}=
\DRR_{R'}^{R'{}^\dag}$ and so can't be proportional to $\al_i$.
\end{proof}

\begin{cor}
\label{ssectEqualRankBasicLemmaCor}
Let $P$ be a parabolic $\QQ$-subgroup and let $D\subset \lsb\RR
\Dstar_\sigma$ be a Satake compactification for which we assume all
real boundary components $D_{R',h}\ge D_{P^\dag,h}$ are equal-rank.
Let $V$ be an irreducible constituent of $H^\l(\n_P;E)$ for which
$V|_{M_P}$ is conjugate self-contragredient.  Then
\begin{align*}
(\xi_V+\hsr)|_{\sa_P} \in -\sa^{*+}_P &\implies \l \ge \frac12(\dim
\n_P+\dim \sa_P^G + \dim D_{P,\l}(V)) \\ \intertext{and}
(\xi_V+\hsr)|_{\sa_P} \in \sa^{*+}_P &\implies \l \le \frac12(\dim
\n_P-\dim \sa_P^G - \dim D_{P,\l}(V)).
\end{align*}
\end{cor}
\begin{proof}
Combine Lemma ~\ref{ssectEqualRankLemma} and Corollary
~\ref{ssectBasicLemmaCor}.
\end{proof}
\begin{rem*}
For the Baily-Borel-Satake compactification in the Hermitian case, the
second inequality (for $\xi_V$ in a larger cone) is the content of
\cite[Prop.~ 11.1]{refnSaperSternTwo}.  Shortly after
\cite{refnSaperSternTwo} appeared, Saper and Stern noted that the result
generalized to Satake compactifications of equal-rank symmetric spaces
where all real boundary components were equal-rank---a proof replaced
Lemmas ~11.6 and 11.7 from \cite{refnSaperSternTwo} with a case-by-case
analysis based on the list appearing in \cite[(A2)]{refnZuckerLtwoIHTwo}.
The proof here via Corollary~\ref{ssectBasicLemmaCor} and Lemma~
\ref{ssectEqualRankLemma} is independent of classification theory.
\end{rem*}

\section{Restriction to Fibers of $\pi\colon \Xhat \to \Xstar_\sigma$}
\label{sectRestrictMicroSupportToFiber}

\begin{thm}
\label{ssectRestrictMicroSupportToFiber}
Let $x\in F_R$ be a point on a stratum of a Satake compactification $\Xstar_\sigma$ for which
all the real boundary components $D_{R',h}\ge D_{R,h}$ of the
associated Satake compactification $\lsb\RR\Dstar_\sigma$ are equal-rank.  Set $W=\pi^{-1}(U(F_R))$ and  consider the inclusions
\[
\pi^{-1}(x)=\Xhat_{R,\l}\overset{\ihat_{R,\l}}\longrightarrow X_R(L_{R,\l})=
\Xhat_R\cap W \overset{\ihat_R}\longrightarrow W \overset{i_W}\longrightarrow \Xhat\ .
\]
Let $\M$ be an $\L_{\Xhat}$-module and let $\mathscr S= \coprod_{P\in\Pl(W)}\IrrRep(L_P)$.
Then
\begin{equation}
  \begin{split}
    \mS(\ihat_{R,\l}^*\ihat_R^* i_W^*\M) &\subset \Res_\l^* \ihat_R^*(\mS(\M)\cap \mathscr S)
    \qquad\text{and} \\
    \mS(\ihat_{R,\l}^!\ihat_R^! i_W^!\M) &\subset \Res_\l^! \ihat_R^!(\mS(\M)\cap \mathscr S)\ .
  \end{split}
  \label{eqnMicrosupportRestrictionToPiInversePoint}
\end{equation}

We have the estimates
\begin{equation}
\begin{split}
d(\ihat_{R,\l}^*\ihat_R^* i_W^*\M) \le \smash[b]{\sup_{\substack{\tilde V \in \mS(\M)\\ 
 F_{\tilde R}\ge F_R,\ \tilde R= \tilde P^\dag}}}\bigl(&\, d(\tilde V;\M) + \frac12(\dim
 D_{\tilde P,\l} + \dim D_{\tilde P,\l}(\tilde V)) \\
 &\qquad + \frac12\codim_{F_{\tilde R}^*} F_R - \dim \sa_{\tilde P\cap
 R}^{\tilde P} \,\bigr)
\end{split}
\label{eqnDREllStarEstimate}
\end{equation}
and 
\begin{equation}
\begin{split}
  c(\ihat_{R,\l}^!\ihat_R^! i_W^!\M) \ge \smash[b]{\inf_{\substack{\tilde V\in \mS(\M)\\ F_{ \tilde R}\ge F_R, \ \tilde R= \tilde P^\dag}}}
  \bigl(&\, c(\tilde V;\M) + \frac12(\dim D_{\tilde P,\l} -
 \dim D_{\tilde P,\l}(\tilde V)) \\
 &\qquad + \frac12\codim_{F_{\tilde R}^*} F_R + \dim F_R
 \bigr).
\end{split}
\label{eqnDREllShriekEstimate}
\end{equation}
In these estimates, $\tilde V$ is an irreducible $L_{\tilde P}$-module.  Furthermore the inequality in \eqref{eqnDREllShriekEstimate} is strict if $F_{\tilde R} > F_R$ for all $\tilde V\in \mS(\M)$.
\end{thm}

\begin{proof}
Equation \eqref{eqnMicrosupportRestrictionToPiInversePoint} is simply
Propositions ~\ref{ssectMicroSupportOpenInverseImage}, ~\ref{ssectMicroSupportInverseImages}, and
\ref{ssectMicroSupportInclusionFiber}.  These propositions further imply
that if $V_\l\in \mS(\ihat_{R,\l}^*\ihat_R^* i_W^*\M)$, then there exists an $L_P$-module $V\in
\mS(\ihat_R^* i_W^*\M)$ and an $L_{\tilde P}$-module $\tilde V\in \mS(i_W^*\M)$ such that  $\Res_{L_{P,\l}}^{L_P}V$ is isotypic of type $V_\l$ and there exist irreducible $L_{P_i}$-modules
\begin{equation*}
V=V_0\prec_+ V_1 \prec_+
\dots \prec_+ V_N=\tilde V
\end{equation*}
satisfying \eqref{eqnPversusR}--\eqref{eqnNumberOfInductionSteps}, \eqref{eqnConjugateSelfContragredient}.  By \eqref{eqndMDefn},
$d(\ihat_{R,\l}^*\ihat_R^* i_W^*\M)$ is bounded by the supremum of
\begin{equation}
  \label{eqnZerothDBound}
 \frac12(\dim D_{P,\l}+\dim D_{P,\l}(V_\l)) + d(V_\l;\ihat_{R,\l}^*\ihat_R^* i_W^*\M)
\end{equation}
over all $V_\l\in \mS(\ihat_{R,\l}^*\ihat_R^* i_W^*\M)$.  Now by 
\eqref{eqnDegreeFiberRestrictionEstimates} and
\eqref{eqnWeakMSInverseImageDegree} we see that \eqref{eqnZerothDBound} is bounded by the supremum of
\begin{equation}
  \label{eqnFirstDBound}
\frac12(\dim D_{P,\l}+\dim D_{P,\l}(V)) + d(\tilde V; \M) + [\tilde V; V]\ .
\end{equation}
over all $V$ and $\tilde V$ (note $D_{P,\l}(V_\l) = D_{P,\l}(V)$ by \S\ref{ssectVcentralizerRationalFactors}).

The term $[\tilde V; V]$ is the sum $\sum_{i=0}^{N-1} [V_{i+1},V_i]$.  In order to bound each summand recall by Lemma ~\ref{lemTwoBoundary}\itemref{itemTwoBoundaryLeviFactors} that $N_{P_i}^{P_{i+1}} \subset \widetilde{L_{P_{i+1},h}}$ and that there is an almost direct product factorization $L_{P_i}=\widetilde{L_{P_i,h}}L_{P_i,\l} = \widetilde{L_{P_i,h}}
L_{P_i,\l}^{P_{i+1},h}L_{P_{i+1},\l}$.
Let $D_{P_{i},\l}^{P_{i+1},h}$ denote the symmetric space associated to
$L_{P_i,\l}^{P_{i+1},h}$.  Then by \eqref{eqnStrictDominantCone} we can apply Corollary
\ref{ssectEqualRankBasicLemmaCor} to the parabolic $(P_i/N_{P_{i+1}})\cap
L_{P_{i+1},h} \subset L_{P_{i+1},h}$ and the induced Satake
compactification of $D_{P_{i+1},h}$ in order to conclude that
\begin{equation}
[V_{i+1}:V_i] \le \frac12 ( \dim \n_{P_i}^{P_{i+1}} - \dim
\sa_{P_i}^{P_{i+1}} - \dim D_{P_{i},\l}^{P_{i+1},h}(V_i) )\ .
\label{eqnOneStepDegreeBound}
\end{equation}

Since $D_{P_{i+1},\l}(V_{i}) =  D_{P_{i+1},\l}(V_{i+1})$  by Lemma~\ref{lemTransferLemma}\itemref{itemRestrictionVToPtildeEll}, we have a factorization $D_{P_i,\l}(V_i) = D_{P_{i},\l}^{P_{i+1},h}(V_i) \times D_{P_{i+1},\l}(V_{i+1})$ and hence $\dim D_{P_{i},\l}^{P_{i+1},h}(V_i) = D_{P_i,\l}(V_i) -  D_{P_{i+1},\l}(V_{i+1})$.
Consequently if we sum
\eqref{eqnOneStepDegreeBound} for $i=0,\dots,N-1$ we obtain
\begin{equation*}
[\tilde V: V] \le \frac12 ( \dim \n_{P}^{\tilde P} - \dim \sa_{P}^{\tilde P}
+ \dim D_{\tilde P,\l}(\tilde V) -  \dim D_{P,\l}(V) )\ .
\label{eqnReprDegreeEstimate}
\end{equation*}
Combined with \eqref{eqnFirstDBound} this shows $d(\ihat_{R,\l}^*\ihat_R^* i_W^*\M)$ is bounded by the supremum of
\begin{equation}
  \label{eqnSecondDBound}
 d(\tilde V; \M) +  \frac12 ( \dim D_{P,\l} + \dim D_{\tilde P,\l}(\tilde V)) + \frac12 ( \dim \n_{P}^{\tilde P} - \dim \sa_{P}^{\tilde P}  ) \ .
\end{equation}
Since however
\begin{equation*}
  \codim_{F_{\tilde R}^*} F_R = \dim D_{\tilde P,h} - \dim D_{P,h}
= \dim D_{P,\l} - \dim D_{\tilde P,\l} + \dim \n_P^{\tilde P} + \dim \sa_P^{\tilde P}
\end{equation*}
equation \eqref{eqnSecondDBound} can be re-expressed as
\begin{equation*}
  d(\tilde V; \M) +  \frac12 (  \dim D_{\tilde P,\l} + \dim D_{\tilde P,\l}(\tilde V))  + \frac12  \codim_{F_{\tilde R}^*} F_R - \dim \sa_{P}^{\tilde P}
\end{equation*}
which proves \eqref{eqnDREllStarEstimate} since $ \dim \sa_{P}^{\tilde P} \ge
\dim \sa_{\tilde P\cap R}^{\tilde P}$.

The proof of \eqref{eqnDREllShriekEstimate} is analogous but observe the shifts by $\dim D_{R,h}$ from \eqref{eqnDegreeFiberRestrictionEstimatesShriek} and by $N\ge0$ from \eqref{eqnWeakMSInverseImageProperSupportDegree}.  To establish the final remark note if $F_{\tilde R}>F_R$ then $\tilde P> P$ and thus $N\ge 1$.
\end{proof}

\begin{cor}
\label{ssectRestrictMicroSupportToFiberCorollary}
In the above situation, assume further that $\emS(\M)=\{E\}$ for
$E\in\IrrRep(G)$ with $c(E;\M)=d(E;\M)=0$.  Then for a proper stratum $F_R$ of $\Xstar_\sigma$,
\begin{align*}
d(\ihat_{R,\l}^*\ihat_R^* i_W^*\M) & \le \frac12\codim F_R - \dim \sa_{R}^{G}
\quad\text{and}\quad \\
c(\ihat_{R,\l}^!\ihat_R^! i_W^!\M) &> \frac12\codim F_R + \dim F_R \ .
  \end{align*}
\end{cor}

\begin{proof}
Any $L_P$-module $V\in \mS(\M)$ satisfies $V\preccurlyeq_0 E$ by Corollary
~\ref{ssectPartialOrderingOnMicroSupport} and therefore has the form
$H^{\l(w)}(\n_{P};E)_w$ where $w$ is fundamental for $P$ in $G$ by Lemma
~\ref{lemPartialOrdering}.  Thus $P$ contains a fundamental parabolic $\RR$-subgroup by Proposition ~\ref{ssectBorelCasselman}.  However since $D$ is equal-rank the only
fundamental parabolic $\RR$-subgroup is $G$.  Hence $P=G$ and $\mS(\M)=\{E\}$.  Now
apply \eqref{eqnDREllStarEstimate}, \eqref{eqnDREllShriekEstimate}, and the remark following with $\tilde V=E$ and $\tilde P=\tilde R = G$.
\end{proof}

\section{Application: the Conjecture of Rapoport and Goresky-MacPherson}
\label{sectRapoportConjecture}

Recall the middle perversities are $m(k) =\lfloor(k-2)/2\rfloor$ and $n(k) =
\lfloor(k-1)/2\rfloor$.
\begin{thm}
\label{ssectRapoportConjecture}
Let $\Xstar$ be a Satake compactification of $X$ for which all real
boundary components $D_{R,h}$ are equal-rank, and let $\pi\colon \Xhat \to
\Xstar$ be the projection from the reductive Borel-Serre compactification.
Let $p$ be a middle perversity.  Then there is a natural quasi-isomorphism
$\pi_*\IpC(\Xhat;\EE) \cong \IpC(\Xstar;\EE)$.
\end{thm}
\begin{rem*}
When $D$ is Hermitian and $\Xstar$ is the Baily-Borel-Satake
compactification, the equal-rank hypothesis is automatic.  In this case the
theorem was conjectured independently by Rapoport
\cite{refnRapoportLetterBorel} and Goresky and MacPherson
\cite{refnGoreskyMacPhersonWeighted}; a proof in this case for $\QQrank
G=1$ was given by Saper and Stern \cite[Appendix]{refnRapoport}.
\end{rem*}

\begin{proof}
For $x\in F_R$, a proper stratum of $\Xstar$, let $\i_x\colon
\{x\}\hookrightarrow \Xstar$, $\ihat_{R,\l}\colon \pi^{-1}(x) =
\Xhat_{R,\l}\hookrightarrow \Xhat_R$, and $\ihat_R\colon \Xhat_R \hookrightarrow \Xhat$ be the inclusion maps.  Let $q(k)= k -2
- p(k)$ be the dual middle perversity.  By the local characterization of
intersection cohomology \cite{refnGoreskyMacPhersonIHTwo}, \cite[V,
4.2]{refnBorelIntersectionCohomology} on $\Xstar$ we need to verify that
\begin{enumerate}
\item $\pi_*\IpC(\Xhat;\EE)$ is $\X$-clc,
\label{itemIHCLC}
\item $\pi_*\IpC(\Xhat;\EE)|_X$ is quasi-isomorphic to $\EE$,
\label{itemIHOpenStratum}
\item $H^i(\i_x^*\pi_*\IpC(\Xhat;\EE)) = 0$ for $x\in F_R$, $i> p(\codim
F_R)$, and
\label{itemIHVanishing}
\item $H^i(\i_x^!\pi_*\IpC(\Xhat;\EE)) = 0$ for $x\in F_R$, $i< \dim X -
q(\codim F_R)$.
\label{itemIHCovanishing}
\end{enumerate}
Condition ~\itemref{itemIHCLC} follows from \cite[V,
10.16]{refnBorelIntersectionCohomology} and condition
~\itemref{itemIHOpenStratum} is obvious.  For the others, calculate
\begin{alignat*}{2}
H(\i_x^*\pi_*\IpC(\Xhat;\EE)) &\cong
H(\Xhat_{R,\l};\ihat_{R,\l}^*\ihat_R^*\IpC(\Xhat;\EE)) &\qquad
&\text{by \cite[V, 10.7]{refnBorelIntersectionCohomology},} \\
&\cong H(\Xhat_{R,\l};\ihat_{R,\l}^*\ihat_R^*\IpC(E)) &\qquad&\text{by Thm.
~\ref{ssectRealizationLModules} and Prop.
~\ref{ssectIntersectionCohomology}}
\end{alignat*}
and similarly $H(\i_x^!\pi_*\IpC(\Xhat;\EE)) \cong
H(\Xhat_{R,\l};\ihat_{R,\l}^!\ihat_R^!\IpC(E))$.  Now since all strata of $\Xstar$
have even dimension by the equal-rank hypothesis, $p(\codim F_R) =
\frac12\codim F_R - 1$ and $\dim X - q(\codim F_R)= \frac12\codim F_R + \dim F_R + 1$.
Thus we need to demonstrate that
\begin{equation}
\begin{alignedat}{2}
H^i(\Xhat_{R,\l};\ihat_{R,\l}^*\ihat_R^*\IpC(E))&=0  &\qquad&\text{for $i\ge
  \frac12\codim F_R$,}  \\ 
H^i(\Xhat_{R,\l};\ihat_{R,\l}^!\ihat_R^!\IpC(E))&=0  &\qquad&\text{for $i\le
  \frac12\codim F_R + \dim F_R$.}
\end{alignedat}
\label{eqnVanishingCoVanishing}
\end{equation}

Under the equal-rank hypotheses, the possible types of irreducible
$\QQ$-root systems that can occur in $G$ are $B_n$, $BC_n$, $C_n$, and
$G_2$ \cite[(A.2)]{refnZuckerLtwoIHTwo}.  Thus we can apply Corollary
~\ref{ssectIHMicroPurityCorollary} to see that the essential micro-support
of $\IpC(E)$ is simply $\{E\}$; Corollary
~\ref{ssectRestrictMicroSupportToFiberCorollary} then implies that
$d(\ihat_{R,\l}^*\ihat_R^*\IpC(E)) < \frac12\codim F_R$ and $c(\ihat_{R,\l}^!\ihat_R^!\IpC(E)) >
\frac12\codim F_R + \dim F_R$.  The proof of \eqref{eqnVanishingCoVanishing} is concluded
by applying our vanishing theorem for the cohomology of $\L$-modules,
Theorem ~\ref{ssectGlobalVanishing}.
\end{proof}

\subsection{}
\label{ssectGoreskyHarderMacPhersonTheorem}
If we replace the use of Corollary ~\ref{ssectIHMicroPurityCorollary} by
Theorem ~\ref{ssectWeightCohomologyMicroSupport} in the preceding proof we
obtain the following theorem.

\begin{thm*}
Let $\Xstar$ be a Satake compactification of $X$ for which all real
boundary components $D_{R,h}$ are equal-rank, and let $\pi\colon \Xhat \to
\Xstar$ be the projection from the reductive Borel-Serre compactification.
Let $\eta$ be a middle weight profile and let $p$ be a middle perversity.
Then there is a natural quasi-isomorphism $\pi_*\WnC(\Xhat;\EE) \cong
\IpC(\Xstar;\EE)$.
\end{thm*}
\begin{rem*}
When $D$ is Hermitian and $\Xstar$ is the Baily-Borel-Satake
compactification, this theorem was the main result of
\cite{refnGoreskyHarderMacPherson}.
\end{rem*}


\bibliographystyle{amsplain}
\bibliography{references}
\end{document}